\documentclass[a4paper]{book}
\usepackage{xcolor}
\definecolor{rtg_red}{RGB}{235,0,9}
\definecolor{rtg_darkblue}{RGB}{0,77,128}
\definecolor{rtg_lightblue}{RGB}{0,175,192}
\usepackage{hyperref}

\AtEndPreamble{%
    \RequirePackage{hyperref}
    \hypersetup{%
        colorlinks    = true,
        linkcolor     = rtg_red,
        citecolor     = rtg_darkblue,
        urlcolor      = rtg_lightblue,
    }
}

\usepackage{emptypage} 
    
\usepackage{subfiles}
\usepackage{amsmath, amssymb, amsthm}

\usepackage{enumitem}   

\usepackage{mathtools}   
\mathtoolsset{showonlyrefs} 

\usepackage{multirow}   

\usepackage[style=authoryear, url=false, dashed=false, uniquelist=false]{biblatex}
\AtEveryBibitem{\clearfield{pages}}
\usepackage{xurl}   

\usepackage{graphics}   

\usepackage[color=teal!50!white]{todonotes}

\newtheoremstyle{myTheoremStyle}
    {}
    {}
    {\itshape}
    {}
    {\bfseries}
    {}
    {\newline}
    {}

\theoremstyle{myTheoremStyle}
\newcounter{theoremCounter}
\newtheorem{theorem}[theoremCounter]{Theorem}
\newtheorem{definition}[theoremCounter]{Definition}
\newtheorem{lemma}[theoremCounter]{Lemma}
\newtheorem{remark}[theoremCounter]{Remark}

\makeatletter
\renewenvironment{proof}[1][\proofname\leavevmode\\]{\par
\pushQED{\qed}%
\normalfont \topsep6\p@\@plus6\p@\relax
\trivlist
\item\relax{\itshape
{#1}}\hspace\labelsep\ignorespaces
}{%
\popQED\endtrivlist\@endpefalse
}
\makeatother

\newcommand{\expliciteCalculation}[1]{\textcolor{gray}{
    \ifx\showExplicitCalculations\undefined%
    \else%
        #1%
    \fi%
    }%
}
  
\newcommand{\D}{\mathbb{D}}

\newcommand{\Nu}{\textnormal{Nu}}
\newcommand{\Ra}{\textnormal{Ra}}
\renewcommand{\Pr}{\textnormal{Pr}}
\newcommand{\np}{{n^+}}
\newcommand{\nm}{{n^-}}
\newcommand{\npm}{{n^\pm}}
\newcommand{\taupm}{{\tau^\pm}}
\newcommand{\hp}{{h^+}}
\newcommand{\hm}{{h^-}}
\newcommand{\kappap}{\kappa_+}
\newcommand{\kappam}{\kappa_-}
\newcommand{\alphap}{\alpha_+}
\newcommand{\alpham}{\alpha_-}
\newcommand{\up}{u_+}
\newcommand{\um}{u_-}
\newcommand{\gammap}{{\gamma^+}}
\newcommand{\gammam}{{\gamma^-}}

\newcommand{\Omegatilde}{{\Omega^{\approx}}}
\newcommand{\Omegabar}{{\Omega^{\eqsim}}}
\newcommand{\Omegadelta}{{\Omega^{\delta}}}
\newcommand{\Omegatildedelta}{{\tilde\Omega^\delta}}
\newcommand{\ngammamp}{{n_{\gammam}^{+}}}
\newcommand{\ngammapp}{{n_{\gammap}^{+}}}
\newcommand{\Omegadash}{{\Omega^{\_}}}
\newcommand{\Omegadashdash}{{\Omega^{\_\hspace{2pt}\_}}}
\newcommand{\Omegavert}{{\Omega^\shortmid}}
\newcommand{\Omegavertvert}{{\Omega^\shortmid_\shortmid}}
\newcommand{\Omegadeltavert}{{\Omega^\delta_\shortmid}}
\newcommand{\Omegastraight}{{\Omega^{=}}}
\newcommand{\omegatildepm}{{\tilde \omega^{\pm}}}
\newcommand{\omegatildem}{{\tilde \omega^{-}}}
\newcommand{\omegatildep}{{\tilde \omega^{+}}}
\newcommand{\omegabarpm}{{\bar \omega^{\pm}}}
\newcommand{\omegabarm}{{\bar \omega^{-}}}
\newcommand{\omegabarp}{{\bar \omega^{+}}}
\newcommand{\omegahatm}{{\hat \omega^{-}}}
\newcommand{\omegahatp}{{\hat \omega^{+}}}
\newcommand{\omegahatpm}{{\hat \omega^{\pm}}}
\newcommand{\omegahatExpl}{{\hat \omega}}
\newcommand{\tmax}{T}

\DeclareMathOperator*{\esssup}{ess\ sup}
\DeclareMathOperator*{\essinf}{ess\ inf}

\def\Xint#1{\mathchoice
   {\XXint\displaystyle\textstyle{#1}}%
   {\XXint\textstyle\scriptstyle{#1}}%
   {\XXint\scriptstyle\scriptscriptstyle{#1}}%
   {\XXint\scriptscriptstyle\scriptscriptstyle{#1}}%
   \!\int}
\def\XXint#1#2#3{{\setbox0=\hbox{$#1{#2#3}{\int}$}
     \vcenter{\hbox{$#2#3$}}\kern-.5\wd0}}

\def\dashint{\Xint-}

\let\temp\phi
\let\phi\varphi
\let\varphi\temp

\let\temp\theta
\let\theta\vartheta
\let\vartheta\temp

\let\temp\epsilon
\let\epsilon\varepsilon
\let\varepsilon\temp

\usepackage[ddmmyyyy]{datetime}

\usepackage[pagestyles]{titlesec}
\newpagestyle{blank}{%
\sethead[][][]{}{}{} 
\setfoot{}{\thepage}{}
}
\newpagestyle{myPageStyle}{%
\sethead[\thechapter. \chaptertitle][][]{}{}{\thesection. \sectiontitle} 
\setfoot{}{\thepage}{}
}

\begin{document}
\pagestyle{blank}
\frontmatter
\begin{titlepage}
  \begin{center}
     \addtolength{\baselineskip}{15pt}
     \textbf{\huge Buoyancy-Driven Flows With Navier-Slip Boundary Conditions}\\
     \bigskip

     \addtolength{\baselineskip}{-15pt}

    \vspace{12ex}
    A\\
    \bigskip
    \begin{large}Dissertation\end{large}\\
    \bigskip
    submitted to the Department of Mathematics\\
    \medskip
    of the Faculty of Mathematics, Informatics and Natural Sciences\\
    \medskip
    of the University of Hamburg
    \bigskip


    \vspace{7ex}
    by
    \bigskip

    \begin{large}Fabian Bleitner\end{large}\medskip

    \vspace{20ex}
    \begin{tabular}[h]{ll}

        Reviewers\hspace{2cm}   & Dr.\ Camilla\ Nobili\\
                                & Prof.\ Dr.\ Jens\ Rademacher\\
                                & Prof.\ Dr.\ Jared\ Whitehead\\

    \end{tabular}

    \vspace{20ex}
    Hamburg 2024
    \smallskip

  \end{center}
\end{titlepage}

\hfill Date of disputation\hspace{1cm} 10.07.2024
\tableofcontents
\markboth{}{}
\newpage
\chapter*{Abstract}
\markboth{}{}

In this dissertation two-dimensional buoyancy-driven flows are investigated. While usually the Navier-Stokes equations are equipped with no-slip boundary conditions here we focus on the Navier-slip conditions that, depending on the system at hand, better reflect the physical behavior. In particular, we study two systems, Rayleigh-Bénard convection and a closely related problem without thermal diffusion. In the former, bounds on the vertical heat transfer, given by the Nusselt number, with respect to the strength of the buoyancy force, characterized by the Rayleigh number, are derived. These bounds hold for a broad range of applications, allowing for non-flat boundaries, any sufficiently smooth positive slip coefficient, and are valid over all ranges of the Prandtl number, a system parameter determined by the fluid. For the thermally non-diffusive system, regularity estimates are proven. Up to a certain order, these bounds hold uniformly in time, which, combined with estimates for their growth, provide insight into the long-time behavior. In particular, solutions converge to the hydrostatic equilibrium, where the fluid's velocity vanishes and the buoyancy force is balanced by the pressure gradient.

\newpage
{\let\cleardoublepage\relax\chapter*{Zusammenfassung}}

In dieser Dissertation werden zweidimensionale auftriebsgetriebene Flüsse untersucht. Während die Navier-Stokes-Gleichungen normalerweise mit Haftrandbedingungen versehen sind, fokussieren wir uns hier auf Navier-Randbedingungen, die abhängig vom betrachteten System, das physikalische Verhalten besser widerspiegeln. Insbesondere untersuchen wir zwei Systeme, Rayleigh-Bénard Konvektion und ein eng verwandtes Problem ohne Wärmediffusion. Im ersten Modell werden Grenzen für den vertikalen Wärmetransport, welcher durch die Nußelt-Zahl gegeben ist, bezüglich der Stärke des Auftriebskraft, charakterisiert durch die Rayleigh-Zahl, hergeleitet. Diese Abschätzungen gelten für einen großen Anwendungsbereich, der gekrümmte Ränder und beliebige, ausreichend glatte, positive Gleitkoeffizienten zulässt, und sind für alle Prandtl-Zahlen, einem durch das Fluid bestimmten Systemparameter, gültig. Für das System ohne Wär\-me\-dif\-fu\-si\-on werden Regularitätsabschätzungen bewiesen. Diese halten bis zu einer gewissen Ordnung gleichmäßig bezüglich der Zeit, was zusammen mit Abschätzungen für deren Wachstum Einsicht in das Langzeitverhalten gibt. Insbesondere konvergieren Lösungen zum hydrostatischen Gleichgewicht, in dem das Geschwindigkeitsfeld verschwindet und die Auftriebskraft durch den Druckgradienten ausgeglichen wird.

\chapter*{List of Publications}
The doctoral studies led to the following publications.
\begin{itemize}
    \item 
        In \cite{bleitnerNobili2024Bounds} Rayleigh-Bénard convection on a domain with potentially curved boundaries is studied. The results are given in Theorem~\ref{theorem_nonlinear_1_2_bound} and Theorem~\ref{theorem_nonlinear_main_theorem} and proven in Section \ref{section:nonlinearity_paper}.
    \item 
        The submitted preprint \cite{bleitner2024Scaling} investigates Rayleigh-Bénard convection with flat boundary profiles. The slightly adjusted results are given in Theorem~\ref{theorem_flat} and Section \ref{section:flat_system} is devoted to their proofs.
    \item 
        The submitted preprint \cite{bleitnerCarlsonNobili2023Large} analyzes the thermally non-diffusive system. Chapter \ref{chapter:thermally_nd_system} is dedicated to these findings and results from this article are stated in Theorems~\ref{theorem_nd_regularity} and \ref{theorem_nd_convergence}.
    \item
        In \cite{bleitner2023Lower} lower bounds on mixing norms for the ad\-vection-hyperdiffusion equation are proven. Due to their disparity to buoyancy-driven fluids, the results are not discussed in this dissertation.
\end{itemize}

{\let\cleardoublepage\relax\chapter*{Declaration of Contribution}}
\begin{itemize}
    \item
        Lemmas \ref{lemma_grad_ids} and \ref{lemma_elliptic_regularity_periodic_domain} are generalizations of findings in \cite{bleitnerNobili2024Bounds} and Theorems~\ref{theorem_nonlinear_1_2_bound} and \ref{theorem_nonlinear_main_theorem}, as well as the proofs provided in Section \ref{section:nonlinearity_paper} correspond to the results given in the publication. Both authors contributed equally to the paper. In particular, the candidate contributed to rigorous proofs and the analytical framework.
    \item
        Lemmas \ref{lemma_elliptic_regularity_bounded_domain}, \ref{lemma_int_Delta_u_v}, \ref{lemma_coercivity}, \ref{lemma_u_cdot_nabla_p_on_boundary} as well as Theorems~\ref{theorem_nd_regularity} and \ref{theorem_nd_convergence}, together with the their proofs in Chapter \ref{chapter:thermally_nd_system} are published in \cite{bleitnerCarlsonNobili2023Large}. The candidate provided rigorous proofs and technical results. In particular the crucial Lemma \ref{lemma_u_cdot_nabla_p_on_boundary} is due to the candidate.
    \item 
        Theorem~\ref{theorem_flat}, together with the proofs in Section \ref{section:flat_system} are slight variations of the findings in \cite{bleitner2024Scaling}. The authors contributed equally.
    \item
        The statements given in Theorems~\ref{theorem_1_2_bound} and \ref{theorem_main_theorem_rb_curved}, proven in Sections \ref{Section:nusselt_number_representations} and \ref{section:general_system}, are new results due to the candidate. However, the majority of the approaches are combinations of arguments in \cite{bleitnerNobili2024Bounds}, \cite{bleitnerCarlsonNobili2023Large} and \cite{bleitner2024Scaling}.
    \item
        All graphics, including the simulations, were created by the candidate. The Tables \ref{Table:flat_result_Ls_depends_on_Pr_Ra} and \ref{Table:flat_result_Pr_depends_on_Ls_Ra} are a joint effort by Camilla Nobili and the candidate.
\end{itemize}

{\let\cleardoublepage\relax \chapter*{Declaration of Authorship}}

\noindent
Hiermit versichere ich an Eides statt, dass ich die vorliegende Dissertation selbst verfasst
und keine anderen als die angegebenen Quellen und Hilfsmittel benutzt habe.
\vspace{12ex}\\

\chapter*{Acknowledgement}

I am extremely grateful to Camilla Nobili, who not only encouraged and enabled me to pursue this academic journey but also guided me with exceptional supervision. Her immense support and advice on mathematics, academia, and beyond are invaluable to me.

I would like to thank the many helpful colleagues and friends in Hamburg for the discussions, feedback, and time spent inside and outside the University.

I am deeply thankful to my parents who always supported me and gave me the opportunity to dive into mathematics.

Finally a special thank you to my friends at home for countless moments of joy and personal support.

\newpage
\mainmatter
\pagestyle{myPageStyle}
\chapter{Introduction}
\vspace{-\baselineskip}
\section{Motivation}

Buoyancy describes the force arising from density variations in a gravitational field, usually due to differences in heat, salinity, or composition. Typical examples range from everyday life, such as heating water in a pot, to engineering challenges like designing a cooling system, and geophysical phenomena such as atmospheric convection, oceanic, or earth mantle currents to large scales as in solar convection layers.

In equilibrium, the net energy that is induced by the forcing has to be balanced by the dissipation in the fluid and on the boundary. It is apparent that the latter plays a significant role. In fact, imagine a horizontal periodic channel where the fluid can freely move along the boundary. If there is no forcing, a constant horizontal flow is expected to move indefinitely. In contrast, if the fluid exhibits friction on the walls then this flow will be slowed down and kinematic viscosity induces a decay in energy.

Naturally, also the geometry of the system plays a crucial role in the dynamics of the fluid. Heat might get trapped in a pocket, and walls constrain the flow to follow a certain path.

Due to the variety of applications, it is of immense interest to understand how these system parameters change the dynamics of the flow. While small-scale experiments and simulations might provide insight, large-scale problems and in particular extreme conditions require mathematical theory to answer these questions.

In this thesis, we want to investigate two systems. The first one is Rayleigh-Bénard convection, where a fluid is trapped in-between a heated bottom and a cooled top plate. The main focus lies on rigorously deriving bounds for the vertical heat transport, measured by the Nusselt number, with respect to the strength of the buoyancy forcing, given by the Rayleigh number. In particular, we allow non-flat boundaries and capture the influence of the geometry and the friction at the walls. The second problem is concerned with a thermally non-diffusive system, which can be interpreted as a specific limit of the former problem. Here the objective is to capture the dynamics by deriving regularity estimates and showing convergence to the hydrostatic equilibrium, where the buoyancy force is balanced by the pressure.

In what follows we will always work in a two-dimensional domain $\Omega$.

\section{Rayleigh-Bénard Convection}
\label{section:introduction_rb}

In the case of Rayleigh-Bénard convection the top and bottom boundary are given by $\gammap$ and $\gammam$, which are additionally restricted to be functions of the horizontal variable. In the horizontal direction, we assume periodic boundary conditions. The domain and vertical boundaries are given by
\begin{align}
    \Omega &= \left\lbrace (x_1,x_2)\in \mathbb{R}^2 \ \middle \vert\ 0\leq x_1\leq \Gamma, \hm(x_1)\leq x_2\leq \hp(x_1) \right\rbrace
    \label{intro_domain_rb}
    \\
    \gammap &= \left\lbrace (x_1,x_2)\in \mathbb{R}^2 \ \middle \vert\ 0\leq x_1\leq \Gamma, x_2= \hp(x_1) \right\rbrace
    \label{gfbkjbgllf}
    \\
    \gammam &= \left\lbrace (x_1,x_2)\in \mathbb{R}^2 \ \middle \vert\ 0\leq x_1\leq \Gamma, x_2= \hm(x_1) \right\rbrace.
    \label{knkbnjlkbvnj}
\end{align}\noeqref{gfbkjbgllf}\noeqref{knkbnjlkbvnj}%
By that, we also introduce the domain width $\Gamma>0$. Additionally, we will assume that $\hp(x_1)>\hm(x_1)$ so that the top and bottom boundary are separated from each other and that the domain size is given by $|\Omega|=\Gamma$. Hence, the average height is set to $1$. In the classical Rayleigh-Bénard problem, the boundaries are flat and the non-dimensionalization leads to a domain height of $1$, so the domain considered here is a generalization thereof.

In Rayleigh-Bénard convection the buoyancy force is a consequence of temperature differences. Specifically, there is a temperature gap between the hotter lower and the colder upper boundary. Accordingly, the hot fluid near the bottom expands and becomes less dense than the cold fluid at the top. Due to gravity, the hot fluid experiences an upward force, resulting in a dynamical system.

The main feature of such fluids is the forcing due to density variations, which are therefore described by the compressible Navier-Stokes equations. Under the assumptions that the density of the fluid varies linearly with the temperature and that density variations except for the forcing due to buoyancy can be neglected, the Boussinesq approximation leads to the following incompressible system (\cite{goluskin2015Internally}). After non-dimensionalizing the velocity $u=(u_1,u_2)(x_1,x_2,t)$, the scalar pressure and temperature fields $p=p(x_1,x_2,t)$ and $\theta=\theta(x_1,x_2,t)$, satisfy
\begin{align}
    \Pr^{-1}(u_t + u\cdot\nabla u) + \nabla p - \Delta u &= \Ra \theta e_2
    \label{intro_rb_navier_stokes}
    \\
    \nabla \cdot u &= 0
    \label{intro_rb_incompressible}
    \\
    \theta_t + u\cdot\nabla \theta - \Delta \theta &= 0,%
    \label{intro_rb_advection_diffusion}%
\end{align}\noeqref{intro_rb_incompressible}%
where $e_2 = (0,1)$ is the unit vector in upward direction. Here, we also introduced the Prandtl number $\Pr$ and the Rayleigh number $\Ra$, defined by
\begin{align}
    \Pr &= \frac{\nu}{\varkappa}
    \label{renjterjt}
    \\
    \Ra &= \frac{g \varsigma d^3(T^--T^+)}{\varkappa \nu},
    \label{definition_Ra}
\end{align}\noeqref{renjterjt}%
where $\nu$ is the kinematic viscosity, $\varkappa$ is the thermal diffusivity, $g$ is the gravitational constant, $\varsigma$ is the thermal expansion coefficient, $d$ is the height gap of the boundaries, and $T^-$ and $T^+$ is the temperature on the respective boundaries in the original system. This non-dimensionalization results in a temperature gap of $1$ at the boundaries, i.e.
\begin{alignat}{2}
    \theta &= 0 \qquad &\textnormal{ on }&\gammap
    \\
    \theta &= 1 \qquad &\textnormal{ on }&\gammam.
\end{alignat}

The Rayleigh number describes the strength of the buoyancy forcing and depends on the underlying setup at hand due to its dependency on $d$, $T^-$, and $T^+$, while the Prandtl number is an intrinsic property of the fluid. Typical values for these parameters are given in Table \ref{table:pr_and_ra_numbers}.

\begin{table}
    \hrule
    \vspace{\baselineskip}
    \begin{minipage}{\textwidth}
        \begin{center}
            \begin{tabular}{lc | l c}
                Fluid
                \phantom{\footnote{\label{citeSchumacher2020Colloquium}\cite{schumacher2020Colloquium}}\footnote{\label{citeNiemela2000Turbulent}\cite{niemela2000Turbulent}}\footnote{\label{citeTropea2007Springer}\cite{tropea2007Springer}}\footnote{\label{citeSchubert2001Mantle}\cite{schubert2001Mantle}}\footnote{\label{citeWolstencroft2009Nusselt}\cite{wolstencroft2009Nusselt}}} & $\Pr$ & System &$\Ra$
                \\
                \hline
                Solar Convection Plasma\footref{citeSchumacher2020Colloquium} & $10^{-6}$ & Solar Convection Zone\footref{citeSchumacher2020Colloquium}\footref{citeNiemela2000Turbulent} & $10^{20}-10^{23}$
                \\
                Air\footref{citeTropea2007Springer} & $0.71$ & Experiments\footref{citeNiemela2000Turbulent} & up to $10^{17}$
                \\
                Water\footref{citeTropea2007Springer} & $7.0$ & Ocean\footref{citeNiemela2000Turbulent} & $10^{20}$
                \\
                Earth Mantle Rock\footref{citeSchubert2001Mantle}
                & $10^{23}$ & Earth Mantle\footref{citeSchubert2001Mantle}\footref{citeWolstencroft2009Nusselt}
                & $10^{7}-10^{9}$
            \end{tabular}
        \end{center}
    \end{minipage}
    \vspace{-\baselineskip}
    \caption{Typical values for the Prandtl and Rayleigh number in selected systems.}
    \label{table:pr_and_ra_numbers}
    \vspace{\baselineskip}
    \hrule
\end{table}

\begin{figure}
    \hrule
    \definecolor{color_hot}{rgb}{0.623529, 0.219608, 0.211765}
    \begin{center}
        \includegraphics[width=\textwidth]{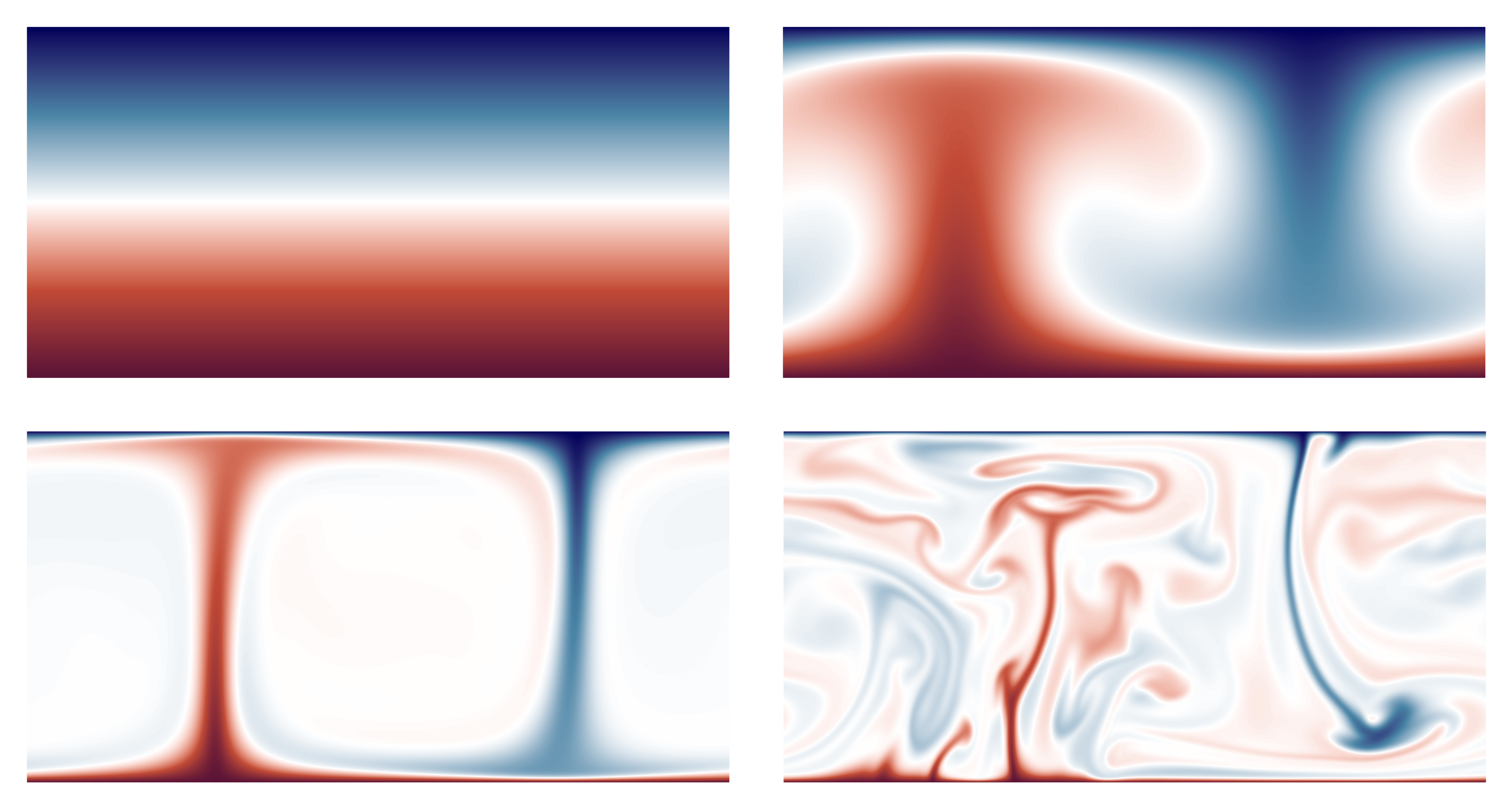}
        \vspace*{-9mm}
        \caption{Snapshots of finite element simulations of the temperature field for Rayleigh-Bénard convection with $\Pr=1$, $\alpha = 1$, $\kappa=0$, and $\Gamma=2$. \textcolor{color_hot}{Red} corresponds to hot regions, cold regions are colored in \textcolor{rtg_darkblue}{blue} and the Rayleigh numbers are $10^2$ (top left), $10^4$ (top right), $10^6$ (bottom left) and $10^8$ (bottom right). The simulations were created with Firedrake (\cite{firedrake}) and visualized in Paraview (\cite{paraview}).}
        \label{fig:simulation_rb}
    \end{center}
    \hrule
\end{figure}

The third non-dimensional number, the Nusselt number $\Nu$, is of particular interest in the realm of Rayleigh-Bénard convection. It measures the excess of upwards heat transport over the purely conducting state and we will define it later in \eqref{nu_definition}. For increasing Rayleigh numbers, the fluid is expected to become more turbulent increasing the energy transfer from the bottom plate to the top one. This is illustrated in Figure \ref{fig:simulation_rb}, where for $\Ra=10^2$ the solution is purely conductive, and for the subsequent increasing Rayleigh numbers it becomes more and more advection dominated. Many works, experimental, numerical, and of theoretic nature, are dedicated to showing the scaling of this number with respect to $\Pr$ and $\Ra$, see \cite{plumley2019Scaling} for a review. While for moderate Rayleigh numbers experiments and simulations might solve this question, extreme cases demand rigorous mathematical bounds. In fact, Table \ref{table:pr_and_ra_numbers} indicates that the solar convection zone exceeds the capabilities of experiments. With higher temperature differences and larger height gaps in other stars \eqref{definition_Ra} implies more extreme $\Ra$ values, demanding theoretical results. In Section \ref{section:rb_results} we will discuss some of the results and explain the findings of this thesis in the context of these works.

Here, we only want to mention some results that show the influence of the boundary conditions and Prandtl number for the fluid in this limit. The top and bottom boundaries are expected to be solid, implying the no-penetration boundary condition in the normal direction, i.e.
\begin{align}
    n\cdot u=0,
\end{align}
where $n$ is the unit normal vector. In the tangential direction, the situation is less clear. While usually the equations are equipped with no-slip boundary conditions, where the fluid sticks to the wall, there has been a long debate including Bernoulli, Couette, Coulomb, Helmholtz, Navier, Poisson, Stokes (\cite{priezjev2006Influence}) whether or not slip occurs. Physical experiments show that depending on the materials slip is expected (\cite{uthe2022optical,neto2005Boundary,tropea2007Springer}). Additionally, theoretical results (\cite{miksis1994Slip,bolanos2017Derivation}) show that imperfections and roughness on the boundary for fluids subject to no-slip boundary conditions lead to slip. Here we consider the Navier-slip boundary conditions, first proposed in \cite{navier1823Memoire} 
and given by
\begin{align}
    \tau \cdot (\D u\ n+\alpha u) =0,
    \label{intro_navier_slip}
\end{align}
where $\tau$ is the unit tangent vector, $\D u = \frac{1}{2}(\nabla u + \nabla u^T)$ is the symmetric gradient and $\alpha=\alpha(x)>0$ is the slip coefficient. Note that in the $\alpha\to \infty$ limit these conditions resemble the no-slip case. 
In fact \cite{amrouche2020Semigroup,kelliher2006Navier} show that also solutions converge to those with no-slip boundary conditions as $\alpha \to \infty$. On the other hand, setting $\alpha = 0$ yields free-slip boundary conditions, implying that Navier-slip boundary conditions interpolate between the two extreme cases.

To see how the boundary conditions influence the scaling laws we want to discuss some results. For no-slip, flat boundaries \cite{doering1996VariationalConvection} showed $\Nu\lesssim \Ra^\frac{1}{2}$ uniform in $\Pr$, even in three spatial dimensions. In contrast, in the two-dimensional free-slip setting \cite{whitehead2011Ultimate} proved $\Nu \lesssim \Ra^\frac{5}{12}$, again uniform in $\Pr$. For the flat infinite Prandtl setup with no-slip boundary conditions \cite{constantin1999Infinite} proved $\Nu\lesssim \Ra^\frac{1}{3}(1+\ln \Ra)^\frac{2}{3}$, for which the logarithmic exponent has been improved (\cite{doering2006Bounds,otto2011Rayleigh}) since then. The gap between the free- and no-slip results was studied in \cite{drivas2022Bounds}, where the authors showed $\Nu \lesssim \alpha^2 \Ra^\frac{1}{2} + \Ra^\frac{5}{12}$ for flat boundaries with constant slip coefficient in the high $\Pr$ regime.

For rough boundaries, results are more limited. \cite{goluskin2016Bounds} proved $\Nu\lesssim \Ra^\frac{1}{2}$ uniformly in $\Pr$ for no-slip boundary conditions when the profile functions $\hp,\hm$ are in $H^1$.

These findings inspired the study of the problem with full Navier-slip boundary conditions on curved domains. Nobili and the author of this thesis proved (\cite{bleitnerNobili2024Bounds}) 
\begin{align}
    \Nu\lesssim \|\alpha+\kappa\|_{W^{1,\infty}}^2 \Ra^\frac{1}{2} + \Ra^\frac{5}{12},    
\end{align}
where $\kappa$ is the curvature of the boundary, generalizing the result of \cite{drivas2022Bounds}, and $\Nu\lesssim_{\alpha,\kappa} \Ra^\frac{3}{7}$, where the implicit constant hides a complex dependency on $\alpha$ and $\kappa$. Both results only hold in the case of sufficiently big $\Pr$ and small $\alpha,\kappa\in W^{1,\infty}$. The proofs of these results will be given in Section \ref{section:nonlinearity_paper}.

A refined pressure estimate yields
\begin{align}
    \Nu \lessapprox \Pr^{-\frac{1}{6}} \Ra^\frac{1}{2} + \Ra^\frac{5}{12}.
\end{align}
The exact statement is given in Theorem \ref{theorem_main_theorem_rb_curved} and proven in Chapter \ref{chapter:rayleigh-benard-convection}. In the case of flat boundaries, constant slip coefficient, and sufficiently large $\Ra$, the bound is given by
\begin{align}
    \Nu
    &\lesssim \Ra^\frac{5}{12} + \alpha^{\frac{1}{12}}\Pr^{-\frac{1}{6}}\Ra^\frac{1}{2}
\end{align}
for $\alpha\leq 1$ and
\begin{align}
    \Nu
    &\leq  \alpha^{\frac{1}{3}} \Ra^\frac{1}{3} + \alpha^{\frac{2}{13}} \Ra^{\frac{5}{13}} +  \Ra^\frac{5}{12} + \alpha^{\frac{1}{2}}\Pr^{-\frac{1}{6}}\Ra^\frac{1}{2}
\end{align}
for $\alpha\geq 1$, similar to the results of \cite{bleitner2024Scaling} that only vary in the $\alpha$ exponent in the $\Pr$ terms. These bounds significantly improve the previous findings. Apart from their improved estimates, these bounds hold in a much broader range of physically relevant settings. They allow any Prandtl number, even showing a crossover at $\Pr=\Ra^\frac{1}{2}$. Additionally, they hold for any, sufficiently smooth, slip coefficient, which in particular allows close to no-slip setups, that might seem physical more realistic. However in Section \ref{section:scaling_kappa_alpha_with_Ra} we will provide an argument showing that $\alpha$ and $\kappa$ can scale with respect to $\Ra$. In Section \ref{section:rb_results} we will discuss the Nusselt number scaling in more detail.

Finally, we want to discuss the different approaches that lead to these bounds. The Constantin and Doering background field method (\cite{doering1994VariationalShear,doering1996VariationalConvection}) led to numerous of the previously mentioned results. The main strategy here is to decompose the temperature field into a steady profile, approximating the expected long-time boundary layer and bulk behavior of the fluid, and fluctuations around it. This method is illustrated in Section \ref{section:nonlinearity_paper}. Contrary in Sections \ref{section:general_system} and \ref{section:flat_system} the direct method (\cite{seis2015Scaling}) is employed, which solely relies on a localization principle of the Nusselt number. 

\section{Thermally Non-Diffusive System}
\label{Section:intro_nd_system}

Further, we want to investigate a closely related system without thermal diffusion, given by
\begin{align}
    u_t + u\cdot \nabla u + \nabla p -\Delta u &= \theta e_2
    \label{intro_nd_navier_stokes}
    \\
    \nabla \cdot u &= 0
    \label{intro_nd_incompressiblitiy}
    \\
    \theta_t + u\cdot \nabla \theta &= 0.
    \label{intro_nd_transport}
\end{align}\noeqref{intro_nd_incompressiblitiy}%
This set of equations is an immediate consequence of \eqref{intro_rb_navier_stokes}-  \eqref{intro_rb_advection_diffusion} when setting $\Ra=\Pr=1$ and disregarding the diffusive term in the advection-diffusion equation. However, it can also be seen as a limit of these equations. In fact consider solutions $\tilde u(x,\tilde t),\tilde p(x,\tilde t), \tilde \theta(x,\tilde t)$ of
\begin{align}
    \Pr^{-1}(\tilde u_{\tilde t} + \tilde u \cdot \nabla \tilde u) + \nabla \tilde p - \Delta \tilde u&= \Ra\tilde\theta e_2
    \\
    \nabla \cdot \tilde u &= 0
    \\
    \nabla \tilde\theta_{\tilde t} + \tilde u \cdot \nabla \tilde \theta - \Delta \tilde \theta &= 0
\end{align}
and rescale them according to
\begin{align}
    u(x,t) &=\Pr^{-\frac{1}{2}} \Ra^{-\frac{1}{2}}\tilde u(x,\tilde t)
    &
    p(x,t) &= \Ra^{-1} \tilde p(x,\tilde t)
    \\
    \theta(x,t) &= \tilde \theta(x,\tilde t)
    &
    t &= \Pr^{\frac{1}{2}} \Ra^{\frac{1}{2}}\tilde t.
\end{align}
Then $u,p,\theta$ solve
\begin{align}
    u_t + u\cdot\nabla u + \nabla p - \Pr^\frac{1}{2}\Ra^{-\frac{1}{2}} \Delta u  &= \theta e_2
    \\
    \nabla \cdot u &= 0
    \\
    \theta_t + u\cdot\nabla \theta - \Pr^{-\frac{1}{2}}\Ra^{-\frac{1}{2}}\Delta \theta &= 0,
\end{align}
which are the equations, rescaled according to free-fall time (\cite{schneide2018Probing}). The previous discussion showed that the limit $\Pr,\Ra\to \infty$ is of particular interest for this problem. Setting $Pr = \nu^2 \Ra$ for some $\nu>0$ and taking the limit $\Ra\to \infty$ yields \eqref{intro_nd_navier_stokes}-\eqref{intro_nd_transport} with an additional viscosity parameter $\nu$, the system studied in \cite{bleitnerCarlsonNobili2023Large}. Here we set $\nu=1$ in order to simplify the notation but remark that all results hold for any $\nu > 0$.

Note that due to the absence of thermal diffusion, the governing equations are only equipped with the velocity boundary conditions
\begin{align}
    n\cdot u&=0
    \\
    \tau\cdot (\D u \ n +\alpha u) &= 0
\end{align}
and instead of a horizontally periodic strip, we assume a bounded Lipschitz domain $\Omega\subset \mathbb{R}^2$. Note that \cite{hu2018Partially} showed global well-posedness for the problem.

\begin{figure}
    \hrule
    \definecolor{color_hot}{rgb}{0.623529, 0.219608, 0.211765}
    \vspace*{1.27mm}
    \begin{center}
        \includegraphics[width=\textwidth]{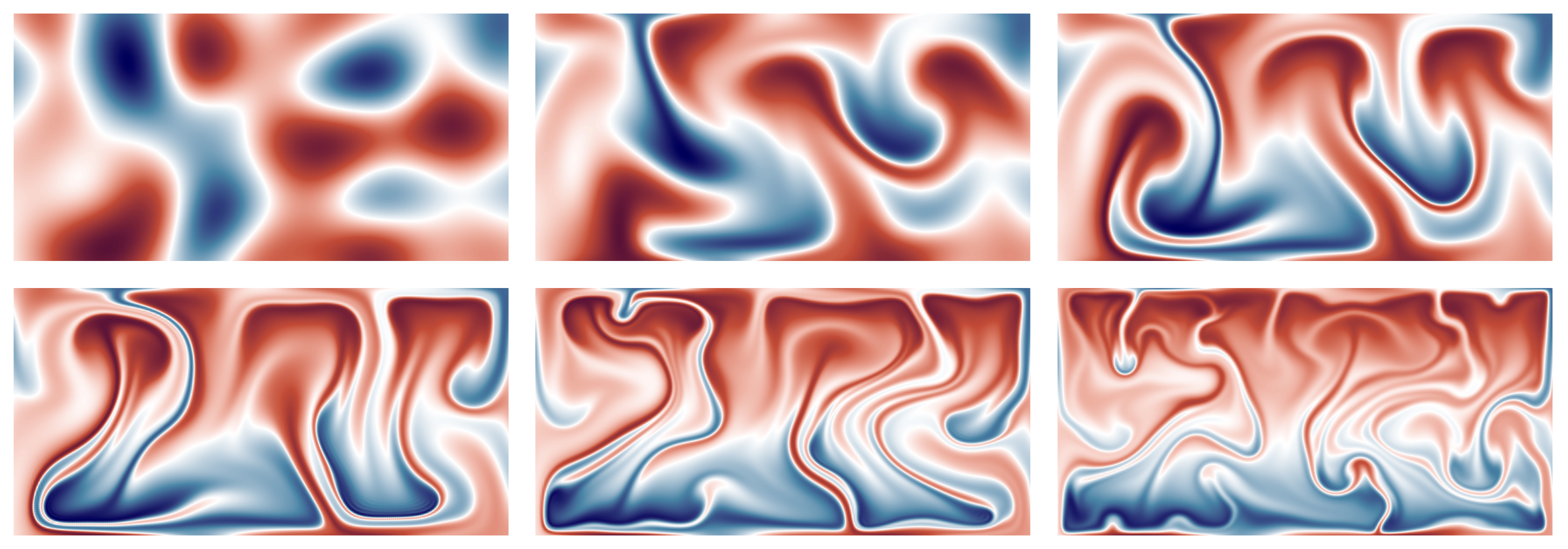}
        \vspace*{-8mm}
        \caption{Snapshots of a simulation of the thermally non-diffusive system with $\alpha = 10^4$, where hot regions are colored in \textcolor{color_hot}{red}, and cold regions in \textcolor{rtg_darkblue}{blue}. The individual pictures correspond to times $0$, $12$, and $25$ in the top row and $50$, $100$, and $200$ in the bottom row. The simulation was created with Firedrake (\cite{firedrake}) and visualized in Paraview (\cite{paraview}).}
        \label{fig:simulation_nd}
    \end{center}
    \hrule
\end{figure}

As a consequence of the absence of thermal boundary conditions, no mechanism in this system provides an influx of energy. Therefore, over time, the velocity is expected to decay due to the viscosity term, resulting in a steady state
\begin{align}
    u &= 0
    \\
    \nabla p &= \theta e_2.
\end{align}
The state, where the velocity vanishes and the buoyancy force is balanced by the pressure gradient is called the hydrostatic equilibrium.

Figure \ref{fig:simulation_nd} shows this behavior. Over time the temperature field becomes vertically stratified and the velocity field decays. Note though that the rectangular domain of the simulation does not satisfy the regularity assumptions of Theorem \ref{theorem_nd_convergence}.

\cite{doering2018longTime} studied the system with stress-free boundary conditions
\begin{align}
    n\cdot u &= 0
    \\
    \omega &= 0,
\end{align}
where $\omega = -\partial_2 u_1 + \partial_1 u_2$ is the vorticity. After proving global well-posedness with regularity estimates
\begin{align}
    u&\in L^\infty\left((0,T);H^3(\Omega)\right) \cap L^2\left((0,T);H^4(\Omega)\right)
    \\
    \theta&\in L^\infty\left((0,T);H^3(\Omega)\right)
\end{align}
for any $T>0$, provided that the boundary and initial data are sufficiently smooth, they showed
\begin{align}
    \|u(t)\|_{H^2}&\leq C
    \\
    \int_0^t \|u(s)\|_{H^1}^2 ds &\leq C,
\end{align}
for constants $C>0$ independent of time. With a slight abuse of notation, we write the bounds as
\begin{align}
    u&\in L^\infty\left((0,\infty);H^2(\Omega)\right) \cap L^2\left((0,\infty);H^1(\Omega)\right).
\end{align}
These uniform in time estimates imply that the velocity decays and solutions converge to the hydrostatic equilibrium. In particular, their analysis shows
\begin{align}
    \|u(t)\|_{H^1} &\to 0
    \\
    \|(\nabla p - \theta e_2)(t)\|_{H^{-1}} &\to 0
\end{align}
for $t\to \infty$. Additionally, \cite{doering2018longTime} studied the linear stability of the system.

Inspired by these findings Elizabeth Carlson, Camilla Nobili and the author of this thesis (\cite{bleitnerCarlsonNobili2023Large}) proved the regularity estimates
\begin{align}
    u&\in L^2\left((0,T);H^3(\Omega)\right) \cap L^\frac{2p}{p-2}\left((0,T);W^{1,p}(\Omega)\right)
    \\
    \theta &\in L^\infty\left((0,T);W^{1,q}(\Omega)\right)
    \\
    u&\in L^\infty\left((0,\infty);H^2(\Omega)\right) \cap L^p\left((0,\infty);W^{1,p}(\Omega)\right)    
\end{align}
for any $T>0$, $2\leq p < \infty$ and $1\leq q<\infty$, provided that the domain and initial data are sufficiently smooth. Again the uniform in time bounds result in convergence to the hydrostatic equilibrium in $H^1 \times H^{-1}$. The explicit statement of the results, the proof thereof, as well as a further discussion, will be given in Chapter \ref{chapter:thermally_nd_system}.

\section{Well-Posedness}

In what follows we will focus on a-priori estimates. The analysis will provide sufficient regularity estimates for the solutions, which subsequently allow Galerkin approximation techniques to prove the existence of solutions. For general flows with Navier-slip boundary conditions, the reader is referred to \cite{clopeau1998On,kelliher2006Navier}. The well-posedness of the Rayleigh-Bénard problem with no-slip boundary conditions is covered in \cite{foias1987Attractors}, which also includes remarks for other boundary conditions. In particular, a combination of their techniques with the functional analysis properties of Navier-slip boundary conditions discussed in \cite{amrouche2020Semigroup} yields the desired result. Similarly, the well-posedness of the thermally non-diffusive system with no-slip boundary conditions is studied in \cite{hu2013On}, and again combining their approach with the results of \cite{amrouche2020Semigroup} leads to existence results. The uniqueness of solutions to the thermally non-diffusive system with Navier-slip boundary conditions is studied in \cite{hu2018Partially}.

We remark that although the approximation techniques only yield the existence of solutions on bounded time intervals we will use $L^p\left((0,\infty);W^{k,q}(\Omega)\right)$ to indicate that the corresponding bounds hold uniformly in time.

\section{Structure of the Thesis}%

In Chapter \ref{chapter:navier_slip} general results and technical properties of fluids subject to Navier-slip boundary conditions are proven. In particular, estimates regarding gradients on the boundary are derived in Section \ref{section_navSlip_boundary_gradients}.

Chapter \ref{chapter:rayleigh-benard-convection} is dedicated to Rayleigh-Bénard convection, where in Section \ref{section:rb_results} the main results are provided, that are proven in Section \ref{section:general_system}, Section \ref{section:flat_system} and Section \ref{section:nonlinearity_paper} for the general system, the system with flat boundaries and constant slip coefficient, and the system with the same boundary profile function, respectively.

Chapter \ref{chapter:thermally_nd_system} is devoted to the non-diffusive equations, where the main findings regarding regularity and convergence are given in Section \ref{section:nd_main_results}, which are proven in Section \ref{section:nd_regularity_estimates} and Section \ref{section:nd_convergence}, respectively.

\section{Notation}
Throughout the thesis, the following notation is used.
\vspace{0.5\baselineskip}
\\
\begin{tabular}[h]{p{0.17\linewidth}p{0.75\linewidth}}
    $\Omega$ & $\Omega$ is a two-dimensional at least Lipschitz domain, either bounded as throughout Chapter \ref{chapter:thermally_nd_system} or a potentially curved periodic channel as specified in \eqref{intro_domain_rb} with top boundary $\gammap$, described by the height function $\hp$ and bottom boundary $\gammam$ described by the height function $\hm$.
    \\[5pt]
    $u$ & The fluid's velocity and its components are denoted by $u=u(x,t)=(u_1,u_2)(x_1,x_2,t)$.
    \\[5pt]
    $a^\perp$ & For a vector $a=(a_1,a_2)$ we use the convention of defining the perpendicular direction as $a^\perp = (-a_2,a_1)$.
    \\[5pt]
    $n,\npm$ & The unit normal vector in outward pointing direction. In case no outward direction is specified, $\np$ is the upward-pointing unit vector, while $\nm$ is the downward-pointing one.
    \\[5pt]
    $\tau$ & The unit tangent vector is defined by $\tau = n^\perp$.
    \\[5pt]
    $u_\tau$ & The tangential velocity $u_\tau = u\cdot \tau$.
    \\[5pt]
    $\omega$ & The vorticity is defined as $\omega = \nabla^\perp\cdot u = -\partial_2 u_1+\partial_1 u_2$. 
    \\[5pt]
    $\D u$ & The symmetric gradient $\D u = \frac{1}{2}(\nabla u+\nabla u^T)$, which in components is given by $(\D u)_{ij} = \frac{1}{2} (\partial_i u_j + \partial_j u_i)$
    \\[5pt]
    $\alpha$ & The slip coefficient $\alpha$, which can vary in space.
    \\[5pt]
    $\kappa$ & The boundary curvature is defined by $\kappa = n\cdot (\tau\cdot \nabla)\tau$.
    \\[5pt]
    $\|\cdot\|_p$, $\|\cdot \|_{W^{k,p}}$ & These norms denote the spatial $L^p$, respectively $W^{k,p}$-norms defined by $\|f\|_p^p = \int_\Omega |f|^p$, respectively $\|f\|_{W^{k,p}}^p=\sum_{|\alpha|\leq k} \|\nabla^\alpha f\|_p^p$.
    \\[5pt]
    $\langle\ \cdot\ \rangle$ & For any domain $X$ the large time and spatial average is defined by $\langle f\rangle_X = \limsup_{T\to\infty}\frac{1}{T} \int_0^T \frac{1}{|\Omega|}\int_X f(x,t) \ dx\ dt$. Note that the convention is used to always divide by $|\Omega|$ instead of $|X|$ and we omit the index if $X=\Omega$, i.e. $\langle \ \cdot\ \rangle = \langle \ \cdot\ \rangle_\Omega$.
    \\[5pt]
    $\lesssim$, $\lessapprox$ & We write $f\lesssim g$ if there exists a constant $C>0$, potentially depending on $|\Omega|$, the Lipschitz constant of $\Omega$ and the Lebesgue norm parameter such that $f\leq c g$. Similarly $f\lessapprox g$, where the constant potentially also depends on $\alpha$, $\kappa$ and $\Omega$ in general.
\end{tabular}
\leavevmode\\
\vspace{0.5\baselineskip}

Additionally, the Einstein summation convention is used, implying summation over identical indices appearing twice, i.e. $u_i v_i = \sum_{i} u_iv_i $.

In the partially periodic domain, the cancellation of boundary terms is used without explicitly mentioning it. With a slight abuse of notation in Chapter \ref{chapter:navier_slip} the relevant boundaries are always referred to as $\partial\Omega$, while in the case of the Rayleigh-Bénard domain only the top and bottom boundaries $\gammam\cup\gammap$ are implied.

\newpage
\chapter{Navier-Slip and Curved Domains}
\label{chapter:navier_slip}

In this chapter, we illustrate how boundary roughness, measured by the curvature, influences fluids with slip boundary conditions. As we will see the curvature and slip coefficient are closely related to each other. Additionally, the estimates derived here will later be used to prove the corresponding results.

While most of this chapter's results can be adapted to other two-dimensional domains we will always assume that the domain is either the partially periodic one given in \eqref{intro_domain_rb} or a bounded Lipschitz domain and only specify the needed regularity of the boundary in the statements. A Lipschitz domain is of class $C^{k,1}$ if locally there exist bijections to the half plane, whose derivatives up to order $k$ are Lipschitz continuous. For the partially periodic domain \eqref{intro_domain_rb} these bijections are given by $\hp$ and $\hm$.

The curvature is defined as
\begin{align}
    \kappa = n \cdot (\tau \cdot \nabla) \tau
    \label{definition_curvature}
\end{align}
and if $\Omega$ is a $C^{k,1}$ domain then the curvature satisfies $\kappa\in W^{k-1,\infty}(\partial\Omega)$.

While for the fluid we have the solution of either Rayleigh-Bénard convection or the thermally non-diffusive system in mind, in this chapter the fluid only needs to satisfy
\begin{alignat}{2}
    \nabla \cdot u &= 0 \qquad & \textnormal{ in } &\Omega
    \label{nav_slip_incompressible}
    \\
    n \cdot u &= 0 \qquad & \textnormal{ on } &\partial\Omega
    \label{nav_slip_no_penetration}
    \\
    \tau\cdot(\D u \ n+\alpha u) &= 0 \qquad & \textnormal{ on } &\partial\Omega.
    \label{nav_slip_nav_slip_bc}
\end{alignat}

\section{Gradients and Vorticity}

This section is devoted to finding estimates that allow us to exchange the gradient $\nabla u$, the symmetric gradient $\frac{1}{2}(\nabla u + \nabla u^T)$ and the vorticity $\omega = \nabla^\perp \cdot u = -\partial_2 u_1 + \partial_1 u_2$.

The first estimate in that direction yields a result for the $L^2$-norm of those quantities, in exchange for a boundary integral. Note that the boundary terms vanish and the bulk integrals coincide for flat domains.

\begin{lemma}
    Let $\Omega$ be $C^{1,1}$, $u\in H^1(\Omega)$ satisfy \eqref{nav_slip_incompressible}-\eqref{nav_slip_nav_slip_bc} and $v\in H^1(\Omega)$ satisfy $n\cdot v = 0$. Then
    \begin{align}
        2 \int_{\Omega} (\D u)_{ij} (\D v)_{ij} - \int_{\partial\Omega} \kappa u\cdot v = \int_{\Omega} \partial_i u_j \partial_i v_j = \int_{\Omega} \nabla^\perp \cdot u \nabla^\perp \cdot v + \int_{\partial\Omega} \kappa u\cdot v.
        \\
        \label{grad_ids_eq}
    \end{align}
    and in particular
    \begin{align}
        2 \|\D u\|_2^2 - \int_{\partial\Omega} \kappa u_\tau^2 = \|\nabla u\|_2^2 = \|\omega\|_2^2 + \int_{\partial\Omega} \kappa u_\tau^2.
    \end{align}
    \label{lemma_grad_ids}
\end{lemma}

\begin{proof}
    Assume at first $u$ is smooth and note that
    \begin{align}
        2 \int_{\Omega} (\D u)_{ij} (\D v)_{ij} &= \int_{\Omega} \partial_i u_j \partial_i v_j + \int_{\Omega} \partial_i u_j \partial_j v_i.
        \label{adfajbg}
    \end{align}
    For the second term on the right-hand side of \eqref{adfajbg} integration by parts yields
    \begin{align}
        \int_{\Omega} \partial_i u_j \partial_j v_i = \int_{\partial\Omega} n \cdot (v \cdot \nabla) u  + \int_{\Omega} \partial_i \partial_j u_j v_i = \int_{\partial\Omega} n \cdot (v \cdot \nabla) u,
        \label{bnisbb}
    \end{align}
    where the second-order term vanished due to the incompressibility condition. In order to calculate the first term on the right-hand side of \eqref{bnisbb}, notice that $u = u_\tau \tau$ as $u\cdot n = 0$ on $\partial\Omega$ and therefore
    \begin{align}
        n\cdot (\tau\cdot \nabla) u = n\cdot (\tau\cdot \nabla) (u_\tau\tau) = u_\tau n\cdot (\tau\cdot \nabla) \tau + n\cdot \tau (\tau \cdot \nabla) u_\tau = \kappa u_\tau,
        \label{n_tau_grad_u}
    \end{align}
    where in the last identity we used \eqref{definition_curvature}. As analogously $v=v_\tau \tau$ one finds
    \begin{align}
        \int_{\partial\Omega} n \cdot (v \cdot \nabla) u = \int_{\partial\Omega} v_\tau n \cdot (\tau \cdot \nabla) u = \int_{\partial\Omega} \kappa u_\tau v_\tau.\label{nbgubingf}
    \end{align}
    The second term on the right-hand side of \eqref{bnisbb} vanishes by the divergence-free condition. Therefore, combining \eqref{adfajbg}, \eqref{bnisbb}, and \eqref{nbgubingf} yields
    \begin{align}
        2 \int_{\Omega} (\D u)_{ij} (\D v)_{ij} = \int_{\Omega} \partial_i u_j \partial_i v_j + \int_{\partial\Omega} \kappa u_\tau v_\tau,
        \label{grad_ids_first_id}
    \end{align}
    proving the first identity of \eqref{grad_ids_eq}.
    
    In order to show the second identity of \eqref{grad_ids_eq}, notice that
    \begin{align}
        2(\D u)_{ij} (\D v)_{ij} &+ (\nabla^\perp\cdot u)(\nabla^\perp\cdot v)
        \\
        &= \partial_i u_j \partial_i v_j + \partial_i u_j \partial_j v_i + (-\partial_2 u_1+\partial_1 u_2)(-\partial_2 v_1+\partial_1 v_2)
        \\
        &= \partial_i u_j \partial_i v_j + \partial_1 u_1 \partial_1 v_1 + \partial_1 u_2 \partial_2 v_1 + \partial_2 u_1 \partial_1 v_2 + \partial_2 u_2 \partial_2 v_2
        \\
        &\qquad + \partial_2 u_1 \partial_2 v_1 - \partial_2 u_1 \partial_1 v_2 - \partial_1 u_2 \partial_2 v_1 + \partial_1 u_2 \partial_1 v_2
        \\
        &= 2 \partial_i u_j \partial_i v_j
    \end{align}
    and therefore integrating and using \eqref{grad_ids_first_id} yields
    \begin{align}
        \int_{\Omega} (\nabla^\perp\cdot u)(\nabla^\perp\cdot v) &= 2 \int_{\Omega} \partial_i u_j \partial_i v_j - 2 \int_{\Omega} (\D u)_{ij} (\D v)_{ij}
        \\
        &= \int_{\Omega} \partial_i u_j \partial_i v_j - \int_{\partial\Omega} \kappa u_\tau v_\tau.
    \end{align}
    By approximation, the assumption of $u$ being smooth can be dropped.
\end{proof}

For second-order derivatives, one has the direct identity
\begin{align}
    \Delta u = \nabla^\perp \omega\label{id_Delta_u_is_nablaPerp_omega}
\end{align}
if $u$ satisfies \eqref{nav_slip_incompressible}. This follows from the immediate computation
\begin{align}
    \Delta u =
    \begin{pmatrix}
        \partial_1^2u_1+\partial_2^2u_1\\\partial_1^2u_2+\partial_2^2u_2
    \end{pmatrix}
    =
    \begin{pmatrix}
        \partial_2(-\partial_1 u_2+ \partial_2u_1)\\\partial_1(\partial_1 u_2- \partial_2 u_1)
    \end{pmatrix}
    = \nabla^\perp \omega.
\end{align}

The next estimate generalizes the previous two results to arbitrary Sobolev norms of order one and two. It is a slight deviation of Lemma 3.6 in \cite{bleitnerNobili2024Bounds}.
\begin{lemma}
    \label{lemma_elliptic_regularity_periodic_domain}
    Let $\Omega$ be the partially periodic domain defined in \eqref{intro_domain_rb}. Assume it is of class $C^{1,1}$ and let $u\in W^{1,p}(\Omega)$ for $2 \leq p < \infty$ satisfy \eqref{nav_slip_incompressible} and \eqref{nav_slip_no_penetration}. Then
    \begin{align}
        \label{elliptic_regularity_periodic_lemma_first_order}
        \|u\|_{W^{1,p}}\lesssim \|\omega\|_p + \left(1+\|\kappa\|_\infty^{1+\frac{2}{q}-\frac{2}{p}}\right)\|u\|_q,
    \end{align}
    where the implicit constant only depends on $p$, $|\Omega|$ and the Lipschitz constant of the boundary. If additionally $\Omega$ is $C^{2,1}$ and $u\in W^{2,q}(\Omega)$ for $1<q<\infty$, then
    \begin{align}
        \label{elliptic_regularity_periodic_lemma_second_order}
        \|u\|_{W^{2,p}} \lesssim \|\omega\|_{W^{1,p}} + \|\kappa\|_\infty \|\omega\|_p + (1+\|\kappa\|_{W^{1,\infty}}+\|\kappa\|_\infty^2) \|u\|_p.
    \end{align}
\end{lemma}

\begin{proof}
    For most of the proof, we follow the approach of Section 6.3.2, Theorem 4 in \cite{evans1998Partial} and provide it in detail here to capture the explicit dependency on $\kappa$. We first provide an overview of the strategy for the proof of \eqref{elliptic_regularity_periodic_lemma_first_order} and \eqref{elliptic_regularity_periodic_lemma_second_order} and prove these individual steps afterward.
    \begin{enumerate}
        \item First order estimate
            \begin{enumerate}[label=(\alph*)]
                \item \hyperlink{sfa_item_streamFunctionPde}{Derive the PDE for the stream function}
                \item \hyperlink{sfa_item_removeBdry}{Remove the boundary conditions by redefining the stream function}
                \item \hyperlink{sfa_item_straighten}{Establish a change in variables that straighten the boundaries}
                \item \hyperlink{sfa_item_straightenedPDE}{Derive the PDE in the straightened variables}
                \item \hyperlink{sfa_item_straightenedBounds}{Derive the estimates for the straightened system}
                \item \hyperlink{sfa_item_streamFunctionBounds}{Translate the bounds back to the original system}
                \item \hyperlink{sfa_item_generalizeLowerOrderTerm}{Generalize the lower order term}
            \end{enumerate}
        \item Second order estimate
            \begin{enumerate}[label=(\alph*)]
                \item \hyperlink{sfa_item_higherOrderPDE}{Derive a PDE for the horizontal derivative of the straightened system}
                \item \hyperlink{sfa_item_higherOrderHorizontalBound}{Derive bounds for the horizontal derivative}
                \item \hyperlink{sfa_item_higherOrderVerticalEquation}{Derive an equation for the remaining derivative}
                \item \hyperlink{sfa_item_higherOrderFullBound}{Get bounds for the full norm}
            \end{enumerate}
    \end{enumerate}
    Now we will prove the steps outlined before.
    
    \begin{enumerate}[label=(1\alph*)]
        \item \hypertarget{sfa_item_streamFunctionPde}{Derive the PDE for the stream function}
        
            Let $\phi$ be the stream function of $u$, i.e. $u=\nabla^\perp \phi$. Then
            \begin{align}
                \Delta \phi = \nabla^\perp\cdot\nabla^\perp \phi = \nabla^\perp \cdot u = \omega.
            \end{align}
            Additionally, the stream function is constant along the individual boundaries $\gammam$ and $\gammap$ as
            \begin{align}
                \tau\cdot\nabla \phi = n^\perp \cdot\nabla \phi = -n\cdot\nabla^\perp \phi = -n\cdot u = 0,
            \end{align}
            where we used the non-penetration boundary condition \eqref{no_penetration_bc}. As $\phi$ is only defined up to a constant we choose it such that $\phi\vert_{\gammam}=0$. Stokes theorem implies
            \begin{align}
                \dashint_{\Omega} u_1 = -\dashint_{\Omega} \partial_2 \phi = -\frac{1}{\Gamma} \int_{\gammap} n_2 \phi -\frac{1}{\Gamma} \int_{\gammam} n_2 \phi = - \phi\vert_{\gammap},
            \end{align}
            where $\dashint_X = \frac{1}{|X|}\int_{X} $ is the averaged integral and we used that since $n(x) = \frac{1}{\sqrt{1+(\hp'(x_1))^2}}\begin{pmatrix}-\hp'(x_1)\\1\end{pmatrix}$ on $\gammap$
            \begin{align}
                \int_{\gammap} n_2 \ dx = \int_{\gammap} \frac{1}{\sqrt{1+(\hp')^2}} \ dx = \int_0^\Gamma \frac{\sqrt{1+(\hp')^2}}{\sqrt{1+(\hp')^2}} \ dx_1 = \Gamma.
            \end{align}
            Therefore, the stream function satisfies
            \begin{alignat}{2}
                \Delta \phi &= \omega & \quad \textnormal{ in } &\Omega
                \\
                \phi &= - \dashint_{\Omega} u_1 & \quad \textnormal{ on } &\gammap
                \\
                \phi &= 0 & \quad \textnormal{ on } &\gammam.
            \end{alignat}

        \item \hypertarget{sfa_item_removeBdry}{Remove the boundary conditions by redefining the stream function}
        
            We define
            \begin{align}
                \tilde \phi = \phi + \frac{x_2-\hm(x_1)}{\hp(x_1)-\hm(x_1)}\dashint_{\Omega} u_1
                \label{elliptic_reg_def_tilde_phi}
            \end{align}
            and omit the argument of the height functions to improve the readability. Then $\tilde\phi$ fulfills
            \begin{align}
                \Delta \tilde\phi &= \Delta \phi + \Delta\left( \frac{x_2-\hm}{\hp-\hm}\right)\dashint_{\Omega} u_1
            \end{align}
            and as on the top boundary
            \begin{align}
                \tilde \phi = \phi + \dashint_{\Omega} u_1 = 0
            \end{align}
            and on the bottom boundary
            \begin{align}
                \tilde \phi = \phi = 0
            \end{align}
            its system can be written as
            \begin{alignat}{2}
                \Delta \tilde \phi &= \tilde f & \quad \textnormal{ in } &\Omega\label{elliptic_reg_no_bc_pde}
                \\
                \phi &= 0 & \quad \textnormal{ on } &\gammam\cup\gammap,\label{elliptic_reg_no_bc_bc}
            \end{alignat}
            where 
            \begin{align}
                \tilde f&= \omega + \Delta\left( \frac{x_2-\hm}{\hp-\hm}\right)\dashint_{\Omega} u_1.
                \label{elliptic_reg_def_tilde_f}
            \end{align}
            Also note that
            \begin{align}
                \left\|\frac{x_2-\hm}{\hp-\hm}\right\|_{W^{1,p}}&\lesssim 1
                \\
                \left\|\frac{x_2-\hm}{\hp-\hm}\right\|_{W^{2,p}}&\lesssim 1+\|\kappa\|_\infty
                \\
                \left\|\frac{x_2-\hm}{\hp-\hm}\right\|_{W^{3,p}}&\lesssim 1+\|\kappa\|_{W^{1,\infty}}.
            \end{align}
            Additionally by Poincaré's inequality as $\tilde \phi$ vanishes on the boundary
            \begin{align}
                \|\phi\|_{p} &\leq \|\tilde \phi\|_p + \|u\|_1 \leq \|\nabla\tilde\phi\|_p+\|u\|_1 \lesssim \|\nabla \phi\|_p + \|u\|_1 \lesssim \|u\|_p, \qquad \label{elliptic_reg_phi_zero_order}
            \end{align}
            where in the last estimate we used Hölder's inequality.

        \item \hypertarget{sfa_item_straighten}{Establish a change in variables that straighten the boundaries}
        
            Next, we establish the change of variables to a system with straightened boundaries. Let $x$ be the coordinates of the original domain $\Omega$ and $y$ be the ones of the straightened system, defined by
            \begin{alignat}{3}
                x&\in \Omega& &= \left\lbrace (x_1,x_2)\ \middle \vert\ 0 < x_1 < \Gamma, \hm(x_1)<x_2<\hp(x_1) \right\rbrace
                \\
                y&\in \Omegastraight& &= \left\lbrace (y_1,y_2)\ \middle \vert\ 0 < y_1 < \Gamma, 0<y_2<1 \right\rbrace.
            \end{alignat}
            We denote the Sobolev spaces in the different coordinate systems by
            \begin{align}
                L^{p}_x &= L^{p}(\Omega) & L^{p}_y &= L^{p}(\Omegastraight)
                \\
                W^{k,p}_x &= W^{k,p}(\Omega) & W^{k,p}_y &= W^{k,p}(\Omegastraight)
            \end{align}
            and the corresponding change of variables by
            \begin{align}
                y&= \Phi(x) = \begin{pmatrix}x_1\\ \frac{x_2-\hm}{\hp-\hm}\end{pmatrix}
                \\
                x&= \Psi(y) = \begin{pmatrix}y_1\\ \hm+(\hp-\hm)y_2\end{pmatrix}
            \end{align}
            with derivatives
            \begin{align}
                \nabla_x \Phi(x)
                &= \begin{pmatrix}1 & 0
                \\
                -\frac{\hm'}{\hp-\hm}-\frac{(x_2-\hm)(\hp'-\hm')}{(\hp-\hm)^2} & \frac{1}{\hp-\hm}\end{pmatrix}
                \\
                \nabla_y \Psi(y)
                &= \begin{pmatrix} 1 & 0
                \\
                \hm'+(\hp'-\hm')y_2 & \hp-\hm \end{pmatrix}
                \\
                \nabla_x \Phi(\Psi(y))
                &= \begin{pmatrix}1 & 0
                \\
                -\frac{\hm'+(\hp'-\hm')y_2}{\hp-\hm} & \frac{1}{\hp-\hm}\end{pmatrix}
                \\
                \nabla_y \Psi(\Phi(x))
                &= \begin{pmatrix} 1 & 0
                \\
                \hm'+\frac{(\hp'-\hm')(x_2-\hm)}{\hp-\hm} & \hp-\hm \end{pmatrix}
            \end{align}
            and in particular
            \begin{align}
                \nabla_y \Psi(y)\cdot \nabla_x \Phi(\Psi(y)) &=\begin{pmatrix}1 & 0 \\ 0 & 1\end{pmatrix}
                \\
                \nabla_x \Phi(x)\cdot \nabla_y \Psi(\Phi(x)) &=\begin{pmatrix}1 & 0 \\ 0 & 1\end{pmatrix}.\label{nibabnifsdnbnb}
            \end{align}
            Note that
            \begin{align}
                \|\Phi\|_{L^\infty_x}&\lesssim 1
                &
                \|\Psi\|_{L^\infty_y}&\lesssim 1
                \\
                \|\nabla_x\Phi\|_{L^\infty_x}&\lesssim 1
                &
                \|\nabla_y\Psi\|_{L^\infty_y}&\lesssim 1
                \\
                \|\nabla_x^2\Phi\|_{L^\infty_x}&\lesssim 1+\|\kappa\|_\infty
                &
                \|\nabla_y^2\Psi\|_{L^\infty_y}&\lesssim 1+\|\kappa\|_\infty
                \\
                \|\nabla_x^3\Phi\|_{L^\infty_x}&\lesssim 1+\|\kappa\|_{W^{1,\infty}}
                &
                \|\nabla_y^3\Psi\|_{L^\infty_y}&\lesssim 1+\|\kappa\|_{W^{1,\infty}}
            \end{align}
            and by chain rule
            \begin{align}
                \|\nabla_x \tilde \rho\|_{L^p_x} &\lesssim \|\nabla_y \bar \rho\|_{L^p_y}
                \\
                \|\nabla_y \bar \rho\|_{L^p_y} &\lesssim \|\nabla_x \tilde \rho\|_{L^p_x}
                \\
                \|\nabla_x^2 \tilde \rho\|_{L^p_x} &\lesssim \|\nabla_y^2 \bar \rho\|_{L^p_y}+(1+\|\kappa\|_\infty)\|\nabla_y \bar \rho\|_{L^p_y}
                \\
                \|\nabla_y^2 \bar \rho\|_{L^p_y} &\lesssim \|\nabla_x^2 \tilde \rho\|_{L^p_x}+(1+\|\kappa\|_\infty)\|\nabla_x \tilde \rho\|_{L^p_x}
                \label{elliptic_reg_chain_rule_change_of_variables}
                \\
                \|\nabla_x^3 \tilde \rho\|_{L^p_x} &\lesssim \|\nabla_y^3 \bar \rho\|_{L^p_y}+(1+\|\kappa\|_\infty)\|\nabla_y^2 \bar \rho\|_{L^p_y}+(1+\|\kappa\|_{W^{1,\infty}})\|\nabla_y \bar \rho\|_{L^p_y}
                \\
                \|\nabla_y^3 \bar \rho\|_{L^p_y} &\lesssim \|\nabla_x^3 \tilde \rho\|_{L^p_x}+(1+\|\kappa\|_\infty)\|\nabla_x^2 \tilde \rho\|_{L^p_x}+(1+\|\kappa\|_{W^{1,\infty}})\|\nabla_x \tilde \rho\|_{L^p_x}
            \end{align}
            for any $\tilde \rho\in W^{1,p}_x$, respectively $\tilde \rho\in W^{2,p}_x$ and $\tilde \rho\in W^{3,p}_x$ and $\bar\rho(y)=\tilde \rho(\Psi(y))=\tilde\rho(x)$.
        
        \item \hypertarget{sfa_item_straightenedPDE}{Derive the PDE in the straightened variables}
    
            Next, we want to derive the straightened system corresponding to \eqref{elliptic_reg_no_bc_pde}, \eqref{elliptic_reg_no_bc_bc}. We define the corresponding stream function in the straightened system by
            \begin{align}
                \bar \phi (y) = \tilde \phi(\Psi(y)) = \tilde \phi(x)\label{elliptic_reg_def_bar_phi}
            \end{align}
            and the operator by
            \begin{align}
                \bar L \bar \phi = \partial_{y_l}(\bar a_{kl}\partial_{y_k} \bar\phi),
            \end{align}
            where
            \begin{align}
                \bar a_{k,l}(y) = |\det \nabla_x \Phi(\Psi(y))|^{-1}\delta_{r,s} \partial_{x_r} \Phi_k(\Psi(y))\partial_{x_s}\Phi_l(\Psi(y)),
                \label{elliptic_regularity_def_a_kl}
            \end{align}
            i.e.
            \begin{align}
                \bar a_{11} &= \hp-\hm, & \bar a_{12} &= -\hm'+(\hp'-\hm')y_2,
                \\
                \bar a_{21} &= -\hm'+(\hp'-\hm')y_2, & \bar a_{22} &= \frac{1+\left(\hm'+(\hp'-\hm')y_2\right)^2}{\hp-\hm}.
            \end{align}
    
            Similarly to \eqref{elliptic_reg_def_bar_phi} we define $\bar f(y)=\tilde f(\Psi(y))=\tilde f(x)$ and for any $\bar \rho \in H^1_0(\Omegastraight)$ let $\tilde\rho(\Psi(y))=\tilde\rho(x)=\bar \rho(y)=\bar \rho(\Phi(x))$. To derive the PDE notice that by partial integration, the previous definitions and chain rule, it holds
            \begin{align}
                - &\int_{\Omegastraight} \bar L(y) \bar\phi(y) \bar\rho(y) \ dy
                \\
                &= \int_{\Omegastraight} \bar a_{kl}(y) \partial_{y_k} \bar \phi(y) \partial_{y_l}\bar\rho(y) \ dy
                \\
                &= \int_{\Omegastraight} |\det \nabla_x \Phi(\Psi(y))|^{-1}\partial_{x_r} \Phi_k(\Psi(y))\partial_{x_r}\Phi_l(\Psi(y)) \partial_{y_k} \bar \phi(y) \partial_{y_l}\bar\rho(y) \ dy
                \\
                &= \int_{\Omegastraight} |\det \nabla_x \Phi(\Psi(y))|^{-1}\partial_{x_r} \Phi_k(\Psi(y))\partial_{x_r}\Phi_l(\Psi(y)) \\
                &\qquad \cdot\partial_{x_i} \tilde \phi(\Psi(y)) \partial_{y_k}\Psi_i(y)\partial_{x_j}\tilde\rho(\Psi(y))\partial_{y_l}\Psi_j(y) \ dy,
            \end{align}
            which according to \eqref{nibabnifsdnbnb} simplifies to
            \begin{align}
                - \int_{\Omegastraight} &\bar L(y) \bar\phi(y) \bar\rho(y) \ dy
                \\
                &= \int_{\Omegastraight} |\det \nabla_x \Phi(\Psi(y))|^{-1} \delta_{r,i} \delta_{r,j}  \partial_{x_i} \tilde \phi(\Psi(y)) \partial_{x_j}\tilde\rho(\Psi(y)) \ dy
                \\
                &= \int_{\Omegastraight} |\det \nabla_x \Phi(\Psi(y))|^{-1} \partial_{x_i} \tilde \phi(\Psi(y)) \partial_{x_i}\tilde\rho(\Psi(y)) \ dy
                \\
                &= \int_{\Omega}  \partial_{x_i} \tilde \phi(x) \partial_{x_i}\tilde\rho(x) \ dx.
                \label{snbuinb}
            \end{align}
            Integration by parts and \eqref{elliptic_reg_no_bc_pde} yields
            \begin{align}
                \int_{\Omega}  \partial_{x_i} \tilde \phi(x) \partial_{x_i}\tilde\rho(x) \ dx
                &= -\int_{\Omega}  \Delta \tilde \phi(x) \tilde\rho(x) \ dx
                \\
                &= -\int_{\Omega}  \tilde f(x) \tilde\rho(x) \ dx
                \\
                &= -\int_{\Omegastraight}  (\hp-\hm)\bar f(y) \bar\rho(y) \ dy.\label{bnsifbngfbi}
            \end{align}
            Combining \eqref{snbuinb} and \eqref{bnsifbngfbi} one finds that $\bar\phi$ fulfills
            \begin{align}
                \bar L \bar\phi = (\hp-\hm)\bar f
            \end{align}
            and on the boundaries
            \begin{align}
                \bar \phi(y_1,0) = \tilde\phi(y_1,\hm(y_1)) = \tilde\phi(x_1,\hm(x_1)) = 0
                \\
                \bar \phi(y_1,1) = \tilde\phi(y_1,\hp(y_1)) = \tilde\phi(x_1,\hp(x_1)) = 0
            \end{align}
            i.e. $\bar \phi$ solves
            \begin{alignat}{2}
                \bar L \bar \phi &= (\hp-\hm) \bar f & \quad \textnormal{ in } &\Omegastraight\label{elliptic_reg_straight_pde}
                \\
                \bar\phi &= 0 & \quad \textnormal{ on } &\partial\Omegastraight.\label{elliptic_reg_straight_bc}
            \end{alignat}

        \item \hypertarget{sfa_item_straightenedBounds}{Derive the estimates for the straightened system}

            $\bar L$ is an elliptic operator as according to \eqref{elliptic_regularity_def_a_kl} for any $\xi\in \mathbb{R}^2$
            \begin{align}
                \bar a_{ij} \xi_i\xi_j &= (\hp-\hm) \partial_{x_r}\Phi_i(\Psi(y))\partial_{x_r}\Phi_j(\Psi(y)) \xi_i\xi_j
                \\
                &= (\hp-\hm) \frac{|\partial_{x_r}\Phi(\Psi(y))\cdot\xi|^2\|\nabla_y\Psi(y)\|^2}{\|\nabla_y\Psi(y)\|^2}
                \\
                &\geq (\hp-\hm) \frac{|\nabla_y\Psi(y)\partial_{x_r}\Phi(\Psi(y))\cdot\xi|^2}{\|\nabla_y\Psi(y)\|^2}
                \\
                &= \frac{\hp-\hm}{1+(\hp-\hm)^2+\left(\hm'+(\hp'-\hm')y_2\right)^2}|\xi|^2
                \\
                &\geq C|\xi|^2
            \end{align}
            for some constant $C>0$, depending on $\min(\hp-\hm)$, $\max(\hp-\hm)$, the Lipschitz constant of $\Omega$ and $\Gamma$. Therefore elliptic regularity for the system \eqref{elliptic_reg_straight_pde}, \eqref{elliptic_reg_straight_bc} implies
            \begin{align}
                \label{elliptic_reg_bar_phi_estimate}
                \|\bar\phi\|_{W^{2,p}_y}\lesssim \|(\hp-\hm)\bar f\|_{L^p_y} \lesssim \|\bar f\|_{L^p_y},
            \end{align}
            where the implicit constant only depends on $\min(\hp-\hm)$, $\max(\hp-\hm)$ and the Lipschitz constant of $\Omega$ and $\Gamma$.
        
        \item \hypertarget{sfa_item_streamFunctionBounds}{Translate the bounds back to the original system}
    
            By the definition of $\tilde\phi$ and $\bar\phi$, i.e. \eqref{elliptic_reg_def_tilde_phi} and \eqref{elliptic_reg_def_bar_phi}, \eqref{elliptic_reg_bar_phi_estimate} implies
            \begin{align}
                \|\phi\|_{W^{2,p}_x}
                &\lesssim \|\tilde\phi\|_{W^{2,p}_x}+(1+\|\kappa\|_\infty)\|u_1\|_{L^1_x}
                \\
                &\lesssim \|\bar\phi\|_{W^{2,p}_y}+(1+\|\kappa\|_\infty)\|\bar\phi\|_{W^{1,p}_y}+(1+\|\kappa\|_\infty)\|u_1\|_{L^1_x}
                \\
                &\lesssim \|\bar f\|_{L^p_y}+(1+\|\kappa\|_\infty)\|\bar\phi\|_{W^{1,p}_y}+(1+\|\kappa\|_\infty)\|u_1\|_{L^1_x}
                \\
                &\lesssim \|\tilde f\|_{L^p_x}+(1+\|\kappa\|_\infty)\|\tilde\phi\|_{W^{1,p}_x}+(1+\|\kappa\|_\infty)\|u_1\|_{L^1_x}
                \\
                &\lesssim \|\tilde f\|_{L^p_x}+(1+\|\kappa\|_\infty)\|\phi\|_{W^{1,p}_x}+(1+\|\kappa\|_\infty)\|u_1\|_{L^1_x}
                \\
                &\lesssim \|\omega\|_{L^p_x}+(1+\|\kappa\|_\infty)\|\phi\|_{W^{1,p}_x}+(1+\|\kappa\|_\infty)\|u_1\|_{L^1_x}\label{elliptic_reg_phi_estimate}
            \end{align}
            where we used \eqref{elliptic_reg_chain_rule_change_of_variables} and the definitions of $\bar f$ and $\tilde f$.
            By the definition of $\phi$ \eqref{elliptic_reg_phi_estimate} yields
            \begin{align}
                \|u\|_{W^{1,p}_x} &\lesssim \|\phi\|_{W^{2,p}_x}
                \\
                &\lesssim \|\omega\|_{L^p_x}+(1+\|\kappa\|_\infty)\|\phi\|_{W^{1,p}_x}+(1+\|\kappa\|_\infty)\|u_1\|_{L^1_x}
                \\
                &\lesssim \|\omega\|_{L^p_x}+(1+\|\kappa\|_\infty)\|u\|_{L^p_x},
                \label{elliptic_reg_u_estimate_only_p}
            \end{align}
            where in the last estimate we used \eqref{elliptic_reg_phi_zero_order}.
        
        \item \hypertarget{sfa_item_generalizeLowerOrderTerm}{Generalize the lower order term}

            Gagliardo-Nirenberg interpolation for the straightened domain 
            and Young's inequality imply
            \begin{align}
                \|u\|_{L^p_y}&\lesssim \|\nabla u\|_{L^p_y}^{\mu}\|u\|_{L^q_y}^{1-\mu}+\|u\|_{L^q_y}
                \\
                &\lesssim \epsilon \|\nabla u\|_{L^p_y} + \left(1+\epsilon^{-\frac{\mu}{1-\mu}}\right)\|u\|_{L^q_y}
                \label{elliptic_reg_gag_niren_and_youngs}
            \end{align}
            for any $\epsilon>0$, where $0<\mu=\frac{2(p-q)}{2(p-q)+pq}<1$ for $1\leq q<p<\infty$. Combining \eqref{elliptic_reg_u_estimate_only_p} and \eqref{elliptic_reg_gag_niren_and_youngs}
            \begin{align}
                \|u\|_{W^{1,p}_x} &\lesssim \|\omega\|_{L^p_x}+(1+\|\kappa\|_\infty)\|u\|_{L^p_x}
                \\
                &\lesssim \|\omega\|_{L^p_x}+(1+\|\kappa\|_\infty)\|u\|_{L^p_y}
                \\
                &\lesssim \|\omega\|_{L^p_x}+\epsilon(1+\|\kappa\|_\infty)\|\nabla u\|_{L^p_y}
                \\
                &\qquad + \left(1+\epsilon^{-\frac{\mu}{1-\mu}}\right)(1+\|\kappa\|_\infty)\|u\|_{L^q_y}
                \\
                &\lesssim \|\omega\|_{L^p_x}+\delta\|\nabla u\|_{L^p_y}
                \\
                &\qquad + \left(1+(1+\|\kappa\|_\infty)^{\frac{\mu}{1-\mu}}\right)(1+\|\kappa\|_\infty)\|u\|_{L^q_y}
                \\
                &\lesssim \|\omega\|_{L^p_x}+\delta\|\nabla u\|_{L^p_x}
                \\
                &\qquad + \left(1+(1+\|\kappa\|_\infty)^{\frac{\mu}{1-\mu}}\right)(1+\|\kappa\|_\infty)\|u\|_{L^q_x},
                \label{ksafjlaskdf}
            \end{align}
            where we chose $\epsilon = (1+\|\kappa\|_\infty)^{-1}\delta$ for any $\delta>0$. Note that since $p > q \geq 1$ and $0<\mu=\frac{2(p-q)}{2(p-q)+pq}<1$ Young's inequality implies
            \begin{align}
                \left(1+(1+\|\kappa\|_\infty)^{\frac{\mu}{1-\mu}}\right)(1+\|\kappa\|_\infty)
                &= 1+\|\kappa\|_\infty + (1+\|\kappa\|_\infty)^\frac{1}{1-\mu}
                \\
                &\lesssim 1 + \|\kappa\|_\infty^{\frac{1}{1-\mu}}.
                \label{nbsjifbaassdfsf}
            \end{align}
            Choosing $\delta$ sufficiently small in order to compensate the gradient term in \eqref{ksafjlaskdf} and using \eqref{nbsjifbaassdfsf} results in
            \begin{align}
                \|u\|_{W^{1,p}_x}
                &\lesssim \|\omega\|_{L^p_x} + \left(1+(1+\|\kappa\|_\infty)^{\frac{\mu}{1-\mu}}\right)(1+\|\kappa\|_\infty)\|u\|_{L^q_x}
                \\
                &\lesssim \|\omega\|_{L^p_x} + \left(1+\|\kappa\|_\infty^{1+\frac{2}{q}-\frac{2}{p}}\right)\|u\|_{L^q_x},
            \end{align}
            proving \eqref{elliptic_regularity_periodic_lemma_first_order}.
    \end{enumerate}
    Next, we focus on \eqref{elliptic_regularity_periodic_lemma_second_order}.    
    \begin{enumerate}[label=(2\alph*)]
        \item \hypertarget{sfa_item_higherOrderPDE}{Derive a PDE for the horizontal derivative of the straightened system}
        
            Calculating the horizontal derivative of \eqref{elliptic_reg_straight_pde} one finds
            \begin{align}
                \partial_{y_1}\left((\hp-\hm)\bar f\right) &= \partial_{y_1}\left(\bar L\bar\phi\right)
                \\
                &= \partial_{y_l}(\partial_{y_1}\bar a_{kl}\partial_{y_k}\bar \phi) + \partial_{y_l}(\bar a_{kl}\partial_{y_k}\partial_{y_1}\bar \phi)
                \\
                &= \partial_{y_l}(\partial_{y_1}\bar a_{kl}\partial_{y_k}\bar \phi) + \bar L\partial_{y_1}\bar \phi
            \end{align}
            and on the boundaries by \eqref{elliptic_reg_straight_bc}
            \begin{align}
                \partial_{y_1} \bar\phi = 0.
            \end{align}
            Therefore $\hat \phi = \partial_{y_1}\bar\phi$ solves
            \begin{alignat}{2}
                \bar L \hat \phi &= \hat f & \quad \textnormal{ in } &\Omegastraight
                \\
                \hat \phi &= 0 & \quad \textnormal{ on } &\lbrace y_2=0\rbrace\cup\lbrace y_2=1\rbrace,
            \end{alignat}
            where $\hat f = \partial_{y_1}\left((\hp-\hm)\bar f\right) - \partial_{y_l}(\partial_{y_1}\bar a_{kl}\partial_{y_k}\bar \phi)$.

        \item \hypertarget{sfa_item_higherOrderHorizontalBound}{Derive bounds for the horizontal derivative}
        
            As in \eqref{elliptic_reg_bar_phi_estimate} elliptic regularity implies
            \begin{align}
                \|\partial_{y_1}\bar \phi\|_{W^{2,p}_y}=\|\hat \phi\|_{W^{2,p}_y} \lesssim \|\hat f\|_{L^p_y},\label{nsbijnsb}
            \end{align}
            where the implicit constant only depends on $\min(\hp-\hm)$, $\max(\hp-\hm)$ and the Lipschitz constant of $\Omega$ and $\Gamma$. Note that this implies bounds for every third order derivative of $\bar\phi$ except for $\partial_{y_2}^3\bar\phi$.

        \item \hypertarget{sfa_item_higherOrderVerticalEquation}{Derive an equation for the remaining derivative}
            
            In order to derive bounds for $\partial_{y_2}^3 \bar\phi$ notice that
            \begin{align}
                \partial_{y_2}\bar f &= \partial_{y_2} \bar L\bar\phi
                \\
                &= \partial_{y_2} \partial_{y_l} (\bar a_{kl} \partial_{y_k} \bar\phi)
                \\
                &= \partial_{y_2} \partial_{y_1} (\bar a_{11} \partial_{y_1} \bar\phi) + \partial_{y_2} \partial_{y_2} (\bar a_{12} \partial_{y_1} \bar\phi) + \partial_{y_2} \partial_{y_1} (\bar a_{21} \partial_{y_2} \bar\phi)
                \\
                &\qquad + \partial_{y_2}^2 \bar a_{22}\partial_{y_2} \bar\phi + 2\partial_{y_2} \bar a_{22}\partial_{y_2}^2 \bar\phi + \bar a_{22}\partial_{y_2}^3 \bar\phi,
            \end{align}
            implying
            \begin{align}
                \partial_{y_2}^3 \bar\phi &= \bar a_{22}^{-1}\left(\partial_{y_2}\bar f-\partial_{y_2}^2 \bar a_{22}\partial_{y_2} \bar\phi - 2\partial_{y_2} \bar a_{22}\partial_{y_2}^2 \bar\phi\right.
                \label{sbnikgsbnf}
                \\
                &\qquad\left.-\partial_{y_2} \partial_{y_1} (\bar a_{11} \partial_{y_1} \bar\phi) - \partial_{y_2} \partial_{y_2} (\bar a_{12} \partial_{y_1} \bar\phi) - \partial_{y_2} \partial_{y_1} (\bar a_{21} \partial_{y_2} \bar\phi)\right)
            \end{align}
            as $\bar a_{22}\neq 0$.
            
        \item \hypertarget{sfa_item_higherOrderFullBound}{Get bounds for the full norm}
        
            Taking the norm of \eqref{sbnikgsbnf}
            \begin{align}
                \|\partial_{y_2}^3\bar\phi\|_{L^p_y}
                &\lesssim \|\bar a_{22}^{-1}\|_{L^\infty_y}\left(\|\partial_{y_2}\bar f\|_{L^p_y} + \|\nabla_y^2 \bar a\|_{L^\infty_y} \|\nabla_y\bar\phi\|_{L^p_y} \right.
                \\
                &\qquad \left.+ \|\nabla_y \bar a\|_{L^\infty_y} \|\nabla_y^2\bar\phi\|_{L^p_y} + \|\bar a\|_{L^\infty_y} \|\nabla_y^2\partial_{y_1}\bar\phi\|_{L^p_y}\right)
                \\
                &\lesssim \|\partial_{y_2}\bar f\|_{L^p_y} + \|\nabla_y^2 \bar a\|_{L^\infty_y} \|\nabla_y\bar\phi\|_{L^p_y}
                \\
                &\qquad + \|\nabla_y \bar a\|_{L^\infty_y} \|\nabla_y^2\bar\phi\|_{L^p_y} +\|\nabla_y^2\partial_{y_1}\bar\phi\|_{L^p_y},
                \label{snjbnsbsjbnfgjb}
            \end{align}
            where we used
            \begin{align}
                \|\bar a_{22}^{-1}\|_{L^\infty_y}&\lesssim 1
                \\
                \|\bar a\|_{L^\infty_y} &\lesssim 1.
            \end{align}
            
            Combining \eqref{snjbnsbsjbnfgjb}, \eqref{nsbijnsb} and the definition of $\hat f$, we obtain
            \begin{align}
                \|\bar \phi\|_{W^{3,p}_y}
                &\lesssim \|\partial_{y_1}\bar\phi\|_{W^{2,p}_y} +\|\partial_{y_2}^3\bar\phi\|_{L^p_y} + \|\bar\phi\|_{W^{2,p}_y}
                \\
                &\lesssim \|\partial_{y_1}\bar\phi\|_{W^{2,p}_y}+ \|\bar\phi\|_{W^{2,p}_y} + \|\partial_{y_2}\bar f\|_{L^p_y} 
                \\
                &\qquad + \|\nabla_y^2 \bar a\|_{L^\infty_y} \|\nabla_y\bar\phi\|_{L^p_y} + \|\nabla_y \bar a\|_{L^\infty_y} \|\nabla_y^2\bar\phi\|_{L^p_y}
                \\
                &\lesssim \|\hat f\|_{L^p_y} + \|\bar\phi\|_{W^{2,p}_y} + \|\partial_{y_2}\bar f\|_{L^p_y} 
                \\
                &\qquad + \|\nabla_y^2 \bar a\|_{L^\infty_y} \|\nabla_y\bar\phi\|_{L^p_y} + \|\nabla_y \bar a\|_{L^\infty_y} \|\nabla_y^2\bar\phi\|_{L^p_y}
                \\
                &\lesssim \|\bar f\|_{W^{1,p}_y} + \|\nabla_y^2 \bar a\|_{L^\infty_y} \|\nabla_y\bar\phi\|_{L^p_y} + (1+ \|\nabla_y \bar a\|_{L^\infty_y})  \|\bar\phi\|_{W^{2,p}_y}
                \\
                &\lesssim \|\bar f\|_{W^{1,p}_y} + (1+\|\kappa\|_{W^{1,\infty}} + \|\kappa\|_\infty^2) \|\bar\phi\|_{W^{1,p}_y} 
                \\
                &\qquad + (1+\|\kappa\|_\infty)  \|\bar\phi\|_{W^{2,p}_y},
                \label{aonfdsfasndf}
            \end{align}            
            where in the last estimate we used
            \begin{align}
                \|\bar a\|_{W^{1,\infty}_y} &\lesssim 1+\|\kappa\|_\infty
                \\
                \|\bar a\|_{W^{2,\infty}_y} &\lesssim 1+\|\kappa\|_{W^{1,\infty}} + \|\kappa\|_\infty^2.
            \end{align}
            Using \eqref{elliptic_reg_chain_rule_change_of_variables} for \eqref{aonfdsfasndf} we find
            \begin{align}
                \|\tilde\phi\|_{W^{3,p}_x} &\lesssim \|\bar\phi\|_{W^{3,p}_y} + (1+\|\kappa\|_\infty) \|\nabla^2_y \bar\phi\|_{L^{p}_y} + (1+\|\kappa\|_{W^{1,\infty}})\|\nabla_y \bar \phi\|_{L^p_y}
                \\
                &\lesssim \|\bar f\|_{W^{1,p}_y} + (1+\|\kappa\|_\infty)  \|\bar\phi\|_{W^{2,p}_y} 
                \\
                &\qquad +  (1+\|\kappa\|_{W^{1,\infty}} + \|\kappa\|_\infty^2) \|\bar\phi\|_{W^{1,p}_y}
                \\
                &\lesssim \|\tilde f\|_{W^{1,p}_x} + (1+\|\kappa\|_\infty)  \|\tilde\phi\|_{W^{2,p}_x} 
                \\
                &\qquad +  (1+\|\kappa\|_{W^{1,\infty}} + \|\kappa\|_\infty^2) \|\tilde\phi\|_{W^{1,p}_x}
            \end{align}
            and by the definition of $\tilde \phi$ and $\tilde f$, i.e. \eqref{elliptic_reg_def_tilde_phi} and \eqref{elliptic_reg_def_tilde_f}
            \begin{align}
                \|\nabla^2 u\|_{W^{2,p}_x}
                &\lesssim \|\phi\|_{W^{3,p}_x}
                \\
                &\lesssim \|\tilde \phi\|_{W^{3,p}_x} + (1+\|\kappa\|_{W^{1,\infty}})\|u\|_{L^1_x}
                \\
                &\lesssim \|\tilde f\|_{W^{1,p}_x} + (1+\|\kappa\|_{W^{1,\infty}})\|u\|_{L^1_x} + (1+\|\kappa\|_\infty)  \|\tilde\phi\|_{W^{2,p}_x}
                \\
                &\qquad + (1+\|\kappa\|_{W^{1,\infty}} + \|\kappa\|_\infty^2) \|\tilde\phi\|_{W^{1,p}_x}
                \\
                &\lesssim \|\omega\|_{W^{1,p}_x} + (1+\|\kappa\|_{W^{1,\infty}})\|u\|_{L^1_x} + (1+\|\kappa\|_\infty)  \|\tilde\phi\|_{W^{2,p}_x}
                \\
                &\qquad + (1+\|\kappa\|_{W^{1,\infty}} + \|\kappa\|_\infty^2) \|\tilde\phi\|_{W^{1,p}_x}
                \\
                &\lesssim \|\omega\|_{W^{1,p}_x} + (1+\|\kappa\|_\infty)  \|\nabla^2_x\phi\|_{L^p_x}
                \\
                &\qquad + (1+\|\kappa\|_{W^{1,\infty}} + \|\kappa\|_\infty^2) \left(\|\phi\|_{W^{1,p}_x}+\|u\|_{L^1_x}\right)
                \\
                &\lesssim \|\omega\|_{W^{1,p}_x} + \|\kappa\|_\infty  \|\omega\|_{L^p_x}
                \\
                &\qquad + (1+\|\kappa\|_{W^{1,\infty}} + \|\kappa\|_\infty^2) \left(\|\phi\|_{W^{1,p}_x}+\|u\|_{L^1_x}\right),
                \label{afbjugbngiidaofg}
            \end{align}
            where in the last inequality we used \eqref{elliptic_reg_phi_estimate}. Finally, $u=\nabla^\perp\phi$, \eqref{afbjugbngiidaofg}, \eqref{elliptic_reg_phi_zero_order} and Hölder's inequality imply
            \begin{align}
                \|\nabla^2 u\|_{W^{2,p}_x} &\lesssim \|\omega\|_{W^{1,p}_x} + \|\kappa\|_\infty \|\omega\|_{L^p_x} + (1+\|\kappa\|_{W^{1,\infty}} + \|\kappa\|_\infty^2) \|u\|_{L^p_x}
            \end{align}
            proving \eqref{elliptic_regularity_periodic_lemma_second_order}.
            \hfill\qed
    \end{enumerate}\renewcommand{\qedsymbol}{}
\end{proof}
\vspace{-2\baselineskip} 
We need a similar estimate for the standard domain. Since in this case, the results are qualitative, we can use a more abstract approach allowing for a simpler proof and a more general result. The following Lemma corresponds to Lemma 5.3 in \cite{bleitnerCarlsonNobili2023Large}.
\begin{lemma}
    \label{lemma_elliptic_regularity_bounded_domain}
    Let $1<q<\infty$, $1\leq r\leq \infty$, $k \in \mathbb{N}_0$, $\Omega$ be a bounded $C^{k+1,1}$-domain and $u\in W^{k+1,q}(\Omega)$ satisfy \eqref{nav_slip_incompressible} and \eqref{nav_slip_no_penetration}. Then there exists a constant $C>0$ depending on $\Omega$, $q$, $k$, $r$ such that
    \begin{align}
        \|\nabla u\|_{W^{k,q}} \leq C (\|\omega\|_{W^{k,q}}+\|u\|_r).
    \end{align}
    If additionally $u\in W^{k+3,q}(\Omega)$ satisfies \eqref{nav_slip_nav_slip_bc} with $\alpha \in W^{k+2,\infty}(\partial\Omega)$ and $\Omega$ is a $C^{k+3,1}$-domain, then
    \begin{align}
        \|\nabla u\|_{W^{k+2,q}}\leq C \left(\|\Delta \omega\|_{W^{k,q}} + (1+ \|\alpha\|_{W^{k+2,\infty}})\|u\|_{W^{k+2,q}}\right).
    \end{align}
\end{lemma}

\begin{proof}
    The proof follows a similar strategy as the one of Lemma \ref{lemma_elliptic_regularity_periodic_domain}. As here we are not interested in the explicit dependency on $\kappa$ we are able to use general results, significantly simplifying the proof.

    As before, one first needs to find the PDE for the stream function $\phi$, which, up to a constant, is given by $\nabla^\perp\phi=u$. Taking the curl, we find $\Delta \phi = \nabla^\perp \cdot u = \omega$. To derive boundary conditions, let $\lambda$ be the parametrization of a portion $\Gamma_i$ of the boundary by arc length. Then
    \begin{align}
        \frac{d}{d\lambda} \phi(x_1(\lambda),x_2(\lambda)) = \frac{d}{d\lambda}x(\lambda)\cdot \nabla \phi = \tau \cdot \nabla\phi = \tau^\perp\cdot\nabla^\perp \phi = - n\cdot u = 0
    \end{align}
    and therefore $\phi$ is constant along connected components of the boundaries. Combining these one finds
    \begin{alignat}{2}
        \Delta \phi &= \omega \qquad &\textnormal{ in }&\Omega
        \\
        \phi &= \psi_i \qquad &\textnormal{ on }&\Gamma_i
    \end{alignat}
    with constants $\psi_i$, where $\Gamma_i$ are the connected components of $\partial\Omega$.
    By elliptic regularity $\phi\mapsto (-\Delta \phi,\phi\vert_{\partial\Omega})$ is an isomorphism from $W^{k+2,q}(\Omega)$ onto $W^{k,q}(\Omega)\times W^{k+2-\frac{1}{q},q}(\partial\Omega)$. For details see Remark 2.5.1.2 in \cite{grisvard1985}. This implies that
    \begin{align}
        \|\phi\|_{W^{k+2,q}} \lessapprox \|\omega\|_{W^{k,q}} + \|\phi\|_{W^{k+2-\frac{1}{q},q}(\partial\Omega)},\label{xziubvzhbajoif}
    \end{align}
    where the implicit constant depends on $\Omega$, $q$ and $k$. In order to estimate the boundary term note that for any $s\geq 0$
    \begin{align}
        \|\phi\|_{W^{s,q}(\partial\Omega)}^q = \sum_i \|\psi_i\|_{W^{s,q}(\Gamma_i)}^q = \sum_i \|\psi_i\|_{L^{q}(\Gamma_i)}^q = \|\phi\|_{L^q(\partial\Omega)}^q.
        \label{hvbjckhbxc}
    \end{align}
    Since $\phi$ is only defined up to a constant we can choose it such that $\phi$ has vanishing average in $\Omega$ and Poincaré's inequality holds. Therefore, \eqref{hvbjckhbxc}, trace theorem, Poincaré's inequality and the definition of $\phi$ yield
    \begin{align}
        \|\phi\|_{W^{s,q}(\partial\Omega)}^q =\|\phi\|_{L^q(\partial\Omega)}^q
        \leq \|\phi\|_{W^{1,q}(\Omega)}^q \lesssim \|\nabla \phi\|_{L^q(\Omega)}^q = \|u\|_{L^q(\Omega)}^q.
        \label{uvbhgiuvbhuvb}
    \end{align}
    Combining \eqref{xziubvzhbajoif} and \eqref{uvbhgiuvbhuvb} and using the definition of $\phi$ one gets
    \begin{align}
        \|\nabla u\|_{W^{k,q}} \leq \|\phi\|_{W^{k+2,q}} \lessapprox \|\omega\|_{W^{k,q}} + \|\phi\|_{W^{k+2-\frac{1}{q},q}(\partial\Omega)} \lesssim  \|\omega\|_{W^{k,q}} + \|u\|_q.
        \\\label{uhvbihuvcb}
    \end{align}
    Finally, in order to change the norm on the zeroth order term on the right-hand side, Gagliardo-Nirenberg interpolation and Young's inequality imply
    \begin{align}
        \|u\|_q \lesssim \|\nabla u\|_q^\rho \|u\|_1^{1-\rho} + \|u\|_1 \lesssim \epsilon \rho \|\nabla u\|_q + \left(1+(1-\rho)\epsilon^{-\frac{\rho}{1-\rho}}\right)\|u\|_1\quad
        \label{jretklrejtl}
    \end{align}
    for any $\epsilon>0$ where $\rho = \frac{2q-2}{3q-2}$. Combining \eqref{uhvbihuvcb} and \eqref{jretklrejtl} we obtain
    \begin{align}
        \|\nabla u\|_{W^{k,q}} \lessapprox  \|\omega\|_{W^{k,q}} + \|u\|_q \lesssim \|\omega\|_{W^{k,q}} + \epsilon \|\nabla u\|_q + \left(1+\epsilon^{-\frac{\rho}{1-\rho}}\right)\|u\|_1
    \end{align}
    and choosing $\epsilon$ sufficiently small one can compensate the third term on the right-hand side implying
    \begin{align}
        \|\nabla u\|_{W^{k,q}} \lessapprox \|\omega\|_{W^{k,q}} + \|u\|_1 \lesssim \|\omega\|_{W^{k,q}} + \|u\|_{r},
        \label{oincvnibocnb}
    \end{align}
    where the last inequality is due to Hölder's inequality, proving the first statement.

    The proof of the second statement follows a similar strategy. In order to get boundary conditions for $\omega$, note that by \eqref{nav_slip_incompressible} and \eqref{nav_slip_no_penetration}
    \begin{align}
        -2(\alpha+\kappa) u_\tau &= 2 \tau \cdot \D u\ n - 2 n\cdot (\tau\cdot\nabla) u
        = \tau \cdot (n\cdot \nabla) u - n \cdot (\tau\cdot \nabla) u
        \\
        &= (\tau_i n_j - \tau_j n_i) \partial_j u_i
        = (\tau_1 n_2 - \tau_2 n_1)\partial_2 u_1 + (\tau_2 n_1 - \tau_1 n_2) \partial_1 u_2
        \\
        &= (-n_2^2-n_1^2) \partial_2 u_1 + (n_1^2 + n_2^2) \partial_1 u_2
        = \omega.
    \end{align}
    Again by elliptic regularity, details can be found in Remark 2.5.1.2 of \cite{grisvard1985}, $\omega\mapsto (-\Delta \omega, \omega\vert_{\partial\Omega})$ is an isomorphism form $W^{k+2,q}(\Omega)$ onto $W^{k,q}(\Omega)\times W^{k+2-\frac{1}{q},q}(\partial\Omega)$, implying
    \begin{align}
        \|\omega\|_{W^{k+2,q}}
        &\lessapprox \|\Delta \omega\|_{W^{k,q}} + \|\omega\|_{W^{k+2-\frac{1}{q},q}(\partial\Omega)}
        \\
        &= \|\Delta \omega\|_{W^{k,q}} + 2\|(\alpha+\kappa)u_\tau\|_{W^{k+2-\frac{1}{q},q}(\partial\Omega)}.
        \label{iugdsiugfdigb}
    \end{align}
    The boundary term can be estimated using Hölder's inequality and trace theorem, for which details can be found in Theorem 1.5.1.2 in \cite{grisvard1985}, by
    \begin{align}
        \|(\alpha+\kappa)u_\tau\|_{W^{k+2-\frac{1}{q},q}(\partial\Omega)} \lessapprox (1+\|\alpha\|_{W^{k+2,\infty}(\partial\Omega)}) \|u\|_{W^{k+2,q}(\Omega)}.
        \label{vbuhibuvhnh}
    \end{align}
    Combining \eqref{iugdsiugfdigb} and \eqref{vbuhibuvhnh} yields
    \begin{align}
        \|\omega\|_{W^{k+2,q}}
        &\lessapprox \|\Delta \omega\|_{W^{k,q}} + 2\|(\alpha+\kappa)u_\tau\|_{W^{k+2-\frac{1}{q},q}(\partial\Omega)}
        \\
        &\lessapprox \|\Delta \omega\|_{W^{k,q}} + (1+\|\alpha\|_{W^{k+2,\infty}}) \|u\|_{W^{k+2,q}}
    \end{align}
    and using \eqref{oincvnibocnb} we arrive at
    \begin{align}
        \|\nabla u\|_{W^{k+2,q}} &\lessapprox \|\omega\|_{W^{k+2,q}} + \|u\|_{r}
        \\
        &\lessapprox \|\Delta \omega\|_{W^{k,q}} + (1+\|\alpha\|_{W^{k+2,\infty}}) \|u\|_{W^{k+2,q}},
    \end{align}
    proving the second statement.
\end{proof}

\section{Diffusion and Navier-Slip}
In this section, we discuss the interplay between the diffusion term and Navier-Slip boundary conditions. Naturally, the boundary conditions enter the system when testing the equations with the corresponding velocity field. Lemma \ref{lemma_int_Delta_u_v} will show that the symmetric gradient and the boundary conditions arise organically when doing so. In order to prove the statements, we first need estimates for terms on the boundary. In particular, we need to extend functions that are only defined on the boundary to the whole domain.

The following inverse trace result will provide this. We will not state a proof here and refer the reader to Theorem 18.40 of \cite{leoni2017A}. 
More details can also be found in Theorem 1.5.1.2 of \cite{grisvard1985}.

\begin{lemma}[Inverse Trace Estimate]
    Let $\Omega$ be $C^{1,1}$ and $1\leq q<\infty$. Then for every $g\in W^{1,\infty}(\partial\Omega)$ there exists some $\zeta \in W^{1,q}(\Omega)$ such that
    \begin{align}
        \zeta\vert_{\partial\Omega}=g, \qquad \|\zeta\|_{W^{1,q}(\Omega)} \leq C \|g\|_{W^{1,\infty}(\partial\Omega)},
    \end{align}
    where the constant $C>0$ only depends on $q$ and $\Omega$.
    \label{lemma_inverse_trace_theorem}
\end{lemma}

With this extension at hand, we are able to derive the identity and bounds that arise when testing the diffusion term with a vector field.

\begin{lemma}
    Assume $\Omega$ is $C^{1,1}$, $u,v\in H^1(\Omega)$. Let $u$ satisfy $\nabla \cdot u=0$ and the boundary conditions $n \cdot u=0$ and $\tau \cdot (\D u\ n +\alpha u)=0$ with $\alpha\in L^\infty(\partial\Omega)$. Then
    \begin{align}
        |\langle \Delta u,v\rangle| \leq \|u\|_{H^1}\|v\|_{H^1}\label{bnjkretbjr}
    \end{align}
    and if additionally  $n \cdot v=0$ on $\partial\Omega$, then
    \begin{align}
        - \langle \Delta u, v \rangle = 2\int_{\Omega} (\D u)_{ij}(\D v)_{ij} + 2\int_{\partial\Omega} \alpha u_\tau v_\tau.
    \end{align}
    \label{lemma_int_Delta_u_v}
\end{lemma}

\begin{proof}
    Assume at first $u\in H^2$, then integration by parts and projecting $v$ onto its tangential and normal part on the boundary, i.e. $v=(v\cdot n) n + (v\cdot\tau)\tau$ yield
    \begin{align}
        -\int_{\Omega} \Delta u\cdot v &= - \int_{\partial\Omega} v_j n_i \partial_i u_j + \int_{\Omega} \partial_i u_j \partial_i v_j
        \\
        &= - \int_{\partial\Omega} v_k \tau_k \tau_j n_i \partial_i u_j - \int_{\partial\Omega} v_k n_k n_j n_i \partial_i u_j + \int_{\Omega} \partial_i u_j \partial_i v_j 
        \\
        &= - \int_{\partial\Omega} v_k \tau_k \tau_j n_i (\partial_i u_j + \partial_j u_i) + \int_{\partial\Omega} v_k \tau_k \tau_j n_i \partial_j u_i
        \\
        &\qquad - \int_{\partial\Omega} v_k n_k n_j n_i \partial_i u_j + \int_{\Omega} \partial_i u_j \partial_i v_j .
        \label{nerjtkerj}
    \end{align}
    Note that, using the boundary condition $\tau \cdot (\D u\ n +\alpha u)=0$ and $u=u\cdot\tau\tau + u\cdot n n = u\cdot\tau\tau$ since $u\cdot n = 0$, the first term on the right-hand side of \eqref{nerjtkerj} satisfies
    \begin{align}
        - \int_{\partial\Omega} v_k \tau_k \tau_j n_i (\partial_i u_j + \partial_j u_i) &= - 2\int_{\partial\Omega} (v\cdot \tau)(\tau \cdot \D u \ n)
        = 2 \int_{\partial\Omega} \alpha (v\cdot \tau)(u\cdot \tau)
        \\
        &
        = 2 \int_{\partial\Omega} \alpha u\cdot v.
        \label{lwjerkwer}
    \end{align}
    For the second term on the right-hand side of \eqref{nerjtkerj} \eqref{n_tau_grad_u}, i.e $n\cdot(\tau\cdot\nabla)u = \kappa u\cdot \tau$ yields
    \begin{align}
        \int_{\partial\Omega} v_k \tau_k \tau_j n_i \partial_j u_i &= \int_{\partial\Omega} (v\cdot \tau) n\cdot (\tau\cdot\nabla) u
        = \int_{\partial\Omega} \kappa (v\cdot \tau) (u\cdot\tau)
        = \int_{\partial\Omega} \kappa u\cdot v.
        \\
        \label{bnktrjbh}
    \end{align}
    In order to estimate the third term on the right-hand side of \eqref{nerjtkerj} we need to extend $n$ to a function in $\Omega$. By Lemma \ref{lemma_inverse_trace_theorem}, there exists $\zeta\in W^{1,4}(\Omega)$ satisfying
    \begin{align}
        \zeta\vert_{\partial\Omega} = n, \qquad \|\zeta\|_{W^{1,4}(\Omega)}\lesssim \|n\|_{W^{1,\infty}(\partial\Omega)} \lesssim (1+\|\kappa\|_\infty)
        \label{nbejrktbre}
    \end{align}
    and therefore Stokes theorem yields
    \begin{align}
        -\int_{\partial\Omega} v_k n_k n_j n_i \partial_i u_j &= -\int_{\partial\Omega} v_k \zeta_k n_j \zeta_i \partial_i u_j 
        \\
        &= - \int_{\Omega} \partial_j (v_k \zeta_k \zeta_i \partial_i u_j)
        \\
        &= - \int_{\Omega} \partial_j v_k \zeta_k \zeta_i \partial_i u_j - \int_{\Omega} v_k \partial_j \zeta_k \zeta_i \partial_i u_j - \int_{\Omega} v_k \zeta_k \partial_j \zeta_i \partial_i u_j
        \\
        &\qquad - \int_{\Omega} v_k \zeta_k \zeta_i \partial_i \partial_j u_j
        \\
        &= - \int_{\Omega} \partial_j v_k \zeta_k \zeta_i \partial_i u_j - \int_{\Omega} v_k \partial_j \zeta_k \zeta_i \partial_i u_j - \int_{\Omega} v_k \zeta_k \partial_j \zeta_i \partial_i u_j,
    \end{align}
    where in the last identity we used that $\nabla \cdot u=0$. Consequently, using Hölder's inequality and Sobolev embedding
    \begin{align}
        &\left|\int_{\partial\Omega} v_k n_k n_j n_i \partial_i u_j\right|
        \\
        &\qquad \qquad \leq \int_{\Omega} |\partial_j v_k \zeta_k \zeta_i \partial_i u_j| + \int_{\Omega} |v_k \partial_j \zeta_k \zeta_i \partial_i u_j| + \int_{\Omega} |v_k \zeta_k \partial_j \zeta_i \partial_i u_j|
        \\
        &\qquad \qquad \leq \|v\|_{H^1} \|\zeta\|_\infty^2 \|u\|_{H^1} + 2\|v\|_4 \|\zeta\|_{W^{1,4}} \|\zeta\|_\infty \|u\|_{H^1}
        \\
        &\qquad \qquad \lesssim \|\zeta\|_{W^{1,4}}^2 \|v\|_{H^1}\|u\|_{H^1}
        \\
        &\qquad \qquad \lesssim (1+\|\kappa\|_\infty) \|u\|_{H^1}\|v\|_{H^1},
        \label{bjtrhrhj}
    \end{align}
    where in the last estimate we used \eqref{nbejrktbre}.
    Combining \eqref{nerjtkerj}, \eqref{lwjerkwer}, \eqref{bnktrjbh}, and \eqref{bjtrhrhj} and using Hölder's inequality and trace theorem
    \begin{align}
        \left|\int_{\Omega} \Delta u\cdot v\right| &\lesssim \int_{\partial\Omega} (|\alpha|+|\kappa|) |u\cdot v| + (1+\|\kappa\|_\infty) \|u\|_{H^1}\|v\|_{H^1}
        \\
        &\lesssim (1+\|\alpha\|_\infty+\|\kappa\|_\infty) \|u\|_{H^1}\|v\|_{H^1}.
    \end{align}
    By approximation, this holds for $u\in H^1$, proving \eqref{bnjkretbjr}.
    If $v\cdot n =0$ on $\partial\Omega$, then the third term on the right-hand side of \eqref{nerjtkerj} vanishes, which according to \eqref{lwjerkwer} and \eqref{bnktrjbh} implies
    \begin{align}
        -\int_{\Omega} \Delta u\cdot v
        &= - \int_{\partial\Omega} v_k \tau_k \tau_j n_i (\partial_i u_j + \partial_j u_i) + \int_{\partial\Omega} v_k \tau_k \tau_j n_i \partial_j u_i + \int_{\Omega} \partial_i u_j \partial_i v_j
        \\
        &= \int_{\partial\Omega} (2\alpha+\kappa) u\cdot v + \int_{\Omega} \partial_i u_j\partial_i v_j.
    \end{align}
    Finally, using Lemma \ref{lemma_grad_ids} one finds
    \begin{align}
        -\int_{\Omega} \Delta u\cdot v
        &= \int_{\partial\Omega} (2\alpha+\kappa) u\cdot v + \int_{\Omega} \partial_i u_j\partial_i v_j
        \\
        &= 2 \int_{\Omega} (\D u)_{ij}(\D v)_{ij} + 2\int_{\partial\Omega} \alpha u\cdot v.
    \end{align}
\end{proof}

Note that by Lemma \ref{lemma_grad_ids} and Lemma \ref{lemma_int_Delta_u_v}
\begin{align}
    -\langle \Delta u,u\rangle &= 2\|\D u\|_2^2 + 2 \int_{\partial\Omega} \alpha u_\tau^2
    \\
    &= \|\nabla u\|_2^2 + \int_{\partial\Omega} (2\alpha+\kappa) u_\tau^2
    \\
    &= \|\omega\|_2^2 + 2\int_{\partial\Omega} (\alpha+\kappa) u_\tau^2,
\end{align}
so $\alpha$ and $\kappa$ are strongly related to each other. In fact, the stress-free boundary conditions
\begin{align}
    \omega = 0,
\end{align}
as investigated for example in \cite{doering2018longTime} could be translated to
\begin{align}
    \alpha = -\kappa
    \label{lkwrjklwwer}
\end{align}
as can be seen in the boundary conditions for the vorticity that we will derive in \eqref{vorticity_bc} and \eqref{nd_vorticity_bc}. However, throughout this thesis, we will assume $\alpha>0$, and therefore \eqref{lkwrjklwwer} is only reasonable in our scenario if $\Omega$ is a convex domain, since then $\kappa\leq 0$. In the case of Rayleigh-Bénard convection, we are working in a domain that is periodic in $x_1$, implying that convexity can only be achieved if $\kappa=0$, i.e. the boundary is flat.

The next estimate, together with Lemma \ref{lemma_int_Delta_u_v}, shows $-\langle \Delta u,u\rangle $ is actually comparable to $\|u\|_{H^1}^2$.

\begin{lemma}[Coercivity]
    Let $\Omega$ be $C^{1,1}$ and $u\in H^1(\Omega)$ satisfy \eqref{nav_slip_incompressible}-\eqref{nav_slip_nav_slip_bc} with $0<\alpha\in L^\infty(\partial\Omega)$. Then
    \begin{align}
        \|\D u\|_{L^2}^2 + \int_{\partial\Omega} \alpha u_\tau^2 &\gtrsim (1+\|\alpha^{-1}(1+|\kappa|)\|_\infty)^{-1} \|u\|_{H^1}^2
        \label{coercivity_lemma_H1}
        \\
        \|\D u\|_{L^2}^2 + \int_{\partial\Omega} \alpha u_\tau^2 &\gtrsim (1+\|\alpha^{-1}\|_\infty)^{-1} \|u\|_{L^2}^2,
        \label{coercivity_lemma_zero_order}
    \end{align}
    where the implicit constant only depends on $|\Omega|$ and the Lipschitz constant of the domain.
    \label{lemma_coercivity}
\end{lemma}

\begin{proof}
    Although the proofs are similar, we split them depending on the considered domain due to technical reasons.
    \begin{itemize}
        \item The Lipschitz domain\\
        For $x\in \Omega$ let $(\tilde x_1(x_1,x_2),x_2)$ be a point on $\partial\Omega$ such that the horizontal line $\Omegadash$ connecting $x$ with $(\tilde x_1,x_2)$ lies in $\Omega$, i.e.        \begin{align}
            \tilde x_1(x_1,x_2) &= \max_{y_1<x_1,(y_1,x_2)\in\partial\Omega}y_1,
            \\
            \Omegadash(x_1,x_2) &= \lbrace (y_1,y_2)\in \Omega \ \vert \ \tilde x_1(x_1,x_2) < y_1 < x_1, y_2=x_2\rbrace
        \end{align}
        and let $\Omegadashdash$ be all the points in $\Omega$ at height $x_2$, i.e.
        \begin{align}
            \Omegadashdash(x_2) = \lbrace (y_1,y_2)\in \Omega \ \vert \ y_2=x_2\rbrace.
        \end{align}
        The definition of these sets is illustrated in Figure \ref{fig:def_tilde_x1_omegadash_omegadashdash_rb} for the partially periodic domain.
        Then the fundamental theorem of calculus yields
        \begin{align}
            u_1^2(x) &=\left|u_1(\tilde x_1, x_2) + \int_{\Omegadash} \partial_1 u_1(y_1,x_2)\ dy_1\right|^2
            \\
            &\leq 2 u_\tau^2(\tilde x_1, x_2) + 2 \|\partial_1 u_1(\cdot,x_2)\|_{L^2(\Omegadashdash(x_2))}^2
            \label{boundary_poincare_u1}
        \end{align}
        and similarly for
        \begin{align}
            \tilde x_2(x_1,x_2) &= \max_{y_2<x_2,(x_1, y_2)\in\partial\Omega}y_2,
            \\
            \Omegavert(x_1,x_2) &= \lbrace (y_1,y_2)\in \Omega \ \vert \ y_1=x_1, \tilde x_2(x_1,x_2) < y_2 < x_2\rbrace
            \\
            \Omegavertvert(x_1) &= \lbrace (y_1,y_2)\in \Omega \ \vert \ y_1=x_1\rbrace
        \end{align}
        one has
        \begin{align}
            u_2^2(x) &=\left|u_2(x_1, \tilde x_2) + \int_{\Omegavert} \partial_2 u_2(x_1,y_2)\ dy_2\right|^2
            \\
            &\leq 2 u_\tau^2(x_1, \tilde x_2) + 2 \|\partial_2 u_2(x_1,\cdot)\|_{L^2(\Omegavertvert(x_1))}^2
            \label{boundary_poincare_u2}
        \end{align}
        and integrating \eqref{boundary_poincare_u1} and \eqref{boundary_poincare_u2} over $\Omega$ implies
        \begin{align}
            \|u\|_{L^2}^2 \lesssim \int_{\partial\Omega} u_\tau^2 + \|\D u\|_{L^2}^2.
            \label{boundary_poincare_u}
        \end{align}
        
        \item The partially periodic domain\\
        Fix some $x\in \Omega$. We first focus on $u_1$ and divide our analysis depending on the value of $x_2$.
        
        \begin{figure}
            \hrule
            \vspace{\baselineskip}
            \begin{center}
                \includegraphics[width=\textwidth]{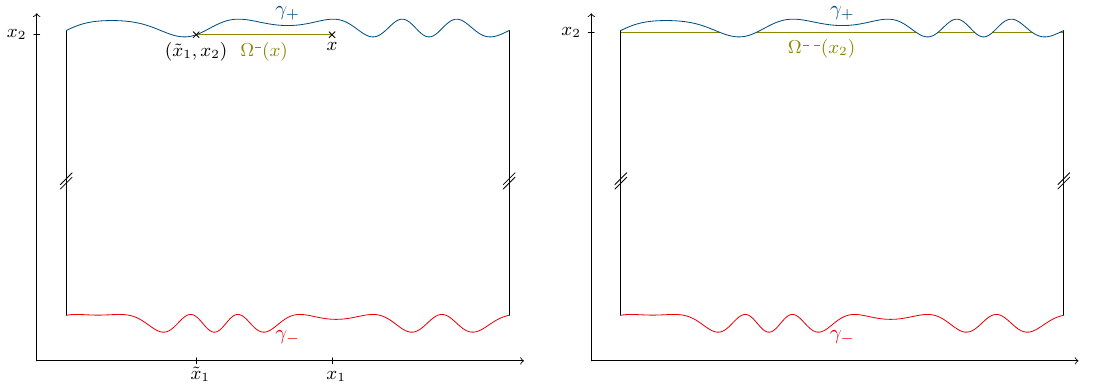}
                \vspace*{-5mm}
                \caption{Illustration of the definitions of $\tilde x_1$, $\Omegadash(x)$ and $\Omegadashdash(x_2)$.}
                \label{fig:def_tilde_x1_omegadash_omegadashdash_rb}
            \end{center}
            \hrule
        \end{figure}
        \begin{itemize}
            \item[-]
                At first assume that either $x_2\leq \max \hm$ or $x_2\geq \min \hp$. If it exists, we define
                \begin{align}
                    \tilde x_1(x_1,x_2) &= \max_{y_1<x_1,(y_1,x_2)\in\partial\Omega}y_1
                \end{align}
                as before and otherwise
                \begin{align}
                    \tilde x_1(x_1,x_2) &= \min_{x_1<y_1,(y_1,x_2)\in\partial\Omega}y_1
                \end{align}
                and again set
                \begin{align}
                    \Omegadashdash(x_2) &= \lbrace (y_1,y_2)\in \Omega \ \vert \ y_2=x_2\rbrace,
                \end{align}
                as illustrated in Figure \ref{fig:def_tilde_x1_omegadash_omegadashdash_rb}.
                Then the fundamental theorem of calculus yields
                \begin{align}
                    |u_1|^2(x)
                    &\lesssim |u_1|^2(\tilde x_1,x_2) + \int_{\tilde x_1}^{x_1} |\partial_1 u_1|^2(s,x_2) \ ds
                    \\
                    &\leq |u_\tau|^2(\tilde x_1,x_2) + \|\D u(\cdot, x_2)\|_{L^2(\Omegadashdash)}^2\label{coerc_est_u1_1}
                \end{align}
            \item[-]
                If $\max \hm <x_2<\min\hp$, define $\Omegabar(x_2) = \lbrace (y_1,y_2) \in \Omega \ \vert \ y_2<x_2 \rbrace$ and notice that
                \begin{align}
                    &\int_0^\Gamma u_1 (y_1,x_2) \ dy_1
                    \\
                    &= \int_0^\Gamma \left(u_1(y_1,\hm(y_1)) + \int_{\hm(y_1)}^{x_2} \partial_2 u_1(y_1,y_2) \ dy_2\right) \ dy_1 
                    \\
                    &=  \int_0^\Gamma u_1(y_1,\hm(y_1)) \ dy_1
                    \\
                    &\qquad + \int_0^\Gamma \left(\int_{\hm(y_1)}^{x_2} (\partial_2 u_1+\partial_1 u_2 - \partial_1 u_2)(y_1,y_2) \ dy_2\right) \ dy_1
                    \\
                    &=  \int_0^\Gamma \left(u_1(y_1,\hm(y_1)) + \int_{\hm(y_1)}^{x_2} (\partial_2 u_1+\partial_1 u_2)\ dy_2\right) \ dy_1
                    \\
                    &\qquad - \int_{\Omegabar(x_2)} \partial_1 u_2.
                    \label{nisfbnsfbgnfs}
                \end{align}
                Note that by periodicity the $\Omegabar$-integral on the right-hand side of \eqref{nisfbnsfbgnfs} can be estimated by
                \begin{align}
                    - \int_{\Omegabar(x_2)} \partial_1 u_2 \ dy &= - \int_{\Omegabar(\max\hm)} \partial_1 u_2 \ dy
                    \\
                    &= - \int_{\gammam} n_1 u_2
                    \\
                    &\leq \int_{\gammam} |u_\tau|.\label{nsbuibsfbd}
                \end{align}
                Combining \eqref{nisfbnsfbgnfs} and \eqref{nsbuibsfbd} one has
                \begin{align}
                    \left|\int_0^\Gamma u_1 (s,x_2) \ ds\right| \lesssim \int_{\gammam} |u_\tau| + \int_{\Omegabar(x_2)} |\partial_2 u_1+\partial_1 u_2|\ dy.\label{snbjgbf}
                \end{align}
                Again, using the fundamental theorem of calculus, we find
                \begin{align}
                    |u_1|^2(x) 
                    &= \left|u_1(x_1,x_2) - \dashint_0^\Gamma u_1(s,x_2)\ ds + \dashint_0^\Gamma u_1(s,x_2)\ ds \right|^2
                    \\
                    &\lesssim \left|\dashint_0^\Gamma (u_1(x_1,x_2) - u_1(s,x_2))\ ds \right|^2 + \left|\dashint_0^\Gamma u_1(s,x_2)\ ds \right|^2
                    \\
                    &\leq \left|\int_{0}^\Gamma |\partial_1 u_1|(l,x_2)\ dl\right|^2 + \left|\dashint_0^\Gamma u_1(s,x_2)\ ds \right|^2
                    \\
                    &\lesssim \|\partial_1 u_1(\cdot,x_2)\|_{L^2(0,\Gamma)}^2 + \left|\dashint_0^\Gamma u_1(s,x_2)\ ds \right|^2,\label{bnnfbfiu}
                \end{align}
                where in the last estimate we used Hölder's inequality. Combining \eqref{snbjgbf} and \eqref{bnnfbfiu}
                \begin{align}
                    |u_1|^2(x) &\lesssim \|\partial_1 u_1(\cdot,x_2)\|_{L^2(0,\Gamma)}^2 + \int_{\gammam} |u_\tau|^2
                    \\
                    &\qquad + \int_{\Omegabar(x_2)} |\partial_2 u_1+\partial_1 u_2|^2\ dy
                    \\
                    &\leq \|\D u(\cdot,x_2)\|_{L^2(0,\Gamma)}^2 + \int_{\gammam} |u_\tau|^2 + \|\D u\|_{L^2(\Omegabar(x_2))}^2\label{coerc_est_u1_3}.\qquad
                \end{align}
        \end{itemize}
        So regardless of the value of $x_2$, we can bound $u_1$ in terms of the symmetric gradient and the boundary values. To get bounds for $u_2$, notice that by the fundamental theorem of calculus
        \begin{align}
            |u_2|^2(x) &\lesssim \left|u_2(x_1,\hm(x_1)) + \int_{\hm(x_1)}^{x_2} \partial_2 u_2 (x_1,z)\ dz\right|^2
            \\
            &\lesssim |u_\tau|^2(x_1,\hm(x_1)) + \left(\int_{\hm(x_1)}^{\hp(x_1)} |\partial_2 u_2|(x_1,z)\ dz\right)^2
            \\
            &\leq |u_\tau|^2(x_1,\hm(x_1)) + \|\D u(x_1,\cdot)\|_{L^2(\hm(x_1),\hp(x_1))}^2,\label{coerc_est_u2}
        \end{align}
        where in the last estimate we used Hölder's inequality. Finally, combining \eqref{coerc_est_u1_1}, \eqref{coerc_est_u1_3} and \eqref{coerc_est_u2}
        \begin{align}
            |u|^2(x)
            &\lesssim |u_\tau|^2(\tilde x_1,x_2) + \|\D u(\cdot,x_2)\|_{L^2(\Omega \cap \lbrace y_2=x_2\rbrace)}^2 + \|\D u(\cdot,x_2)\|_{L^2(0,\Gamma)}^2
            \\
            &\qquad + \int_{\gammam} |u_\tau|^2 + \|\D u\|_{L^2(\Omega\cap \lbrace y_2 \leq x_2\rbrace)}^2 + |u_\tau|^2(x_1,\hm(x_1))
            \\
            &\qquad + \|\D u(x_1,\cdot)\|_{L^2(\hm(x_1),\hp(x_1))}^2\label{coerc_est_ges_rb}
        \end{align}
        and integrating \eqref{coerc_est_ges_rb} over $\Omega$ yields
        \begin{align}
            \|u\|_{L^2}^2 &\lesssim \|\D u\|_{L^2}^2 + \int_{\partial\Omega}u_\tau^2.
            \label{coercivity_rb_domain_zero_order}
        \end{align}
    \end{itemize}
    Smuggling in a factor of $\alpha$ in \eqref{boundary_poincare_u} and \eqref{coercivity_rb_domain_zero_order}, one gets
    \begin{align}
        \|u\|_{L^2}^2 &\lesssim \|\D u\|_{L^2}^2 + \int_{\partial\Omega}u_\tau^2
        \leq \|\D u\|_{L^2}^2 + \|\alpha^{-1}\|_\infty \int_{\partial\Omega}\alpha u_\tau^2,
    \end{align}
    which implies \eqref{coercivity_lemma_zero_order}.
    In order to estimate the full $H^1$-norm Lemma \ref{lemma_grad_ids} and \eqref{coercivity_rb_domain_zero_order} yield
    \begin{align}
        \|u\|_{H^1}^2 &= \|\nabla u\|_{L^2}^2 + \|u\|_{L^2}^2
        \\
        &= 2 \|\D u\|_{L^2}^2 - \int_{\partial\Omega} \kappa u_\tau^2 + \|u\|_{L^2}^2
        \\
        &\lesssim \|\D u\|_{L^2}^2 + \int_{\partial\Omega} \alpha \frac{1+\max\lbrace 0, -\kappa\rbrace}{\alpha} u_\tau^2
        \\
        &\lesssim (1+\|\alpha^{-1}(1+|\kappa|)\|_\infty) \left(\|\D u\|_{L^2}^2 + \int_{\partial\Omega} \alpha u_\tau^2\right).
    \end{align}
\end{proof}

Note that we estimated $\max\lbrace 0,-\kappa\rbrace\leq |\kappa|$ in order to simplify the terms. This is suboptimal. Imagine for example the domain to be a rectangle with an elliptic hole on the inside. Then $\kappa = 0$ almost everywhere on the outer boundary and $\kappa>0$ at the boundary where the hole is. So, in that case the curvature would not influence the estimate, which could therefore be improved to
\begin{align}
    \|\D u\|_{L^2}^2 + \int_{\partial\Omega}\alpha u_\tau^2\gtrsim (1+\|\alpha^{-1}\|_\infty)^{-1}\|u\|_{H^1}^2.
\end{align}

\section{Boundary Gradients}
\label{section_navSlip_boundary_gradients}

Naturally, when deriving $H^2$ estimates for the velocity, terms of type
\begin{align}
    \int_{\partial\Omega} u\cdot \nabla p
\end{align}
will arise as in the proof of Lemma \ref{lemma_long_time_H2_bound} for example. Note that if this was a bulk integral, then integration by parts would imply that it vanished as
\begin{align}
    \int_{\Omega} u\cdot \nabla p = - \int_{\partial\Omega} u\cdot n p - \int_{\Omega} p \nabla \cdot u = 0
\end{align}
due to \eqref{nav_slip_no_penetration} and \eqref{nav_slip_incompressible}. In \cite{drivas2022Bounds} the authors bound this term using the periodicity of the boundary and trace theorem in the case of flat boundaries as
\begin{align}
    \int_{\partial\Omega} u\cdot \nabla p = \int_{\partial\Omega} u_1 \partial_1 p = - \int_{\partial\Omega} \partial_1 u_1 p \lesssim \|u\|_{H^2}\|p\|_{H^1}.
\end{align}
We will use an analogous estimate in the case of curved boundaries in \eqref{bsdjgdfsgdfsghj}. The following Lemma will show that this term can be estimated by
\begin{align}
    \int_{\partial\Omega} u\cdot \nabla p \lesssim \|u\|_{H^1}\|p\|_{H^1},
\end{align}
reducing the norm on the right-hand side. This is the crucial improvement of \cite{bleitner2024Scaling} over \cite{bleitnerNobili2024Bounds} and \cite{drivas2022Bounds}. The Lemma is a slight deviation of Lemma 5.4 in \cite{bleitnerCarlsonNobili2023Large}.

\begin{lemma}
    \label{lemma_u_cdot_nabla_p_on_boundary}
    Let $\Omega$ be $C^{1,1}$, $f\in W^{1,\infty}(\partial\Omega)$, $\rho\in W^{1,q}(\Omega)$ and $u\in W^{1,q'}(\Omega)$ with $\frac{1}{q}+\frac{1}{q'}=1$ satisfy \eqref{nav_slip_incompressible} and \eqref{nav_slip_no_penetration}. Then
    \begin{align}
        \left|\int_{\partial\Omega} f u\cdot \nabla \rho\right| \lesssim
        \|fn\|_{W^{1,\infty}} \| u \|_{W^{1,q'}} \|\rho\|_{W^{1,q}},
    \end{align}
    where the implicit constant only depends on $|\Omega|$ and the Lipschitz constant of the domain.
\end{lemma}

\begin{proof}
    Note that connected components of the boundary are closed curves and $\tau\cdot \nabla$ is the derivative along the boundary parameterized by arc length. Therefore, due to the periodicity of these boundary components, we find the integration by parts formula
    \begin{align}
        \int_{\partial\Omega} f u\cdot \nabla \rho
        = \int_{\partial\Omega} f u_\tau \tau\cdot \nabla \rho
        = - \int_{\partial\Omega} \rho \tau\cdot \nabla (f u_\tau),
    \end{align}
    which by product rule yields
    \begin{align}
        \int_{\partial\Omega} f u\cdot \nabla \rho = - \int_{\partial\Omega} f \rho \tau\cdot \nabla u_\tau - \int_{\partial\Omega} \rho u_\tau \tau\cdot \nabla f.\label{nifgbnsiufgbn}
    \end{align}
    The second term on the right-hand side of \eqref{nifgbnsiufgbn} can be estimated by trace theorem and Hölder's inequality as
    \begin{align}
        \left|\int_{\partial\Omega} \rho u_\tau \tau\cdot \nabla f\right| \lesssim \|f\|_ {W^{1,\infty}} \|\rho u\|_{W^{1,1}(\Omega)} \lesssim \|f\|_{W^{1,\infty}} \|\rho\|_{W^{1,q}} \|u\|_{W^{1,q'}}.\qquad
        \label{ycxvyxcv}
    \end{align}
    In order to estimate the first term on the right-hand side of \eqref{nifgbnsiufgbn}, notice that since $\tau\cdot (\tau\cdot \nabla) \tau = \frac{1}{2}\tau\cdot\nabla (\tau^2)=0$ due to $\tau^2=1$ one has
    \begin{align}
        \tau\cdot \nabla u_\tau = \tau\cdot \nabla (u\cdot \tau) = \tau\cdot (\tau\cdot \nabla) u+u_\tau \tau\cdot (\tau\cdot \nabla) \tau = \tau\cdot (\tau\cdot \nabla) u.\quad\label{qwerqwer}
    \end{align}
    As $\tau = n^\perp$ one has $\tau_i\tau_j+n_i n_j = \delta_{ij}$ and therefore \eqref{qwerqwer} implies
    \begin{align}
        \tau\cdot \nabla u_\tau = \tau_i \tau_j \partial_j u_i = \delta_{ij}\partial_j u_i - n_i n_j \partial_j u_i = \nabla\cdot u - n\cdot (n\cdot \nabla) u = - n\cdot (n\cdot \nabla) u,
    \end{align}
    where in the last inequality we used \eqref{incompressibility}. Accordingly, the first term on the right-hand side of \eqref{nifgbnsiufgbn} can be written as
    \begin{align}
        - \int_{\partial\Omega} f \rho \tau\cdot \nabla u_\tau = \int_{\partial\Omega} f \rho n\cdot (n\cdot \nabla)u.\label{yxcvyyxcv}
    \end{align}
    To use Stokes theorem we need to first extend $fn$ to a function defined on $\Omega$. By Lemma \ref{lemma_inverse_trace_theorem} there exists $\zeta\in W^{1,\tilde q}(\Omega)$ such that
    \begin{align}
        \zeta\vert_{\partial\Omega}=fn, \qquad \|\zeta\|_{W^{1,\tilde q}(\Omega)} \leq C \|fn\|_{W^{1,\infty}(\partial\Omega)}
        \label{gnsigb}
    \end{align}
    and therefore Stokes theorem for \eqref{yxcvyyxcv} yields
    \begin{align}
        - \int_{\partial\Omega} f \rho \tau\cdot \nabla u_\tau = \int_{\partial\Omega} f \rho n\cdot (n\cdot \nabla) u = \int_{\partial\Omega} \rho n\cdot (\zeta\cdot \nabla) u = \int_{\Omega} \nabla \cdot \left(\rho  (\zeta\cdot \nabla) u\right).
        \label{yxcvyx}
    \end{align}
    In the case of the partially periodic domain, an explicit example of such an extension is given by
    \begin{align}
        \zeta(x_1,x_2) = \frac{\hp(x_1) - x_2}{\hp(x_1) - \hm(x_1)} f^-(x_1) n^-(x_1) + \frac{x_2 - \hm(x_1)}{\hp(x_1) - \hm(x_1)} f^+(x_1) n^+(x_1),
    \end{align}
    where $f^-(x_1)=f(x_1,\hm(x_1))$ and $f^+(x_1)=f(x_1,\hp(x_1))$, which fulfills
    \begin{align}
        \|\zeta\|_{L^\infty(\Omega)}\leq 2 \|f\|_{L^\infty(\partial\Omega)}, \qquad \|\zeta\|_{W^{1,\infty}(\Omega)}\leq C \|fn\|_{W^{1,\infty}(\partial\Omega)},
    \end{align}
    where $C>0$ depends only on the Lipschitz constant of $\Omega$ and the minimal distance of the boundaries in the vertical direction, i.e. $\min_{1\leq x_1\leq \Gamma}h^+(x_1)-h^-(x_1)$.
    \expliciteCalculation{
    \begin{align}
        \tilde h:&= \frac{\hp - x_2}{\hp - \hm}
        \\
        \zeta &= \tilde h f^- n^- + (1-\tilde h) f^+n^+
        \\
        \nabla\zeta &=
        \begin{pmatrix}
            \tilde h'f^- n^- +\tilde h{f^-}'n^-+\tilde hf^-{n^-}'+\tilde h'f^+ n^+ +(1-\tilde h){f^+}'n^+ + (1-\tilde h)f^+{n^+}'
            \\
            \frac{f^+n^+-f^-n^-}{h^+-h^-}
        \end{pmatrix}
        \\
        \|\zeta\|_\infty &\leq 2\|f\|_\infty
        \\
        \|\nabla \zeta\|_\infty &\leq 2\|\tilde h'\|_\infty \|f\|_\infty + 2\|f'\|_\infty + 2\|fn'\|_\infty
        \\
        \tilde h' &= \frac{\hp'(\hp-\hm)-(\hp-x_2)(\hp'-\hm')}{(\hp-\hm)^2}
        \\
        \|\tilde h'\|_\infty &\leq 2\frac{\|h'\|_\infty}{d} \leq C(\textnormal{lip},d)
        \\
        \|\zeta\|_\infty &\lesssim \|f\|_\infty
        \\
        \|\zeta\|_{W^{1,\infty}} &\lesssim_{\textnormal{lip},d} \|f\|_\infty + \|f\|_{W^{1,\infty}} + \|f n'\|_\infty \lesssim \|fn\|_{W^{1,\infty}}
    \end{align}
    }
    Using product rule \eqref{yxcvyx} yields
    \begin{align}
        -\int_{\partial\Omega} f \rho \tau\cdot \nabla u_\tau  &= \int_{\Omega} \nabla \cdot \left(\rho  (\zeta\cdot \nabla) u\right)
        \\
        &= \int_{\Omega}  (\zeta\cdot \nabla) u \cdot\nabla  \rho  + \int_{\Omega} \rho  \partial_i \zeta_j\partial_j u_i + \int_{\Omega} \rho  (\zeta\cdot \nabla)\nabla \cdot u
        \\
        &= \int_{\Omega}  (\zeta\cdot \nabla) u \cdot\nabla  \rho  + \int_{\Omega} \rho  \partial_i \zeta_j\partial_j u_i,
        \label{yxcvyxcv}
    \end{align}
    where the last term vanished due to \eqref{incompressibility}. Hölder's inequality with
    \begin{align}
            \begin{aligned}
                r&=\tfrac{2q}{2-q},&\hspace{-7pt} s&=2 & \quad\textnormal{ if }q&<2
                \\
                r&=4,&\hspace{-7pt} s&=4 & \quad\textnormal{ if }q&=2
                \\
                r&=2q,&\hspace{-7pt} s&=2q & \quad\textnormal{ if }q&>2
            \end{aligned}
    \end{align}
    we can estimate the right-hand side of \eqref{yxcvyxcv} by
    \begin{align}
        \left|\int_{\partial\Omega} f \rho \tau\cdot \nabla u_\tau\right|
        &= \left|\int_{\Omega}  (\zeta\cdot \nabla) u \cdot\nabla  \rho  + \int_{\Omega} \rho  \partial_i \zeta_j\partial_j u_i\right|
        \\
        &\lesssim \|\zeta\|_\infty \| u \|_{W^{1,q'}} \|\rho\|_{W^{1,q\phantom{'}}} + \|\zeta\|_{W^{1,s\phantom{'}}}\|\rho\|_{r}\|u\|_{W^{1,q'}}
    \end{align}
    and with these values of $r$ and $s$ the Sobolev embedding implies
    \begin{align}
        \left|\int_{\partial\Omega} f \rho \tau\cdot \nabla u_\tau\right|
        &\lesssim \|\zeta\|_\infty \| u \|_{W^{1,q'}} \|\rho\|_{W^{1,q\phantom{'}}} + \|\zeta\|_{W^{1,s\phantom{'}}}\|\rho\|_{r}\|u\|_{W^{1,q'}}
        \\
        &\lesssim \|\zeta\|_{W^{1,s+1}} \| u \|_{W^{1,q'}} \|\rho\|_{W^{1,q\phantom{'}}}
        \\
        &\lesssim \|fn\|_{W^{1,\infty}(\partial\Omega)} \| u \|_{W^{1,q'}} \|\rho\|_{W^{1,q}},
        \label{ycvxvxyv}
    \end{align}
    where in the last estimate we used \eqref{gnsigb}. Combining \eqref{nifgbnsiufgbn}, \eqref{ycxvyxcv} and \eqref{ycvxvxyv} yields
    \begin{align}
        \left|\int_{\partial\Omega} f u\cdot \nabla \rho\right| &\lesssim \left|\int_{\partial\Omega} f \rho \tau\cdot \nabla u_\tau\right| + \left|\int_{\partial\Omega} \rho u_\tau \tau\cdot \nabla f\right|
        \\
        &\lesssim (\|f\|_{W^{1,\infty}}+\|fn\|_{W^{1,\infty}}) \| u \|_{W^{1,q'}} \|\rho\|_{W^{1,q}}.
        \label{yxcvyv}
    \end{align}
    Finally, to unify the boundary terms notice that as $|n|=1$ and $2 n\cdot (\tau \cdot \nabla) n = \tau \cdot \nabla (n\cdot n) = 0$ one has
    \begin{align}
        \|f\|_\infty &= \|fn\cdot n\|_\infty \leq \|fn\|_\infty
        \\
        \|\tau\cdot \nabla f\|_\infty &= \|\tau\cdot \nabla (fn\cdot n)\|_\infty \leq \|n\cdot (\tau\cdot \nabla) (fn) \|_\infty + \|fn\cdot (\tau\cdot \nabla) n \|_\infty
        \\
        &\leq \|\tau\cdot\nabla (fn)\|_\infty
    \end{align}
    implying
    \begin{align}
        \|f\|_{W^{1,\infty}}\leq \|fn\|_{W^{1,\infty}},
    \end{align}
    which combined with \eqref{yxcvyv} yields the claim.

\end{proof}

\newpage
\chapter{Rayleigh-Bénard Convection}
\label{chapter:rayleigh-benard-convection}

\section{The Model}

\begin{figure}
    \hrule
    \vspace{0.5\baselineskip}
    \begin{center}
        \includegraphics[width=\textwidth]{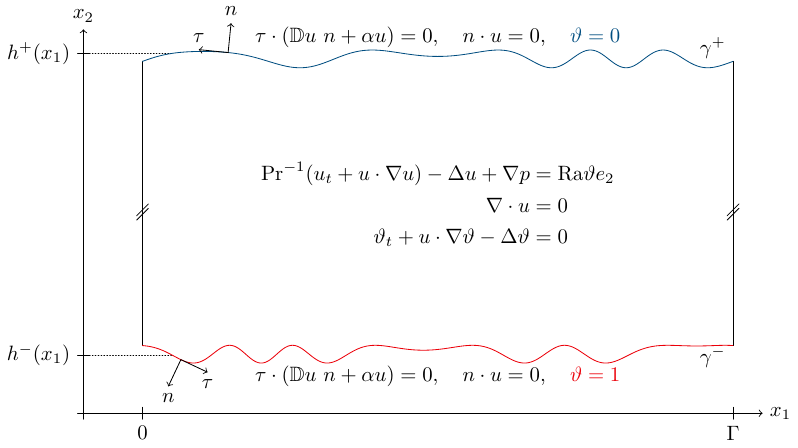}
        \vspace*{-7mm}
        \caption{Overview of the system considered in Section \ref{section:general_system}.}
        \label{fig:rb_overview_curved}
    \end{center}
    \hrule
\end{figure}

Here we study the following system, motivated in Section \ref{section:introduction_rb}.
\begin{alignat}{2}
    \frac{1}{\Pr}(u_t + u\cdot \nabla u) + \nabla p - \Delta u &= \Ra \theta e_2 \qquad & \textnormal{ in }&\Omega
    \label{navier_stokes}
    \\
    \label{incompressibility}
    \nabla \cdot u &= 0 \qquad & \textnormal{ in }&\Omega
    \\
    \label{advection_diffusion_equation}
    \theta_t + u\cdot \nabla \theta - \Delta \theta &= 0 \qquad & \textnormal{ in }&\Omega
    \\
    \label{no_penetration_bc}
    n\cdot u &= 0 \qquad & \textnormal{ on }&\gammam\cup\gammap
    \\
    \label{nav_slip_bc}
    \tau \cdot (\D u\ n +\alpha u)&=0 \qquad & \textnormal{ on }&\gammam\cup\gammap
    \\
    \label{heat_bc}
    \theta &= \bigg\lbrace \begin{matrix}0\\ 1 \end{matrix}  \qquad& \begin{matrix}\textnormal{ on }\\\textnormal{ on }\end{matrix}&\begin{matrix}\gammap\\\gammam\end{matrix}
    \\
    (u,\theta)(\cdot,0)&= (u_0,\theta_0) \qquad & \textnormal{ in }&\Omega
    \label{rb_ic}
\end{alignat}
with periodic boundary conditions in the horizontal direction, where\noeqref{rb_ic}
\begin{align}
    \Omega &= \left\lbrace (x_1,x_2)\in \mathbb{R}^2 \ \middle \vert\ 0\leq x_1\leq \Gamma, \hm(x_1)\leq x_2\leq \hp(x_1) \right\rbrace
    \\
    \gammap &= \left\lbrace (x_1,x_2)\in \mathbb{R}^2 \ \middle \vert\ 0\leq x_1\leq \Gamma, x_2= \hp(x_1) \right\rbrace
    \\
    \gammam &= \left\lbrace (x_1,x_2)\in \mathbb{R}^2 \ \middle \vert\ 0\leq x_1\leq \Gamma, x_2= \hm(x_1) \right\rbrace
    \label{gammam_def}
\end{align}
for sufficiently smooth $\hp$, $\hm$, which we also assume to not intersect and on average to be separated by distance $1$, such that the domain size is given by $|\Omega|=\Gamma$. An overview of the general system is given in Figure \ref{fig:rb_overview_curved}.

Additionally, we impose that $0\leq \theta_0\leq 1$ almost everywhere in $\Omega$, and therefore, the maximum principle implies that
\begin{align}
    \label{maximum_principle}
    0\leq \theta(x,t) &\leq 1
\end{align}
for almost every $x\in \Omega$ and $t\geq 0$ and $\theta(t)\in H^2(\Omega)$ for $t>0$ (\cite{foias1987Attractors}), which from now on we will always assume without explicitly mentioning.

The following vorticity formulation is extensively employed during the upcoming analysis. The vorticity $\omega = \nabla^\perp \cdot u = -\partial_2 u_1 + \partial_1 u_2$ satisfies
\begin{alignat}{2}
    \frac{1}{\Pr} (\omega_t + u\cdot \nabla\omega) - \Delta \omega &= \Ra \partial_1 \theta \qquad &\textnormal{ in }&\Omega
    \label{vorticity_equation}
    \\
    \omega &= -2 (\alpha+\kappa) u_\tau  \qquad &\textnormal{ on }&\gammap\cup\gammam.
    \label{vorticity_bc}
\end{alignat}
\eqref{vorticity_equation} is an immediate consequence of applying $\nabla^\perp \cdot$ to \eqref{navier_stokes} and using 
\begin{align}
    \nabla^\perp \cdot (u\cdot\nabla u) = u\cdot \nabla \omega + \omega \nabla \cdot u = u\cdot\nabla \omega
\end{align}
due to \eqref{incompressibility}. The boundary condition follows from \eqref{nav_slip_bc} and \eqref{n_tau_grad_u} as
\begin{align}
    -2(\alpha+\kappa) u_\tau &= 2 \tau \cdot \D u\ n - 2 n\cdot (\tau\cdot\nabla) u
    = \tau \cdot (n\cdot \nabla) u - n \cdot (\tau\cdot \nabla) u
    \\
    &= (\tau_i n_j - \tau_j n_i) \partial_j u_i
    = (\tau_1 n_2 - \tau_2 n_1)\partial_2 u_1 + (\tau_2 n_1 - \tau_1 n_2) \partial_1 u_2
    \\
    &= (-n_2^2-n_1^2) \partial_2 u_1 + (n_1^2 + n_2^2) \partial_1 u_2
    = \omega
\end{align}

In particular for this system, one is interested in bounds for the Nusselt number, which is defined as follows.

\begin{definition}[Nusselt Number]
    The Nusselt number is defined as
    \begin{align}
        \label{nu_definition}
        \Nu &= \langle n \cdot \nabla \theta \rangle_{\gammam},
    \end{align}
    where $\langle f \rangle_X = \limsup_{T\to\infty}\frac{1}{T}\int_0^T \frac{1}{\Gamma}\int_X f(x,t)\ dx \ dt$.
\end{definition}

\newpage
\section{Main Results}
\label{section:rb_results}

The first, most general result is the following Theorem.

\begin{theorem}\label{theorem_1_2_bound}
    Let $\Omega$ be $C^{1,1}$, $\Ra\geq 1$, $d=\min_{0\leq x_1\leq\Gamma} (\hp(x_1)-\hm(x_1)) >0$, $u_0\in L^2(\Omega)$ and $0<\alpha\in L^\infty(\gammam\cup\gammap)$. Then the Nusselt number is bounded by
    \begin{align}
        \Nu &\leq C_0 C_2^\frac{1}{2} \Ra^{\frac{1}{2}}
        \label{main_result_1_2_bound}
    \end{align}
    uniformly in $\Pr$, where $C_2(\alpha,\kappa)= 1+\|\alpha^{-1}(1+|\kappa|)\|_\infty$ and $C_0>0$ only depends on $d$, the Lipschitz constant of the domain and $\Gamma$. If $\Ra$ is sufficiently large, i.e. $\Ra\geq d^{-2}$, $C_0$ is independent of $d$.
    If additionally $|\kappa|\leq 2\alpha$, the bound improves to
    \begin{align}
        \Nu &\leq C_0 \Ra^{\frac{1}{2}}.
        \label{main_result_1_2_bound_kappa_leq_2alpha}
    \end{align}
\end{theorem}

The \hyperlink{proof_theorem_1_2_bound}{proof} of Theorem \ref{theorem_1_2_bound} is given in Section \ref{section:direct_method_application}.

We first remark that Theorem \ref{theorem_1_2_bound} holds in a broad class of applications. It also shows the same scaling as the results of \cite{doering1996VariationalConvection}, who proved
$\Nu\lesssim \Ra^\frac{1}{2}$ for no-slip boundary conditions in three spatial dimensions. In the two-dimensional setting \cite{goluskin2016Bounds} generalized this result to rough boundaries, where the height functions solely satisfy $\hm,\hp \in H^1(0,\Gamma)$. Although, in the assumptions of Theorem \ref{theorem_1_2_bound} more regularity on the boundary functions is assumed, the bound resembles the one of \cite{goluskin2016Bounds} if the slip coefficient is sufficiently large, i.e. close to no-slip.

If one further assumes more regularity, the following Theorem provides a stricter bound.

\begin{theorem}[Main Result]\label{theorem_main_theorem_rb_curved}
    Let $\Omega$ be $C^{2,1}$, $\Ra\geq 1$, $d=\min_{0\leq x_1\leq\Gamma} (\hp(x_1)-\hm(x_1)) >0$, $u_0\in W^{1,4}(\Omega)$, $0<\alpha\in W^{1,\infty}(\gammam\cup\gammap)$. Then
    \begin{align}
        \Nu
        &\leq C_0 C_1^\frac{1}{3}C_2^\frac{1}{3}\left(1+\Pr^{-\frac{1}{6}}\|u_0\|_{W^{1,4}}^\frac{1}{6} \right) \Ra^\frac{1}{3} + C_0C_2^\frac{3}{13}C_3^\frac{2}{13}\Ra^\frac{5}{13}
        \\
        &\qquad + C_0C_2^\frac{1}{4} \Ra^\frac{5}{12} + C_0 C_1^\frac{1}{3}C_2^\frac{1}{3} \Pr^{-\frac{1}{6}}\Ra^\frac{1}{2},
        \label{main_result_5_12_bound}
    \end{align}
    where
    \begin{align}
        C_1(\alpha,\kappa)&=1+\|\alpha\|_{W^{1,\infty}}+ \|\kappa\|_{W^{1,\infty}}+\|\alpha\|_\infty^3+\|\kappa\|_\infty^3
        \\
        C_2(\alpha,\kappa)&= 1+\|\alpha^{-1}(1+|\kappa|)\|_\infty
        \\
        C_3(\alpha,\kappa)&= \|\alpha+\kappa\|_\infty
    \end{align}
     and $C_0>0$ only depends on $d$, the Lipschitz constant of the domain and $\Gamma$.        
\end{theorem}

The \hyperlink{proof_theorem_main_theorem_rb_curved}{proof} of Theorem \ref{theorem_main_theorem_rb_curved} is given in Section \ref{section:direct_method_application}.

Because of the more tailored approach used in deriving the bound, i.e. using the $H^2$-norm of $u$ instead of the $H^1$-norm, one needs slightly more regularity assumptions, which in turn provide a stricter bound. These bounds yield a rich description of different scaling regimes. 

If $\Ra$ is sufficiently large, i.e. $\Ra > \|u_0\|_{W^{1,4}}$, the bound holds independent of the initial values. Due to the energy dissipation and the fact that we are investigating long-time averages, this is not surprising. The reason for the initial value to be present in the bound is that to bound the $H^2$-norm one needs bounds on the pressure gradient for the boundary conditions at hand, which in turn requires estimates on the nonlinearity, where one uses uniform in time bounds for the velocity. In fact, if one could directly obtain bounds on the long-time averages of the term arising from the nonlinearity, either in the pressure bound or when testing the gradient of \eqref{navier_stokes} with $\nabla u$, one could potentially improve the $\Ra$ exponent in the $\Pr^{-\frac{1}{6}}\Ra^\frac{1}{2}$ term.

If $\alpha\to\infty$ or $\kappa\to\infty$, the coefficients $C_1$ and $C_3$ blow up, indicating that similar bounds can not be expected to hold for no-slip boundary conditions. Similarly, if $\alpha\to 0$, the coefficient $C_2$ diverges. In contrast for $\kappa= 0$, i.e. straight boundaries, the bounds hold true. However, they are suboptimal as can be seen in the bounds derived for the flat system in Theorem \ref{theorem_flat}. There, a more detailed description of scaling laws with respect to $\alpha$ is given. Additionally, in Section \ref{section:scaling_kappa_alpha_with_Ra}, we will give an argument as to why the slip coefficient and the curvature can scale with respect to the Rayleigh number.

Assuming constant $\alpha$ and $\kappa$, the bound can be simplified to
\begin{align}
    \Nu \lessapprox \Pr^{-\frac{1}{6}}\Ra^\frac{1}{2} + \Ra^\frac{5}{12}
\end{align}
and for $\Pr\geq \Ra^\frac{1}{2}$ this recaptures the $\Nu\lesssim\Ra^\frac{5}{12}$ bound derived by \cite{whitehead2012Rigid} for the infinite Prandtl number, free-slip setting in three dimensions. Note that the authors also proved the same bound in the two-dimensional free-slip setup, where the bound holds uniform in $\Pr$ (\cite{whitehead2011Ultimate}).

The next result covers the case of flat boundaries, i.e. $\kappa=0$ and $\hm=0$, $\hp=1$. In that case, the analysis and also the bounds simplify significantly.

\begin{theorem}[Flat System]
    Let $\Ra\geq 1$, $\kappa=0$, $\alpha(x) = \frac{1}{2L_s}>0$ be constant on the boundaries and $u_0\in L^2$. Then there exists a constant $C_0>0$, only depending on $\Gamma$, such that
    \label{theorem_flat}
    \begin{align}
        \Nu &\leq C_0 \Ra^{\frac{1}{2}}
        \label{flat_theorem_bound_1_2}
    \end{align}
    holds uniformly in $\Pr$.
    
    Assume additionally $u_0\in W^{1,4}$, then
    \begin{align}
        \Nu
        &\leq C_0 L_s^{-\frac{1}{6}}\|u_0\|_{W^{1,4}}^\frac{1}{6} \Pr^{-\frac{1}{6}}\Ra^\frac{1}{3} + C_0\Ra^\frac{5}{12} +C_0 L_s^{-\frac{1}{12}}\Pr^{-\frac{1}{6}}\Ra^\frac{1}{2}
        \label{rjentnerj}
    \end{align}
    if $L_s\geq 1$ and
    \begin{align}
        \Nu
        &\leq  C_0 L_s^{-\frac{1}{3}} \Ra^\frac{1}{3} + C_0 L_s^{-\frac{1}{2}}\|u_0\|_{W^{1,4}}^\frac{1}{6}\Pr^{-\frac{1}{6}}\Ra^\frac{1}{3} + C_0 L_s^{-\frac{2}{13}} \Ra^{\frac{5}{13}}
        \\
        &\qquad + C_0 \Ra^\frac{5}{12} + C_0 L_s^{-\frac{1}{2}}\Pr^{-\frac{1}{6}}\Ra^\frac{1}{2}
        \label{jnrtkjwerj}
    \end{align}
    if $L_s\leq 1$.
\end{theorem}

The \hyperlink{proof_theorem_flat}{proof} of Theorem \ref{theorem_flat} is given in Section \ref{section:flat_theorem_proof}.
The bound in \eqref{flat_theorem_bound_1_2} is already included in Theorem \ref{theorem_1_2_bound} and we state it here for the sake of completeness.
Similarly to before, if $\Ra> \left(1+L_s^\frac{1}{2}\right) \|u_0\|_{W^{1,4}}$ the initial value terms can be absorbed in the other terms.

\begin{table}
    \hrule
    \vspace{\baselineskip}
    \begin{minipage}{\textwidth}
        \begin{center}
            \begin{tabular}{lcl}
                \multicolumn{2}{c}{Assumptions}  & \multicolumn{1}{c}{Bound}
                \\[3pt]\hline
                \multirow{2}{*}{$1\leq L_s$} & $L_s^{-\frac{1}{2}}\Ra^\frac{1}{2} \leq \Pr$ & $\Nu\lesssim\Ra^\frac{5}{12}$
                \\
                & $\Pr \leq L_s^{-\frac{1}{2}}\Ra^\frac{1}{2}$ & $\Nu\lesssim L_s^{-\frac{1}{12}}\Pr^{-\frac{1}{6}}\Ra^\frac{1}{2}$
                \\\hline
                \multirow{2}{*}{$\Ra^{-\frac{5}{24}} \leq L_s \leq 1$} & $L_s^{-3}\Ra^\frac{1}{2}\leq\Pr$ & $\Nu\lesssim \Ra^{\frac{5}{12}}$
                \\
                & $\Pr\leq L_s^{-3}\Ra^\frac{1}{2}$ & $\Nu\lesssim L_s^{-\frac{1}{2}}\Pr^{-\frac{1}{6}}\Ra^\frac{1}{2}$
                \\\hline
                \multirow{2}{*}{$\Ra^{-\frac{2}{7}} \leq L_s \leq \Ra^{-\frac{5}{24}}$} & $L_s^{-\frac{27}{13}}\Ra^\frac{9}{13}\leq\Pr$ & $\Nu\lesssim L_s^{-\frac{2}{13}}\Ra^\frac{5}{13}$
                \\
                & $\Pr\leq L_s^{-\frac{27}{13}}\Ra^\frac{9}{13}$ & $\Nu\lesssim L_s^{-\frac{1}{2}}\Pr^{-\frac{1}{6}}\Ra^\frac{1}{2}$
                \\\hline
                \multirow{2}{*}{$L_s\leq \Ra^{-\frac{2}{7}}$} & $L_s^{-1}\Ra\leq\Pr$ & $\Nu\lesssim L_s^{-\frac{1}{3}}\Ra^\frac{1}{3}$
                \\
                & $\Pr\leq L_s^{-1}\Ra$ & $\Nu\lesssim L_s^{-\frac{1}{2}}\Pr^{-\frac{1}{6}}\Ra^\frac{1}{2}$
            \end{tabular}
        \end{center}
    \end{minipage}
    \begin{center}
        \caption{Overview of the scaling laws in Theorem \ref{theorem_flat} for $\Ra\geq (1+L_s^\frac{1}{2}) \|u_0\|_{W^{1,4}}$ in the regimes $L_s(\Ra)$, $\Pr(L_s,\Ra)$.}
        \label{Table:flat_result_Pr_depends_on_Ls_Ra}
        \vspace{\baselineskip}
        \hrule
    \end{center}
\end{table}

Note that the results in \cite{bleitner2024Scaling} read
\begin{align}
    \Nu \lesssim \Ra^\frac{5}{12} + L_s^{-\frac{1}{6}} \Pr^{-\frac{1}{6}}\Ra^\frac{1}{2}
\end{align}
if $L_s\geq 1$ and
\begin{align}
    \Nu \lesssim L_s^{-\frac{1}{3}} \Ra^\frac{1}{3} + L_s^{-\frac{2}{13}}\Ra^\frac{5}{13} + \Ra^\frac{5}{12} + L_s^{-\frac{2}{3}} \Pr^{-\frac{1}{6}} \Ra^\frac{1}{2}
\end{align}
if $L_s\leq 1$. The results only differ slightly in the $L_s$ exponent in the $\Pr$ term. The reason for the different estimates will be discussed in Remark \ref{Remark:difference_pressure_bound}. Accordingly, using $p>4$ in \eqref{difference_to_arxiv_referd_in_explanation} would further alter this exponent.

If $L_s\geq 1$, \eqref{rjentnerj} shows that if either $\Pr\to\infty$ or $L_s\to\infty$ one recaptures the previously mentioned $\Nu\lesssim \Ra^\frac{5}{12}$ results of \cite{whitehead2011Ultimate} and \cite{whitehead2012Rigid}. 

Assuming $L_s$ scales with $\Ra$, i.e. $L_s = \Ra^{\delta}$ for some $\delta\geq 0$ in order to match $L_s\geq 1$, (an argument why this is reasonable will be given in Section \ref{section:scaling_kappa_alpha_with_Ra}) the bound reads 
\begin{alignat}{2}
    \Nu &\lesssim \Ra^{\frac{5}{12}} & \textnormal{ if }& \Pr\geq \Ra^{\frac{1}{2}-\delta}
    \\
    \Nu &\lesssim \Pr^{-\frac{1}{6}}\Ra^{\frac{1}{2}-\frac{\delta}{6}} & \textnormal{ if }& \Pr\leq \Ra^{\frac{1}{2}-\delta}
\end{alignat}
A similar crossover behavior was shown in \cite{choffrut2016Upper} for the no-slip setup, where, up to logarithmic corrections, the bounds change from $\Nu\lesssim \Ra^\frac{1}{3}$ when $\Pr\geq \Ra^\frac{1}{3}$ to $\Nu\lesssim \Pr^{-\frac{1}{2}}\Ra^\frac{1}{2}$ when $\Pr\leq \Ra^\frac{1}{3}$.

If instead $L_s \leq 1$, the scaling with respect to $L_s$ is not that obvious. Formally setting $\Pr = \infty$ and looking at the no-slip limit \eqref{jnrtkjwerj} indicates a $\Nu\sim \Ra^\frac{1}{3}$ scaling, which would, again up to logarithmic corrections, match results of \cite{constantin1999Infinite,doering2006Bounds,otto2011Rayleigh}. Note though, that this argument is not rigorous and one could make any term dominate in that limit by multiplying it with $L_s^{-\frac{1}{3}}$.

The different regimes of both \eqref{rjentnerj} and \eqref{jnrtkjwerj} are shown in Tables \ref{Table:flat_result_Pr_depends_on_Ls_Ra} and \ref{Table:flat_result_Ls_depends_on_Pr_Ra}. Note that Table \ref{Table:flat_result_Ls_depends_on_Pr_Ra} corresponds to defining the ranges of $L_s$ with respect to $\Pr$ while Table \ref{Table:flat_result_Pr_depends_on_Ls_Ra} does the opposite.

Assuming again that $L_s$ scales with $\Ra$, as argued in Section \ref{section:scaling_kappa_alpha_with_Ra}, then Table \ref{Table:flat_result_Pr_depends_on_Ls_Ra} shows the corresponding scaling law. Table \ref{Table:flat_result_Ls_depends_on_Pr_Ra} corresponds to the more direct approach of having a specific fluid in mind, i.e. fixing $\Pr$ and then investigating how the scaling law changes when varying the slip length by for example adding a lubricant or changing the boundary material.

\begin{table}
    \newcommand{\colorTableBig}{rtg_darkblue}
    \newcommand{\colorTableBoth}{black}
    \newcommand{\colorTableSmall}{rtg_red}
    \hrule
    \vspace{\baselineskip}
    \begin{minipage}{\textwidth}
        \begin{center}
            \begin{tabular}{lcl}
                \multicolumn{2}{c}{Assumptions}  & \multicolumn{1}{c}{Bound}
                \\[3pt]\hline
                & \textcolor{\colorTableBoth}{$\Ra^{-\frac{5}{24}}\leq L_s$} &\textcolor{\colorTableBoth}{$\Nu \lesssim \Ra^{\frac{5}{12}}$}
                \\
                \multirow{2}{*}{$\Ra^\frac{9}{7}\leq \Pr$} & \textcolor{\colorTableSmall}{$\Ra^{-\frac{2}{7}}\leq L_s\leq \Ra^{-\frac{5}{24}}$} & \textcolor{\colorTableSmall}{$\Nu \lesssim L_s^{-\frac{2}{13}}\Ra^{\frac{5}{13}}$}
                \\
                & \textcolor{\colorTableSmall}{$\Pr^{-1}\Ra\leq L_s\leq \Ra^{-\frac{2}{7}}$} &\textcolor{\colorTableSmall}{$\Nu \lesssim L_s^{-\frac{1}{3}}\Ra^{\frac{1}{3}}$}
                \\
                & \textcolor{\colorTableSmall}{$L_s\leq \Pr^{-1} \Ra$} & \textcolor{\colorTableSmall}{$\Nu \lesssim L_s^{-\frac{1}{2}}\Pr^{-\frac{1}{6}}\Ra^\frac{1}{2}$}
                \\\hline
                & \textcolor{\colorTableBoth}{$\Ra^{-\frac{5}{24}}\leq L_s$} & \textcolor{\colorTableBoth}{$\Nu \lesssim \Ra^{\frac{5}{12}}$}
                \\
                $\Ra^{\frac{8}{9}}\leq \Pr\leq \Ra^\frac{9}{7}$ & \textcolor{\colorTableSmall}{$\Pr^{-\frac{13}{27}}\Ra^{\frac{1}{3}}\leq L_s\leq \Ra^{-\frac{5}{24}}$} & \textcolor{\colorTableSmall}{$\Nu \lesssim L_s^{-\frac{2}{13}}\Ra^{\frac{5}{13}}$}
                \\
                & \textcolor{\colorTableSmall}{$L_s\leq \Pr^{-\frac{13}{27}} \Ra ^{\frac{1}{3}}$} & \textcolor{\colorTableSmall}{$\Nu \lesssim L_s^{-\frac{1}{2}}\Pr^{-\frac{1}{6}}\Ra^\frac{1}{2}$}
                \\\hline
                \multirow{2}{*}{$\Ra^\frac{1}{2}\leq \Pr\leq \Ra^\frac{9}{8}$} & \textcolor{\colorTableBoth}{$\Pr^{-\frac{1}{3}}\Ra^{\frac{1}{6}}\leq L_s$} & \textcolor{\colorTableBoth}{$\Nu \lesssim \Ra^{\frac{5}{12}}$}
                \\
                & \textcolor{\colorTableSmall}{$L_s\leq \Pr^{-\frac{1}{3}} \Ra ^{\frac{1}{6}}$} & \textcolor{\colorTableSmall}{$\Nu \lesssim L_s^{-\frac{1}{2}}\Pr^{-\frac{1}{6}}\Ra^\frac{1}{2}$}
                \\\hline
                & \textcolor{\colorTableBig}{$\Pr^{-2}\Ra\leq L_s$} & \textcolor{\colorTableBig}{$\Nu \lesssim \Ra^{\frac{5}{12}}$}
                \\
                $\Pr\leq \Ra^\frac{1}{2}$& \textcolor{\colorTableBig}{$1\leq L_s \leq \Pr^{-2}\Ra$} & \textcolor{\colorTableBig}{$\Nu \lesssim L_s^{-\frac{1}{12}}\Pr^{-\frac{1}{6}}\Ra^{\frac{1}{2}}$}
                \\
                & \textcolor{\colorTableSmall}{$L_s\leq 1$} & \textcolor{\colorTableSmall}{$\Nu \lesssim L_s^{-\frac{1}{2}}\Pr^{-\frac{1}{6}}\Ra^\frac{1}{2}$}
            \end{tabular}
        \end{center}
    \end{minipage}
    \begin{center}
        \vspace{-\baselineskip}
        \caption{Overview of the scaling laws in Theorem \ref{theorem_flat} for $\Ra\geq (1+L_s^\frac{1}{2}) \|u_0\|_{W^{1,4}}$ in the regimes $\Pr(\Ra)$, $L_s(\Pr,\Ra)$. The color reflects the slip length cases, i.e. \textcolor{\colorTableSmall}{$L_s\leq 1$}, \textcolor{\colorTableBig}{$1\leq L_s$} and uncolored if there exist $L_s\leq 1$ and $L_s\geq 1$ satisfying the assumptions.}
        \label{Table:flat_result_Ls_depends_on_Pr_Ra}
        \vspace{\baselineskip}
        \hrule
    \end{center}
\end{table}

The final results cover the system with identical, potentially curved boundaries. These findings are published in \cite{bleitnerNobili2024Bounds}. For this system, the most general result is given in the following Theorem.

\begin{theorem}
    Let $\Omega$ be $C^{1,1}$ with $\hp=\hm+1$, $u_0\in L^2(\Omega)$, $0<\alpha\in L^\infty(\gammam\cup\gammap)$. Then
    \label{theorem_nonlinear_1_2_bound}
    \begin{align}
        \Nu \lesssim \Ra^\frac{1}{2} + \|\kappa\|_\infty
    \end{align}
    if $|\kappa|\leq 2\alpha$ and
    \begin{align}
        \Nu \lesssim \left(1+\|\alpha^{-1}\|_\infty^\frac{1}{2}\right)\Ra^\frac{1}{2} + \|\kappa\|_\infty
    \end{align}
    if $|\kappa|\leq 2\alpha + \frac{1}{4\sqrt{1+(h'(x_1))^2}}\min \left\{1,\sqrt{\alpha}\right\}$. In both cases, the implicit constant only depends on $\Gamma$ and the Lipschitz constant of the boundary.
\end{theorem}

The \hyperlink{proof_theorem_nonlinear_1_2_bound}{proof} of Theorem \ref{theorem_nonlinear_1_2_bound} is given in Section \ref{section:nonlin_regularity_estimates}.

Similarly to before, assuming more regularity one can prove a refined statement given below.

\begin{theorem}
    \label{theorem_nonlinear_main_theorem}
    Let $\Omega$ be $C^{2,1}$, $\hp=1+\hm$, $u_0\in W^{1,4}(\Omega)$ and $0<\alpha \in W^{1,\infty}(\gammam\cup\gammap)$. Then there exists a constant $0<\bar C<1$ such that for all $\alpha$ and $\kappa$ satisfying
    \begin{align}
        \|\alpha+\kappa\|_\infty \leq \bar C
    \end{align}
    one has
    \begin{itemize}
        \item
            for $|\kappa|\leq \alpha$ on $\gammam\cup\gammap$, $\Pr\geq \|\alpha^{-1}\|_\infty^\frac{3}{2}\Ra^\frac{3}{4}$ and $\|\alpha^{-1}\|_\infty \leq \Ra$
            \begin{align}
                \Nu \lesssim \|\alpha+\kappa\|_{W^{1,\infty}}^2 \Ra^\frac{1}{2}+ C_1 \Ra^\frac{5}{12},
            \end{align}
        \item
            for $|\kappa|\leq 2\alpha +\frac{1}{4}\sqrt{\frac{\alpha}{1+(h')^2}}$, $\Pr\geq \|\alpha^{-1}\|_\infty^\frac{5}{4}\Ra^\frac{3}{4}$ and $\|\alpha^{-1}\|_\infty\leq\Ra$
            \begin{align}
                \Nu \lesssim \|\alpha^{-1}\|_\infty^{-\frac{1}{2}}\|\alpha+\kappa\|_{W^{1,\infty}}^2 \Ra^\frac{1}{2} + C_1 \|\alpha^{-1}\|_\infty^\frac{1}{12}\Ra^\frac{5}{12}
            \end{align}
        \item
            for $|\kappa|\leq 2\alpha +\frac{1}{4}\sqrt{\frac{\alpha}{1+(h')^2}}$ and $\Pr\geq \Ra^\frac{5}{7}$
            \begin{align}
                \Nu \lesssim C_2 \Ra^\frac{3}{7}
            \end{align}
    \end{itemize}
    where 
    \begin{align}
        C_1(u_0,\alpha,\kappa) &= C\left(\|u_0\|_{W^{1,4}}^\frac{1}{3} + \|\alpha\|_{W^{1,\infty}}^\frac{1}{3}+\|\kappa\|_{W^{1,\infty}}^\frac{2}{3} + 1\right),
        \\
        C_2(u_0,\alpha,\kappa) &= \|\alpha^{-1}\|_\infty^{-3}\|\alpha+\kappa\|_{W^{1,\infty}}^2 + \|\alpha^{-1}\|_\infty^\frac{1}{2} + C_1(u_0,\alpha,\kappa)\|\alpha^{-1}\|_\infty^\frac{1}{6}.
    \end{align}
    and the implicit constant only depends on $\Gamma$, the Lipschitz constant of the domain.
\end{theorem}

The \hyperlink{proof_theorem_nonlinear_main_theorem}{proof} of Theorem \ref{theorem_nonlinear_main_theorem} is given in Section \ref{section:background_field_method_application}.

The bounds of Theorem \ref{theorem_nonlinear_1_2_bound} and \ref{theorem_nonlinear_main_theorem} demand more assumptions and yield less strict bounds than the ones given in Theorem \ref{theorem_1_2_bound} and \ref{theorem_main_theorem_rb_curved}, with the one exception that the Theorem \ref{theorem_nonlinear_main_theorem} yields a better result if $\alpha\geq |\kappa|$ in the limit $\|\alpha+\kappa\|_{W^{1,\infty}}\to 0$. As the curvature vanishes in this limit, Theorem \ref{theorem_flat} yields a more strict description of the scaling behavior. Therefore, we refrain from interpreting the results further.

An overview of the results is given in Table \ref{Table:overview_nusselt_bounds}. As indicated in the table, to the best of the author's knowledge no result for free-slip boundary conditions on curved domains is known. The approach given in Theorem \ref{theorem_main_theorem_rb_curved} fails in the limit $\alpha \to 0$, due to the lack of control of the energy. This can be seen in the energy balance
\begin{align}
    \frac{1}{2\Pr}\frac{d}{dt} \|u\|_{L^2}^2 &+ \|\nabla u\|_{L^2}^2 + \int_{\gammam\cup\gammap} (2\alpha+\kappa) u_\tau^2
    \\
    &= \frac{1}{2\Pr}\frac{d}{dt} \|u\|_{L^2}^2 + 2 \|\D u\|_{L^2}^2 + 2 \int_{\gammam\cup\gammap} \alpha u_\tau^2
    \\
    &= \Ra \int_{\Omega} u_2 \theta,
    \label{rjwekrjnwe}
\end{align}
proven in \eqref{energy_balance}. If $\alpha = 0$, then the convex portions of the boundaries imply $\kappa<0$, and therefore, no bound on the energy can be inferred from \eqref{rjwekrjnwe}, at least directly. Note, that if the problem was set on a non-periodic domain $\Omega$, then the computation
\begin{align}
    \int_\Omega u = \int_{\Omega} 1\!\!1 \ u = \int_{\Omega}  u \cdot \nabla x = \int_{\partial\Omega} u\cdot n x - \int_{\Omega} x \nabla \cdot u = 0,
    \label{huweriwehui}
\end{align}
where $1\!\!1$ is the identity matrix, shows that the incompressibility and no-pene\-tration conditions are sufficient to show that $u$ is average free. Therefore Poincaré's inequality together with \eqref{rjwekrjnwe} would yield control over the energy. In the partially periodic domain, the computation \eqref{huweriwehui} only works for the second component, i.e. showing $u_2$ has vanishing average.

Finally, we want to discuss the different methods used to prove the above results. Theorems \ref{theorem_1_2_bound}, \ref{theorem_main_theorem_rb_curved}, \ref{theorem_flat} and \ref{theorem_nonlinear_1_2_bound} are proved using the direct method introduced in \cite{seis2015Scaling}. In this method, one localizes the Nusselt number in a strip of width $\delta$ at the boundary. Using the boundary conditions and the fundamental theorem of calculus, together with the long time bounds for the corresponding norms of the solution yields estimates depending on $\delta$. The desired result is then obtained after optimizing in $\delta$.

In contrast, the Background field method, described in \cite{doering1994VariationalShear,doering1996VariationalConvection}, is used to prove Theorem \ref{theorem_nonlinear_main_theorem}. In this method, one splits the temperature into a steady profile and fluctuations around it. Here this profile is a piecewise linear function in $x_2$ with slope $\delta^{-1}$ near the boundary, while in $x_1$ it matches the boundary height functions. Defining a variational quadratic functional, that is motivated by the long time bounds, results in a variational problem that when solved determines $\delta$ and therefore the scaling law for the Nusselt number.

As indicated by the same choice of parameter $\delta$, the methods are strongly connected. In fact \cite{chernyshenko2022Relationship} showed that the background field method is a special case of the auxiliary functional method. The study shows that bounds derived by the auxiliary functional method can be proven using the direct method, while under certain assumptions the converse is also true.

\begin{table}
    \hrule
    \vspace{\baselineskip}
    \begin{minipage}{\textwidth}
        \begin{center}
            \resizebox{0.97\columnwidth}{!}{%
                \begin{tabular}{l|lll}
                    \phantom{
                        \footnote{\label{citeChoffrut}\cite{choffrut2016Upper}}\footnote{\label{citeTheorem_flat}Theorem \ref{theorem_flat}, \cite{bleitner2024Scaling}}
                        \footnote{\label{citeWhiteheadUltimate}\cite{whitehead2011Ultimate}}
                        \footnote{\label{citeDoering1996VariationalConvection}\cite{doering1996VariationalConvection}}
                    }
                    & no-slip & Navier-slip & free-slip
                    \\[3pt]\hline
                    Flat & \multirow{2}{*}{$\Nu\lesssim \Ra^\frac{1}{3}(\ln \Ra)^\frac{1}{3}$\footref{citeChoffrut}} & \multirow{2}{*}{$\Nu\lesssim\Ra^\frac{5}{12}$\footref{citeTheorem_flat}} & \multirow{2}{*}{$\Nu\lesssim \Ra^\frac{5}{12}$\footref{citeWhiteheadUltimate}}
                    \\
                    $\Pr\geq \Ra^\mu$
                    \\\hline
                    \multirow{2}{*}{Flat} & $\Nu\lesssim \Pr^{-\frac{1}{2}}\Ra^\frac{1}{2}(\ln \Ra)^\frac{1}{2}$\footref{citeChoffrut} & $\Nu\lesssim \Pr^{-\frac{1}{6}}\Ra^{\frac{1}{2}}+\Ra^\frac{5}{12}$\footref{citeTheorem_flat} & \multirow{2}{*}{$\Nu\lesssim \Ra^\frac{5}{12}$\footref{citeWhiteheadUltimate}}
                    \\
                    & $\Nu\lesssim \Ra^\frac{1}{2}$\footref{citeDoering1996VariationalConvection} & $\Nu\lesssim \Ra^{\frac{1}{2}}$\footref{citeTheorem_flat}\footref{citeDrivas} &
                    \\\hline
                    \multirow{2}{*}{Cuved} & \multirow{2}{*}{$\Nu\lesssim \Ra^\frac{1}{2}$\footref{citeGoluskinBounds}} & $\Nu\lesssim \Pr^{-\frac{1}{6}}\Ra^{\frac{1}{2}}+\Ra^\frac{5}{12}$\footref{citeTheorem_curved} & 
                    \\
                    & & $\Nu\lesssim \Ra^{\frac{1}{2}}$\footref{citeTheorem_12bound} & 
                    \phantom{
                        \footnote{\label{citeDrivas}\cite{drivas2022Bounds}}
                        \footnote{\label{citeGoluskinBounds}\cite{goluskin2016Bounds}}
                        \footnote{\label{citeTheorem_curved}Theorem \ref{theorem_main_theorem_rb_curved}}
                        \footnote{\label{citeTheorem_12bound}Theorem \ref{theorem_1_2_bound}, \cite{bleitnerNobili2024Bounds}}
                    }
                \end{tabular}
            }
        \end{center}
    \end{minipage}
    \vspace{-\baselineskip}
    \begin{center}
        \caption{Overview of Nusselt number bounds for the specific systems. For Navier-slip the results are given for some fixed slip coefficient.}
        \label{Table:overview_nusselt_bounds}
        \vspace{\baselineskip}
        \hrule
    \end{center}
\end{table}

\section{Scaling of the Curvature and Slip Coefficient}
\label{section:scaling_kappa_alpha_with_Ra}

This section is devoted to providing an argument (see \cite[A.2]{bleitnerNobili2024Bounds}) as to why the curvature and slip coefficient might scale with respect to the Rayleigh number. 

At first, we focus on the curvature. The original system, i.e. the system that is not non-dimensionalized yet, can be given by
\begin{align}
    u_t + u\cdot\nabla u + \nabla \frac{p}{\rho_0}-\nu \Delta u &= -\varsigma g (T-T^-)
    \\
    \partial_t (T-T^-) + u\cdot\nabla(T-T^-)-\varkappa \Delta (T-T^-) &= 0
\end{align}
where $\nu$ is the kinematic viscosity, $\varkappa$ is the thermal diffusivity, $g$ is the gravitational constant, $\varsigma$ is the thermal expansion coefficient. Additionally, the average height gap is denoted by $d$ and $T^-$ and $T^+$ is the fixed temperature on the lower, respectively upper boundary in the original system, which are described by the boundary height functions $\hm(x_1)$ and $\hp(x_1)$. Non-dimensionalizing by introducing the variables
\begin{align}
    \hat x &= \frac{1}{d}x, & \hat t &= \frac{\varkappa}{d^2}t,& \hat u &= \frac{d}{\varkappa}u,& \hat \theta &= \frac{T-T^-}{T^+-T^-}, & \hat p&= \frac{d^2}{\kappa^2 \rho_0}p,
\end{align}
the system is given by
\begin{align}
    \Pr^{-1}(\hat u_t + \hat u \cdot\nabla \hat u)+ \nabla\hat p -\Delta \hat u &= \Ra \hat\theta e_2
    \\
    \hat \theta_t + \hat u\cdot\nabla \hat \theta- \Delta \hat \theta &= 0,
\end{align}
where
\begin{align}
    \Pr &= \frac{\nu}{\varkappa} & \Ra &= \frac{g \varsigma d^3(T^--T^+)}{\varkappa \nu}.
    \label{jknkjrtnzhjrtbzrtjh}
\end{align}
Accordingly, the boundary height functions $\hm$ and $\hp$ are rescaled as
\begin{align}
    \hat h^+ &= \frac{1}{d} \hp, & \hat h^- &= \frac{1}{d} \hm
\end{align}
Note, that varying the temperature gap $\delta T=T^+-T^-$ of the original system changes the Rayleigh number, while $\hat h^+$ and $\hat h^-$ stay the same. On the other hand, varying the average boundary gap $d$ leads to a change in the non-dimensional boundary height functions $\hat h^+$ and $\hat h^-$, without changing the profile in the original system. As
\begin{align}
    \hat \kappa \sim \hat h^{\pm\prime \prime}
\end{align}
on the respective boundaries, the scaling yields
\begin{align}
    \hat \kappa &= d \kappa,& \hat \kappa' &= d^2 \kappa, & \|\hat\kappa\|_{W^{1,\infty}}\sim d^2 + d.
    \label{werjkwerjwen}
\end{align}
Comparing two systems with different boundary temperature gaps $\delta T_1$ and $\delta T_2$ and boundary height gaps $d_1$ and $d_2$, that otherwise coincide, one can achieve any desired scaling of $\|\kappa\|_{W^{1,\infty}}$ with respect to $\Ra$, i.e.
\begin{align}
    \frac{\|\kappa_2\|_{W^{1,\infty}}}{\|\kappa_1\|_{W^{1,\infty}}} \approx \left(\frac{\Ra_2}{\Ra_1}\right)^\rho
    \label{njbbngkj}
\end{align}
for any $\rho \in \mathbb{R}$. In fact, setting
\begin{align}
    \frac{d_2}{d_1} = \left(\frac{\delta T_2}{\delta T_1}\right)^{\frac{\rho}{2-3\rho}}
\end{align}
and using \eqref{werjkwerjwen} one has
\begin{align}
    \frac{\|\kappa_2\|_{W^{1,\infty}}}{\|\kappa_1\|_{W^{1,\infty}}} &\approx \frac{d_2^2}{d_1^2} = \left(\frac{d_2}{d_1}\right)^{3\rho} \left(\frac{d_2}{d_1}\right)^{2-3\rho} = \left(\frac{d_2}{d_1}\right)^{3\rho} \left(\frac{\delta T_2}{\delta T_1}\right)^{\rho}
    \\
    &= \left(\frac{d_2^3\delta T_2}{d_1^3\delta T_1}\right)^{\rho} = \left(\frac{\Ra_2}{\Ra_1}\right)^\rho
\end{align}
for sufficiently large boundary height gaps $d_i\gg 1$, where in the last inequality we used \eqref{jknkjrtnzhjrtbzrtjh}. The cases $\rho =0$ and $\rho = \frac{2}{3}$ are covered with  $d_2=d_1$, respectively $\delta T_2=\delta T_1$, concluding the argument for \eqref{njbbngkj}.

In order to show the scaling of $\alpha$ with respect to $\Ra$, we note that the slip coefficient can be derived as a highly oscillating, small amplitude limit of the boundary roughness (\cite{miksis1994Slip}), where it is proportional to the average height function of the roughness, i.e. $\hat \alpha \sim \hat h^\pm$. Therefore, an analogous argument for the curvature shows the scaling.

\section{Nusselt Number Representations}
\label{Section:nusselt_number_representations}

\begin{lemma}
    \label{lemma_nusselt_number_representations}
    Let $\Omega$ be $C^{1,1}$ and $u_0\in L^2(\Omega)$. Then for any $0\leq z < \min \hp - \max \hm$
    \begin{align}
        \Nu &= \langle \np \cdot (u-\nabla)\theta\rangle_{\gammam + z} \label{nu_representation_shifted}
        \\
        &= \langle \|\nabla \theta\|_2^2\rangle_{\Omega} \label{nu_representation_gradient}
        \\
        &\geq (\max \hp - \min \hm)^{-1} \langle (u_2-\partial_2)\theta \rangle_{\Omega}, \label{nu_representation_straight}
    \end{align}
    where $\gammam + z = \lbrace (x_1,x_2)\ \vert \ 0\leq x_1\leq \Gamma, x_2 = \hm(x_1) + z\rbrace$ is the shifted bottom boundary.
\end{lemma}

\begin{proof}
    \vspace{-\baselineskip}
    \begin{itemize}
        \item Argument for \eqref{nu_representation_shifted}.

            Without loss of generality let $z>0$ as for $z=0$ it matches the definition \eqref{nu_definition} because of \eqref{no_penetration_bc}. Define $\Omegatilde = \lbrace (x_1,x_2)\ \vert \ 0 < x_1 < \Gamma, \hm(x_1) < x_2 < \hm(x_1) + z\rbrace$ as illustrated Figure \ref{fig:nu_representation}. Then by \eqref{advection_diffusion_equation}, \eqref{incompressibility} and Stokes' theorem
            \begin{align}
                \frac{d}{dt} \int_{\Omegatilde} \theta &= - \int_{\Omegatilde} (u\cdot \nabla \theta - \Delta \theta)
                = - \int_{\Omegatilde} \nabla \cdot ((u-\nabla)\theta)
                \\
                &= - \int_{\partial\Omegatilde} n\cdot (u-\nabla)\theta
                = - \int_{\gammam + z} \np \cdot (u-\nabla)\theta + \int_{\gammam} n\cdot \nabla \theta,
            \end{align}
            where in the last equality we used \eqref{no_penetration_bc}. Taking the long time average and using \eqref{maximum_principle} we find
            \begin{align}
                0 &= \limsup_{T\to \infty}T^{-1}|\Omega|^{-1} \left( \int_\Omegatilde \theta(x,T)  \ dx - \int_\Omegatilde \theta_0(x) \ dx\right)
                \\
                &= \limsup_{T\to \infty}\dashint_0^T |\Omega|^{-1} \frac{d}{dt} \int_\Omegatilde \theta \ dx \ dt
                \\
                &= -\limsup_{T\to \infty}\dashint_0^T |\Omega|^{-1} \left( \int_{\gammam + z} \np \cdot (u-\nabla)\theta \ dx + \int_{\gammam} n\cdot \nabla \theta \ dx\right) \ dt
                \\
                &= - \langle \np \cdot (u-\nabla)\theta \rangle_{\gammam + z} + \Nu.
            \end{align}
            
    \begin{figure}
        \hrule
        \vspace{0.5\baselineskip}
        \begin{center}
            \includegraphics[width=\textwidth]{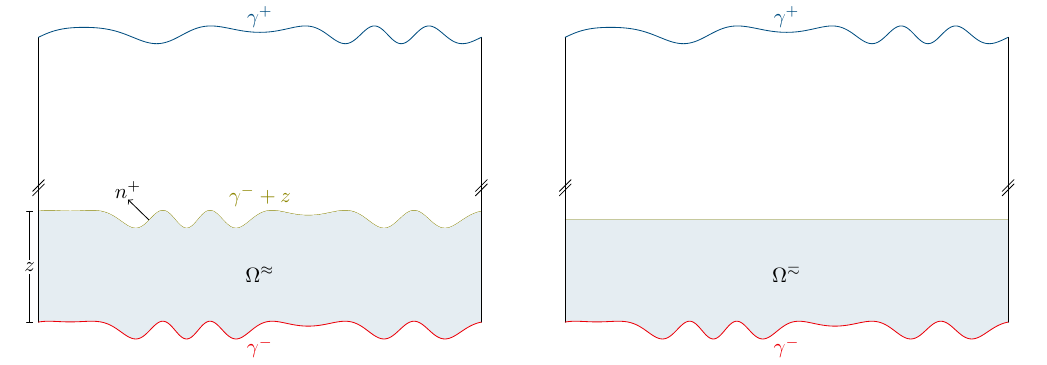}
            \vspace*{-7mm}
            \caption{Illustration of the definitions of $\Omegatilde$ and $\Omegabar$.}
            \label{fig:nu_representation}
        \end{center}
        \hrule
    \end{figure}
        \item Argument for \eqref{nu_representation_gradient}

            Testing \eqref{advection_diffusion_equation} with $\theta$ and using integration by parts we find
            \begin{align}
                \frac{1}{2}\frac{d}{dt}\|\theta\|_{L^2}^2 &= - \int_\Omega \theta u\cdot \nabla \theta + \int_\Omega \theta \Delta \theta
                \\
                &= - \int_\Omega \theta u\cdot \nabla \theta + \int_{\gammam\cup\gammap} \theta n\cdot \nabla \theta - \|\nabla \theta\|_{L^2}^2 \label{nbsidb}
            \end{align}
            Note that the first term on the right-hand side of \eqref{nbsidb} vanishes by \eqref{incompressibility}, \eqref{no_penetration_bc} and integration by parts as
            \begin{align}
                - 2 \int_\Omega \theta u\cdot \nabla \theta &= - \int_\Omega u\cdot \nabla \theta^2 = - \int_ {\gammam\cup\gammap} u\cdot n \theta^2 + \int_\Omega \nabla \cdot u \theta^2 = 0.\\[-10pt]\label{alnbjgns}
            \end{align}
            Taking the long time average of \eqref{nbsidb}, using \eqref{maximum_principle}, \eqref{heat_bc} and \eqref{alnbjgns}, yields
            \begin{align}
                0 &= \limsup_{T\to\infty} (2T|\Omega|)^{-1} \left( \|\theta(T)\|_{L^2}^2-\|\theta_0\|_{L^2}^2\right)
                \\
                &= \limsup_{T\to\infty} |\Omega|^{-1} \dashint_0^T \frac{1}{2}\frac{d}{dt} \|\theta\|_{L^2}^2 \ dt
                \\
                &= \limsup_{T\to\infty} |\Omega|^{-1} \dashint_0^T \left(\int_{\gammam\cup\gammap} \theta n\cdot \nabla \theta - \|\nabla \theta\|_{L^2}^2 \right)\ dt
                \\
                &= \langle n\cdot \nabla \theta\rangle_{\gammam} - \langle\|\nabla \theta\|_{L^2}^2 \rangle_\Omega,
            \end{align}
            which according to definition \eqref{nu_definition} proves the claim.
        \item Argument for \eqref{nu_representation_straight}.
            For $\min \hm \leq y \leq \max \hp$ define $\Omegabar(y) = \lbrace (x_1,x_2) \in \Omega \ \vert \ x_2<y \rbrace$. Then similar to the proof of \eqref{nu_representation_shifted} one finds
            \begin{align}
                \frac{d}{dt}\int_{\Omegabar(y)} \theta &= - \int_{\Omegabar(y)} (u\cdot \nabla \theta -\Delta \theta)
                = -\int_{\Omegabar(y)} \nabla \cdot ((u-\nabla)\theta)
                \\
                &= -\int_{\partial\Omegabar(y)} n\cdot (u-\nabla)\theta
                \\
                &= \int_{\gammap \cap \lbrace x_2 \leq y \rbrace} \np \cdot \nabla\theta + \int_{\gammam \cap \lbrace x_2 \leq y \rbrace} \nm \cdot \nabla\theta
                \\
                &\qquad - \int_{\Omega \cap \lbrace x_2 = y \rbrace} (u_2-\partial_2)\theta,\label{bsniadsf}
            \end{align}
            where we split the non-canceling parts of the integral into three relevant boundary terms and used \eqref{no_penetration_bc}. As for $x\in \Omega$ we have $0\leq \theta(x) \leq 1$ by \eqref{maximum_principle} and $\theta = 0$ on $\gammap$ one has $\np\cdot \nabla \theta \leq 0$ on $\gammap$ and therefore
            \begin{align}
                \int_{\gammap \cap \lbrace x_2 \leq y \rbrace} \np \cdot \nabla\theta \leq 0,
            \end{align}
            implying
            \begin{align}
                \frac{d}{dt}\int_{\Omegabar(y)} \theta &\leq \int_{\gammam \cap \lbrace x_2 \leq y \rbrace} \nm \cdot \nabla\theta - \int_{\Omega \cap \lbrace x_2 = y \rbrace} (u_2-\partial_2)\theta.\label{bnugfifgbn}
            \end{align}
            Integrating \eqref{bnugfifgbn} in $y$ yields
            \begin{align}
                \frac{d}{dt}\int_{\min \hm}^{\max \hp}&\int_{\Omegabar(y)} \theta \ dx \ dy
                \\
                &\leq \int_{\min \hm}^{\max \hp}\int_{\gammam \cap \lbrace x_2 \leq y \rbrace} \nm \cdot \nabla\theta \ dx \ dy
                \\
                &\qquad - \int_{\min \hm}^{\max \hp}\int_{\Omega \cap \lbrace x_2 = y \rbrace} (u_2-\partial_2)\theta \ dx \ dy
                \\
                &= \int_{\min \hm}^{\max \hp}\int_{\gammam} \nm \cdot \nabla\theta \ dx \ dy
                \\
                &\qquad - \int_{\min \hm}^{\max \hp}\int_{\gammam \cap \lbrace x_2 > y \rbrace} \nm \cdot \nabla\theta \ dx \ dy 
                \\
                &\qquad - \int_{\min \hm}^{\max \hp}\int_{\Omega \cap \lbrace x_2 = y \rbrace} (u_2-\partial_2)\theta \ dx \ dy
                \\
                &= (\max\hp-\min\hm)\int_{\gammam} \nm \cdot \nabla\theta \ dS
                \\
                &\qquad - \int_{\min \hm}^{\max \hp}\int_{\gammam \cap \lbrace x_2 > y \rbrace} \nm \cdot \nabla\theta \ dx \ dy 
                \\
                &\qquad - \int_{\Omega} (u_2-\partial_2)\theta \ dx \label{nsiubfgn}
            \end{align}
            Similar to the top boundary, since $\theta = 1$ on $\gammam$ and $0\leq \theta(x) \leq 1$ for $x\in \Omega$, one has $\nm\cdot \nabla \theta\geq 0$ as $\nm$ is the outwards pointing unit normal vector. Therefore \eqref{nsiubfgn} yields
            \begin{align}
                \frac{d}{dt}\int_{\min \hm}^{\max \hp}&\int_{\Omegabar(y)} \theta \ dx \ dy
                \\
                &\leq (\max\hp-\min\hm)\int_{\gammam} \nm \cdot \nabla\theta - \int_{\Omega} (u_2-\partial_2)\theta \ dx.
            \end{align}
            Again taking the long time average and and using \eqref{maximum_principle}
            \begin{align}
                0 &= \limsup_{T\to\infty} (2T|\Omega|)^{-1} \int_{\min \hm}^{\max \hp}\int_{\Omegabar(y)} (\theta(x,T)-\theta_0(x)) \ dx \ dy
                \\
                &= \limsup_{T\to\infty} |\Omega|^{-1}  \dashint_0^T \frac{d}{dt}\int_{\min \hm}^{\max \hp}\int_{\Omegabar(y)} \theta \ dx \ dy\ dt
                \\
                &\leq \limsup_{T\to\infty} |\Omega|^{-1}  \dashint_0^T \bigg((\max\hp-\min\hm)\int_{\gammam} \nm \cdot \nabla\theta \ dS
                \\
                &\qquad - \int_{\Omega} (u_2-\partial_2)\theta \ dx\bigg) \ dt
                \\
                &= (\max\hp-\min\hm) \langle \nm \cdot \nabla \theta\rangle_{\gamma^-} - \langle (u_2-\partial_2)\theta\rangle_{\Omega},
            \end{align}
            which by the definition \eqref{nu_definition} yields the claim.        
    \end{itemize}
\end{proof}

\section{General System}
\label{section:general_system}

This section is devoted to proving the results given in Theorem \ref{theorem_1_2_bound} and Theorem \ref{theorem_main_theorem_rb_curved}. In Section \ref{section:general_system_direct_method}, estimates for the Nusselt number with respect to long-time bounds for the fluid velocity are provided. These long-time averages are estimated with respect to the system parameters in Section \ref{section:general_aprior}, which in Section \ref{section:direct_method_application} yields the final proofs by optimizing the boundary layer depth.

\subsection{The Direct Method}
\label{section:general_system_direct_method}

The following Lemma together with the regularity estimates proven later describes the direct method, where we use the localization principle of the Nusselt number close to the boundary. For the second order Poincaré type estimate a generalization of the argument given in \cite[Lemma 3.5]{drivas2022Bounds} is used.

\begin{lemma}
    Let $\Omega$ be $C^{1,1}$, $d = \min (\hp - \hm)>0$ and $u_0\in L^2(\Omega)$. Then
    \begin{align}
        \Nu
        &\lesssim \delta^\frac{1}{2} \langle \|\nabla u\|_2^2 \rangle^\frac{1}{2} + \delta^{-\frac{1}{2}}\Nu^\frac{1}{2}
        \label{nusselt_bound_H1_in_curved_lemma}
    \end{align}
    for every $\delta<d$ and if $\Omega$ is $C^{2,1}$ and $u_0\in W^{1,4}(\Omega)$
    \begin{align}
        \Nu
        &\lesssim \delta (1+d^{-2}) \langle \|u\|_{H^1}^2 \rangle^\frac{1}{4} \langle \|u\|_{H^2}^2\rangle^\frac{1}{4} + \delta^{-\frac{1}{2}}\Nu^\frac{1}{2}
    \end{align}
    for every $\delta <\frac{d}{2}$.
    \label{lemma_nusselt_bounded_by_Hx}
\end{lemma}

\begin{figure}
    \hrule
    \vspace{\baselineskip}
    \begin{center}
        \includegraphics[width=0.7\textwidth]{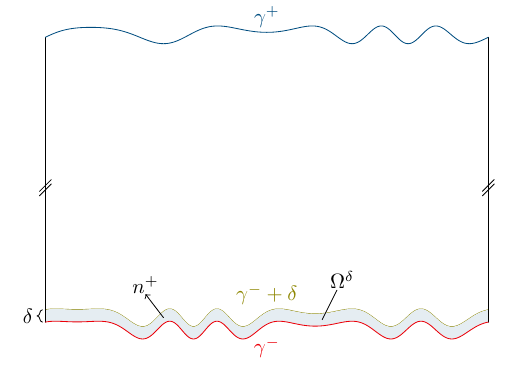}
        \vspace*{-5mm}
        \caption{Illustration of the definitions of $\Omegadelta$ and $\np$. Note that $\np$ is a vector field defined in the entire $\Omega$.}
        \label{fig:Omega_delta}
    \end{center}
    \hrule
\end{figure}
\begin{proof}
    We define $\Omegadelta=\lbrace (x_1,x_2)\in \Omega \ \vert \ x_2 < \gammam(x_1) + \delta \rbrace$ as illustrated in Figure \ref{fig:Omega_delta}. Then averaging \eqref{nu_representation_shifted} over $z\in [0,\delta]$ one finds
    \begin{align}
        \Nu
        &= \dashint_0^\delta \langle \np \cdot (u-\nabla) \theta\rangle_{\gammam + z} \ dz
        \\
        &= \dashint_0^\delta \langle \np \cdot u \theta\rangle_{\gammam + z} \ dz - \dashint_0^\delta \langle \np \cdot \nabla \theta\rangle_{\gammam + z} \ dz,
        \label{direct_method_nu_1}
    \end{align}
    where $\np(x_1,x_2)$ is the upwards pointing unit normal vector of $\gammam$ shifted in to $x_2$, i.e. $\np(x_1,x_2)=-n(x_1)$ where $n$ is the unit outward normal vector of $\Omega$ on $\gammam$, see Figure \ref{fig:Omega_delta}. Using Hölder's inequality the second term on the right-hand side of \eqref{direct_method_nu_1} can be estimated as
    \begin{align}
        - \dashint_0^\delta \langle \np \cdot \nabla \theta\rangle_{\gammam + z} \ dz
        &\leq \delta^{-1} \int_{\Omegadelta} |\nabla\theta|
        \leq \delta^{-\frac{1}{2}} \langle \|\nabla \theta\|_{L^2(\Omegadelta)} \rangle
        \\
        &\leq \delta^{-\frac{1}{2}} \langle \|\nabla \theta\|_{L^2(\Omegadelta)}^2 \rangle^\frac{1}{2}
        \leq \delta^{-\frac{1}{2}} \Nu^\frac{1}{2}.
        \label{direct_method_gradTheta}
    \end{align}
    By \eqref{maximum_principle} the first term on the right-hand side of \eqref{direct_method_nu_1} satisfies
    \begin{align}
        \dashint_0^\delta \langle \np \cdot u \theta\rangle_{\gammam + z} \ dz
        &\leq \delta^{-1} \langle |\np \cdot u|\rangle_{\Omegadelta}.
        \label{direct_method_un_1}
    \end{align}
    Note that since $\np\cdot u=0$ on $\gammam$ and $\np$ is constant in $x_2$ in $\Omegadelta$ we find for $x\in \Omegadelta$ by the fundamental theorem of calculus
    \begin{align}
        |\np\cdot u|(x_1,x_2) &\leq |\np\cdot u|(x_1,\hm(x_1)) + \int_0^\delta |\np\cdot \partial_2 u|(x_1,z)\ dz
        \\
        &\leq \delta \|(\np\cdot \partial_2 u)(x_1,\cdot)\|_{L^\infty(\Omegadeltavert)},
        \label{direct_method_un_2}
    \end{align}
    where $\Omegadeltavert(x_1)=(\hm(x_1),\hm(x_1)+\delta)$.
    Next, we extend $\np$ to the whole domain $\Omega$ by a smooth function such that it matches the normal vector on $\gammap$. In order to provide such an explicit extension let $\ngammamp(x_1)$ be the upwards pointing normal vector of $\gammam$ in $x_1$, $\ngammapp(x_1)$ be the upwards pointing normal vector of $\gammap$ in $x_1$ and define the transition function $G$ as
    \begin{align}
        G(x_2)=\int_{\hm(x_1)+\delta}^{x_2} e^{\frac{1}{\left(\frac{2 z-\hp(x_1)-\hm(x_1)-\delta}{\hp(x_1)-\hm(x_1)-\delta}\right)^2-1}}\ dz   
    \end{align}
    and the extension $\zeta$ as
    \begin{align}
        \zeta(x) =
        \begin{cases}
                \ngammamp & \text{ if } x_2<\hm(x_1)+\delta
                \\
                \left(1-\frac{G(x_2)}{G(\hp(x_1))}\right)\ngammamp
                + \frac{G(x_2)}{G(\hp(x_1))}\ngammapp & \text{ if } \hm(x_1)+\delta<x_2<\hp(x_1)
                \\
                \ngammapp & \text{ if } x_2=\hp(x_1),
        \end{cases}
    \end{align}
    where $\ngammamp=\ngammamp(x_1)$ and $\ngammapp=\ngammapp(x_1)$.
    \expliciteCalculation{
        \begin{align}
            G(\hp(x_1)) &= \frac{\hp(x_1)-\hm(x_1)-\delta}{2}\underbrace{\int_{-1}^1 e^{\frac{1}{s^2-1}}\ ds}_{\sim 0,44}\begin{cases}
                < \frac{\hp(x_1)-\hm(x_1)-\delta}{4}
                \\
                > \frac{\hp(x_1)-\hm(x_1)-\delta}{5}
            \end{cases}
            \\
            \frac{4}{\hp(x_1)-\hm(x_1)-\delta}&<\frac{1}{G(\hp(x_1))}<\frac{5}{\hp(x_1)-\hm(x_1)-\delta}
            \\
            \left|\partial_2 G(x_2)\right| &= e^{\frac{1}{\left(\frac{2 x_2-\hp(x_1)-\hm(x_1)-\delta}{\hp(x_1)-\hm(x_1)-\delta}\right)^2-1}} \leq 1
            \\
            g(s)&:= e^{\frac{1}{s^2-1}},\qquad \partial_2G(x_2)=g\left(\frac{2 x_2-\hp(x_1)-\hm(x_1)-\delta}{\hp(x_1)-\hm(x_1)-\delta}\right)
            \\
            \implies \left|\partial_2^2 G(x_2)\right| &= \left| g'(s) \frac{d}{dx_2} \left(\frac{2 x_2-\hp(x_1)-\hm(x_1)-\delta}{\hp(x_1)-\hm(x_1)-\delta}\right)\right|
            \overset{|g'(s)|\leq 1}{\leq} \frac{2}{\hp(x_1)-\hm(x_1)-\delta}
        \end{align}
        and therefore
        \begin{align}
            \|\zeta\|_{L^\infty(\Omega)} &\leq 2
            \\
            \|\partial_2\zeta\|_{L^\infty(\Omega)} &\leq \frac{5}{\hp(x_1)-\hm(x_1)-\delta}\|G'\|_\infty
            \leq \frac{5}{\hp(x_1)-\hm(x_1)-\delta}
            \\
            \|\partial_2^2\zeta\|_{L^\infty(\Omega)} &\leq \frac{5}{\hp(x_1)-\hm(x_1)-\delta}\|G''\|_\infty
            \leq \frac{10}{(\hp(x_1)-\hm(x_1)-\delta)^2}
        \end{align}
    }
    Note that
    \begin{align}
        \|\zeta\|_{L^\infty(\Omega)} &\lesssim 1
        \label{direct_method_un_3a}
        \\
        \|\partial_2\zeta\|_{L^\infty(\Omega)} &\lesssim \|(\hp(x_1)-\hm(x_1)-\delta)^{-1}\|_\infty \leq d^{-1}
        \label{direct_method_un_3b}
        \\
        \|\partial_2^2\zeta\|_{L^\infty(\Omega)} &\lesssim \|(\hp(x_1)-\hm(x_1)-\delta)^{-2}\|_\infty\leq d^{-2}.
        \label{direct_method_un_3c}
    \end{align}\noeqref{direct_method_un_3b}%
    Additionally, since for every $x_1\in (0,\Gamma)$
    \begin{align}
        \int_{\hm(x_1)}^{\hp(x_1)} \partial_2(\zeta\cdot u)\ dz = \ngammapp \cdot u - \ngammamp\cdot u = 0
    \end{align}
    there exists $\tilde x_2 (x_1)$ such that
    \begin{align}
        \partial_2(\zeta\cdot u)(x_1,\tilde x_2(x_1)) = 0.
    \end{align}
    By the fundamental theorem of calculus
    \begin{align}
        |\partial_2(\zeta\cdot u)|^2 (x) &= \left| \partial_2(\zeta\cdot u)(x_1,\tilde x_2(x_1)) + \int_{\tilde x_2(x_1)}^{x_2} \partial_2 \left((\partial_2(\zeta\cdot u))^2\right)\ dz\right|
        \\
        &\leq 2 \int_{\tilde x_2(x_1)}^{x_2} |\partial_2(\zeta\cdot u) \partial_2^2(\zeta\cdot u)| \ dz
        \\
        &\leq 2 \|\partial_2(\zeta\cdot u)\|_{L^2(\Omegavertvert)}\|\partial_2^2(\zeta\cdot u)\|_{L^2(\Omegavertvert)}
        \label{direct_method_un_4}
    \end{align}
    for any $x\in \Omega$, where $\Omegavertvert(x_1)=(\hm(x_1),\hp(x_1))$. Combining \eqref{direct_method_un_2} and \eqref{direct_method_un_4} one has for $x\in \Omegadelta$
    \begin{align}
        |\np \cdot u|(x)
        &\leq \delta\|\np \cdot \partial_2 u\|_{L^\infty(\Omegadeltavert)}
        \\
        &= \delta \|\partial_2 (\np\cdot  u)\|_{L^\infty(\Omegadeltavert)}
        \\
        &\leq \delta \|\partial_2 (\zeta \cdot  u)\|_{L^\infty(\Omegavertvert)}
        \\
        &\leq 2\delta \|\partial_2(\zeta\cdot u)\|_{L^2(\Omegavertvert)}^\frac{1}{2}\|\partial_2^2(\zeta\cdot u)\|_{L^2(\Omegavertvert)}^\frac{1}{2}
        \label{direct_method_un_H2}
    \end{align}
    and therefore for the spatial and long time average, using Hölder's inequality, it follows that
    \expliciteCalculation{
    \begin{align}
        \delta^{-1}&\langle |\np \cdot u|\rangle_{\Omegadelta}
        \\
        &\lesssim \langle\|\partial_2(\zeta\cdot u)\|_{L^2(\hm(x_1),\hp(x_1))}^\frac{1}{2}\|\partial_2^2(\zeta\cdot u)\|_{L^2(\hm(x_1),\hp(x_1))}^\frac{1}{2}\rangle_{\Omegadelta}
        \\
        &= \delta |\Omega|^{-1} \limsup_{T\to\infty} \dashint_0^T \int_0^\Gamma \|\partial_2(\zeta\cdot u)\|_{L^2(\hm(x_1),\hp(x_1))}^\frac{1}{2}\|\partial_2^2(\zeta\cdot u)\|_{L^2(\hm(x_1),\hp(x_1))}^\frac{1}{2} \ dx_1 \ dt 
        \\
        &\leq \delta |\Omega|^{-1} \limsup_{T\to\infty} \dashint_0^T \left(\int_0^\Gamma \|\partial_2(\zeta\cdot u)\|_{L^2(\hm(x_1),\hp(x_1))} \ dx_1\right)^\frac{1}{2}
        \\
        &\qquad \cdot \left(\int_0^\Gamma \|\partial_2^2(\zeta\cdot u)\|_{L^2(\hm(x_1),\hp(x_1))} \ dx_1\right)^\frac{1}{2} \ dt 
        \\
        &\leq \delta |\Omega|^{-\frac{1}{2}} \limsup_{T\to\infty} \dashint_0^T \left(\int_0^\Gamma \|\partial_2(\zeta\cdot u)\|_{L^2(\hm(x_1),\hp(x_1))}^2 \ dx_1\right)^\frac{1}{4}
        \\
        &\qquad \cdot \left(\int_0^\Gamma \|\partial_2^2(\zeta\cdot u)\|_{L^2(\hm(x_1),\hp(x_1))}^2 \ dx_1\right)^\frac{1}{4} \ dt 
        \\
        &= \delta|\Omega|^{-\frac{1}{2}} \limsup_{T\to\infty} \dashint_0^T \|\partial_2(\zeta\cdot u)\|_{L^2}^\frac{1}{2}\|\partial_2^2(\zeta\cdot u)\|_{L^2}^\frac{1}{2}\ dt
        \\
        &= \delta |\Omega|^{\frac{1}{2}} \langle\|\partial_2(\zeta\cdot u)\|_{L^2}^\frac{1}{2}\|\partial_2^2(\zeta\cdot u)\|_{L^2}^\frac{1}{2}\rangle
        \\
        &\leq \delta |\Omega|^{\frac{1}{2}}\langle\|\partial_2(\zeta\cdot u)\|_{L^2}^2\rangle^\frac{1}{4}\langle\|\partial_2^2(\zeta\cdot u)\|_{L^2}^2\rangle^\frac{1}{4}
        \\
        &\lesssim \delta |\Omega|^{\frac{1}{2}}(d^{-2}\langle\|u\|_{L^2}^2\rangle + \langle\|\partial_2u\|_{L^2}^2\rangle)^\frac{1}{4}(d^{-4}\langle\|u\|_{L^2}^2\rangle+d^{-2}\langle\|\partial_2 u\|_{L^2}^2\rangle+\langle\|\partial_2^2 u\|_{L^2}^2\rangle)^\frac{1}{4}
        \\
        &\lesssim \delta |\Omega|^{\frac{1}{2}}\left( d^{-\frac{3}{2}} \langle\|u\|_{L^2}^2\rangle^\frac{1}{2}
        + d^{-1}\langle\|u\|_{L^2}^2\rangle^\frac{1}{4}\langle\|\partial_2 u\|_{L^2}^2\rangle^\frac{1}{4}
        + d^{-\frac{1}{2}}\langle\|u\|_{L^2}^2\rangle^\frac{1}{4}\langle\|\partial_2^2 u\|_{L^2}^2\rangle^\frac{1}{4}
        \right.
        \\
        &\qquad \left.
        +d^{-\frac{1}{2}}\langle\|\partial_2 u\|_{L^2}^2\rangle^\frac{1}{2}
        +\langle\|\partial_2 u\|_{L^2}^2\rangle^\frac{1}{4}\langle\|\partial_2^2 u\|_{L^2}^2\rangle^\frac{1}{4}\right)
    \end{align}
    }
    \begin{align}
        \delta^{-1}\langle |\np \cdot u|\rangle_{\Omegadelta} &\leq \delta |\Omega|^{\frac{1}{2}} \langle\|\partial_2(\zeta\cdot u)\|_{L^2}^\frac{1}{2}\|\partial_2^2(\zeta\cdot u)\|_{L^2}^\frac{1}{2}\rangle
        \\
        &\lesssim \delta (1+d^{-2})\langle \|u\|_{H^1}^2 \rangle^\frac{1}{4} \langle \|u\|_{H^2}^2\rangle^\frac{1}{4},
        \label{direct_method_un_5}
    \end{align}
    where we used \eqref{direct_method_un_3a}-\eqref{direct_method_un_3c} and Young's inequality in the last estimate. Combining \eqref{direct_method_nu_1}, \eqref{direct_method_gradTheta}, \eqref{direct_method_un_1} and \eqref{direct_method_un_5}
    \begin{align}
        \Nu
        &= \dashint_0^\delta \langle \np \cdot u \theta\rangle_{\gammam + z} \ dz - \dashint_0^\delta \langle \np \cdot \nabla \theta\rangle_{\gammam + z} \ dz
        \\
        &\leq \dashint_0^\delta \langle \np \cdot u \theta\rangle_{\gammam + z} \ dz + \delta^{-\frac{1}{2}}\Nu^\frac{1}{2}
        \\
        &\leq \delta^{-1} \langle |\np \cdot u|\rangle_{\Omegadelta} + \delta^{-\frac{1}{2}}\Nu^\frac{1}{2}
        \\
        &\lesssim \delta (1+d^{-2}) \langle \|u\|_{H^1}^2 \rangle^\frac{1}{4} \langle \|u\|_{H^2}^2\rangle^\frac{1}{4} + \delta^{-\frac{1}{2}}\Nu^\frac{1}{2}.
    \end{align}
    Analogously for $x\in \Omegadelta$ one could have estimated
    \begin{align}
        |\np \cdot u|(x)
        &= \left|(\np\cdot u)(x_1,\hm(x_1)) + \int_{\hm(x_1)}^{x_2}\partial_2 (\np\cdot u)\ dz\right|
        \\
        &\leq \int_{\hm(x_1)}^{\hm(x_1)+\delta} |\partial_2u|\ dz \leq \delta^{\frac{1}{2}} \|\partial_2 u\|_{L^2(\hm(x_1),\hm(x_1)+\delta)}
    \end{align}
    instead of \eqref{direct_method_un_H2}, which would have resulted in
    \expliciteCalculation{
    \begin{align}
        \delta^{-1} \langle |\np \cdot u|\rangle_{\Omegadelta} 
        &= \delta^{-1}\limsup_{T\to\infty} |\Omega|^{-1} \dashint_0^T \int_0^\Gamma \int_{\hm(x_1)}^{\hm(x_1)+\delta} |\np \cdot u| \ dx_2 \ dx_1 \ dt
        \\
        &\leq \delta^{-\frac{1}{2}}\limsup_{T\to\infty} |\Omega|^{-1} \dashint_0^T \int_0^\Gamma \int_{\hm(x_1)}^{\hm(x_1)+\delta} \|\partial_2 u\|_{L^2(\hm(x_1),\hm(x_1)+\delta)} \ dx_2 \ dx_1 \ dt
        \\
        &\leq \delta^{\frac{1}{2}}\limsup_{T\to\infty} |\Omega|^{-1} \dashint_0^T \int_0^\Gamma \|\partial_2 u\|_{L^2(\hm(x_1),\hm(x_1)+\delta)} \ dx_1 \ dt
        \\
        &\leq \delta^{\frac{1}{2}}\limsup_{T\to\infty} |\Omega|^{-\frac{1}{2}} \dashint_0^T  \|\partial_2 u\|_{L^2(\Omegadelta)} \ dt
        \\
        &\leq \delta^{\frac{1}{2}}\limsup_{T\to\infty} |\Omega|^{-\frac{1}{2}} \dashint_0^T  \|\partial_2 u\|_{L^2} \ dt
        \\
        &= \delta^{\frac{1}{2}} |\Omega|^{\frac{1}{2}} \langle  \|\partial_2 u\|_{L^2} \rangle
        \\
        &\leq \delta^{\frac{1}{2}} |\Omega|^{\frac{1}{2}} \langle  \|\partial_2 u\|_{L^2}^2 \rangle^\frac{1}{2}
    \end{align}
    }
    \begin{align}
        \Nu
        &= \dashint_0^\delta \langle \np \cdot u \theta\rangle_{\gammam + z} \ dz - \dashint_0^\delta \langle \np \cdot \nabla \theta\rangle_{\gammam + z} \ dz
        \\
        &\leq \dashint_0^\delta \langle \np \cdot u \theta\rangle_{\gammam + z} \ dz + \delta^{-\frac{1}{2}}\Nu^\frac{1}{2}
        \\
        &\leq \delta^{-1} \langle |\np \cdot u|\rangle_{\Omegadelta} + \delta^{-\frac{1}{2}}\Nu^\frac{1}{2}
        \\
        &\lesssim \delta^\frac{1}{2} |\Omega|^\frac{1}{2} \langle \|\nabla u\|_2^2 \rangle^\frac{1}{2} + \delta^{-\frac{1}{2}}\Nu^\frac{1}{2},
    \end{align}
    proving \eqref{nusselt_bound_H1_in_curved_lemma}.    
\end{proof}

\subsection{A-Priori Estimates}
\label{section:general_aprior}

With Lemma \ref{lemma_nusselt_bounded_by_Hx} at hand one needs regularity estimates for the velocity to get bounds on the Nusselt number. The first such estimate is an energy bound, which in turn implies a long-time bound for the corresponding $H^1$-norm. In particular, using the symmetric gradient instead of the full gradient removes the assumption of small curvature with respect to the slip coefficient in the corresponding estimates (see Lemma 3.1 - Corollary 3.4 in \cite{bleitnerNobili2024Bounds}), which we will handle in Section \ref{section:nonlinearity_paper}.

\begin{lemma}[Energy Bound]
    Let $\Omega$ be $C^{1,1}$, $0<\alpha\in L^\infty(\gammam\cup\gammap)$ and $u_0\in L^2$. Then
    \label{lemma_energy_bound}
    \begin{align}
        \|u\|_{2}^2&\lesssim \|u_0\|_{2}^2 + (1+\|\alpha^{-1}\|_\infty^2) \Ra^2 \label{energy_bound_lemma_L2_bounded}
        \\
        \langle\|u\|_{H^1}^2\rangle &\lesssim (1+\|\alpha^{-1}(1+|\kappa|)\|_\infty) \Nu \Ra,
        \label{energy_bound_lemma_H1_bound}
    \end{align}
    where the implicit constant only depends on $|\Gamma|$ and the Lipschitz constant of the boundary.
    Additionally, if $|\kappa|\leq 2\alpha$
    \begin{align}
        \langle \|\nabla u\|_2^2\rangle \lesssim \Nu \Ra. \label{energy_bound_lemma_kappa_leq_2alpha_case}
    \end{align}
\end{lemma}

\begin{proof}
    Testing \eqref{navier_stokes} with $u$ yields
    \begin{align}
        \frac{1}{2\Pr}\frac{d}{dt} \|u\|_{L^2}^2 &= - \frac{1}{\Pr}\int_\Omega u\cdot (u\cdot \nabla) u - \int_{\Omega} u\cdot \nabla p + \int_{\Omega} u\cdot \Delta u + \Ra \int_{\Omega} u_2 \theta.
        \\
        \label{nbjsfibnfb}
    \end{align}
    Note that the first term on the right-hand side of \eqref{nbjsfibnfb} vanishes by \eqref{incompressibility} and \eqref{no_penetration_bc} as
    \begin{align}
        \int_{\Omega} u\cdot (u\cdot \nabla) u &= \frac{1}{2}\int_{\Omega} u\cdot \nabla |u|^2
        =-\frac{1}{2}\int_{\Omega} \nabla \cdot u |u|^2 + \frac{1}{2}\int_{\gammam\cup\gammap} u\cdot n |u|^2
        = 0,
        \\
        \label{ndfuindf}
    \end{align}
    where we use integration by parts in the second identity. Similarly, the second term on the right-hand side of \eqref{nbjsfibnfb} vanishes as
    \begin{align}
        \label{bnuiadfasd}
        \int_\Omega u\cdot \nabla p = \int_{\gammam\cup\gammap} p u\cdot n - \int_\Omega p\nabla \cdot u = 0.
    \end{align}
    Applying \eqref{ndfuindf}, \eqref{bnuiadfasd} and Lemma \ref{lemma_int_Delta_u_v} in \eqref{nbjsfibnfb} yields
    \begin{align}
        \label{jalksdfjakls}
        \frac{1}{2\Pr}\frac{d}{dt} \|u\|_{L^2}^2 + 2 \|\D u\|_{L^2}^2 + 2 \int_{\gammam\cup\gammap} \alpha u_\tau^2 &= \Ra \int_{\Omega} u_2 \theta.
    \end{align}
    Using Hölder's inequality in order to estimate the right-hand side of \eqref{jalksdfjakls} and Lemma \ref{lemma_coercivity}, i.e. \eqref{coercivity_lemma_zero_order}, one finds
    \begin{align}
        \frac{1}{2\Pr} \frac{d}{dt} \|u\|_{L^2}^2 + c (1+\|\alpha^{-1}\|_\infty)^{-1} \|u\|_{L^2}^2
        &\leq \Ra \|u_2\|_{L^2} \|\theta\|_{L^2}
        \\
        &\leq \frac{\epsilon}{4} \|u\|_{L^2}^2 + \epsilon^{-1} \Ra^2.
    \end{align}
    for some $c>0$ depending only on $\Gamma$ and any $\epsilon>0$, where in the last inequality we used Young's inequality and \eqref{maximum_principle}. Choosing $\epsilon = c(1+\|\alpha^{-1}\|_\infty)^{-1}$ in order to absorb the first term on the right-hand side yields
    \begin{align}
        \frac{d}{dt} \|u\|_{L^2}^2 \leq -\frac{c}{2}(1+\|\alpha^{-1}\|_\infty)^{-1} \Pr \|u\|_{L^2}^2 + (1+\|\alpha^{-1}\|_\infty) \Pr \Ra^2
    \end{align}
    and Grönwall's inequality results in
    \begin{align}
        \|u\|_{L^2}^2 
        &\lesssim \|u_0\|_{L^2}^2 + (1+\|\alpha^{-1}\|_\infty^2) \Ra^2,
        \label{skdlabndisb}
    \end{align}
    proving \eqref{energy_bound_lemma_L2_bounded}.

    Applying Lemma \ref{lemma_coercivity}, i.e. \eqref{coercivity_lemma_H1}, to \eqref{jalksdfjakls} implies that there exists a constant $c>0$ depending only on $\Gamma$ such that
    \begin{align}
        \frac{1}{2\Pr}\frac{d}{dt} \|u\|_{L^2}^2 &+ c(1+\|\alpha^{-1}(1+|\kappa|)\|_\infty)^{-1} \|u\|_{H^1}^2
        \\
        &\leq \frac{1}{2\Pr}\frac{d}{dt} \|u\|_{L^2}^2 + 2 \|\D u\|_{L^2}^2 + 2 \int_{\gammam\cup\gammap} \alpha u_\tau^2
        \\
        &= \Ra \int_{\Omega} u_2 \theta.
        \label{nufsbifb}
    \end{align}
    The integral on the right-hand side of \eqref{nufsbifb} can be rewritten as
    \begin{align}
        \int_{\Omega} u_2 \theta &= \int_{\Omega} (u_2-\partial_2) \theta + \int_{\Omega} \partial_2 \theta
        \\
        &= \int_{\Omega} (u_2-\partial_2) \theta + \int_0^\Gamma \int_{\hm(x_1)}^{\hp(x_1)} \partial_2 \theta \ dx_2 \ dx_1
        \\
        &= \int_{\Omega} (u_2-\partial_2) \theta + \int_0^\Gamma \theta(x_1,\hp(x_1))\ dx_1 - \int_0^\Gamma \theta(x_1,\hm(x_1))\ dx_1
        \\
        &= \int_{\Omega} (u_2-\partial_2) \theta - \Gamma,
        \label{gbnsifbusndb}
    \end{align}
    where we used the boundary condition \eqref{heat_bc}, which combined with \eqref{nufsbifb} yields
    \begin{align}
        \frac{1}{2\Pr}\frac{d}{dt} \|u\|_{L^2}^2+c(1+\|\alpha^{-1}(1+|\kappa|)\|_\infty)^{-1}\|u\|_{H^1}^2 &\lesssim \Ra \int_{\Omega} (u_2-\partial_2) \theta\qquad
        \label{nsifgbfginb}
    \end{align}
    Note that the first term on the left-hand side of \eqref{nsifgbfginb} vanishes in the long-time average as
    \begin{align}
        \limsup_{T\to\infty} \dashint_0^T \frac{d}{dt} \|u\|_{L^2}^2 \ dt &= \limsup_{T\to\infty} \frac{\|u\|_{L^2}^2-\|u_0\|_{L^2}^2}{T} = 0,
        \label{sdbnisdfbsujb}
    \end{align}
    where we used that $\|u\|_{L^2}^2$ is uniformly bounded in time due to \eqref{skdlabndisb}. Taking the long-time and spatial average of \eqref{nsifgbfginb} and using \eqref{sdbnisdfbsujb} results in
    \begin{align}
        \langle\|u\|_{H^1}^2\rangle
        &\lesssim (1+\|\alpha^{-1}(1+|\kappa|)\|_\infty) \langle(u_2-\partial_2) \theta\rangle_{\Omega} \Ra
        \\
        &\leq (1+\|\alpha^{-1}(1+|\kappa|)\|_\infty) (\max \hp - \min \hm) \Nu \Ra
        \\
        &\lesssim (1+\|\alpha^{-1}(1+|\kappa|)\|_\infty) \Nu \Ra
    \end{align}
    due to \eqref{nu_representation_straight} in Lemma \ref{lemma_nusselt_number_representations}.

    In order to prove \eqref{energy_bound_lemma_kappa_leq_2alpha_case} note that if $|\kappa|\leq 2\alpha$ Lemma \ref{lemma_grad_ids} and \eqref{nufsbifb} imply
    \begin{align}
        \frac{1}{2\Pr}\frac{d}{dt} \|u\|_{L^2}^2 + \|\nabla u\|_{L^2}^2
        &\leq \frac{1}{2\Pr}\frac{d}{dt} \|u\|_{L^2}^2 + \|\nabla u\|_{L^2}^2 + \int_{\gammam\cup\gammap} (2\alpha+\kappa) u_\tau^2
        \\
        &= \frac{1}{2\Pr}\frac{d}{dt} \|u\|_{L^2}^2 + 2 \|\D u\|_{L^2}^2 + 2 \int_{\gammam\cup\gammap} \alpha u_\tau^2
        \\
        &= \Ra \int_{\Omega} u_2 \theta
        \label{energy_balance}
    \end{align}
    and estimating the right-hand side similar to before and taking the long-time average yields
    \begin{align}
        \langle \|\nabla u\|_2^2\rangle \lesssim \Nu \Ra.
    \end{align}
\end{proof}

These bounds are already sufficient to establish the bounds of Theorem \ref{theorem_1_2_bound}. In contrast the more tailored results of Theorem \ref{theorem_main_theorem_rb_curved} request for second-order bounds on the velocity. Therefore we first prove bounds for the vorticity that will afterwards be used to bound the excessive term arising from the nonlinearity.
The following estimates can be found in \cite[Section 3.2]{bleitnerNobili2024Bounds} and we reprove them here for the convenience of the reader.
\begin{lemma}[Vorticity Bound]
    \label{lemma_vorticity_bound}
    Let $2<p\in 2\mathbb{N}$, $\Omega$ be $C^{1,1}$, $u_0\in W^{1,p}$ and $0<\alpha\in L^\infty(\gammam\cup\gammap)$. Then
    \begin{align}
        \|\omega\|_{p}
        &\lesssim \|\omega_0\|_p + C_{\alpha,\kappa,p} \|u_0\|_2 + (1+C_{\alpha,\kappa,p}) (1+\|\alpha^{-1}\|_\infty)\Ra,
    \end{align}
    where $C_{\alpha,\kappa,p} = \left(1+\|\kappa\|_\infty+\|\alpha+\kappa\|_\infty^{\frac{p}{p-2}}\right)\|\alpha+\kappa\|_\infty$ and the implicit constant depends only on $p$, $\Gamma$ and the Lipschitz constant of the boundary.
\end{lemma}

\begin{remark}[\texorpdfstring{$\|u\|_{W^{1,p}}$ is Uniformly Bounded in Time}{Uniform in Time Vorticity} Bound]
    \label{remark_Linfty_W1p_bound}
    Note that by \eqref{elliptic_regularity_periodic_lemma_first_order} and Hölder's inequality, Lemma \ref{lemma_vorticity_bound} and Lemma \ref{lemma_energy_bound} yield a uniform in time bound for the full $W^{1,p}$-norm of $u$.
\end{remark}

\begin{proof}
    Fix some arbitrary time $\tmax>0$, define $\Lambda = 2 \|(\alpha+\kappa)u_\tau\|_{L^\infty([0,\tmax]\times \lbrace\gammam\cup\gammap\rbrace)}$ and let $\omegatildepm$ solve    
    \begin{alignat}{2}
        \frac{1}{\Pr}\left(\omegatildepm_t+u\cdot\nabla \omegatildepm \right) - \Delta \omegatildepm &= \Ra \partial_1 \theta & \quad \textnormal{ in } &\Omega
        \\
        \omegatildepm_0 &= \pm |\omega_0| & \quad \textnormal{ in } &\Omega
        \\
        \omegatildepm &= \pm \Lambda & \quad \textnormal{ on } &\gammam\cup\gammap.
    \end{alignat}
    Then by \eqref{vorticity_equation}, \eqref{vorticity_bc} $\omegabarpm=\omega-\omegatildepm$ solves
    \begin{alignat}{2}
        \frac{1}{\Pr}(\omegabarpm_t + u\cdot \nabla \omegabarpm) - \Delta \omegabarpm &= 0 & \quad \textnormal{ in } &\Omega
        \\
        \omegabarpm_0 &= \omega_0 \mp |\omega_0| & \quad \textnormal{ in } &\Omega
        \\
        \omegabarpm &= -2(\alpha+\kappa)u_\tau \mp \Lambda & \quad \textnormal{ on } &\gammam\cup\gammap.
    \end{alignat}
    As the initial and boundary values of $\omegabarp$ and $\omegabarm$ are non-positive, respectively non-negative, the maximum principle shows $\omegabarp(x,t)\leq 0$ and $\omegabarm(x,t)\geq 0$, implying
    \begin{align}
        |\omega|\leq |\omegatildem|+|\omegatildep|.
        \label{bngfbzhnfbhfb}
    \end{align}
    Therefore it is sufficient to get bounds on $\omegatildepm$. Next, we define $\omegahatpm = \omegatildepm \mp \Lambda$ to remove the boundary condition. The following estimates work analogously for $\omegahatp$ and $\omegahatm$. Hence we will focus on $\omegahatp$ and omit the $^+$ to simplify notation. With this definition $\omegahatExpl$ satisfies
    \begin{alignat}{2}
        \frac{1}{\Pr} \left(\omegahatExpl_t + u\cdot \nabla \omegahatExpl\right)- \Delta \omegahatExpl &= \Ra \partial_1 \theta & \quad \textnormal{ in } &\Omega\label{omegahat_pde}
        \\
        \omegahatExpl_0 &= |\omega_0| - \Lambda & \quad \textnormal{ in } &\Omega
        \label{omegahat_ic}
        \\
        \omegahatExpl &= 0 & \quad \textnormal{ on } &\gammam\cup\gammap.\label{omegahat_bc}
    \end{alignat}
    Testing \eqref{omegahat_pde} with $\omegahatExpl^{p-1}$, where $2<p\in 2\mathbb{N}$, yields
    \begin{align}
        \frac{1}{p\Pr} \frac{d}{dt} \|\omegahatExpl\|_p^p = -\frac{1}{\Pr} \int_{\Omega} \omegahatExpl^{p-1} u\cdot\nabla\omegahatExpl + \int_{\Omega} \omegahatExpl^{p-1}\Delta \omegahatExpl + \Ra \int_{\Omega} \omegahatExpl^{p-1} \partial_1 \theta .\quad\label{bnfusibf}
    \end{align}
    Note that when integrating by parts the first term on the right-hand side of \eqref{bnfusibf} vanishes as
    \begin{align}
        p\int_{\Omega} \omegahatExpl^{p-1}u\cdot \nabla \omegahatExpl = \int_{\Omega} u\cdot\nabla(\omegahatExpl^p) = \int_{\gammam\cup\gammap} \omegahatExpl^p n\cdot u - \int_{\Omega} \omegahatExpl^p \nabla\cdot u = 0
        \quad\label{dsakfjsadlf}
    \end{align}
    by \eqref{no_penetration_bc} and \eqref{incompressibility}. For the second term on the right-hand side of \eqref{bnfusibf} integration by parts and the boundary condition \eqref{omegahat_bc} yield
    \begin{align}
        \int_{\Omega} \omegahatExpl^{p-1} \Delta\omegahatExpl &= \int_{\gammam\cup\gammap} \omegahatExpl^{p-1} n\cdot\nabla\omegahatExpl - \int_{\Omega} \nabla\omegahatExpl\cdot \nabla (\omegahatExpl^{p-1})
        \\
        &= - (p-1)\int_{\Omega} \omegahatExpl^{p-2} |\nabla\omegahatExpl|^2
        \label{nbugifgnb}
    \end{align}
    and similarly for the third term on the right-hand side of \eqref{bnfusibf}
    \begin{align}
        \Ra \int_{\Omega} \omegahatExpl^{p-1} \partial_1 \theta &= \Ra \int_{\gammam\cup\gammap} n_1 \omegahatExpl^{p-1} \theta - \Ra \int_{\Omega} \theta \partial_1 \left( \omegahatExpl^{p-1} \right)
        \\
        &= - (p-1) \Ra \int_{\Omega} \theta \omegahatExpl^{p-2} \partial_1 \omegahatExpl
        \\
        &\leq (p-1) \Ra \|\theta\|_\infty \|\omegahatExpl^{\frac{p-2}{2}} \|_2\|\omegahatExpl^{\frac{p-2}{2}}|\nabla\omegahatExpl|\|_2
        \\
        &\leq \frac{p-1}{4\epsilon} \Ra^2 \|\theta\|_\infty^2 \|\omegahatExpl^{\frac{p}{2}} \|_2^2 + \epsilon(p-1)\|\omegahatExpl^{\frac{p-2}{2}}|\nabla\omegahatExpl|\|_2^2
        \\
        &\leq \frac{p-1}{2} \Ra^2 \int_{\Omega}\omegahatExpl^{p-2} + \frac{p-1}{2} \int_{\Omega} \omegahatExpl^{p-2} |\nabla\omegahatExpl|^2,
    \end{align}
    where in the second to last estimate we used Hölder's inequality and in the last estimate we used Young's inequality and \eqref{maximum_principle}.
    Therefore, combining \eqref{bnfusibf}, \eqref{dsakfjsadlf} and \eqref{nbugifgnb} yields
    \begin{align}
        \frac{1}{p\Pr} \frac{d}{dt} \|\omegahatExpl\|_p^p
        &\leq - \frac{p-1}{2}\int_{\Omega} \omegahatExpl^{p-2} |\nabla\omegahatExpl|^2 + \frac{p-1}{2} \Ra^2 \int_{\Omega}\omegahatExpl^{p-2}
        \\
        &= - \frac{2(p-1)}{p^2} \|\nabla(\omegahatExpl^\frac{p}{2})\|_2^2 + \frac{p-1}{2} \Ra^2 \int_{\Omega}\omegahatExpl^{p-2}.
    \end{align}
    As $\omegahatExpl$ vanishes on $\gammam\cup\gammap$, Poincaré's inequality implies the existence of a constant $C>0$, depending on $p,|\Omega|$ and $\|h\|_\infty$ such that
    \begin{align}
        \frac{1}{2\Pr} \|\omegahatExpl\|_p^2 \frac{d}{dt} \|\omegahatExpl\|_p^2
        &= \frac{1}{p\Pr} \frac{d}{dt} \|\omegahatExpl\|_p^p
        \\
        &\leq - C \|\omegahatExpl^\frac{p}{2}\|_2^2 + \frac{p-1}{2} \Ra^2 \int_{\Omega}\omegahatExpl^{p-2}
        \\
        &\lesssim - C \|\omegahatExpl\|_p^p + \Ra^2 \|\omegahatExpl\|_p^{p-2},
        \label{aosidfjaosif}
    \end{align}
    where we estimated the second term on the right-hand side by Hölder's inequality. Dividing \eqref{aosidfjaosif} by $\|\omegahatExpl\|_p^2$ we find
    \begin{align}
        \frac{d}{dt}\|\omegahatExpl\|_p^2 + C\Pr \|\omegahatExpl\|_p^2 \lesssim \Pr\Ra^2
    \end{align}
    and therefore Grönwall's inequality and \eqref{omegahat_ic} results in
    \begin{align}
        \|\omegahatExpl\|_p^2 \lesssim \|\omegahatExpl_0\|_p^2 + \Ra^2 \lesssim \|\omega_0\|_p^2 + \Lambda^2 + \Ra^2.
    \end{align}
    Applying the analogous steps for $\omegahatm$ implies
    \begin{align}
        \|\omegahatpm\|_p^2 \lesssim \|\omega_0\|_p^2 + \Lambda^2 + \Ra^2.\label{nbufbgnbg}
    \end{align}
    In order to estimate $\Lambda$ notice that 
    \begin{align}
        \Lambda &= 2\|(\alpha+\kappa)u_\tau\|_{L^\infty([0,\tmax]\times \lbrace\gammam\cup\gammap\rbrace)}
        \\
        &\lesssim \|\alpha+\kappa\|_\infty \|u_\tau\|_{L^\infty\left(0,\tmax;L^\infty(\gammam\cup\gammap)\right)}
        \\
        &\lesssim \|\alpha+\kappa\|_\infty \|u\|_{L^\infty(0,\tmax;L^\infty(\Omega))},\label{adsnbiugntb}
    \end{align}
    (see for instance \cite[Satz 7.1.26]{emmrichGewoehnlicheUndOperatorDifferentialgleichungen}, 
    respectively \cite[Exercise 4.1]{troeltzschOptimaleSteuerungpartiellerDifferentialgleichungen}
    )
    and as $p>2$, \eqref{adsnbiugntb} and Gagliardo-Nirenberg interpolation yield
    \begin{align}
        \Lambda &\lesssim \|\alpha+\kappa\|_\infty \|u\|_{L^\infty_t;L^\infty_x}
        \\
        &\lesssim \|\alpha+\kappa\|_\infty \|u\|_{L^\infty_t;W^{1,p}_x}^{\frac{p}{2(p-1)}}\|u\|_{L^\infty_t;L^2_x}^{\frac{p-2}{2(p-1)}}
        \\
        &\lesssim \|\alpha+\kappa\|_\infty\|\omega\|_{L^\infty_t;L^p_x}^{\frac{p}{2(p-1)}}\|u\|_{L^\infty_t;L^2_x}^{\frac{p-2}{2(p-1)}}+\left(1+\|\kappa\|_\infty\right)\|\alpha+\kappa\|_\infty\|u\|_{L^\infty_t;L^2_x}
        ,
    \end{align}
    with $L^\infty;L^q_x=L^\infty(0,\tmax;L^q(\Omega))$ and $L^\infty;W^{1,p}_x=L^\infty(0,\tmax;W^{1,p}(\Omega))$, where in the last inequality we used \eqref{elliptic_regularity_periodic_lemma_first_order}. Next, Young's inequality yields
    \begin{align}
        \Lambda &\lesssim \|\alpha+\kappa\|_\infty\|\omega\|_{L^\infty_t;L^p_x}^{\frac{p}{2(p-1)}}\|u\|_{L^\infty_t;L^2_x}^{\frac{p-2}{2(p-1)}}+\left(1+\|\kappa\|_\infty\right)\|\alpha+\kappa\|_\infty\|u\|_{L^\infty_t;L^2_x}
        \\
        &\leq \epsilon \|\omega\|_{L^\infty_t;L^p_x} + \left(1+\|\kappa\|_\infty+\epsilon^{-\frac{p}{p-2}}\|\alpha+\kappa\|_\infty^{\frac{p}{p-2}}\right)\|\alpha+\kappa\|_\infty\|u\|_{L^\infty_t;L^2_x}
        \qquad
        \label{bnufhbfnbh}
    \end{align}
    for any $\epsilon>0$. With this estimate at hand, we can proceed bounding the vorticity. Combining \eqref{bngfbzhnfbhfb}, the definition of $\omegahatExpl$, i.e. $\omegahatExpl = \omegatildep-\Lambda$, \eqref{nbufbgnbg} and \eqref{bnufhbfnbh} and the analogous estimates for $\omegatildem$ we find
    \begin{align}
        \|\omega\|_{L^\infty_t;L^p_x}^2 &\lesssim \|\omegatildem\|_{L^\infty_t;L^p_x}^2 + \|\omegatildep\|_{L^\infty_t;L^p_x}^2
        \\
        &\lesssim \|\omegahatm\|_{L^\infty_t;L^p_x}^2 + \|\omegahatp\|_{L^\infty_t;L^p_x}^2 + \Lambda^2
        \\
        &\lesssim \|\omega_0\|_{p}^2 + \Ra^2 + \Lambda^2
        \\
        &\lesssim \|\omega_0\|_{p}^2 + \Ra^2 + \epsilon^2 \|\omega\|_{L^\infty_t;L^p_x}^2
        \label{bnisubnisfng}
        \\
        &\qquad + \left(1+\|\kappa\|_\infty+\epsilon^{-\frac{p}{p-2}}\|\alpha+\kappa\|_\infty^{\frac{p}{p-2}}\right)^2\|\alpha+\kappa\|_\infty^2\|u\|_{L^\infty_t;L^2_x}^2.
    \end{align}
    Choosing $\epsilon$ sufficiently small, we can compensate the third term on the right-hand side of \eqref{bnisubnisfng} implying
    \begin{align}
        \|\omega\|_{L^\infty_t;L^p_x}^2 
        &\lesssim \|\omega_0\|_{p}^2 + \Ra^2 + \left(1+\|\kappa\|_\infty+\|\alpha+\kappa\|_\infty^{\frac{p}{p-2}}\right)^2\|\alpha+\kappa\|_\infty^2\|u\|_{L^\infty_t;L^2_x}^2
        \\
        &\lesssim \|\omega_0\|_{p}^2 + \Ra^2 + C_{\alpha,\kappa,p}^2\|u\|_{L^\infty_t;L^2_x}^2,
    \end{align}
    where $C_{\alpha,\kappa,p} = \left(1+\|\kappa\|_\infty+\|\alpha+\kappa\|_\infty^{\frac{p}{p-2}}\right)\|\alpha+\kappa\|_\infty$ and finally using \eqref{energy_bound_lemma_L2_bounded} yields
    \begin{align}
        \|\omega\|_{L^\infty(0,\tmax;L^p(\Omega))}
        &\lesssim \|\omega_0\|_p + \Ra + C_{\alpha,\kappa,p} \|u\|_{L^\infty_t;L^2_x}
        \\
        &\lesssim \|\omega_0\|_p + C_{\alpha,\kappa,p} \|u_0\|_2 + (1+C_{\alpha,\kappa,p}) (1+\|\alpha^{-1}\|_\infty)\Ra.
    \end{align}
    Note that this bound holds independent of $\tmax$ concluding the proof.
\end{proof}

As remarked in Section \ref{section_navSlip_boundary_gradients}, boundary terms of type $\int_{\gammam\cup\gammap} u \cdot \nabla p$, naturally arise when deriving second order estimates for fluids with Navier-slip boundary conditions. These terms can be bounded using Lemma \ref{lemma_u_cdot_nabla_p_on_boundary}, which demands estimates on the pressure gradient. The following Lemma provides such a result. In contrast to the results in \cite{bleitnerNobili2024Bounds} we are not deriving a Poisson equation and Neumann boundary conditions for the pressure, but rather directly prove the results by testing \eqref{navier_stokes} with $\nabla p$. Both methods yield the same results after integration by parts.

\begin{remark}[Difference to \cite{bleitner2024Scaling}]
\label{Remark:difference_pressure_bound}
The difference in the result of Theorem \ref{theorem_flat} and \cite{bleitner2024Scaling} lies in the variation of estimating the nonlinear term in the following pressure bound.
Here, in essence the argumentation is Hölder's inequality and Gagliardo-Nieren\-berg interpolation imply
\begin{align}
    \int_\Omega \nabla p \cdot (u\cdot\nabla)u &\leq \|p\|_{H^1} \|u\|_\infty \|\nabla u\|_2
    \\
    &\leq \|p\|_{H^1} \left(\|\nabla u\|_p^\rho\|u\|_2^{1-\rho} + \|u\|_2)\right) \|\nabla u\|_2,\label{jntrwekrjwerwjk}
\end{align}
while, for the flat system, the same term could be estimated using integration by parts, Hölder's inequality, and Sobolev embedding as
\begin{align}
    \int_\Omega \nabla p \cdot (u\cdot\nabla)u &= -\int_\Omega p \nabla u\colon \nabla u^T + \int_{\partial\Omega} p n \cdot (u\cdot \nabla )u 
    \\
    &= -\int_\Omega p \nabla u\colon \nabla u^T + \int_{\partial\Omega} p n_2 u_1 \partial_1 u_2
    \\
    &=-\int_\Omega p \nabla u\colon \nabla u^T \leq \|p\|_4 \|\nabla u\|_4 \|\nabla u\|_2
    \\
    &\lesssim \|p\|_{H^1}\|\nabla u\|_4 \|\nabla u\|_2,
\end{align}
where $\colon$ denotes the tensor contraction, i.e. $\nabla u\colon \nabla u^T = \partial_i u_j \partial_j u_i$, which reflects the approach in \cite{bleitner2024Scaling}. The reason for opting for the former estimate is twofold. The first reason is to demonstrate that these estimates are flexible and the analysis can be tuned to the regime of interest. The second reason for using this approach is the exchange of higher-order bounds for lower-order ones, which is better suited in the regime of big $\alpha$. Note that using $p>4$ on the right-hand side results in smaller $\rho$, which is even more suitable for the regime of big slip coefficientsm.
\end{remark}

\begin{lemma}[Pressure Bound]
    \label{lemma_pressure_bound}
    Let $\Omega$ be $C^{2,1}$, $u_0\in W^{1,4}(\Omega)$ and $0<\alpha \in W^{1,\infty}(\gammam\cup\gammap)$. Then
    \begin{align}
        \|p\|_{H^1}\lesssim \left(\frac{1}{\Pr} \|u\|_\infty + \|(\alpha+\kappa)n\|_{W^{1,\infty}}\right) \|u\|_{H^1} + \Ra,
    \end{align}
    where the implicit constant only depends on $\Gamma$ and the Lipschitz constant of the boundary.
\end{lemma}

\begin{proof}
    Testing \eqref{navier_stokes} with $\nabla p$, one gets
    \begin{align}
        \|\nabla p\|_2^2 = -\frac{1}{\Pr} \left(\int_{\Omega} u_t \cdot \nabla p + \int_{\Omega} \nabla p \cdot (u\cdot \nabla) u \right) + \int_{\Omega} \nabla p \cdot \Delta u + \Ra \int \partial_2 p \theta.
        \\
        \label{abniuvab}
    \end{align}
    Note that the first term on the right-hand side of \eqref{abniuvab} vanishes under partial integration as
    \begin{align}
        \int_{\Omega} u_t \cdot \nabla p
        = \int_{\gammam\cup\gammap} u_t\cdot n p - \int_{\Omega} p \nabla \cdot u_t
        = 0,\label{qwerqwe}
    \end{align}
    where in the last equality we used that, after deriving \eqref{incompressibility} and \eqref{no_penetration_bc} with respect to time, $\nabla \cdot u_t=0$ and $u_t\cdot n=0$.
    To estimate the third term on the right-hand side of \eqref{abniuvab} notice that due to \eqref{id_Delta_u_is_nablaPerp_omega} and integration by parts
    \begin{align}
        \int_{\Omega} \nabla p\cdot \Delta u
        &= \int_{\Omega} \nabla p\cdot \nabla^\perp \omega
        \\
        &= \int_{\gammam\cup\gammap} \omega n^\perp \cdot \nabla p- \int_{\Omega} \omega \nabla^\perp\cdot \nabla p
        \\
        &= \int_{\gammam\cup\gammap} \omega \tau \cdot \nabla p,
        \label{bnigdnsijbfg}
    \end{align}
    where the bulk term vanished because of $\nabla^\perp\cdot \nabla p = 0$. Due to \eqref{vorticity_bc} we can further estimate \eqref{bnigdnsijbfg} as
    \begin{align}
        \int_{\Omega} \nabla p\cdot \Delta u
        &= -2\int_{\gammam\cup\gammap} (\alpha+\kappa)  u_\tau \tau \cdot \nabla p
        = -2\int_{\gammam\cup\gammap} (\alpha+\kappa) u \cdot \nabla p \qquad 
        \label{bnxcjicbnc}
    \end{align}
    and by Lemma \ref{lemma_u_cdot_nabla_p_on_boundary} the left-hand side can be bounded by
    \begin{align}
        \left|\int_{\Omega} \nabla p\cdot \Delta u\right| \lesssim \|(\alpha+\kappa)n\|_{W^{1,\infty}} \|u\|_{H^1}\|p\|_{H^1}.\label{ycxvyxvx}
    \end{align}
    Using \eqref{qwerqwe}, \eqref{ycxvyxvx}, and Hölder's inequality in order to estimate the second and fourth term on the right-hand side of \eqref{abniuvab} we find
    \begin{align}
        \|\nabla p\|_2^2 \lesssim \frac{1}{\Pr} \|u\|_\infty \|p\|_{H^1}\|u\|_{H^1} + \|(\alpha+\kappa)n\|_{W^{1,\infty}} \|u\|_{H^1}\|p\|_{H^1} + \Ra \|p\|_{H^1}\|\theta\|_2.\label{yxcvybxb}
    \end{align}
    As $p$ is only defined up to a constant, we choose it such that $p$ is average free and therefore Poincaré's inequality and \eqref{yxcvybxb} imply
    \begin{align}
        \|p\|_{H^1}^2&\lesssim \|\nabla p\|_2^2
        \\
        &\lesssim \frac{1}{\Pr} \|u\|_\infty \|p\|_{H^1}\|u\|_{H^1} + \|(\alpha+\kappa)n\|_{W^{1,\infty}} \|u\|_{H^1}\|p\|_{H^1} + \Ra \|p\|_{H^1}\|\theta\|_2,
    \end{align}
    which after dividing by $\|p\|_{H^1}$ and using the maximum principle \eqref{maximum_principle} for $\theta$, yields
    \begin{align}
        \|p\|_{H^1}
        &\lesssim \left(\frac{1}{\Pr} \|u\|_\infty + \|(\alpha+\kappa)n\|_{W^{1,\infty}}\right) \|u\|_{H^1} + \Ra,
    \end{align}
    proving the claim.
\end{proof}

Having tools at hand that allow us to cope with the nonlinearity and the boundary terms we are able to prove the long-time average bound for the $H^2$ norm of $u$.

\begin{lemma}[Long time \texorpdfstring{$H^2$}{H2} bound]
    \label{lemma_long_time_H2_bound}
    Let $\Omega$ be $C^{2,1}$, $u_0\in W^{1,4}(\Omega)$ and $0<\alpha\in W^{1,\infty}(\gammam\cup\gammap)$. Then
    \begin{align}
        \langle\|u\|_{H^2}^2\rangle
        &\lesssim  C_1^2 C_2\left(1+\Pr^{-1}\left(\|u_0\|_{W^{1,4}} +  C_2\Ra\right)\right) \Nu \Ra
        \\
        &\qquad + C_3 C_2^\frac{1}{2} \Nu^\frac{1}{2} \Ra^\frac{3}{2} + C_2^\frac{1}{2} \Nu \Ra^\frac{3}{2},
    \end{align}
    where
    \begin{align}
        C_1(\alpha,\kappa)&=1+\|\alpha\|_{W^{1,\infty}}+ \|\kappa\|_{W^{1,\infty}}+\|\alpha\|_\infty^3+\|\kappa\|_\infty^3
        \\
        C_2(\alpha,\kappa)&= 1+\|\alpha^{-1}(1+|\kappa|)\|_\infty
        \\
        C_3(\alpha,\kappa)&= \|\alpha+\kappa\|_\infty.
    \end{align}
    and the implicit constant only depends on $\Gamma$ and the Lipschitz constant of the boundary.
\end{lemma}

\begin{proof}
    Testing \eqref{vorticity_equation} with $\omega$ one finds
    \begin{align}
        \frac{1}{2\Pr}\frac{d}{dt}\|\omega\|_2^2 = \frac{1}{\Pr}\int_{\Omega} \omega u\cdot \nabla \omega + \int_\Omega \omega \Delta \omega + \Ra\int_\Omega \omega\partial_1 \theta.
        \label{yvcxvbx}
    \end{align}
    The first term on the right-hand side of \eqref{yvcxvbx} vanishes under partial integration as
    \begin{align}
        2\int_{\Omega} \omega u\cdot\nabla \omega = \int_\Omega u\cdot\nabla \omega^2 = \int_{\gammam\cup\gammap} \omega^2 u\cdot n - \int_\Omega \omega^2 \nabla\cdot u=0,
        \label{bxcvbcv}
    \end{align}
    where in the last equality we exploit \eqref{no_penetration_bc} and \eqref{incompressibility}. For the second term on the right-hand side of \eqref{yvcxvbx} integration by parts yields
    \begin{align}
        \int_{\Omega} \omega\Delta \omega =  - \|\nabla \omega\|_2^2 + \int_{\gammam\cup\gammap} \omega n\cdot\nabla \omega.
        \label{xbxcbv}
    \end{align}
    In order to estimate the boundary term in \eqref{xbxcbv} note that by \eqref{vorticity_bc} and \eqref{id_Delta_u_is_nablaPerp_omega}
    \begin{align}
        \int_{\gammam\cup\gammap} \omega n\cdot\nabla \omega
        &= \int_{\gammam\cup\gammap} \omega n\cdot\nabla \omega
        = \int_{\gammam\cup\gammap} \omega \tau\cdot\nabla^\perp \omega
        \\
        &= -2\int_{\gammam\cup\gammap} (\alpha+\kappa)u \cdot\Delta u.
    \end{align}
    Inserting \eqref{navier_stokes} one finds
    \begin{align}
        \frac{1}{2}\int_{\gammam\cup\gammap} \omega n\cdot\nabla \omega
        &= -\int_{\gammam\cup\gammap} (\alpha+\kappa)u \cdot\Delta u
        \\
        &= -\int_{\gammam\cup\gammap} (\alpha+\kappa)u \cdot \left(\frac{1}{\Pr}(u_t+u\cdot\nabla u)+\nabla p-\Ra \theta e_2\right)
        \\
        &= -\frac{1}{2\Pr}\frac{d}{dt}\int_{\gammam\cup\gammap} (\alpha+\kappa)u_\tau^2-\frac{1}{2\Pr}\int_{\gammam\cup\gammap} (\alpha+\kappa)u\cdot\nabla u_\tau^2
        \\
        &\qquad - \int_{\gammam\cup\gammap} (\alpha+\kappa)u\cdot \nabla p + \Ra \int_{\gammam} (\alpha+\kappa)u_2
        \label{nbvbcnvnbc}
    \end{align}
    Combining \eqref{yvcxvbx}, \eqref{bxcvbcv}, \eqref{xbxcbv} and \eqref{nbvbcnvnbc} yields
    \begin{align}
        \frac{1}{2\Pr}\frac{d}{dt}&\left(\|\omega\|_2^2 + 2\int_{\gammam\cup\gammap} (\alpha+\kappa) u_\tau^2\right) + \|\nabla\omega\|_2^2
        \\
        &= - \frac{1}{\Pr} \int_{\gammam\cup\gammap} (\alpha+\kappa) u\cdot \nabla u_\tau^2 - 2\int_{\gammam\cup\gammap} (\alpha+\kappa) u\cdot\nabla p
        \\
        &\qquad + 2 \Ra \int_{\gammam} (\alpha+\kappa)u_2 + \Ra \int_{\Omega} \omega \partial_1 \theta
        \label{cvjnbjbcvxbni}
    \end{align}
    and using Lemma \ref{lemma_u_cdot_nabla_p_on_boundary}, the trace theorem and Hölder's inequality
    \begin{align}
        \frac{1}{2\Pr}\frac{d}{dt}&\left(\|\omega\|_2^2 + 2\int_{\gammam\cup\gammap} (\alpha+\kappa) u_\tau^2\right) + \|\nabla\omega\|_2^2
        \\
        &\lesssim  \|(\alpha+\kappa)n\|_{W^{1,\infty}} \left(\frac{1}{\Pr}\|u^2\|_{H^1} +\|p\|_{H^1}\right)\|u\|_{H^1} 
        \\
        &\qquad + \|\alpha+\kappa\|_\infty \Ra \|u\|_{H^1} + \Ra \|\omega\|_2 \|\nabla\theta\|_2
    \end{align}
    By Lemma \ref{lemma_grad_ids} and \ref{lemma_pressure_bound}, $\|\omega\|_2\leq\|\nabla u\|_2$ and Hölder's inequality one gets
    \begin{align}
        \frac{1}{\Pr}\frac{d}{dt}&\left(\|\D u\|_2^2 + \int_{\gammam\cup\gammap} \alpha u_\tau^2\right) + \|\nabla\omega\|_2^2
        \\
        &\lesssim  \|(\alpha+\kappa)n\|_{W^{1,\infty}} \left(\frac{1}{\Pr}\|u^2\|_{H^1} +\|p\|_{H^1}\right) \|u\|_{H^1}
        \\
        &\qquad + \left(\|\alpha+\kappa\|_\infty + \|\nabla\theta\|_2\right)\Ra\|u\|_{H^1}
        \\
        &\lesssim  \|(\alpha+\kappa)n\|_{W^{1,\infty}}\left(\|(\alpha+\kappa)n\|_{W^{1,\infty}}+ \frac{1}{\Pr}\|u\|_\infty\right)\|u\|_{H^1}^2
        \\
        &\qquad + \left(\|\alpha+\kappa\|_{\infty} + \|\nabla\theta\|_2\right)\Ra\|u\|_{H^1}.
        \label{vbnccnbbvn}
    \end{align}
    Note that by Lemma \ref{lemma_elliptic_regularity_periodic_domain}
    \begin{align}
        \|u\|_{H^2} &\lesssim \|\omega\|_{H^1} + \|\kappa\|_\infty \|\omega\|_2 + (1+\|\kappa\|_{W^{1,\infty}}+\|\kappa\|_\infty^2) \|u\|_2
        \\
        &\lesssim \|\nabla\omega\|_{2} + (1 + \|\kappa\|_\infty) \|u\|_{H^1} + (1+\|\kappa\|_{W^{1,\infty}}+\|\kappa\|_\infty^2) \|u\|_2
        \\
        &\lesssim \|\nabla\omega\|_{2}+ (1+\|\kappa\|_{W^{1,\infty}}+\|\kappa\|_\infty^2) \|u\|_{H^1},
        \label{bvncnbcb}
    \end{align}
    where in the last estimate we used Young's inequality. Combining \eqref{vbnccnbbvn} and \eqref{bvncnbcb} there exists a constant $C>0$ depending only on $\Gamma$ and the Lipschitz constant of the boundary such that
    \begin{align}
        \frac{1}{\Pr}\frac{d}{dt}&\left(\|\D u\|_2^2 + \int_{\gammam\cup\gammap} \alpha u_\tau^2\right) + C\|u\|_{H^2}^2
        \\
        &\lesssim \left(1+\|\kappa\|_{W^{1,\infty}} + \|\kappa\|_\infty^2+\|(\alpha+\kappa)n\|_{W^{1,\infty}}\right)^2\|u\|_{H^1}^2
        \\
        &\qquad + \frac{1}{\Pr} \|(\alpha+\kappa)n\|_{W^{1,\infty}}\|u\|_\infty\|u\|_{H^1}^2
        \\
        &\qquad + \left(\|\alpha+\kappa\|_{\infty} + \|\nabla\theta\|_2\right)\Ra\|u\|_{H^1}.
        \label{cxbvcxb}
    \end{align}
    For $2<p\in 2\mathbb{N}$ by Sobolev embedding, Lemma \ref{lemma_elliptic_regularity_periodic_domain}, Lemma \ref{lemma_vorticity_bound} and Lemma \ref{lemma_energy_bound}
    \begin{align}
        \|u\|_\infty &\lesssim \|u\|_{W^{1,p}} \lesssim \|\omega\|_p + \left(1+\|\kappa\|^{2-\frac{2}{p}}\right) \|u\|_2
        \\
        &\lesssim \left(1+\|\kappa\|^{\frac{2p-2}{p}}+\|\alpha+\kappa\|_\infty^\frac{2p-2}{p-2}\right) \left(\|u_0\|_{W^{1,p}} +  (1+\|\alpha^{-1}\|_\infty)\Ra\right)\qquad
        \label{cvnbvcnbcvb}
    \end{align}
    where we additionally used Young's and Hölder's inequality to simplify the prefactors and initial data with the cost of slightly worsening the estimate. Combining \eqref{cxbvcxb} and \eqref{cvnbvcnbcvb} yields
    \begin{align}
        \frac{1}{\Pr}\frac{d}{dt}&\left(\|\D u\|_2^2 + \int_{\gammam\cup\gammap} \alpha u_\tau^2\right) + C\|u\|_{H^2}^2
        \\
        &\lesssim \left(1+\|\kappa\|_{W^{1,\infty}} + \|\kappa\|_\infty^2+\|(\alpha+\kappa)n\|_{W^{1,\infty}}\right)^2\|u\|_{H^1}^2
        \\
        &\qquad + \frac{1}{\Pr} \|(\alpha+\kappa)n\|_{W^{1,\infty}}\left(1+\|\kappa\|^{\frac{2p-2}{p}}+\|\alpha+\kappa\|_\infty^\frac{2p-2}{p-2}\right) 
        \\
        &\qquad\qquad \cdot\left(\|u_0\|_{W^{1,p}} +  (1+\|\alpha^{-1}\|_\infty)\Ra\right)\|u\|_{H^1}^2
        \\
        &\qquad + \left(\|\alpha+\kappa\|_{\infty} + \|\nabla\theta\|_2\right)\Ra\|u\|_{H^1}.
        \label{xcvbxcbcv}
    \end{align}
    Note that there are at most quadratic terms depending on time on the right-hand side of \eqref{xcvbxcbcv}, allowing us to take the long-time average under which the time derivative on the left-hand side vanishes as $\|u\|_{H^1}^2\lessapprox \|\D u\|_2^2 + \int_{\gammam\cup\gammap} \alpha u_\tau^2 \lessapprox \|u\|_{H^1}^2$ by Lemma \ref{lemma_coercivity} and $\|u\|_{H^1}^2$ is uniformly bounded in time by Remark \ref{remark_Linfty_W1p_bound}. Therefore, using $\langle fg \rangle\lesssim \langle f^2\rangle^\frac{1}{2}\langle g^2\rangle^\frac{1}{2}$ due to Young's inequality, Lemma \ref{lemma_energy_bound} and Lemma \ref{lemma_nusselt_number_representations} yield
    \begin{align}
        \langle\|u\|_{H^2}^2\rangle
        &\lesssim \left(1+\|\kappa\|_{W^{1,\infty}} + \|\kappa\|_\infty^2+\|(\alpha+\kappa)n\|_{W^{1,\infty}}\right)^2\langle\|u\|_{H^1}^2\rangle
        \\
        &\qquad + \frac{1}{\Pr} \|(\alpha+\kappa)n\|_{W^{1,\infty}}\left(1+\|\kappa\|^{\frac{2p-2}{p}}+\|\alpha+\kappa\|_\infty^\frac{2p-2}{p-2}\right) 
        \\
        &\qquad\qquad \cdot\left(\|u_0\|_{W^{1,p}} +  (1+\|\alpha^{-1}\|_\infty)\Ra\right)\langle\|u\|_{H^1}^2\rangle
        \\
        &\qquad + \left(\|\alpha+\kappa\|_{\infty} + \langle\|\nabla\theta\|_2^2\rangle^\frac{1}{2}\right)\Ra\langle\|u\|_{H^1}^2\rangle^\frac{1}{2}
        \\
        &\lesssim \left(1+\|\kappa\|_{W^{1,\infty}}^2 + \|\kappa\|_\infty^4+\|(\alpha+\kappa)n\|_{W^{1,\infty}}^2\right)\left(1+\left\|\tfrac{1+|\kappa|}{\alpha}\right\|_\infty\right) \Nu \Ra
        \\
        &\qquad + \|(\alpha+\kappa)n\|_{W^{1,\infty}}\left(1+\|\kappa\|^{\frac{2p-2}{p}}+\|\alpha+\kappa\|_\infty^\frac{2p-2}{p-2}\right) \left(1+\left\|\tfrac{1+|\kappa|}{\alpha}\right\|_\infty\right)
        \\
        &\qquad\qquad \cdot\frac{\|u_0\|_{W^{1,p}} +  (1+\|\alpha^{-1}\|_\infty)\Ra}{\Pr} \Nu \Ra
        \\
        &\qquad + \|\alpha+\kappa\|_{\infty} \left(1+\left\|\tfrac{1+|\kappa|}{\alpha}\right\|_\infty\right)^\frac{1}{2} \Nu^\frac{1}{2} \Ra^\frac{3}{2}
        \\
        &\qquad + \left(1+\left\|\tfrac{1+|\kappa|}{\alpha}\right\|_\infty\right)^\frac{1}{2} \Nu \Ra^\frac{3}{2}
    \end{align}
    In order to simplify the estimate we use $p=4$, the lowest value for which the analysis works, Young's inequality in multiple ways and
    \begin{align}
        \|(\alpha+\kappa)n\|_{W^{1,\infty}}&\lesssim \|n\|_\infty \|\alpha+\kappa\|_{W^{1,\infty}} + \|\alpha+\kappa\|_\infty \|n\|_{W^{1,\infty}}
        \\
        &\lesssim\|\alpha+\kappa\|_{W^{1,\infty}} + \|\alpha+\kappa\|_\infty \|\kappa\|_\infty
    \end{align}
    to get
    \begin{align}
        \langle\|u\|_{H^2}^2\rangle
        &\lesssim \left(1+\|\alpha\|_{W^{1,\infty}}+ \|\kappa\|_{W^{1,\infty}}+\|\alpha\|_\infty^3+\|\kappa\|_\infty^3\right)^2\left(1+\left\|\tfrac{1+|\kappa|}{\alpha}\right\|_\infty\right)
        \\
        &\qquad\qquad \cdot\left(1+\frac{\|u_0\|_{W^{1,4}} +  (1+\|\alpha^{-1}\|_\infty)\Ra}{\Pr}\right) \Nu \Ra
        \\
        &\qquad + \|\alpha+\kappa\|_{\infty} \left(1+\left\|\tfrac{1+|\kappa|}{\alpha}\right\|_\infty\right)^\frac{1}{2}\Nu^\frac{1}{2} \Ra^\frac{3}{2}
        \\
        &\qquad + \left(1+\left\|\tfrac{1+|\kappa|}{\alpha}\right\|_\infty\right)^\frac{1}{2} \Nu \Ra^\frac{3}{2}
        \\
        &\lesssim C_1^2 C_2\left(1+\Pr^{-1}\left(\|u_0\|_{W^{1,4}} +  C_2\Ra\right)\right) \Nu \Ra
        \\
        &\qquad + C_3 C_2^\frac{1}{2} \Nu^\frac{1}{2}\Ra^\frac{3}{2} + C_2^\frac{1}{2} \Nu\Ra^\frac{3}{2}
    \end{align}
    concluding the proof.
\end{proof}

\subsection{Proof of the Theorems}
\label{section:direct_method_application}

With all this preparation at hand, we are able to combine the results and prove the Theorms \ref{theorem_1_2_bound} and \ref{theorem_main_theorem_rb_curved}.

\begin{proof}[\hypertarget{proof_theorem_1_2_bound}{Proof of Theorem} \ref{theorem_1_2_bound}]
\\
    Lemma \ref{lemma_nusselt_bounded_by_Hx} and Lemma \ref{lemma_energy_bound} imply
    \begin{align}
        \Nu &\lesssim \delta^\frac{1}{2} \langle \|\nabla u\|_2^2 \rangle^\frac{1}{2} + \delta^{-\frac{1}{2}}\Nu^\frac{1}{2}
        \\
        &\lesssim \delta^\frac{1}{2} \left(1+\left\|\tfrac{1+|\kappa|}{\alpha}\right\|_\infty\right)^\frac{1}{2} \Nu^\frac{1}{2} \Ra^\frac{1}{2} + \delta^{-\frac{1}{2}}\Nu^\frac{1}{2}
    \end{align}
    for any $d>\delta>0$, which after dividing by $\Nu^\frac{1}{2}$ and squaring yields
    \begin{align}
        \Nu
        &\lesssim \delta C_2 \Ra + \delta^{-1},
    \end{align}
    where $C_2=1+\left\|\tfrac{1+|\kappa|}{\alpha}\right\|_\infty$. Choosing
    \begin{align}
        \delta = C_2^{-\frac{1}{2}} \min\left\lbrace \Ra^{-\frac{1}{2}},d\right\rbrace
    \end{align}
    results in
    \begin{align}
        \Nu &\lesssim C_2^\frac{1}{2}\Ra^{\frac{1}{2}},
    \end{align}
    where in the case of $d \leq \Ra^{-\frac{1}{2}}$ we used
    \begin{align}
        \Nu \lesssim C_2^\frac{1}{2} (d\Ra + d^{-1}) \leq C_2^\frac{1}{2} (\Ra^\frac{1}{2} + d^{-1}) \leq C(d) C_2^\frac{1}{2}\Ra^\frac{1}{2}
    \end{align}
    proving the first claim.

    The same proof, using \eqref{energy_bound_lemma_kappa_leq_2alpha_case} instead of \eqref{energy_bound_lemma_H1_bound}, yields \eqref{main_result_1_2_bound_kappa_leq_2alpha}.

\end{proof}

\begin{proof}[\hypertarget{proof_theorem_main_theorem_rb_curved}{Proof of Theorem} \ref{theorem_main_theorem_rb_curved}]
\\
    The proof follows the same strategy as the one of Theorem \ref{theorem_1_2_bound}, but instead uses the higher order bounds. By Lemma \ref{lemma_nusselt_bounded_by_Hx} the Nusselt number is bounded by
    \begin{align}
        \Nu
        &\lesssim \delta \langle \|u\|_{H^1}^2 \rangle^\frac{1}{4} \langle \|u\|_{H^2}^2\rangle^\frac{1}{4} + \delta^{-\frac{1}{2}}\Nu^\frac{1}{2}
    \end{align}
    for any $\delta < \frac{d}{2}$ and using the long time bounds of $\|u\|_{H^1}$ and $\|u\|_{H^2}$, i.e. Lemma \ref{lemma_energy_bound} respectively Lemma \ref{lemma_long_time_H2_bound}, one gets
    \begin{align}
        \Nu
        &\lesssim \delta  \Big[C_1^\frac{1}{2} C_2^\frac{1}{2}\left(1+\Pr^{-\frac{1}{4}}\left(\|u_0\|_{W^{1,4}}^\frac{1}{4} +  C_2^\frac{1}{4}\Ra^\frac{1}{4}\right)\right) \Nu^\frac{1}{2} \Ra^\frac{1}{2}
        \\
        &\qquad + C_2^\frac{3}{8} C_3^\frac{1}{4} \Nu^\frac{3}{8} \Ra^\frac{5}{8}
        + C_2^\frac{3}{8} \Nu^\frac{1}{2} \Ra^\frac{5}{8}\Big]
        + \delta^{-\frac{1}{2}}\Nu^\frac{1}{2},
        \label{ertjnerkjt}
    \end{align}
    where
    \begin{align}
        C_1(\alpha,\kappa)&=1+\|\alpha\|_{W^{1,\infty}}+ \|\kappa\|_{W^{1,\infty}}+\|\alpha\|_\infty^3+\|\kappa\|_\infty^3
        \\
        C_2(\alpha,\kappa)&= 1+\left\|\tfrac{1+|\kappa|}{\alpha}\right\|_\infty
        \\
        C_3(\alpha,\kappa)&= \|\alpha+\kappa\|_\infty.
    \end{align}

    In order to optimize $\delta$ we distinguish between two cases. Note that by similar consideration as in the proof of Theorem \ref{theorem_1_2_bound}, one can without loss of generality assume $\Ra^\frac{5}{12} \geq \frac{2}{d}$ such that the subsequent choices imply $\delta<\frac{d}{2}$. 
    \begin{itemize}
        \item[a)]
            Assume that the third term in the squared brackets is dominating, i.e.
            \begin{align}
                C_3^\frac{1}{4} &\geq C_1^\frac{1}{2} C_2^\frac{1}{8}\left(1+\Pr^{-\frac{1}{4}}\left(\|u_0\|_{W^{1,4}}^\frac{1}{4} +  C_2^\frac{1}{4}\Ra^\frac{1}{4}\right)\right) \Nu^\frac{1}{8} \Ra^{-\frac{1}{8}} +\Nu^\frac{1}{8}.
            \end{align}
            Then \eqref{ertjnerkjt} implies
            \begin{align}
                \Nu
                &\lesssim \delta  C_2^\frac{3}{8} C_3^\frac{1}{4}\Nu^\frac{3}{8} \Ra^\frac{5}{8}
                + \delta^{-\frac{1}{2}}\Nu^\frac{1}{2}
            \end{align}
            and optimizing in $\delta$ by setting
            \begin{align}
                \delta = \Nu^\frac{1}{12}\left(C_2^\frac{3}{8} C_3^\frac{1}{4} \Ra^\frac{5}{8}\right)^{-\frac{2}{3}}
            \end{align}
            yields
            \begin{align}
                \Nu
                \lesssim \Nu^\frac{11}{24}\left( C_2^\frac{3}{8} C_3^\frac{1}{4} \Ra^\frac{5}{8}\right)^\frac{1}{3}
                \lesssim \Nu^\frac{11}{24}  C_2^\frac{3}{24} C_3^\frac{1}{12} \Ra^\frac{5}{24}
            \end{align}
            and after division by $\Nu^\frac{11}{24}$ and exponentiation 
            \begin{align}
                \Nu \lesssim C_2^\frac{3}{13}C_3^\frac{2}{13}\Ra^\frac{5}{13}.
            \end{align}
        \item[b)]
            If instead the third term in the squared bracket is dominated by the others, i.e.
            \begin{align}
                C_3^\frac{1}{4} &\leq C_1^\frac{1}{2} C_2^\frac{1}{8}\left(1+\Pr^{-\frac{1}{4}}\left(\|u_0\|_{W^{1,4}}^\frac{1}{4} +  C_2^\frac{1}{4}\Ra^\frac{1}{4}\right)\right) \Nu^\frac{1}{8} \Ra^{-\frac{1}{8}} +\Nu^\frac{1}{8},
            \end{align}
            \eqref{ertjnerkjt} yields
            \begin{align}
                \Nu
                &\lesssim \delta \Nu^\frac{1}{2}  \Big[C_1^\frac{1}{2} C_2^\frac{1}{2}\left(1+\Pr^{-\frac{1}{4}}\left(\|u_0\|_{W^{1,4}}^\frac{1}{4} +  C_2^\frac{1}{4}\Ra^\frac{1}{4}\right)\right) \Ra^\frac{1}{2} 
                \\
                &\qquad + C_2^\frac{3}{8} \Ra^\frac{5}{8}\Big]
                + \delta^{-\frac{1}{2}}\Nu^\frac{1}{2}.
            \end{align}
            Again division by $\Nu^\frac{1}{2}$ and squaring implies
            \begin{align}
                \Nu \lesssim \delta^2 \Big[C_1^\frac{1}{2} C_2^\frac{1}{2}\left(1+\Pr^{-\frac{1}{4}}\left(\|u_0\|_{W^{1,4}}^\frac{1}{4} +  C_2^\frac{1}{4}\Ra^\frac{1}{4}\right)\right) \Ra^\frac{1}{2} + C_2^\frac{3}{8} \Ra^\frac{5}{8}\Big]^2+\delta^{-1}
            \end{align}
            and optimizing in $\delta$ by setting
            \begin{align}
                \delta=\Big[C_1^\frac{1}{2} C_2^\frac{1}{2}\left(1+\Pr^{-\frac{1}{4}}\left(\|u_0\|_{W^{1,4}}^\frac{1}{4} +  C_2^\frac{1}{4}\Ra^\frac{1}{4}\right)\right) \Ra^\frac{1}{2} + C_2^\frac{3}{8} \Ra^\frac{5}{8}\Big]^{-\frac{2}{3}}
            \end{align}
            results in
            \begin{align}
                \Nu &\lesssim \Big[C_1^\frac{1}{2} C_2^\frac{1}{2}\left(1+\Pr^{-\frac{1}{4}}\left(\|u_0\|_{W^{1,4}}^\frac{1}{4} +  C_2^\frac{1}{4}\Ra^\frac{1}{4}\right)\right) \Ra^\frac{1}{2} + C_2^\frac{3}{8} \Ra^\frac{5}{8}\Big]^\frac{2}{3}
                \\
                &\lesssim C_1^\frac{1}{3}C_2^\frac{1}{3}\left(1+\Pr^{-\frac{1}{6}}\left(\|u_0\|_{W^{1,4}}^\frac{1}{6} +  C_2^\frac{1}{6}\Ra^\frac{1}{6}\right)\right) \Ra^\frac{1}{3} + C_2^\frac{1}{4} \Ra^\frac{5}{12}.
            \end{align}
    \end{itemize}
    Combining the different estimates proves
    \begin{align}
        \Nu
        &\lesssim C_1^\frac{1}{3}C_2^\frac{1}{3}\left(1+\Pr^{-\frac{1}{6}}\|u_0\|_{W^{1,4}}^\frac{1}{6} \right) \Ra^\frac{1}{3} + C_1^\frac{1}{3}C_2^\frac{1}{3} \Pr^{-\frac{1}{6}}\Ra^\frac{1}{2}
        \\
        &\qquad +  C_2^\frac{3}{13}C_3^\frac{2}{13}\Ra^\frac{5}{13} + C_2^\frac{1}{4} \Ra^\frac{5}{12}.
    \end{align}

\end{proof}

\section{Flat System}
\label{section:flat_system}

\begin{figure}
    \hrule
    \vspace{0.5\baselineskip}
    \begin{center}
        \includegraphics[width=\textwidth]{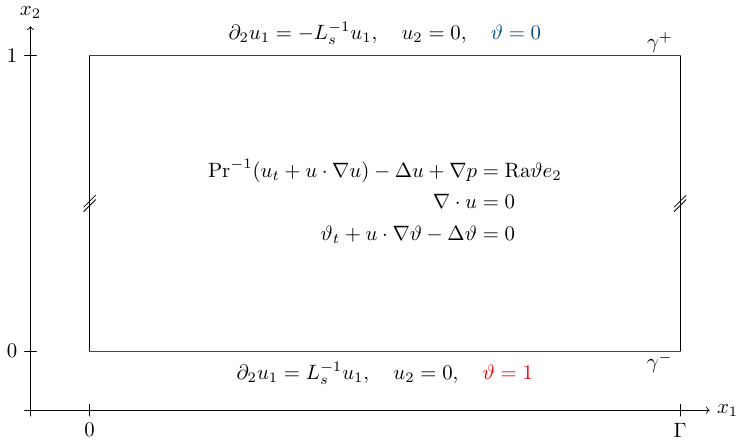}
        \vspace*{-7mm}
        \caption{Overview of the system considered in Section \ref{section:flat_system}, i.e.~with flat boundaries and constant slip coefficient $\alpha = \frac{1}{2L_s}$.}
        \label{fig:rb_overview_flat}
    \end{center}
    \hrule
\end{figure}

In this section, the special case of flat boundaries, i.e. $\hm=0$, $\hp=1$, $\kappa=0$, and constant slip coefficient $\alpha(x)=\frac{1}{2L_s}>0$ is studied. With these choices, the Navier-slip boundary conditions are simplified to
\begin{align}
    0=\tau \cdot (\D u \ n+\alpha u)= \frac{1}{2L_s} \tau_1 u_1 + \frac{1}{2} \tau_1 n_2 (\partial_1 u_2 + \partial_2 u_1)
\end{align}
and using that $u_2=0$ implies $\partial_1 u_2 =0$ on $\gammam\cup\gammap$, they further reduce to
\begin{alignat}{2}
    \partial_2 u_1&=- L_s^{-1} u_1 \qquad & \textnormal{ on }&\gammap
    \\
    \partial_2 u_1&= L_s^{-1} u_1 \qquad & \textnormal{ on }&\gammam.
\end{alignat}
$L_s$ can be viewed as an effective slip length (\cite{miksis1994Slip,bolanos2017Derivation}). An overview of the system is given in Figure \ref{fig:rb_overview_flat}. The subsequent analysis is a slight deviation of \cite{bleitner2024Scaling}.

\subsection{The Direct Method}

The strategy is the same as for Theorem \ref{theorem_1_2_bound} and Theorem \ref{theorem_main_theorem_rb_curved} in Section \ref{section:general_system}. Here we will only state the improvements over the general case and how they affect the result.

Using $\zeta = \np =(0,1)$ in the proof of Lemma \ref{lemma_nusselt_bounded_by_Hx}, the statement reduces to the following.

\begin{lemma}
    Let $u_0\in L^2(\Omega)$. Then
    \begin{align}
        \Nu
        &\lesssim \delta^\frac{1}{2} \langle \|\nabla u\|_2^2 \rangle^\frac{1}{2} + \delta^{-\frac{1}{2}}\Nu^\frac{1}{2}
    \end{align}
    for every $\delta<1$, where the implicit constant only depends on $\Gamma$.
    
    If $u_0\in W^{1,4}(\Omega)$, then
    \begin{align}
        \Nu
        &\lesssim \delta \langle \|\nabla u\|_2^2 \rangle^\frac{1}{4} \langle \|\nabla^2 u\|_2^2\rangle^\frac{1}{4} + \delta^{-\frac{1}{2}}\Nu^\frac{1}{2}
    \end{align}
    for every $\delta <1$, where the implicit constant only depends on $\Gamma$.
    \label{lemma_flat_nusselt_bounded_by_Hx}
\end{lemma}

\subsection{A-Priori Estimates}

Compared to Lemma \ref{lemma_energy_bound}, there is a slight improvement in the prefactor of the $\Ra$ term, as stated in the following.

\begin{lemma}
    Let $u_0\in L^2(\Omega)$. Then
    \begin{align}
        \|u\|_2^2 &\lesssim \|u_0\|_2^2 + (1+L_s)\Ra^2
        \\
        \langle \|\nabla u\|_2^2\rangle &\lesssim \Nu\Ra,
    \end{align}
    where the implicit constant only depends on $\Gamma$.
    \label{lemma_energy_bound_flat}
\end{lemma}

\begin{proof}
    Note that by Lemma \ref{lemma_grad_ids}
    \begin{align}
        2\|\D u\|_2^2 = \|\nabla u\|_2^2,
    \end{align}
    which applied to \eqref{jalksdfjakls} yields
    \begin{align}
        \frac{1}{2\Pr}\frac{d}{dt} \|u\|_{L^2}^2 + \|\nabla u\|_{L^2}^2
        &\leq \frac{1}{2\Pr}\frac{d}{dt} \|u\|_{L^2}^2 + \|\nabla u\|_{L^2}^2 + \frac{1}{L_s} \int_{\gammam\cup\gammap} u_1^2
        \\
        &= \Ra \int_{\Omega} u_2 \theta
        \\
        &\leq \Ra \int_{\Omega} (u_2 -\partial_2) \theta,
        \label{bxjnkbncjkxbn}
    \end{align}
    where in the last estimate we used \eqref{gbnsifbusndb}.
    As in \eqref{coerc_est_u2}, we find
    \begin{align}
        |u|^2(x)\leq 2 |u_\tau|^2(x_1,0) + 2\|\partial_2 u(x_1,\cdot)\|_{L^2(0,1)}^2,
    \end{align}
    which after integrating over $\Omega$ yields
    \begin{align}
        \|u\|_2^2 \lesssim \int_{\gammam\cup\gammap} u_1^2 + \|\nabla u\|_2^2 \leq \max (1,L_s)\left(L_s^{-1}\int_{\gammam\cup\gammap} u_1^2 + \|\nabla u\|_2^2\right).\qquad\label{vycxvjklyv}
    \end{align}
    Combining \eqref{bxjnkbncjkxbn} and \eqref{vycxvjklyv} there exists a constant $C>0$ such that
    \begin{align}
        \frac{1}{2\Pr}\frac{d}{dt} \|u\|_{L^2}^2 + C \min(1,L_s^{-1}) \|u\|_2^2
        &\leq \Ra \int_{\Omega} u_2 \theta\leq \epsilon^{-1}\Ra^2 + \epsilon\|u_2\|_2^2
    \end{align}
    for any $\epsilon>0$, where in the last estimate we used Young's inequality. Choosing $\epsilon$ sufficiently small, we can compensate the $\|u_2\|_2^2$ term and Grönwall's inequality yields
    \begin{align}
        \|u\|_2^2\lesssim \|u_0\|_2^2 + \max(1,L_s)\Ra^2\label{snjngfi}.
    \end{align}
    Similar to before taking the long-time average of \eqref{bxjnkbncjkxbn}, using \eqref{nu_representation_straight} and the fact that $\|u\|_2^2$ is uniformly bounded in time due to \eqref{snjngfi} one finds
    \begin{align}
        \langle\|\nabla u\|_2^2\rangle \leq \Nu\Ra.
        \label{sndifsdg}
    \end{align}    
\end{proof}

Again the strategy is to fall back to the vorticity formulation in order to circumvent the nonlinear term in higher order bounds. Therefore, one needs estimates that allow for exchange between $\nabla u$ and $\omega$, which the following Lemma, corresponding to Lemma \ref{lemma_elliptic_regularity_periodic_domain} for curved domains, provides.

\begin{lemma}
    Let $1<q<\infty$, $k\in \mathbb{N}_0$ and $u\in W^{k+1,q}(\Omega)$. Then
    \begin{align}
        \|\nabla u\|_{W^{k,q}}\lesssim \|\omega\|_{W^{k,q}},
    \end{align}
    where the implicit constant only depends on $\Gamma$ and $q$. Additionally
    \begin{align}
        \|\nabla^2u\|_2 \leq \|\nabla \omega\|_2.
    \end{align}
    \label{lemma_flat_elliptic_regularity}
\end{lemma}

\begin{proof}
    The elliptic regularity estimates, Lemma \ref{lemma_elliptic_regularity_periodic_domain}, are given by the following, where we use a similar approach as \cite[Lemma A.2]{drivas2022Bounds}. As $\phi$ satisfies
    \begin{alignat}{2}
        \Delta \phi &= \omega & \quad \textnormal{ in } &\Omega
        \\
        \phi &= - \dashint_{\Omega} u_1 & \quad \textnormal{ on } &\gammap
        \\
        \phi &= 0 & \quad \textnormal{ on } &\gammam.
    \end{alignat}
    defining $\tilde \phi = \phi-x_2\dashint_{\Omega} u_1$ one has
    \begin{alignat}{2}
        \Delta \tilde \phi &= \omega & \quad \textnormal{ in } &\Omega
        \\
        \tilde \phi &= 0 & \quad \textnormal{ on } &\gammam\cup\gammap.
    \end{alignat}
    Elliptic regularity implies
    \begin{align}
        \|\tilde \phi\|_{W^{k+2,q}} \lesssim \|\omega\|_{W^{k,q}}
    \end{align}
    for any $q\in (1,\infty)$ and therefore
    \begin{align}
        \|\nabla u_2\|_{W^{k,q}} &= \|\nabla \partial_1 \phi\|_{W^{k,q}}
        \leq \|\nabla \partial_1 \tilde \phi\|_{W^{k,q}}
        \lesssim \|\tilde \phi\|_{W^{k+2,q}}
        \lesssim \|\omega\|_{W^{k,q}}.\qquad\label{yxcvynxvm}
    \end{align}
    In order to get an estimate for the first component of the velocity notice that by \eqref{incompressibility} and $\omega = -\partial_2 u_1 + \partial_1 u_2$
    \begin{align}
        \nabla u_1 =
        \begin{pmatrix}
            - \partial_2 u_2\\\partial_1 u_2 - \omega
        \end{pmatrix}.\label{bcxhjv}
    \end{align}
    Combining \eqref{yxcvynxvm} and \eqref{bcxhjv}, we obtain
    \begin{align}
        \|\nabla u\|_{W^{k,q}}
        \lesssim \|\nabla u_1\|_{W^{k,q}} + \|\nabla u_2\|_{W^{k,q}}
        \lesssim \|\nabla u_2\|_{W^{k,q}} + \|\omega\|_{W^{k,q}} \lesssim \|\omega\|_{W^{k,q}}.\label{xcvbncvcihiu}
    \end{align}
    
    Integrating by parts twice, using $n= (0,1)$ on $\gammap$ and $n=(0,-1)$ on $\gammam$, and the cancellation of terms where $i=j=2$, one finds
    \begin{align}
        \|\nabla^2 u\|_2^2&=\int_{\Omega}\partial_i\partial_j u_k \partial_i\partial_j u_k
        \\
        &=\int_{\Omega}\partial_i^2 u_k\partial_j^2 u_k -\int_{\gammam\cup\gammap}\partial_i^2 u_k\partial_j u_k n_j +\int_{\gammam\cup\gammap}\partial_i\partial_j u_k\partial_j u_k n_i
        \\
        &= \|\Delta u\|_2^2 -\int_{\gammam\cup\gammap}\partial_i^2 u_k\partial_2 u_k n_2 +\int_{\gammam\cup\gammap}\partial_2\partial_j u_k\partial_j u_k n_2
        \\
        &= \|\Delta u\|_2^2 -\int_{\gammam\cup\gammap}\partial_1^2 u_1\partial_2 u_1 n_2 +\int_{\gammam\cup\gammap}\partial_2\partial_1 u_1\partial_1 u_1 n_2
        \\
        &\qquad -\int_{\gammam\cup\gammap}\partial_1^2 u_2\partial_2 u_2 n_2 +\int_{\gammam\cup\gammap}\partial_2\partial_1 u_2\partial_1 u_2 n_2.
    \end{align}
    By $u_2=0$ on $\gammam\cup\gammap$, one has $\partial_1u_2=\partial_1^2u_2=0$ and therefore the last two terms vanish. Using the Navier-slip boundary conditions, the identity further simplifies to
    \begin{align}
        \|\nabla^2 u\|_2^2
        &= \|\Delta u\|_2^2 -\int_{\gammam\cup\gammap}\partial_1^2 u_1\partial_2 u_1 n_2 +\int_{\gammam\cup\gammap}\partial_2\partial_1 u_1\partial_1 u_1 n_2
        \\
        &=\|\Delta u\|_2^2 -\int_{\gammam\cup\gammap}\partial_1^2 u_1\partial_2 u_1 n_2 +\int_{\gammam\cup\gammap}\partial_1\partial_2 u_1\partial_1 u_1 n_2
        \\
        &=\|\Delta u\|_2^2 + L_s^{-1} \int_{\gammam\cup\gammap}\partial_1^2 u_1 u_1 - L_s^{-1}\int_{\gammam\cup\gammap}\partial_1u_1\partial_1 u_1.
    \end{align}
    Finally, using the periodicity and $L_s>0$, it holds
    \begin{align}
        \|\nabla^2 u\|_2^2
        &=\|\Delta u\|_2^2 + L_s^{-1} \int_{\gammam\cup\gammap}\partial_1^2 u_1 u_1 - L_s^{-1}\int_{\gammam\cup\gammap}\partial_1u_1\partial_1 u_1  
        \\
        &= \|\Delta u\|_2^2 -2L_s^{-1} \int_{\gammam\cup\gammap} (\partial_1u_1)^2
        \\
        &\leq\|\Delta u\|_2^2
    \end{align}
    and by \eqref{id_Delta_u_is_nablaPerp_omega} one gets
    \begin{align}
        \|\nabla^2 u\|_2^2
        &\leq\|\Delta u\|_2^2
        =\|\nabla \omega\|_2^2.
        \label{bnxcvibncvb}
    \end{align}
\end{proof}

The estimate matching the one of Lemma \ref{lemma_vorticity_bound} is given in the following Lemma.

\begin{lemma}
    Let $2<p\in 2\mathbb{N}$, $u_0\in W^{1,p}(\Omega)$. Then
    \begin{align}
        \|\nabla u\|_p
        \lesssim \|\omega_0\|_{p} + L_s^{-1} \left(1+L_s^{-\frac{p}{p-2}}\right) \|u_0\|_2 + \left(1+L_s^{-\frac{2p-2}{p-2}}\right)\Ra,
    \end{align}
    where the implicit constant only depends on $\Gamma$ and $p$ and
    \begin{align}
        \|u\|_\infty
        &\lesssim \left(1+L_s^{-2}\right)\|u_0\|_{W^{1,4}} + \left(L_s^\frac{1}{2}+L_s^{-2}\right)\Ra,
    \end{align}
    where the implicit constant only depends on $\Gamma$.
    \label{lemma_flat_nablau_and_uinfty_bound}
\end{lemma}

\begin{proof}
    Using the same arguments as in the proof of Lemma \ref{lemma_vorticity_bound} one finds the analogous of \eqref{bnisubnisfng}, i.e.
    \begin{align}
        \|\omega\|_{L^\infty_t;L^p_x}^2 \lesssim \|\omega_0\|_p^2 + \Ra^2 + \Lambda^2,
        \label{jkbnwejkrbwejbwjer}
    \end{align}
    where $\Lambda = L_s^{-1}\|u_1\|_{L^\infty([0,T]\times \lbrace\gammam\cup\gammap\})}\lesssim L_s^{-1}\|u\|_{L^\infty_t;L^\infty_x}$.
    Estimating similar to before one finds, using Gagliardo-Nirenberg interpolation, Hölder's inequality
    \begin{align}
        \Lambda &\lesssim L_s^{-1}\|u\|_{L^\infty_t;L^\infty_x}
        \\
        &
        \lesssim L_s^{-1}\|\nabla u\|_{L^\infty_t;L^p_x}^\frac{p}{2(p-1)}\|u\|_{L^\infty_t;L^2_x}^\frac{p-2}{2(p-1)} + L_s^{-1}\left\|u\right\|_{L^\infty_t;L^2_x}
        \\
        &\lesssim L_s^{-1}\|\omega\|_{L^\infty_t;L^p_x}^\frac{p}{2(p-1)}\|u\|_{L^\infty_t;L^2_x}^\frac{p-2}{2(p-1)} + L_s^{-1}\left\|u\right\|_{L^\infty_t;L^2_x}
        \\
        &\lesssim \epsilon\|\omega\|_{L^\infty_t;L^p_x} + \left(L_s^{-1}+ \epsilon^{-\frac{p}{p-2}} L_s^{-\frac{2(p-1)}{p-2}} \right) \|u\|_{L^\infty_t;L^2_x}
        \label{nbjkxcvbnxcvjk}
    \end{align}
    for any $\epsilon>0$ as $2<p$, where in the last estimate we used Young's inequality. Therefore, combining \eqref{jkbnwejkrbwejbwjer} and \eqref{nbjkxcvbnxcvjk}, one has
    \begin{align}
        \|\omega\|_{L^\infty_t;L^p_x}^2
        &\lesssim \|\omega_0\|_{p}^2 + \Ra^2 +\epsilon^2\|\omega\|_{L^\infty_t;L^p_x}^2 + \left(L_s^{-1}+ \epsilon^{-\frac{p}{p-2}} L_s^{-\frac{2(p-1)}{p-2}} \right)^2 \|u\|_{L^\infty_t;L^2_x}^2
    \end{align}
    such that choosing $\epsilon$ sufficiently small one finds
    \begin{align}
        \|\omega\|_{L^\infty_t;L^p_x}^2
        &\lesssim \|\omega_0\|_{p}^2 + \Ra^2 + \left(L_s^{-1}+ L_s^{-\frac{2(p-1)}{p-2}} \right)^2 \|u\|_{L^\infty_t;L^2_x}^2
        \\
        &\lesssim \|\omega_0\|_{p}^2 + \Ra^2 + L_s^{-2} \left(1+ L_s^{-\frac{p}{p-2}} \right)^2 \left(\|u_0\|_2^2 + (1+L_s)\Ra^2\right),
    \end{align}
    where we used Lemma \ref{lemma_energy_bound_flat}. As this bound holds uniform in time it implies
    \begin{align}
        \|\nabla u\|_p &\lesssim \|\omega\|_{p}
        \lesssim \|\omega_0\|_{p} + L_s^{-1} \left(1+L_s^{-\frac{p}{p-2}}\right) \|u_0\|_2 + \left(1+L_s^{-\frac{2p-2}{p-2}}\right)\Ra,\qquad
        \label{bnxjvncxjb}
    \end{align}
    where we again used Young's inequality in order to unify the exponents.
    Using again Gagliardo-Nirenberg interpolation as in \eqref{nbjkxcvbnxcvjk} one finds
    \begin{align}
        \|u\|_\infty \lesssim \|\nabla u\|_{p}^{\frac{p}{2(p-1)}}\|u\|_2^{\frac{p-2}{2(p-1)}}+\|u\|_2\label{difference_to_arxiv_referd_in_explanation}
    \end{align}
    and choosing $p=4$ and plugging in \eqref{snjngfi} and \eqref{bnxjvncxjb} results in
    \begin{align}
        \|u\|_\infty
        &\lesssim  \|\nabla u\|_{4}^{\frac{2}{3}}\|u\|_2^{\frac{1}{3}}+\|u\|_2
        \\
        &\lesssim \left(\|\omega_0\|_{4} + L_s^{-1} \left(1+L_s^{-2}\right) \|u_0\|_2 + \left(1+L_s^{-3}\right)\Ra\right)^\frac{2}{3}
        \\
        &\qquad \cdot \left(\|u_0\|_2 + \left(1+L_s^\frac{1}{2}\right)\Ra\right)^\frac{1}{3} + \|u_0\|_2 + \left(1+L_s^\frac{1}{2}\right)\Ra
        \\
        &\lesssim \left(1+L_s^{-2}\right)\|u_0\|_{W^{1,4}} + \left(L_s^\frac{1}{2}+L_s^{-2}\right)\Ra,
    \end{align}
    where we again used Young's inequality.
\end{proof}

The proof of Lemma \ref{lemma_u_cdot_nabla_p_on_boundary} simplifies as follows.

\begin{lemma}
    If $u,p\in H^1(\Omega)$ it holds
    \begin{align}
        \left|\int_{\gammam\cup\gammap}u\cdot\nabla p \right|\leq 4 \|p\|_{H^1}\|\nabla u\|_2.
    \end{align}
    \label{lemma_rb_flat_int_u_nabla_p_on_the_boundary}
\end{lemma}

\begin{proof}
    Define $\zeta(x_2)=(2x_2-1)e_2$, then $\zeta(0)=-e_2=n$, $\zeta(1)=e_2=n$ and because of the periodicity in $x_1$-direction, \eqref{no_penetration_bc} and \eqref{incompressibility} and Stokes theorem yield
    \begin{align}
        \int_{\gammam\cup\gammap} u\cdot \nabla p &= \int_{\gammam\cup\gammap} u_1\partial_1 p  = - \int_{\gammam\cup\gammap} p\partial_1 u_1
        \\
        &= \int_{\gammam\cup\gammap} p\partial_2 u_2 = \int_{\gammam\cup\gammap} pn \cdot (\zeta\cdot\nabla) u
        \\
        &= \int_{\Omega} \nabla \cdot \left( p(\zeta\cdot\nabla) u\right)
        \\
        &= \int_{\Omega} \nabla p \cdot (\zeta\cdot\nabla) u + \int_{\Omega} p \partial_j u_i \partial_i \zeta_j + \int_{\Omega} p\zeta \cdot \nabla (\nabla\cdot u)
        \\
        &= \int_{\Omega} (2x_2-1) \nabla p \cdot \partial_2 u + 2 \int_{\Omega} p \partial_2 u_2.
    \end{align}
    Therefore Hölder's inequality yields
    \begin{align}
        \left|\int_{\gammam\cup\gammap} u\cdot \nabla p \right|
        \leq 2(\|\nabla p\|_2+ \|p\|_2) \|\nabla u\|_2
        \leq 4 \|p\|_{H^1}\|\nabla u\|_2,
        \label{cvjxbkbj}
    \end{align}
    where we additionally used Young's inequality to unify the norms of the pressure.
\end{proof}

Next, we derive the result corresponding to the pressure bound of Lemma \ref{lemma_pressure_bound}.

\begin{lemma}
    Let $u_0\in W^{1,4}(\Omega)$. Then
    \begin{align}
        \|p\|_{H^1}
        &\lesssim \left(\Pr^{-1} \|u\|_\infty + L_s^{-1}\right)\|\nabla u\|_2 + \Ra,
    \end{align}
    where the implicit constant only depends on $\Gamma$.
    \label{lemma_flat_pressure_bound}
\end{lemma}

\begin{proof}
    Starting from \eqref{abniuvab}, using \eqref{qwerqwe} and \eqref{bnxcjicbnc}, one has
    \begin{align}
        \|\nabla p\|_2^2 &= -\frac{1}{\Pr} \left(\int_{\Omega} u_t \cdot \nabla p + \int_{\Omega} \nabla p \cdot (u\cdot \nabla) u \right) + \int_{\Omega} \nabla p \cdot \Delta u + \Ra \int \partial_2 p \theta
        \\
        &=  -\Pr^{-1} \int_{\Omega} \nabla p \cdot (u\cdot \nabla) u - L_s^{-1} \int_{\gammam\cup\gammap} u\cdot\nabla p + \Ra \int \partial_2 p \theta
        \\
        &\lesssim \Pr^{-1} \|p\|_{H^1} \|u\|_\infty \|\nabla u\|_{L^2} + L_s^{-1}\|p\|_{H^1}\|\nabla u\|_2 + \Ra \|p\|_{H^1},
        \label{uaidfuidfg}
    \end{align}
    where in the last estimate we used Hölder's inequality, Lemma \ref{lemma_rb_flat_int_u_nabla_p_on_the_boundary} and \eqref{maximum_principle}. As $p$ is only defined up to a constant, we can choose it such that it is average free and Poincaré's inequality holds. Then \eqref{uaidfuidfg} yields
    \begin{align}
        \|p\|_{H^1}^2 &\lesssim \|\nabla p\|_2^2 
        \\
        &\lesssim \Pr^{-1} \|p\|_{H^1} \|u\|_\infty \|\nabla u\|_{L^2} + L_s^{-1}\|p\|_{H^1}\|\nabla u\|_2 + \Ra \|p\|_{H^1}
    \end{align}
    and dividing by $\|p\|_{H^1}$ results in
    \begin{align}
        \|p\|_{H^1}
        &\lesssim \left(\Pr^{-1} \|u\|_\infty + L_s^{-1}\right)\|\nabla u\|_2 + \Ra.
        \label{xbcivnxc}
    \end{align}
\end{proof}

The equivalent of Lemma \ref{lemma_long_time_H2_bound} is given in the following.

\begin{lemma}
    Let $u_0\in W^{1,4}(\Omega)$. Then
    \begin{align}
        \langle \|\nabla^2 u\|_2^2\rangle
        &\lesssim \left(\Pr^{-1} \left(L_s^{-1}+L_s^{-3}\right)\|u_0\|_{W^{1,4}} + L_s^{-2} \right)\Nu\Ra
        \\
        &\qquad + L_s^{-1}\Nu^\frac{1}{2}\Ra^\frac{3}{2} + \Nu\Ra^\frac{3}{2} + \Pr^{-1}\left(L_s^{-\frac{1}{2}}+L_s^{-3}\right)\Nu\Ra^2,
    \end{align}
    where the implicit constant only depends on $\Gamma$.
    \label{lemma_flat_H2_bound}
\end{lemma}

\begin{proof}
    By \eqref{cvjnbjbcvxbni} one has
    \begin{align}
        \frac{1}{2\Pr}\frac{d}{dt}&\left(\|\omega\|_2^2 + L_s^{-1}\int_{\gammam\cup\gammap} u_\tau^2\right) + \|\nabla\omega\|_2^2
        \\
        &= - \frac{1}{2 \Pr L_s} \int_{\gammam\cup\gammap} u_1\partial_1 (u_1^2) - L_s^{-1}\int_{\gammam\cup\gammap} u\cdot\nabla p + \Ra \int_{\Omega} \omega \partial_1 \theta.
    \end{align}
    Note that because of the periodicity of the boundaries
    \begin{align}
        \frac{1}{2}\int_{\gammam\cup\gammap} u_1\partial_1 (u_1^2) = \int_{\gammam\cup\gammap} u_1^2\partial_1 u_1 = \frac{1}{3}\int_{\gammam\cup\gammap} \partial_1 (u_1^3)=0
    \end{align}
    and therefore
    \begin{align}
        \frac{1}{2\Pr}\frac{d}{dt}&\left(\|\omega\|_2^2 + L_s^{-1}\int_{\gammam\cup\gammap} u_\tau^2\right) + \|\nabla\omega\|_2^2
        \\
        &\qquad \qquad = - L_s^{-1}\int_{\gammam\cup\gammap} u\cdot\nabla p + \Ra \int_{\Omega} \omega \partial_1 \theta
        \\
        &\qquad \qquad \lesssim L_s^{-1} \|\nabla u\|_2\|p\|_{H^1} + \Ra \|\omega\|_2 \|\nabla \theta\|_2,
        \label{bxhiuxcb}
    \end{align}
    where we used Lemma \ref{lemma_rb_flat_int_u_nabla_p_on_the_boundary} in the last estimate. By Lemma \ref{lemma_flat_pressure_bound} and Lemma \ref{lemma_grad_ids} \eqref{bxhiuxcb} can be further estimated as
    \begin{align}
        \frac{1}{2\Pr}\frac{d}{dt}&\left(\|\omega\|_2^2 + L_s^{-1}\int_{\gammam\cup\gammap} u_\tau^2\right) + \|\nabla\omega\|_2^2
        \\
        &\lesssim \left(\Pr^{-1}L_s^{-1} \|u\|_\infty + L_s^{-2}\right)\|\nabla u\|_2^2 + \left(L_s^{-1}+\|\nabla \theta\|_2\right) \|\nabla u\|_2 \Ra
        \\
        &\lesssim \Pr^{-1} \left(\left(L_s^{-1}+L_s^{-3}\right)\|u_0\|_{W^{1,4}} + \left(L_s^{-\frac{1}{2}}+L_s^{-3}\right)\Ra\right)\|\nabla u\|_2^2
        \\
        &\qquad + L_s^{-2}\|\nabla u\|_2^2 + \left(L_s^{-1}+\|\nabla \theta\|_2\right) \|\nabla u\|_2 \Ra,
    \end{align}
    where in the last estimate we used Lemma \ref{lemma_flat_nablau_and_uinfty_bound}. Using Lemma \ref{lemma_flat_elliptic_regularity} one gets
    \begin{align}
        \frac{1}{2\Pr}\frac{d}{dt}&\left(\|\omega\|_2^2 + L_s^{-1}\int_{\gammam\cup\gammap} u_\tau^2\right) + \|\nabla^2 u\|_2^2
        \\
        &\lesssim \Pr^{-1} \left(\left(L_s^{-1}+L_s^{-3}\right)\|u_0\|_{W^{1,4}} + \left(L_s^{-\frac{1}{2}}+L_s^{-3}\right)\Ra\right)\|\nabla u\|_2^2
        \\
        &\qquad + L_s^{-2}\|\nabla u\|_2^2 + \left(L_s^{-1}+\|\nabla \theta\|_2\right) \|\nabla u\|_2 \Ra,
        \label{bnisgadfg}
    \end{align}
    and taking the long time average of \eqref{bnisgadfg}, using Lemma \ref{lemma_energy_bound_flat}, Lemma \ref{lemma_nusselt_number_representations} and that $\|u\|_{H^1}^2$ is uniformly bounded in time, we find
    \begin{align}
        \langle \|\nabla^2 u\|_2^2\rangle
        &\lesssim \left(\Pr^{-1} \left(L_s^{-1}+L_s^{-3}\right)\|u_0\|_{W^{1,4}} + L_s^{-2} \right)\Nu\Ra
        \\
        &\qquad + L_s^{-1}\Nu^\frac{1}{2}\Ra^\frac{3}{2} + \Nu\Ra^\frac{3}{2} + \Pr^{-1}\left(L_s^{-\frac{1}{2}}+L_s^{-3}\right)\Nu\Ra^2.
    \end{align}    
\end{proof}

\subsection{Proof of the Theorem}
\label{section:flat_theorem_proof}

\begin{proof}[\hypertarget{proof_theorem_flat}{Proof of Theorem} \ref{theorem_flat}]
\\
    By Lemma \ref{lemma_flat_nusselt_bounded_by_Hx} and Lemma \ref{lemma_energy_bound_flat}
    \begin{align}
        \Nu
        &\lesssim \delta^\frac{1}{2} \langle \|\nabla u\|_2^2 \rangle^\frac{1}{2} + \delta^{-\frac{1}{2}}\Nu^\frac{1}{2}
        \\
        &\lesssim \delta^\frac{1}{2} \Nu^\frac{1}{2}\Ra^\frac{1}{2} + \delta^{-\frac{1}{2}}\Nu^\frac{1}{2}
    \end{align}
    for any $\delta>0$, implying
    \begin{align}
        \Nu \lesssim \delta \Ra + \delta^{-1}
    \end{align}
    and setting $\delta=\Ra^{-\frac{1}{2}}$ yields
    \begin{align}
        \Nu \lesssim \Ra^\frac{1}{2},
    \end{align}
    proving the first claim.
    Similarly by Lemma \ref{lemma_flat_nusselt_bounded_by_Hx}, Lemma \ref{lemma_energy_bound_flat} and Lemma \ref{lemma_flat_H2_bound}
    \begin{align}
        \Nu
        &\lesssim \delta \langle \|\nabla u\|_2^2 \rangle^\frac{1}{4} \langle \|\nabla^2 u\|_2^2\rangle^\frac{1}{4} + \delta^{-\frac{1}{2}}\Nu^\frac{1}{2}
        \\
        &\lesssim \delta \Pr^{-\frac{1}{4}} \left(L_s^{-\frac{1}{4}}+L_s^{-\frac{3}{4}}\right)\|u_0\|_{W^{1,4}}^\frac{1}{4} \Nu^\frac{1}{2}\Ra^\frac{1}{2}
        + \delta L_s^{-\frac{1}{2}}\Nu^\frac{1}{2}\Ra^\frac{1}{2}
        \\
        &\qquad + \delta L_s^{-\frac{1}{4}}\Nu^\frac{3}{8}\Ra^\frac{5}{8} + \delta \Nu^\frac{1}{2}\Ra^\frac{5}{8} + \delta\Pr^{-\frac{1}{4}} \left(L_s^{-\frac{1}{8}}+L_s^{-\frac{3}{4}}\right)\Nu^\frac{1}{2}\Ra^\frac{3}{4}
        \\
        &\qquad + \delta^{-\frac{1}{2}}\Nu^\frac{1}{2}.
        \label{nxcvbnicvb}
    \end{align}
    As in the general domain we distinguish between two cases because of the $\Nu^\frac{3}{8}$ term.
    \begin{itemize}
        \item At first assume that this $\Nu^\frac{3}{8}$ term is dominated by the other terms, i.e.
        \begin{align}
            L_s^{-\frac{1}{4}} &\leq \Pr^{-\frac{1}{4}} \left(L_s^{-\frac{1}{4}}+L_s^{-\frac{3}{4}}\right)\|u_0\|_{W^{1,4}}^\frac{1}{4} \Nu^\frac{1}{8}\Ra^{-\frac{1}{8}}
            + L_s^{-\frac{1}{2}}\Nu^\frac{1}{8}\Ra^{-\frac{1}{8}}
            \\
            &\qquad + \Nu^\frac{1}{8} + \Pr^{-\frac{1}{4}} \left(L_s^{-\frac{1}{8}}+L_s^{-\frac{3}{4}}\right)\Nu^\frac{1}{8}\Ra^\frac{1}{8}.
        \end{align}
        Then \eqref{nxcvbnicvb} implies
        \begin{align}
            \Nu 
            &\lesssim \delta \Pr^{-\frac{1}{4}} \left(L_s^{-\frac{1}{4}}+L_s^{-\frac{3}{4}}\right)\|u_0\|_{W^{1,4}}^\frac{1}{4} \Nu^\frac{1}{2}\Ra^\frac{1}{2}
            + \delta L_s^{-\frac{1}{2}}\Nu^\frac{1}{2}\Ra^\frac{1}{2}
            \\
            &\qquad + \delta \Nu^\frac{1}{2}\Ra^\frac{5}{8} + \delta\Pr^{-\frac{1}{4}} \left(L_s^{-\frac{1}{8}}+L_s^{-\frac{3}{4}}\right)\Nu^\frac{1}{2}\Ra^\frac{3}{4} + \delta^{-\frac{1}{2}}\Nu^\frac{1}{2},
        \end{align}
        which after dividing by $\Nu^\frac{1}{2}$ and squaring yields
        \begin{align}
            \Nu
            &\lesssim \delta^2 \Pr^{-\frac{1}{2}} \left(L_s^{-\frac{1}{2}}+L_s^{-\frac{3}{2}}\right)\|u_0\|_{W^{1,4}}^\frac{1}{2} \Ra
            + \delta^2 L_s^{-1}\Ra
            \\
            &\qquad + \delta^2 \Ra^\frac{5}{4} + \delta^2 \Pr^{-\frac{1}{2}} \left(L_s^{-\frac{1}{4}}+L_s^{-\frac{3}{2}}\right)\Ra^\frac{3}{2} + \delta^{-1}.
        \end{align}
        Setting
        \begin{align}
            \delta &= \Big(\Pr^{-\frac{1}{2}} \left(L_s^{-\frac{1}{2}}+L_s^{-\frac{3}{2}}\right)\|u_0\|_{W^{1,4}}^\frac{1}{2} \Ra
            + L_s^{-1}\Ra
            \\
            &\qquad + \Ra^\frac{5}{4} + \Pr^{-\frac{1}{2}} \left(L_s^{-\frac{1}{4}}+L_s^{-\frac{3}{2}}\right)\Ra^\frac{3}{2}\Big)^{-\frac{1}{3}}
        \end{align}
        results in
        \begin{align}
            \Nu
            &\lesssim \Pr^{-\frac{1}{6}} \left(L_s^{-\frac{1}{6}}+L_s^{-\frac{1}{2}}\right)\|u_0\|_{W^{1,4}}^\frac{1}{6} \Ra^\frac{1}{3}
            + L_s^{-\frac{1}{3}}\Ra^\frac{1}{3}
            \\
            &\qquad + \Ra^\frac{5}{12} + \Pr^{-\frac{1}{6}} \left(L_s^{-\frac{1}{12}}+L_s^{-\frac{1}{2}}\right)\Ra^\frac{1}{2}.
            \label{bnxcvuixcv}
        \end{align}
        
        \item If instead 
        \begin{align}
            L_s^{-\frac{1}{4}} &> \Pr^{-\frac{1}{4}} \left(L_s^{-\frac{1}{4}}+L_s^{-\frac{3}{4}}\right)\|u_0\|_{W^{1,4}}^\frac{1}{4} \Nu^\frac{1}{8}\Ra^{-\frac{1}{8}}
            + L_s^{-\frac{1}{2}}\Nu^\frac{1}{8}\Ra^{-\frac{1}{8}}
            \\
            &\qquad + \Nu^\frac{1}{8} + \Pr^{-\frac{1}{4}} \left(L_s^{-\frac{1}{8}}+L_s^{-\frac{3}{4}}\right)\Nu^\frac{1}{8}\Ra^\frac{1}{8}.
        \end{align}
        Then \eqref{nxcvbnicvb} implies
        \begin{align}
            \Nu 
            &\lesssim \delta L_s^{-\frac{1}{4}}\Nu^\frac{3}{8}\Ra^\frac{5}{8} + \delta^{-\frac{1}{2}}\Nu^\frac{1}{2}
        \end{align}
        and setting $\delta = L_s^{\frac{1}{6}} \Nu^{\frac{1}{12}}\Ra^{-\frac{5}{12}}<\Ra^{-\frac{5}{12}}$ yields
        \begin{align}
            \Nu 
            &\lesssim L_s^{-\frac{1}{12}} \Nu^{\frac{11}{24}}\Ra^{\frac{5}{24}}
        \end{align}
        and therefore division by $\Nu^\frac{11}{24}$ and exponentiation results in
        \begin{align}
            \Nu
            &\lesssim L_s^{-\frac{2}{13}} \Ra^{\frac{5}{13}}.
        \end{align}
    \end{itemize}
    Combining the two cases one has
    \begin{align}
            \Nu
            &\lesssim \left( L_s^{-\frac{1}{3}} + \Pr^{-\frac{1}{6}} \left(L_s^{-\frac{1}{6}}+L_s^{-\frac{1}{2}}\right)\|u_0\|_{W^{1,4}}^\frac{1}{6}\right)\Ra^\frac{1}{3}
            \\
            &\qquad + L_s^{-\frac{2}{13}} \Ra^{\frac{5}{13}} + \Ra^\frac{5}{12} + \Pr^{-\frac{1}{6}} \left(L_s^{-\frac{1}{12}}+L_s^{-\frac{1}{2}}\right)\Ra^\frac{1}{2}
    \end{align}
    and therefore, since $\Ra>1$
    \begin{align}
        \Nu
        &\lesssim L_s^{-\frac{1}{6}}\|u_0\|_{W^{1,4}}^\frac{1}{6} \Pr^{-\frac{1}{6}}\Ra^\frac{1}{3} + \Ra^\frac{5}{12} + L_s^{-\frac{1}{12}}\Pr^{-\frac{1}{6}}\Ra^\frac{1}{2}
    \end{align}
    if $L_s\geq 1$ and
    \begin{align}
        \Nu
        &\lesssim  L_s^{-\frac{1}{3}} \Ra^\frac{1}{3} +  L_s^{-\frac{1}{2}}\|u_0\|_{W^{1,4}}^\frac{1}{6}\Pr^{-\frac{1}{6}}\Ra^\frac{1}{3} + L_s^{-\frac{2}{13}} \Ra^{\frac{5}{13}} + \Ra^\frac{5}{12} + L_s^{-\frac{1}{2}}\Pr^{-\frac{1}{6}}\Ra^\frac{1}{2}
    \end{align}
    if $L_s\leq 1$, concluding the proof.
\end{proof}

\section{System With Identical Boundaries}
\label{section:nonlinearity_paper}

The following findings are published in \cite{bleitnerNobili2024Bounds}. Although the results are suboptimal when compared to Theorem \ref{theorem_main_theorem_rb_curved}, we state the corresponding Lemmas and changes in the proofs in order to discuss the differences and improvements. 

\begin{figure}
    \hrule
    \vspace{0.5\baselineskip}
    \begin{center}
        \includegraphics[width=\textwidth]{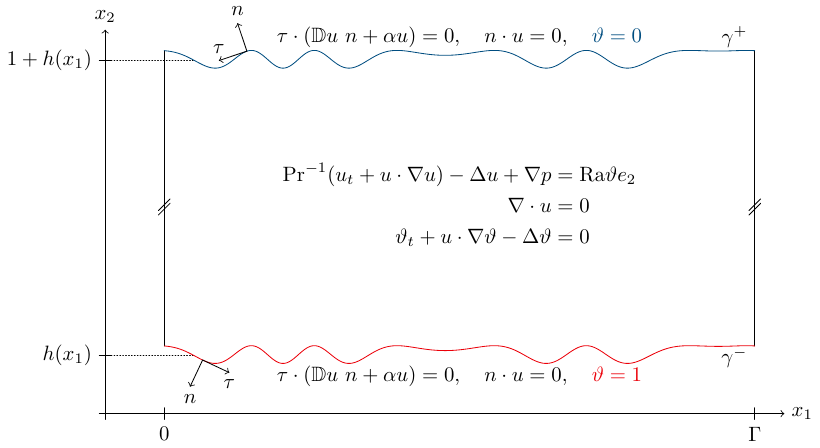}
        \vspace*{-7mm}
        \caption{Overview of the system considered in Section \ref{section:nonlinearity_paper}, where the boundary profile functions coincide.}
        \label{fig:rb_overview_identical}
    \end{center}
    \hrule
\end{figure}

The particular system of interest corresponds to the previously studied general system when the boundary profile functions coincide, i.e.
\begin{align}
    \gammam &= \lbrace (x_1,x_2) \ \vert \ 0 \leq x_1 \leq \Gamma, x_2 = h(x_1)\rbrace
    \\
    \gammap &= \lbrace (x_1,x_2) \ \vert \ 0 \leq x_1 \leq \Gamma, x_2 = 1 + h(x_1)\rbrace
\end{align}
The setup is illustrated in Figure \ref{fig:rb_overview_identical}.

\subsection{Nusselt Number}
First, we state the result corresponding to Lemma \ref{lemma_nusselt_number_representations}.

\begin{lemma}[Nusselt Number Representations]
    \label{lemma_nonlin_nusselt_number_representations}
    Let $\Omega$ be $C^{1,1}$ und $u_0\in L^2(\Omega)$. Then for any $0\leq z \leq 1$
    \begin{align}
        \Nu &= \langle \np \cdot (u-\nabla)\theta\rangle_{\gammam + z}
        \label{nonlin_nu_representation_shifted}
        \\
        &= \langle \|\nabla \theta\|_2^2\rangle_{\Omega} \label{nonlin_nu_representation_gradient}
        \\
        &\geq (1+\max h - \min h)^{-1} \langle (u_2-\partial_2)\theta \rangle_{\Omega}, \label{nonlin_nu_representation_straight}
    \end{align}\noeqref{nonlin_nu_representation_gradient}\noeqref{nu_representation_straight}%
    where $\gammam + z = \lbrace (x_1,x_2)\ \vert \ 0\leq x_1\leq \Gamma, x_2 = h(x_1) + z\rbrace$.
\end{lemma}

The proof of Lemma \ref{lemma_nonlin_nusselt_number_representations} follows the same way as the proof of Lemma \ref{lemma_nusselt_number_representations}.

\subsection{A-Priori Estimates}
\label{section:nonlin_regularity_estimates}

\begin{lemma}[Energy Balance]
    \label{lemma_nonlin_energy_balance}
    Let $\Omega$. Then strong solutions satisfy
    \begin{align}
        \frac{1}{2\Pr}\frac{d}{dt} \|u\|_2^2 + \|\nabla u\|_2^2 + \int_{\gammam\cup\gammap} (2\alpha+\kappa) u_\tau^2 &= \Ra \int_{\Omega} \theta u_2.
        \label{bcnxibknj}
    \end{align}
\end{lemma}

\begin{proof}
    Testing \eqref{navier_stokes} with $u$ one finds
    \begin{align}
        \frac{1}{2\Pr}\frac{d}{dt} \|u\|_2^2 &= - \frac{1}{\Pr}\int_\Omega u\cdot (u\cdot \nabla) u - \int_{\Omega} u\cdot \nabla p + \int_{\Omega} u\cdot \Delta u + \Ra \int_{\Omega} \theta u_2.\qquad\label{uivbucx}
    \end{align}
    Integrating by parts, the boundary conditions \eqref{no_penetration_bc} and \eqref{incompressibility} yield
    \begin{align}
        \int_{\Omega} \nabla p\cdot u = \int_{\gammam\cup\gammap} p n\cdot u - \int_{\Omega} p \nabla \cdot u = 0\label{buviiuxc}
    \end{align}
    and
    \begin{align}
        \int_{\Omega} u\cdot (u\cdot\nabla) u &= \int_{\gammam\cup\gammap} u\cdot (u\cdot n) u - \int_{\Omega} (u\cdot u) \nabla \cdot u - \int_{\Omega} u\cdot (u\cdot \nabla) u
        \\
        &= - \int_{\Omega} u\cdot (u\cdot \nabla) u,
    \end{align}
    implying that the first term on the right-hand side of \eqref{uivbucx} vanishes. Therefore combining \eqref{uivbucx} and \eqref{buviiuxc} 
    \begin{align}
        \frac{1}{2\Pr}\frac{d}{dt} \|u\|_2^2 &= \int_{\Omega} u\cdot \Delta u + \Ra \int_{\Omega} \theta u_2.\label{xcvnbxci}
    \end{align}
    For the viscosity term integration by parts and \eqref{no_penetration_bc}yield
    \begin{align}
        \int_{\Omega} u\cdot \Delta u &= \int_{\gammam\cup\gammap} u \cdot (n\cdot \nabla) u - \|\nabla u\|_2^2
        \\
        &= \int_{\gammam\cup\gammap} u_\tau \tau_i n_j \partial_j u_i - \|\nabla u\|_2^2
        \\
        &= \int_{\gammam\cup\gammap} u_\tau \tau_i n_j (\partial_j u_i+\partial_i u_j) - \int_{\gammam\cup\gammap} u_\tau \tau_i n_j \partial_i u_j - \|\nabla u\|_2^2
        \\
        &= 2 \int_{\gammam\cup\gammap} u_\tau \tau \cdot \D u \ n - \int_{\gammam\cup\gammap} u_\tau n\cdot (\tau\cdot \nabla) u - \|\nabla u\|_2^2
        \\
        &= - 2 \int_{\gammam\cup\gammap} \alpha  u_\tau^2 - \int_{\gammam\cup\gammap} \kappa u_\tau^2 - \|\nabla u\|_2^2,\label{bivic}
    \end{align}
    where in the last identity we used \eqref{nav_slip_bc} and \eqref{n_tau_grad_u}. Plugging \eqref{bivic} into \eqref{xcvnbxci} we find
    \begin{align}
        \frac{1}{2\Pr}\frac{d}{dt} \|u\|_2^2 + \|\nabla u\|_2^2 + \int_{\gammam\cup\gammap} (2\alpha+\kappa) u_\tau^2 &= \Ra \int_{\Omega} \theta u_2.
    \end{align}
\end{proof}

From \eqref{bcnxibknj} it is not directly clear that energy is dissolved, i.e. a coercivity bound holds for the gradient and boundary term. If $2\alpha> |\kappa|$ one can expect such an estimate. If however $2 \alpha<|\kappa|$ the following observation motivates the subsequent lemma, which will provide such a coercivity estimate. Assume $\kappa(x_1,h(x_1))<0$ at some point $x_1$, i.e. the bottom boundary is convex. Then as the boundary profiles are the same the top boundary has to be concave in $x_1$, i.e. $\kappa(x_1,1+h(x_1))>0$, and vice versa. The fundamental theorem of calculus allows us to exchange between these boundary values for $u$.

\begin{lemma}[Coercivity]\label{lemma_nonlin_coercivity}
    Assume $\Omega$ is $C^{1,1}$ with $\hp=1+\hm$, $u\in H^1(\Omega)$, $0<\alpha\in L^\infty(\gammam\cup\gammap)$ and
    \begin{align}
        |\kappa|(x)\leq 2\alpha(x) + \frac{1}{4\sqrt{1+(h'(x_1))^2}}\min \left\{1,\sqrt{\alpha(x)}\right\}
        \label{nonlin_coercivity_assumption}
    \end{align}
    holds for almost every $x\in \gammam\cup\gammap$. Then
    \begin{align}
        \frac{3}{4}\|\nabla u\|_2^2 +\int_{\gammam\cup\gammap} (2\alpha+\kappa)u_\tau^2 \geq \frac{1}{4} \min\left\{1,\|\alpha^{-1}\|_\infty^{-1}\right\} \|u\|_{H^1}^2.
    \end{align}
\end{lemma}

\begin{remark}[Comparison of Lemma \ref{lemma_coercivity} and Lemma \ref{lemma_nonlin_coercivity}]
Note that under assumption \eqref{nonlin_coercivity_assumption} one has
\begin{align}
    1+\|\alpha^{-1}(1+|\kappa|)\|_\infty 
    \lesssim 1+\left\|\tfrac{1+\alpha+\sqrt{\alpha}}{\alpha}\right\|_\infty
    \lesssim 1+\left\|\tfrac{1+\alpha}{\alpha}\right\|_\infty
    \lesssim 1 + \|\alpha^{-1}\|_\infty 
\end{align}
and Lemma \ref{lemma_coercivity} yields
\begin{align}
    \|\nabla u\|_2^2 +\int_{\partial\Omega} (2\alpha+\kappa)u_\tau
    &= 2 \|\D u\|_2^2 + 2\int_{\partial\Omega}\alpha u_\tau 
    \\
    &\gtrsim (1+\|\alpha^{-1}(1+|\kappa|)\|_\infty)^{-1} \|u\|_{H^1}^2
    \\
    &\gtrsim (1+\|\alpha^{-1}\|_\infty)^{-1} \|u\|_{H^1}^2
    \\
    &\gtrsim \min\{1,\|\alpha^{-1}\|_\infty^{-1}\} \|u\|_{H^1}^2.
\end{align}
This shows that Lemma \ref{lemma_coercivity} is indeed a generalization of Lemma \ref{lemma_nonlin_coercivity}. Additionally Lemma \ref{lemma_nonlin_coercivity} relies on the fact that the boundaries have the same profile function $h$ and is therefore not applicable to the domain studied in the general Rayleigh-Bénard system and especially not in the case of an arbitrary Lipschitz domain.
\end{remark}

\begin{proof}
    We will use the following notation
    \begin{alignat}{3}
        \kappam &= \kappa(x_1,h(x_1)),\quad & \alpham &= \alpha(x_1,h(x_1)),\quad & \um &= u_\tau(x_1,h(x_1)),
        \\
        \kappap &= \kappa(x_1,1+h(x_1)),\quad & \alphap &= \alpha(x_1,1+h(x_1)),\quad & \up &= u_\tau(x_1,1+h(x_1)),
    \end{alignat}
    i.e. the evaluation of the functions on the bottom and top boundary. Note that since the boundary profiles are the same and $\np=-\nm$
    \begin{align}
        \kappam=-\kappap.
    \end{align}
    By the fundamental theorem of calculus and Young's inequality
    \begin{align}
        |u|^2(x_1,x_2) &= \left(\um + \int_{h(x_1)}^{x_2} \partial_2 u(x_1,z)\ dz\right)^2
        \\
        &\leq (1+\epsilon) \um^2 + (1+\epsilon^{-1})\left(\int_{h(x_1)}^{x_2} \partial_2 u(y_1,z)\ dz\right)^2
        \\
        &\leq (1+\epsilon) \um^2 + (1+\epsilon^{-1})(x_2-h(x_1))\|\partial_2 u\|_{L^2(\Omegavertvert)}^2
        \label{nxcixbx}
    \end{align}
    for any $\epsilon>0$, where $L^2(\Omegavertvert) = L^2(h(x_1),1+h(x_1))$ and analogously
    \begin{align}
        |u|^2(x_1,x_2)
        &\leq (1+\epsilon) \up^2 + (1+\epsilon^{-1})(1+h(x_1)-x_2)\|\partial_2 u\|_{L^2(\Omegavertvert)}^2.\quad
        \label{jnvcibxjcnb}
    \end{align}
    Integrating \eqref{nxcixbx} and \eqref{jnvcibxjcnb} in $x_2$ one gets
    \begin{align}
        &\|u(x_1,\cdot)\|_{L^2(\Omegavertvert)}^2
        \\
        &\qquad \leq (1+\epsilon)\min\{\um^2,\up^2\} + \frac{1+\epsilon^{-1}}{2}\|\partial_2 u (x_1,\cdot)\|_{L^2(\Omegavertvert)}^2
        \\
        &\qquad \leq (1+\epsilon)\max\{1,\alpham^{-1},\alphap^{-1}\}\Big(\min\{\alpham,\alphap\}\sqrt{1+(h')^2}\min\{\um^2,\up^2\}
        \\
        &\qquad \qquad + (2\epsilon)^{-1}\|\partial_2 u (x_1,\cdot)\|_{L^2(\Omegavertvert)}^2\Big),\label{ncvbjxcn}
    \end{align}
    where we smuggled in the factor $\sqrt{1+(h')^2}>1$. Next, we claim that
    \begin{align}
        \min\{\alpham,\alphap\}\sqrt{1+(h')^2}u_i^2 &\leq \frac{5}{16}\|\partial_2 u\|_{L^2(\Omegavertvert)}^2 + (2\alpham+\kappam)\sqrt{1+(h')^2}\um^2
        \\
        &\qquad  + (2\alphap+\kappap)\sqrt{1+(h')^2}\up^2
        \label{njcvbnxjvb}
    \end{align}
    holds for either $i=-$ or $i=+$. Plugging \eqref{njcvbnxjvb} into \eqref{ncvbjxcn} one gets
    \begin{align}
        &\|u(x_1,\cdot)\|_{L^2(\Omegavertvert)}^2
        \\
        &\qquad\leq (1+\epsilon)\max\{1,\alpham^{-1},\alphap^{-1}\}\Big(\min\{\alpham,\alphap\}\sqrt{1+(h')^2}\min\{\um^2,\up^2\}
        \\
        &\qquad\qquad + (2\epsilon)^{-1}\|\partial_2 u (x_1,\cdot)\|_{L^2(\Omegavertvert)}^2\Big)
        \\
        &\qquad\leq (1+\epsilon)\max\{1,\alpham^{-1},\alphap^{-1}\}\bigg((2\alpham+\kappam)\sqrt{1+(h')^2}\um^2
        \\
        &\qquad\qquad + (2\alphap+\kappap)\sqrt{1+(h')^2}\up^2 + \left(\frac{5}{16}+(2\epsilon)^{-1}\right)\|\partial_2 u (x_1,\cdot)\|_{L^2(\Omegavertvert)}^2\bigg).
    \end{align}
    Integrating with respect to $x_1$ and choosing $\epsilon=3$ yields
    \begin{align}
        \|u\|_2^2 &\leq 4 \max\{1,\|\alpha^{-1}\|_\infty\}\left(\int_{\gammam\cup\gammap}(2\alpha+\kappa)u_\tau^2 + \frac{1}{2}\|\partial_2 u\|_2^2\right).
    \end{align}
    It is left to show that \eqref{njcvbnxjvb} holds. In order to prove the claim we distinguish between two cases.
    \begin{itemize}
        \item Assume $|\kappa|\leq 2\alpha$.
        
        Then $2 \alphap+\kappap\geq 0$ and $2\alpham+\kappam\geq 0$. As $\kappap=-\kappam$, either $\kappap>0$ or $\kappam>0$. If $\kappam>0$ one has
        \begin{align}
            \min\{\alpham,\alphap\}\sqrt{1+(h')^2}u_i^2
            &\leq 2\alpham\sqrt{1+(h')^2}\um^2
            \\
            &\leq (2\alpham+\kappam)\sqrt{1+(h')^2}\um^2
            \\
            &\qquad +(2\alphap+\kappap)\sqrt{1+(h')^2}\up^2
        \end{align}
        and if $\kappap>0$
        \begin{align}
            \min\{\alpham,\alphap\}\sqrt{1+(h')^2}u_i^2
            &\leq 2\alphap\sqrt{1+(h')^2}\up^2
            \\
            &\leq (2\alpham+\kappam)\sqrt{1+(h')^2}\um^2
            \\
            &\qquad +(2\alphap+\kappap)\sqrt{1+(h')^2}\up^2.
        \end{align}
        \item Assume without loss of generality $\kappap<0$ and $|\kappap|>2\alphap$ as the case $\kappam<0$ and $|\kappam|>2\alpham$ follows by exchanging $+$ and $-$.

        Using \eqref{nxcixbx} with $x_2=1+h(x_1)$, respectively \eqref{jnvcibxjcnb} with $x_2=h(x_1)$, and observing that $\kappam=-\kappap>2\alphap$, it holds
        \begin{align}
            &-(2\alpham+\kappam)\sqrt{1+(h')^2}\um^2 -(2\alphap+\kappap)\sqrt{1+(h')^2}\up^2
            \\
            &\qquad = -(2\alpham+\kappam)\sqrt{1+(h')^2}\um^2 +(-2\alphap+\kappam)\sqrt{1+(h')^2}\up^2
            \\
            &\qquad \leq -(2\alpham+\kappam)\sqrt{1+(h')^2}\um^2 +(-2\alphap+\kappam)\sqrt{1+(h')^2}
            \\
            &\qquad\qquad \cdot\left((1+\epsilon) \um^2 + (1+\epsilon^{-1})\|\partial_2 u\|_{L^2(\Omegavertvert)}^2\right)
            \\
            &\qquad \leq -\left(2\alphap+2\alpham-\epsilon(\kappam-2\alphap)\right)\sqrt{1+(h')^2}\um^2
            \\
            &\qquad\qquad  +(\kappam-2\alphap)(1+\epsilon^{-1})\sqrt{1+(h')^2} \|\partial_2 u\|_{L^2(\Omegavertvert)}^2.
        \end{align}
        Note that $\kappam=|\kappap|>2\alphap$. Therefore, for the first bracket to be positive we choose $\epsilon=\frac{\alpham+\alphap}{\kappam-2\alphap}$ to get
        \begin{align}
            &-(2\alpham+\kappam)\sqrt{1+(h')^2}\um^2 -(2\alphap+\kappap)\sqrt{1+(h')^2}\up^2
            \\
            &\qquad \leq -\left(\alphap+\alpham\right)\sqrt{1+(h')^2}\um^2
            \\
            &\qquad\qquad  +(\kappam-2\alphap)(1+\epsilon^{-1})\sqrt{1+(h')^2} \|\partial_2 u\|_{L^2(\Omegavertvert)}^2
            \\
            &\qquad \leq -\left(\alphap+\alpham\right)\sqrt{1+(h')^2}\um^2
            \\
            &\qquad\qquad  +\left(\kappam-2\alphap+\frac{(\kappam-2\alphap)^2}{\alpham+\alphap}\right)\sqrt{1+(h')^2} \|\partial_2 u\|_{L^2(\Omegavertvert)}^2 \quad\qquad
            \label{bnixcunbcix}
        \end{align}
        The assumption \eqref{nonlin_coercivity_assumption} implies
        \begin{align}
            \kappam&=|\kappap|\leq 2\alphap + \frac{1}{4\sqrt{1+(h')^2}}\min\{1,\sqrt{\alphap}\}
            \\
            &\leq 2\alphap + \frac{1}{4}\min\left\{\frac{1}{\sqrt{1+(h')^2}},\frac{\sqrt{\alphap+\alpham}}{\left(1+(h')^2\right)^\frac{1}{4}}\right\},
            \label{nibxcjunbjuc}
        \end{align}
        where $\left(1+(h')^2\right)^\frac{1}{4}\leq \left(1+(h')^2\right)^\frac{1}{2}$ was used. Combining \eqref{bnixcunbcix} and \eqref{nibxcjunbjuc}, one gets
        \begin{align}
            &-(2\alpham+\kappam)\sqrt{1+(h')^2}\um^2 -(2\alphap+\kappap)\sqrt{1+(h')^2}\up^2
            \\
            &\qquad \leq -\left(\alphap+\alpham\right)\sqrt{1+(h')^2}\um^2
            \\
            &\qquad\qquad  +\left(\kappam-2\alphap+\frac{(\kappam-2\alphap)^2}{\alpham+\alphap}\right)\sqrt{1+(h')^2} \|\partial_2 u\|_{L^2(\Omegavertvert)}^2
            \\
            &\qquad \leq -\left(\alphap+\alpham\right)\sqrt{1+(h')^2}\um^2 + \frac{5}{16} \|\partial_2 u\|_{L^2(\Omegavertvert)}^2,
        \end{align}
        proving the claim.
    \end{itemize}
\end{proof}

\begin{lemma}[Energy bound]
    \label{lemma_nonlin_energy_bound}
    Let $\Omega$ be $C^{1,1}$ with $\hp=1+\hm$, $u_0\in L^2(\Omega)$ and $0<\alpha \in L^\infty$ satisfy \eqref{nonlin_coercivity_assumption}. Then
    \begin{align}
        \|u\|_2^2 &\leq \|u_0\|_2^2 + 16(1+\|\alpha^{-1}\|_\infty)^2\Gamma \Ra^2
        \\
        \langle \|\nabla u\|_2^2\rangle + \langle (2\alpha+\kappa)u_\tau^2\rangle_{\gammam\cup\gammap}
        &\leq (1+\max h-\min h) \Nu \Ra
    \end{align}
\end{lemma}

\begin{proof}
    By Lemma \ref{lemma_nonlin_energy_balance}, Lemma \ref{lemma_nonlin_coercivity}, Hölder's and Young's inequality
    \begin{align}
        \frac{1}{2\Pr}\frac{d}{dt} \|u\|_2^2 &+ \frac{1}{4}(1+\|\alpha^{-1}\|)^{-1}\|u\|_{H^1}^2
        \\
        &\leq 
        \frac{1}{2\Pr}\frac{d}{dt} \|u\|_2^2 + \|\nabla u\|_2^2 + \int_{\gammam\cup\gammap} (2\alpha+\kappa) u_\tau^2
        \\
        &= \Ra \int_{\Omega} \theta u_2
        \\
        &\leq \frac{1}{4\epsilon}\Gamma\Ra^2 + \epsilon\|u\|_2^2
        \label{jnfjaksdf}
    \end{align}
    for any $\epsilon>0$, where we also used the maximum principle \eqref{maximum_principle}. Choosing $\epsilon=\frac{1}{8}(1+\|\alpha^{-1}\|_\infty)^{-1}$ one gets
    \begin{align}
        \frac{1}{2\Pr}\frac{d}{dt} \|u\|_2^2 + \frac{1}{8}(1+\|\alpha^{-1}\|_\infty)^{-1}\|u\|_{H^1}^2
        &\leq 2 (1+\|\alpha^{-1}\|_\infty)\Gamma\Ra^2.
    \end{align}
    Now, Grönwall's inequality yields
    \begin{align}
        \|u\|_2^2 \leq \|u_0\|_2^2 + 16(1+\|\alpha^{-1}\|_\infty)^2\Gamma \Ra^2,
        \label{jnibjsdfg}
    \end{align}
    proving the first claim.
    
    In order to prove the long-time bound, notice that by \eqref{lemma_nonlin_energy_balance},
    \begin{align}
        \frac{1}{2\Pr}\frac{d}{dt} \|u\|_2^2 + \|\nabla u\|_2^2 + &\int_{\gammam\cup\gammap} (2\alpha+\kappa) u_\tau^2
        \\
        &= \Ra \int_{\Omega} \theta u_2
        \\
        &= \Ra \left( \int_{\Omega} (u_2-\partial_2)\theta +\int_{\Omega} \partial_2 \theta\right)
        \\
        &= \Ra \left( \int_{\Omega} (u_2-\partial_2)\theta +\int_{\gammam} n_2\right),
        \label{lknvcjkbnsdgfjkdfg}
    \end{align}
    where Stokes theorem and the boundary conditions \eqref{heat_bc} were used. 
    As $n_2 = - \frac{1}{\sqrt{1+(h')^2}}<0$ on $\gammam$ taking the long time average of \eqref{lknvcjkbnsdgfjkdfg}, using that $\|u\|_2^2$ is uniformly bounded in time by \eqref{jnibjsdfg}, we obtain
    \begin{align}
        \langle \|\nabla u\|_2^2\rangle + \langle (2\alpha+\kappa)u_\tau^2\rangle_{\gammam\cup\gammap} &\leq \Ra \langle (u_2-\partial_2)\theta \rangle
        \\
        &\leq (1+\max h-\min h) \Nu \Ra,
    \end{align}
    where we used Lemma \ref{lemma_nonlin_nusselt_number_representations} in the last estimate.
\end{proof}

\begin{lemma}[Vorticity Bound]
    \label{lemma_nonlin_vorticity_bound}
    Let $\Omega$ be $C^{1,1}$ with $\hp=1+\hm$, $2<p\in 2\mathbb{N}$, $u_0\in W^{1,p}(\Omega)$ and $0<\alpha\in L^\infty(\gammam\cup\gammap)$ satisfy \eqref{nonlin_coercivity_assumption}. Then there exists a constant $C_0>0$ depending only on the Lipschitz constant of the boundary, $p$ and $\Gamma$ such that
    \begin{align}
        \|\omega\|_{p}
        &\leq C \left(\|\omega_0\|_p + C_{\alpha,\kappa}\|u_0\|_2 + \left(1+C_{\alpha,\kappa}\right) (1+\|\alpha^{-1}\|_\infty)\Ra\right),
        \label{nauifdnsaf}
    \end{align}
    where $C_{\alpha,\kappa}=1+\|\alpha+\kappa\|_\infty^{\frac{2(p-1)}{p-2}}$. Additionally, one has
    \begin{align}
        0 &= \langle \|\nabla \omega\|_2^2\rangle +2\Pr^{-1} \langle (\alpha+\kappa)u\cdot (u\cdot\nabla)u\rangle_{\gammam\cup\gammap} +2\langle(\alpha+\kappa)u\cdot\nabla p\rangle_{\gammam\cup\gammap}
        \\
        &\qquad - 2\Ra \langle(\alpha+\kappa) u_\tau n_1 \rangle_{\gammam} - \Ra \langle \omega \partial_1\theta\rangle.
        \label{adnbvimndfjg}
    \end{align}
\end{lemma}

\begin{proof}
    The proof of \eqref{nauifdnsaf} follows the same approach as the proof of Lemma \ref{lemma_vorticity_bound}. Additionally, since here the boundary profiles are the same, we can estimate the $\|\kappa\|_\infty$ term as follows. Since $\kappap(x_1)=-\kappam(x_1)$ for every $x_1$ either $\kappap(x_1)$ or $\kappam(x_1)$ is non-negative, and as $\alpha>0$ one has
    \begin{align}
        \|\kappa\|_\infty &= \|\max \lbrace \kappam,\kappap\rbrace\|_\infty
        \leq \|\max \lbrace \alpham+\kappam,\alphap+\kappap\rbrace\|_\infty
        = \|\alpha+\kappa\|_\infty.
    \end{align}
    Therefore, the bound of Lemma \ref{lemma_vorticity_bound} can be estimated as
    \begin{align}
        \|\omega\|_{p}
        &\lesssim \|\omega_0\|_p + \left(1+\|\kappa\|_\infty+\|\alpha+\kappa\|_\infty^{\frac{p}{p-2}}\right)\|\alpha+\kappa\|_\infty \|u_0\|_2
        \\
        &\qquad + \left(1+\left(1+\|\kappa\|_\infty+\|\alpha+\kappa\|_\infty^{\frac{p}{p-2}}\right)\|\alpha+\kappa\|_\infty\right) (1+\|\alpha^{-1}\|_\infty)\Ra
        \\
        &\lesssim \|\omega_0\|_p + \left(1+\|\alpha+\kappa\|_\infty+\|\alpha+\kappa\|_\infty^{\frac{p}{p-2}}\right)\|\alpha+\kappa\|_\infty \|u_0\|_2
        \\
        &\qquad + \left(1+\left(1+\|\alpha+\kappa\|_\infty+\|\alpha+\kappa\|_\infty^{\frac{p}{p-2}}\right)\|\alpha+\kappa\|_\infty\right) (1+\|\alpha^{-1}\|_\infty)\Ra
        \\
        &\lesssim \|\omega_0\|_p + \left(1+\|\alpha+\kappa\|_\infty^{\frac{2(p-1)}{p-2}}\right)\|u_0\|_2
        \\
        &\qquad + \left(1+\|\alpha+\kappa\|_\infty^{\frac{2(p-1)}{p-2}}\right) (1+\|\alpha^{-1}\|_\infty)\Ra,
        \label{njasdfjnasd}
    \end{align}
    where we used Young's inequality in the last estimate.

    In order to prove the long time bound, note that as in \eqref{cvjnbjbcvxbni} one finds
    \begin{align}
        \frac{1}{2\Pr}\frac{d}{dt}&\left(\|\omega\|_2^2 + 2\int_{\gammam\cup\gammap} (\alpha+\kappa) u_\tau^2\right) + \|\nabla\omega\|_2^2
        \\
        &= - \frac{1}{\Pr} \int_{\gammam\cup\gammap} (\alpha+\kappa) u\cdot \nabla u_\tau^2 - 2\int_{\gammam\cup\gammap} (\alpha+\kappa) u\cdot\nabla p
        \\
        &\qquad + 2 \Ra \int_{\gammam} (\alpha+\kappa)u_2 + \Ra \int_{\Omega} \omega \partial_1 \theta
        \\
        &= - \frac{2}{\Pr} \int_{\gammam\cup\gammap} (\alpha+\kappa) u\cdot(u\cdot\nabla)u - 2\int_{\gammam\cup\gammap} (\alpha+\kappa) u\cdot\nabla p
        \\
        &\qquad + 2 \Ra \int_{\gammam} (\alpha+\kappa)u_\tau n_1 + \Ra \int_{\Omega} \omega \partial_1 \theta,
        \label{njasdjfnasdfk}
    \end{align}
    where in the last identity we used that
    \begin{align}
        u\cdot \nabla u_\tau^2 = u\cdot \nabla (u\cdot u) = 2 u\cdot(u\cdot\nabla) u,
    \end{align}
    and
    \begin{align}
        u_\tau n_1 = u_\tau \tau_2 = u_2
    \end{align}
    on $\gammam\cup\gammap$. By Hölder's inequality, trace theorem, and Lemma \ref{lemma_elliptic_regularity_periodic_domain} we find
    \begin{align}
        \|\omega\|_2^2 + 2\int_{\gammam\cup\gammap} (\alpha+\kappa) u_\tau^2
        &\lesssim \|\omega\|_2^2 + \|\alpha+\kappa\|_\infty \|u\|_{H^1}^2
        \\
        &\lesssim (1+\|\alpha+\kappa\|_\infty) \|\omega\|_2^2 + (1+\|\kappa\|_\infty^2)\|\alpha+\kappa\|_\infty \|u\|_2^2
        \\
        &\lesssim (1+\|\alpha+\kappa\|_\infty) \|\omega\|_p^2 + (1+\|\kappa\|_\infty^2)\|\alpha+\kappa\|_\infty \|u\|_2^2,
    \end{align}
    which is uniformly bounded in time by \eqref{njasdfjnasd} and Lemma \ref{lemma_nonlin_energy_bound}. Therefore, taking the long time average of \eqref{njasdjfnasdfk} we obtain
    \begin{align}
        \langle \|\nabla \omega\|_2^2\rangle &= -2\Pr^{-1} \langle (\alpha+\kappa)u\cdot (u\cdot\nabla)u\rangle_{\gammam\cup\gammap} -2\langle(\alpha+\kappa)u\cdot\nabla p\rangle_{\gammam\cup\gammap}
        \\
        &\qquad + 2\Ra \langle(\alpha+\kappa) u_\tau n_1 \rangle_{\gammam} + \Ra \langle \omega \partial_1\theta\rangle,
    \end{align}
    proving \eqref{adnbvimndfjg}.
\end{proof}

The pressure satisfies
\begin{alignat}{2}
    \Delta p &= - \Pr^{-1} \nabla u\colon \nabla u^T + \Ra \partial_2 \theta & \quad \textnormal{ in } &\Omega
    \label{nonlin_pressure_pde}
    \\
    n \cdot \nabla p &= -\Pr^{-1}\kappa u_\tau^2 + 2 \tau \cdot \nabla ((\alpha+\kappa)u_\tau) + n_2\Ra \theta & \quad \textnormal{ on } &\gammam\cup\gammap,\quad
    \label{nonlin_pressure_bc}
\end{alignat}
where $\colon$ denotes the tensor contraction, i.e. $\nabla u\colon \nabla u^T = \partial_i u_j \partial_j u_i$.

\eqref{nonlin_pressure_pde} follows immediately by taking the divergence of \eqref{navier_stokes} and using \eqref{incompressibility} and $\nabla \left((u\cdot\nabla)u\right)=\nabla u\colon\nabla u^T$. To show \eqref{nonlin_pressure_bc}, dotting $n$ into \eqref{navier_stokes} yields
\begin{align}
    \Pr^{-1}n\cdot u_t + \Pr^{-1} n\cdot (u\cdot\nabla) u + n\cdot \nabla p - n\cdot\Delta u = \Ra \theta n_2.
    \label{uisdfnign}
\end{align}
Taking the time derivative of $\eqref{no_penetration_bc}$, one obtains
\begin{align}
    n\cdot u_t = 0.
    \label{nibsjunifndsfa}
\end{align}
The second term on the left-hand side of \eqref{uisdfnign} can be calculated, using \eqref{n_tau_grad_u}, as
\begin{align}
    n \cdot (u\cdot \nabla) u = u_\tau n \cdot (\tau \cdot \nabla) u = \kappa u_\tau^2.
    \label{bnixjubcxnib}
\end{align}
By \eqref{id_Delta_u_is_nablaPerp_omega}, the third term on the left-hand side of \eqref{uisdfnign} can be written as
\begin{align}
    n\cdot \Delta u = n \cdot \nabla^\perp \omega = -\tau\cdot\nabla \omega.
    \label{jinsfijugnsdf}
\end{align}
Taking the derivative of \eqref{vorticity_bc} along the boundary and using \eqref{jinsfijugnsdf} it follows that
\begin{align}
    -2\tau\cdot\nabla\left((\alpha+\kappa) u_\tau\right) = \tau \cdot \nabla \omega = -n\cdot \Delta u.
    \label{naifnsdaf}
\end{align}
Combining \eqref{uisdfnign}, \eqref{nibsjunifndsfa}, \eqref{bnixjubcxnib} and \eqref{naifnsdaf} yields \eqref{nonlin_pressure_bc}.

\begin{lemma}[Pressure Bound]
    \label{lemma_nonlin_pressure_bound}
    Let $\Omega$ be $C^{1,1}$ with $\hp=\hm+1$, $r>2$, $u\in H^2(\Omega)$ and $\alpha \in W^{1,\infty}(\gammam\cup\gammap)$. Then
    \begin{align}
        \|p\|_{H^1} \lesssim \left(\frac{1+\|\kappa\|_\infty}{\Pr} \|u\|_{W^{1,r}}+ \|\dot\alpha+\dot\kappa\|_\infty \right)\|u\|_{H^1} + \|\alpha+\kappa\|_\infty\|u\|_{H^2} + \Ra,
    \end{align}
    where the implicit constant only depends on $r$, $\Gamma$ and the Lipschitz constant of the boundary.
\end{lemma}

\begin{remark}
Note that the bound in Lemma \ref{lemma_nonlin_pressure_bound} is significantly worse than the one in Lemma \ref{lemma_pressure_bound}. The $\|u\|_{H^2}$ term will later result in the smallness condition on $\|\alpha+\kappa\|_\infty$. In particular, this only allows boundaries that are close to free-slip.
\end{remark}

\begin{proof}
    Testing \eqref{nonlin_pressure_pde} with $p$, integrating by parts and using \eqref{nonlin_pressure_bc} it follows that
    \begin{align}
        \int_{\Omega} p\Delta p &= -\|\nabla p\|_2^2 + \int_{\partial\Omega} p n\cdot \nabla p 
        \\
        &= -\|\nabla p\|_2^2 - \Pr^{-1} \int_{\partial\Omega}\kappa p u_\tau^2 + 2 \int_{\partial\Omega} p \tau \cdot \nabla \left((\alpha+\kappa) u_\tau\right) + \Ra \int_{\partial\Omega} p n_2\theta.
        \label{uiadfniasf}
    \end{align}
    By \eqref{nonlin_pressure_pde} and integration by parts, the left-hand side of \eqref{uiadfniasf} also satisfies
    \begin{align}
        \int_{\Omega} p\Delta p &= - \Pr^{-1}\int_{\Omega} p \nabla u \colon \nabla u^T + \Ra \int_{\Omega} p \partial_2 \theta
        \\
        &= - \Pr^{-1}\int_{\Omega} p \nabla u \colon \nabla u^T + \Ra \int_{\Omega} n_2 p \theta - \Ra \int_{\Omega} \theta \partial_2 p
        \label{uianfisdnaf}
    \end{align}
    Subtracting \eqref{uianfisdnaf} from \eqref{uiadfniasf} yields
    \begin{align}
        \|\nabla p\|_2^2 &= - \Pr^{-1} \int_{\partial\Omega}\kappa p u_\tau^2 + 2 \int_{\partial\Omega} p \tau \cdot \nabla \left((\alpha+\kappa) u_\tau\right)
        \\
        &\qquad + \Pr^{-1}\int_{\Omega} p \nabla u \colon \nabla u^T + \Ra \int_{\Omega} \theta \partial_2 p.
        \label{lkanflkasndf}
    \end{align}
    Next, we bound the terms on the right-hand side of \eqref{lkanflkasndf} individually.

    For the first term on the right-hand side of \eqref{lkanflkasndf}, Hölder's inequality, trace theorem, and Sobolev embedding yield
    \begin{align}
        - \int_{\partial\Omega}\kappa p u_\tau^2 &\lesssim \|\kappa\|_\infty \|p u^2\|_{W^{1,1}}
        \\
        &\lesssim \|\kappa\|_\infty \left(\|p u^2\|_1 + \| p u \nabla u\|_1+\|u^2 \nabla p\|_1\right)
        \\
        &\lesssim \|\kappa\|_\infty \left(\|p \|_3 \|u\|_3 \|u\|_3 + \|p\|_4 \|u\|_4 \|\nabla u\|_2 + \|\nabla p\|_2 \|u\|_4 \|u\|_4\right)
        \\
        &\lesssim \|\kappa\|_\infty \|p\|_{H^1}\|u\|_{H^1}^2.
        \label{adjkfnasdf}
    \end{align}

    For the second term on the right-hand side of \eqref{lkanflkasndf}, Hölder's inequality, trace theorem, and Young's inequality yield
    \begin{align}
        \int_{\partial\Omega} |p \tau \cdot \nabla \left((\alpha+\kappa) u_\tau\right)|
        &\leq \|\dot\alpha+\dot\kappa\|_\infty \|p\|_{L^2(\gammam\cup\gammap)} \|u\|_{L^2(\gammam\cup\gammap)}
        \\
        &\qquad
        + \|\alpha+\kappa\|_\infty \|p\|_{L^2(\gammam\cup\gammap)} \|\nabla u\|_{L^2(\gammam\cup\gammap)}
        \\
        &\lesssim \|\dot\alpha+\dot\kappa\|_\infty \|p\|_{H^1} \|u\|_{H^1} + \|\alpha+\kappa\|_\infty \|p\|_{H^1} \|u\|_{H^2}
        \\
        \label{iuauisdbf}
    \end{align}

    In order to estimate the third term on the right-hand side of \eqref{lkanflkasndf}, Hölder's inequality and Sobolev embedding imply
    \begin{align}
        \|p\nabla u\colon\nabla u^T \|_1 \leq \|p\|_q \|\nabla u\|_2 \|\nabla u\|_r \leq \|p\|_{H^1}\|\nabla u\|_2 \|\nabla u\|_r
        \label{bnibhbfsd}
    \end{align}
    for any $2< r,q<\infty$ with $\frac{1}{r}+\frac{1}{q}=\frac{1}{2}$.

     Using Hölder's inequality and \eqref{maximum_principle}, the fourth term on the right-hand side of \eqref{lkanflkasndf} can be bounded by
    \begin{align}
        -\Ra \int \theta \partial_2 p \lesssim \Ra \|p\|_{H^1}.
        \label{uibiubbnnsdaf}
    \end{align}

    Combining \eqref{lkanflkasndf}, \eqref{adjkfnasdf}, \eqref{iuauisdbf}, \eqref{bnibhbfsd}, and \eqref{uibiubbnnsdaf} one gets
    \begin{align}
        \|\nabla p\|_2^2
        &\lesssim \|p\|_{H^1}\bigg[\left(\frac{1+\|\kappa\|_\infty}{\Pr} \|u\|_{W^{1,r}}+ \|\dot\alpha+\dot\kappa\|_\infty \right)\|u\|_{H^1}
        \\
        &\qquad + \|\alpha+\kappa\|_\infty\|u\|_{H^2} + \Ra \bigg]
        \label{lnajksdfnasjkdfn}
    \end{align}
    for any $r>2$.
    As the pressure is only defined up to a constant, we can choose this constant such that $p$ is average free, and therefore, Poincaré's inequality yields $\|p\|_{H^1}^2\leq \|\nabla p\|_2^2 + \|p\|_2^2 \leq (1+C) \|\nabla p\|_2^2$, which together with \eqref{lnajksdfnasjkdfn} implies
    \begin{align}
        \|p\|_{H^1}^2 &\lesssim \|p\|_{H^1}\bigg[\left(\frac{1+\|\kappa\|_\infty}{\Pr} \|u\|_{W^{1,r}}+ \|\dot\alpha+\dot\kappa\|_\infty \right)\|u\|_{H^1}
        \\
        &\qquad + \|\alpha+\kappa\|_\infty\|u\|_{H^2} + \Ra \bigg]
    \end{align}
    and therefore, dividing by $\|p\|_{H^1}$ yields the claim.
\end{proof}

\begin{remark}
We remark that Lemma \ref{lemma_nusselt_bounded_by_Hx} produces the stricter bound. The main difference is the estimate \eqref{direct_method_gradTheta}, respectively \eqref{hvbuxbucvh}, where the latter produces the additional curvature term.
\end{remark}

\begin{proof}[\hypertarget{proof_theorem_nonlinear_1_2_bound}{Proof of Theorem} \ref{theorem_nonlinear_1_2_bound}]
\\
    The proof is a slight modification of \eqref{nusselt_bound_H1_in_curved_lemma} in Lemma \ref{lemma_nusselt_bounded_by_Hx} using the top boundary instead of the bottom one and taking advantage of the same boundary profiles. Averaging \eqref{nonlin_nu_representation_shifted} in Lemma \ref{lemma_nonlin_nusselt_number_representations} over $z\in (1-\delta,1)$ yields
    \begin{align}
        \Nu = \delta^{-1} \langle \np \cdot (u-\nabla) \theta\rangle_{\Omegatildedelta}
        =\delta^{-1}\langle \np \cdot u\theta\rangle_{\Omegatildedelta}-\delta^{-1}\langle \np \cdot \nabla \theta \rangle_{\Omegatildedelta},\label{bjsiobio}
    \end{align}
    where
    \begin{align}
        \Omegatildedelta = \lbrace (x_1,x_2) \ \vert \  0\leq x_1\leq \Gamma, 1+h(x_1)-\delta \leq x_2\leq 1+ h(x_1) \rbrace.
    \end{align}
    In order to estimate the first term on the right-hand side of \eqref{bjsiobio}, the fundamental theorem of calculus implies
    \begin{align}
        |\np\cdot u|(x) &= \left| \np\cdot u\vert_{\gammap} + \int_{1+h(x_1)}^{x_2}\partial_2(\np\cdot u)\ dz \right| \leq \int_{1+h(x_1)-\delta}^{1+h(x_1)}|\np\cdot\partial_2 u|\ dz
        \\
        &\leq \delta^{\frac{1}{2}}\|\nabla u(x_1,\cdot)\|_{L^2(\Omegavertvert)}
        \label{xbuihcvuibxu}
    \end{align}
    for $x\in \Omegatildedelta$, where we used the boundary conditions \eqref{no_penetration_bc}, the fact that $\np$ is constant in $x_2$-direction, Hölder's inequality and
    \begin{align}
        L^2(\Omegavertvert) = L^2(h(x_1),1+h(x_1)).
    \end{align}
    Analogously, for $\theta$ and $x\in \Omegatildedelta$ it holds that
    \begin{align}
        |\theta|(x)\leq \delta^\frac{1}{2} \|\nabla\theta(x_1,\cdot)\|_{L^2(\Omegavertvert)}
        \label{ansdfklsnf}
    \end{align}
    as $\theta=0$ on $\gammap$, and combining \eqref{xbuihcvuibxu} and \eqref{ansdfklsnf}, implies
    \begin{align}
        |\np \cdot u \theta|\leq \delta \|\nabla u\|_{L^2(\Omegavertvert)}\|\nabla\theta\|_{L^2(\Omegavertvert)},
    \end{align}
    which after integration over $\Omegatildedelta$ yields
    \begin{align}
        \int_{\Omegatildedelta} |\np \cdot u \theta| \leq \delta^2 \|\nabla u\|_2 \|\nabla\theta\|_2,\label{buaisdfabdf}
    \end{align}
    where we used Hölder's inequality for the integration with respect to $x_1$.
    
    In order to estimate the second term on the right-hand side of \eqref{bjsiobio}, integration by parts and the boundary condition \eqref{heat_bc} imply
    \begin{align}
        \left|\int_{\Omegatildedelta} \np \cdot \nabla \theta \right| \leq \int_{\gammap} |\theta| + \int_{\gammam+1-\delta} |\np\cdot\nm\theta| + \int_{\Omegatildedelta} |\theta \nabla \cdot \np|.
        \label{hvbuxbucvh}
    \end{align}
    Next, we focus on the divergence of $n$. The bottom boundary can be parameterized by $(x_1,h(x_1))$, implying that the tangential is parallel to $(1,h')$, which yields that the unit tangent and unit normal vectors are given by
    \begin{align}
        \taupm = \mp \left(1+(h')^2\right)^{-\frac{1}{2}} \begin{pmatrix}1\\h'\end{pmatrix}
        \quad\textnormal{and}\quad
        \npm = \pm \left(1+(h')^2\right)^{-\frac{1}{2}} \begin{pmatrix}-h'\\1\end{pmatrix}\quad\label{xjubcvbhi}
    \end{align}
    and therefore,
    \begin{align}
        \frac{d}{dx_1}\taupm = -\frac{h''}{1+(h')^2}\npm.\label{iuabsdfus}
    \end{align}
    Changing the parameterization to arc length $\lambda$, one gets
    \begin{align}
        \lambda(x_1) &= \int_0^{x_1} \sqrt{1+(h'(s))^2} \ ds\textnormal{ on }\gammam,
        \\
        \lambda(x_1) &= \int_{\Gamma}^{x_1} \sqrt{1+(h'(s))^2} \ ds\textnormal{ on }\gammap,
    \end{align}
    implying
    \begin{align}
        \frac{d}{dx_1}\lambda(x_1)=\mp \sqrt{1+(h'(s))^2}
    \end{align}
    and therefore, \eqref{iuabsdfus} implies
    \begin{align}
        \tau\cdot\nabla \taupm = \frac{d}{d\lambda} \taupm = \mp \left(1+(h'(s))^2\right)^{-\frac{1}{2}} \frac{d}{dx_1} \taupm = \pm \left(1+(h'(s))^2\right)^{-\frac{3}{2}} h'' \npm,
    \end{align}
    which by the definition of the curvature, i.e. \eqref{definition_curvature}, yields
    \begin{align}
        \kappa = \pm \left(1+(h'(s))^2\right)^{-\frac{3}{2}} h''.
    \end{align}
    Taking the divergence of \eqref{xjubcvbhi}, it follows that
    \begin{align}
        \nabla \cdot \npm = \mp \left(1+(h'(s))^2\right)^{-\frac{3}{2}} h'' = - \kappa.
        \label{bufidsafbsdf}
    \end{align}
    Combining \eqref{hvbuxbucvh} and \eqref{bufidsafbsdf}, and using that $\theta$ is bounded by \eqref{maximum_principle}
    \begin{align}
        \left|\int_{\Omegatildedelta} \np \cdot \nabla \theta \right| &\leq \int_{\gammap} |\theta| + \int_{\gammam+1-\delta} |\np\cdot\nm\theta| + \int_{\Omegatildedelta} |\theta \nabla \cdot \np|
        \leq C+ \delta \Gamma \|\kappa\|_\infty.
        \label{njuafisdna}
    \end{align}

    By \eqref{bjsiobio}, \eqref{buaisdfabdf} and \eqref{njuafisdna} we obtain
    \begin{align}
        \Nu \leq C\left( \delta \langle \|\nabla u\|_2^2\rangle^\frac{1}{2}\langle \|\nabla \theta\|_2^2\rangle^\frac{1}{2} + \frac{1}{\delta}\right) + \|\kappa\|_\infty,\label{bnaidsfbsadf}
    \end{align}
    where we used $\langle fg\rangle\leq \langle f^2\rangle^\frac{1}{2} \langle g^2\rangle^\frac{1}{2}$ due to Young's inequality. By Lemma \ref{lemma_nonlin_coercivity} and Lemma \ref{lemma_nonlin_energy_bound}, one gets
    \begin{alignat}{1}
        C_\alpha^{-1}\langle\|\nabla u\|_2^2\rangle &\lesssim \langle\|\nabla u\|_2^2\rangle + \langle (2\alpha+\kappa) u_\tau^2\rangle_{\gammam\cup\gammap} \lesssim \Nu\Ra,\label{nasdfjinsajfdan}
    \end{alignat}
    where
    \begin{align}
        C_\alpha = \begin{cases}
            1 & \textnormal{ if } |\kappa|<2\alpha
            \\
            1+\|\alpha^{-1}\|_\infty & \textnormal{ if } |\kappa|\leq 2\alpha + \frac{1}{4\sqrt{1+(h'(x_1))^2}}\min \left\{1,\sqrt{\alpha}\right\}
        \end{cases}
    \end{align}
    and combining \eqref{bnaidsfbsadf}, \eqref{nasdfjinsajfdan} and Lemma \eqref{lemma_nonlin_nusselt_number_representations} it holds
    \begin{align}
        \Nu \leq C \left(\delta C_\alpha^\frac{1}{2}\Nu \Ra^\frac{1}{2}+\delta^{-1}\right) + \|\kappa\|_\infty.
    \end{align}
    Thus, choosing $\delta = C_\alpha^{-\frac{1}{4}}\Nu^{-\frac{1}{2}}\Ra^{-\frac{1}{4}}$ yields
    \begin{align}
        \Nu \lesssim C_\alpha^{\frac{1}{4}} \Nu^\frac{1}{2}\Ra^\frac{1}{4} + \|\kappa\|_\infty \lesssim \epsilon \Nu + \epsilon^{-1}C_\alpha^\frac{1}{2}\Ra^\frac{1}{2} + \|\kappa\|_\infty,\label{ojsdafobnsdufiabsf}
    \end{align}
    for any $\epsilon>0$, where in the last estimate we used Young's inequality. Choosing $\epsilon$ sufficiently small, one can compensate the $\Nu$ term on the right-hand side of \eqref{ojsdafobnsdufiabsf}, concluding the proof.    
\end{proof}

\subsection{Background Field Method Application}
\label{section:background_field_method_application}

\begin{figure}
    \hrule
    \vspace{0.5\baselineskip}
    \begin{center}
        \includegraphics[width=\textwidth]{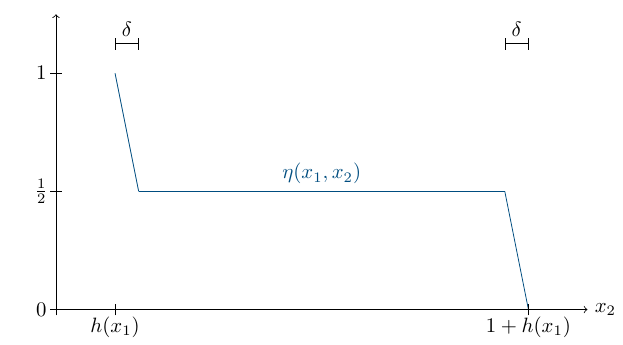}
        \vspace*{-7mm}
        \caption{Illustration of the background field profile. The picture corresponds to a vertical slice of $\Omega$ at $x_1$.}
        \label{fig:rb_background_field_profile}
    \end{center}
    \hrule
\end{figure}

We define the background profile by
\begin{align}
    \eta(x) =
    \begin{cases}
        \begin{aligned}
            &\frac{1+h(x_1)-x_2}{2\delta} & &\textnormal{ for }& 1+h(x_1)-\delta \leq &\ x_2 \leq 1+h(x_1)
            \\
            &\frac{1}{2} & &\textnormal{ for }& h(x_1)+\delta \leq &\ x_2 \leq 1+h(x_1)-\delta
            \\
            &\frac{2\delta + h(x_1)-x_2}{2\delta} & &\textnormal{ for }& h(x_1) \leq &\ x_2 \leq h(x_1)+\delta
        \end{aligned}        
    \end{cases}
\end{align}
for any $\delta\leq\frac{1}{2}$ as illustrated in Figure \ref{fig:rb_background_field_profile}.

Also, note that
\begin{align}
    \nabla\eta
    &=
    \begin{cases}
        \begin{aligned}
            &0 & &\textnormal{ if } h(x_1)+\delta \leq \ x_2 \leq 1+h(x_1)-\delta
            \\
            &\frac{1}{2\delta}
            \begin{pmatrix}
                h'
                \\
                -1
            \end{pmatrix} & &\textnormal{ else }
        \end{aligned}        
    \end{cases}
    \\
    &=
    \begin{cases}
        \begin{aligned}
            &0 & &\textnormal{ if } h(x_1)+\delta \leq \ x_2 \leq 1+h(x_1)-\delta
            \\
            &\frac{\sqrt{1+(h')^2}}{2\delta}\nm  & &\textnormal{ else }
        \end{aligned}
    \end{cases}
    \label{ifsdfasfsgjflg}
\end{align}
for almost every $x\in\Omega$, implying
\begin{align}
    \|\nabla \eta\|_\infty \lesssim \delta^{-1}
\end{align}
and, because of its support,
\begin{align}
    \|\nabla\eta\|_2^2 = \int_0^\Gamma \int_{h(x_1)}^{h(x_1)+\delta} |\nabla \eta|^2\ dx_2+\int_{1+h(x_1)-\delta}^{1+h(x_1)} |\nabla \eta|^2 \ dx_2 \ dx_1 \lesssim \|\nabla\eta\|_\infty \delta \lesssim \delta^{-1}.
    \label{ztvcubtcuzzvub}
\end{align}
Additionally, define
\begin{align}
    \varsigma=\theta-\eta.
    \label{nonlin_def_varsigma}
\end{align}
Then $\varsigma$ fulfills
\begin{alignat}{2}
    \varsigma_t + u\cdot\nabla \eta + u\cdot\nabla \varsigma - \Delta\eta - \Delta\varsigma &=0 & \quad \textnormal{ in } &\Omega,
    \label{nonlin_varsigma_pde}
    \\
    \varsigma &= 0 & \quad \textnormal{ on } &\gammam\cup\gammap.
    \label{nonlin_varsigma_bc}
\end{alignat}
Testing \eqref{nonlin_varsigma_pde} with $\varsigma$ and integrating by parts, one finds
\begin{align}
    0&=\frac{1}{2}\frac{d}{dt}\|\varsigma\|_2^2 + \int_\Omega \varsigma u\cdot\nabla \eta + \int_{\Omega} \varsigma u\cdot\nabla \varsigma - \int_{\Omega} \varsigma\Delta\eta - \int_{\Omega} \varsigma \Delta \varsigma
    \\
    &=\frac{1}{2}\frac{d}{dt}\|\varsigma\|_2^2 + \int_\Omega \varsigma u\cdot\nabla \eta + \int_{\Omega} \nabla \varsigma \cdot \nabla \eta + \|\nabla \varsigma\|_2^2,
\end{align}
where in the second identity the boundary terms vanish because of \eqref{nonlin_varsigma_bc}. Additionally, we used
\begin{align}
    2\int_{\Omega} \varsigma u\cdot \nabla \varsigma = \int_{\Omega} u\cdot\nabla (\varsigma^2) = \int_{\partial\Omega} u\cdot n\varsigma^2 - \int_{\Omega} \varsigma^2 \nabla \cdot u = 0,
\end{align}
which is due to \eqref{no_penetration_bc} and \eqref{incompressibility}. Taking the long-time average, we obtain
\begin{align}
    \langle \nabla \varsigma \cdot \nabla \eta \rangle&=-\langle \varsigma u\cdot\nabla \eta \rangle - \langle\|\nabla \varsigma\|_2^2\rangle\label{basudfbahzsdf}
\end{align}
since $\|\varsigma\|_2^2$ is uniformly bounded in time as both $\|\theta\|_2^2$ and $\|\eta\|_2^2$ are uniformly bounded in time. Because of
\begin{align}
    \|\nabla \theta\|_2^2 = \|\nabla (\varsigma+\eta)\|_2^2 = \|\nabla \varsigma\|_2^2 + 2\int_{\Omega} \nabla\varsigma\cdot\nabla\eta + \|\nabla\eta\|_2^2
\end{align}
Lemma \eqref{lemma_nonlin_nusselt_number_representations} yields
\begin{align}
    \Nu &= \langle \|\nabla\theta\|_2^2 \rangle = \langle\|\nabla \varsigma\|_2^2\rangle + 2\langle \nabla\varsigma\cdot\nabla\eta \rangle + \langle\|\nabla\eta\|_2^2\rangle
    \\
    &= \langle\|\nabla\eta\|_2^2\rangle - 2\langle \varsigma u\cdot\nabla\eta\rangle-\langle\|\nabla \varsigma\|_2^2\rangle,
    \label{bsahfbsahfhj}
\end{align}
where in the last identity we used \eqref{basudfbahzsdf}.

\begin{lemma}
    \label{lemma_nonlin_varsigmaucdotnalbaeta}
    Let $\Omega$ be $C^{1,1}$ with $\hp=1+\hm$ and $u\in H^2(\Omega)$. Then for any $a,\epsilon>0$ it holds
    \begin{align}
        2\left|\langle \varsigma u\cdot\nabla \eta\rangle\right|
        & \leq \delta^6 (a\epsilon)^{-1} C\langle\|\partial_2 u\|_2^2\rangle + a\epsilon\langle\|\partial_2^2 u\|_2^2\rangle + \frac{1}{2}\langle\|\nabla \varsigma\|_2^2\rangle,
    \end{align}
    where $C= 4(1+\|h'\|_\infty)^2$.
\end{lemma}

\begin{remark}
This estimate is the crucial one corresponding to Lemma \ref{lemma_nusselt_bounded_by_Hx}. In fact, the proof uses a similar argument as in the direct method.
\end{remark}

\begin{proof}
    By \eqref{ifsdfasfsgjflg} it holds
    \begin{align}
        2 \int_{\Omega} \varsigma u\cdot\nabla \eta &= \frac{1}{\delta} \int_0^\Gamma\int_{h(x_1)}^{h(x_1)+\delta} \sqrt{1+(h')^2} \varsigma u\cdot\nm 
        \\
        &\qquad + \frac{1}{\delta} \int_0^\Gamma\int_{1+h(x_1)-\delta}^{1+h(x_1)} \sqrt{1+(h')^2} \varsigma u\cdot\nm
        \label{absnduifbasf}
    \end{align}
    and as $u\cdot \nm=0$ on $\gammam\cup\gammap$
    \begin{align}
        |u\cdot \nm|(x) &= \left|u\cdot\nm\vert_{\gammam} + \int_{h(x_1)}^{x_2} \partial_2 (u\cdot\nm)\ dz\right|
        \\
        &\leq \|\partial_2(u\cdot\nm)\|_{L^\infty(\Omegavertvert)} \int_{h(x_1)}^{h(x_1)+\delta} \partial_2 (u\cdot\nm)\ dz
        \\
        &= \delta\|\partial_2(u\cdot\nm)\|_{L^\infty(\Omegavertvert)}\label{basdhfbasdhf}
    \end{align}
    for $x_2\in (h(x_1),h(x_1)+\delta)$, where $\Omegavertvert(x_1)=\lbrace (y_1,y_2)\in \Omega\ \vert\ y_1=x_1\rbrace$. The equivalent estimate yields
    \begin{align}
        |u\cdot \nm|(x)
        &\leq \delta \|\partial_2(u\cdot\nm)\|_{L^\infty(\Omegavertvert)}
        \label{bsdhfbasdufbhasdf}
    \end{align}
    for $x_2\in (1+h(x_1)-\delta,1+h(x_1))$. Similarly as $\varsigma\vert_{\gammam\cup\gammap}=0$
    \begin{align}
        |\varsigma|(x)\leq \delta^\frac{1}{2}\|\partial_2 \varsigma\|_{L^2(\Omegavertvert)}
        \label{bszdfhubzudg}
    \end{align}
    for $x_2\in (h(x_1),h(x_1)+\delta)\cup(1+h(x_1)-\delta,1+h(x_1))$. In order to get a further Poincaré type estimate note that since
    \begin{align}
        \int_{h(x_1)}^{1+h(x_1)} \partial_2 (u\cdot\nm) = n_2 u\cdot\nm\vert_{\gammap}+n_2 u\cdot\nm\vert_{\gammam} =0.
    \end{align}
    Therefore for every $0\leq x_1\leq \Gamma$ there exists $\tilde x_2\in (h(x_1),1+h(x_2))$ with $\partial_2(u\cdot\nm)(x_1,\tilde x_2)=0$. Applying the fundamental theorem of calculus again it holds
    \begin{align}
        (\partial_2(u\cdot\nm))^2(x)
        &\leq \left| (\partial_2(u\cdot\nm))^2(x_1,\tilde x_2)+\int_{\tilde x_2}^{x_2} \partial_2 (\partial_2(u\cdot\nm))^2 \ dz\right|
        \\
        &\leq 2 \int_{h(x_1)}^{1+h(x_1)}|\partial_2 (u\cdot\nm) \partial_2^2 (u\cdot\nm)|\ dz
        \\
        &\leq 2 \|\partial_2 u\|_{L^2(\Omegavertvert)} \|\partial_2^2 u\|_{L^2(\Omegavertvert)}
        \label{auisdfbasidf}
    \end{align}
    for any $x\in \Omega$, where in the last estimate we used Hölder's inequality, that $n$ is constant in $x_2$-direction and $|n|=1$.
    Combining \eqref{absnduifbasf}, \eqref{basdhfbasdhf}, \eqref{bsdhfbasdufbhasdf} and \eqref{bszdfhubzudg}
    \begin{align}
        2\left|\int_{\Omega} \varsigma u\cdot\nabla \eta \right|
        &\leq 2 \delta^\frac{3}{2}\int_0^\Gamma\|\partial_2 (u\cdot\nm)\|_{L^\infty(\Omegavertvert)}\|\nabla \varsigma\|_{L^2(\Omegavertvert)} \sqrt{1+(h')^2} \ dx_1
        \\
        &\leq (2 \delta)^\frac{3}{2}\int_0^\Gamma\|\partial_2 u\|_{L^2(\Omegavertvert)}^\frac{1}{2}\|\partial_2^2 u\|_{L^2(\Omegavertvert)}^\frac{1}{2}\|\nabla \varsigma\|_{L^2(\Omegavertvert)} \sqrt{1+(h')^2} \ dx_1,
        \\[-10pt]\label{sbudfiasbdf}
    \end{align}
    where in the second estimate we used \eqref{auisdfbasidf}. Applying Young's inequality twice \eqref{sbudfiasbdf} can be bounded by
    \begin{align}
        2&\left|\int_{\Omega} \varsigma u\cdot\nabla \eta\right|
        \\
        &\qquad \leq 2\delta^\frac{3}{2}\mu C^\frac{1}{2}\int_0^\Gamma \|\partial_2 u\|_{L^2(\Omegavertvert)}\|\partial_2^2 u\|_{L^2(\Omegavertvert)} \ dx_1 + \delta^\frac{3}{2}\mu^{-1}\|\nabla \varsigma\|_2^2 
        \\
        &\qquad \leq \mu\nu\delta^\frac{3}{2} \frac{C}{4}\|\partial_2 u\|_2^2 + \mu\nu^{-1}\delta^\frac{3}{2}\|\partial_2^2 u\|_2^2 + \mu^{-1}\delta^\frac{3}{2}\|\nabla \varsigma\|_2^2,
    \end{align}
    for any $\mu,\nu >0$ where $C = 4(1+\|h'\|_\infty)^2$. Finally choosing $\mu=2\delta^{-\frac{3}{2}}$ and $\nu = 2(a\epsilon)^{-1}\delta^3$ and taking the long-time average yields the claim.
\end{proof}

\begin{proof}[\hypertarget{proof_theorem_nonlinear_main_theorem}{Proof of Theorem} \ref{theorem_nonlinear_main_theorem}]
\\
    First note that since $\|\alpha+\kappa\|_\infty<1$ and $\alpha>0$ one has
    \begin{align}
        \|\alpha^{-1}\|_\infty^{-1}=\essinf \alpha\leq \essinf_{\kappa>0} (\alpha+\kappa) \leq \|\alpha+\kappa\|_\infty \leq \bar C < 1.
        \label{nbvnbvznuvbtuvzbt}
    \end{align}
    and since $\kappa(x_1,h(x_1))=-\kappa(x_1,1+h(x_1))$ and $\alpha>0$
    \begin{align}
        \|\kappa\|_\infty &= \esssup_{\kappa\geq 0} \kappa \leq \esssup_{\kappa>0} (\alpha + \kappa)\leq \|\alpha+\kappa\|_\infty \leq \bar C <1.
        \label{dabsfkjadfsb_kappa_leq_1}
    \end{align}
    Additionally, we define
    \begin{align}
        \mathcal{A}&=\langle \|\nabla \omega\|_2^2\rangle +2\Pr^{-1} \langle (\alpha+\kappa)u\cdot (u\cdot\nabla)u\rangle_{\gammam\cup\gammap} +2\langle(\alpha+\kappa)u\cdot\nabla p\rangle_{\gammam\cup\gammap}
        \\
        &\qquad - 2\Ra \langle(\alpha+\kappa) u_\tau n_1 \rangle_{\gammam} - \Ra \langle \omega \partial_1\theta\rangle,
        \\
        \mathcal{B}&= \langle \|\nabla u\|_2^2\rangle + \langle (2\alpha+\kappa)u_\tau^2\rangle_{\gammam\cup\gammap} - (1+\max h-\min h) \Nu \Ra,
        \label{vuhixbhcuvbixhc}
        \\
        \mathcal{Q}&=M \Ra^2 + \langle \|\nabla\eta\|_2^2\rangle + \langle \|\nabla\varsigma\|_2^2\rangle + 2\langle\varsigma u\cdot\nabla \eta\rangle
        \label{iuvbhcibhxcvuhb}
        \\
        &\qquad + b\Ra^{-1}\langle\|\nabla u\|_2^2\rangle +b\Ra^{-1} \langle(2\alpha+\kappa)u_\tau^2\rangle_{\gammam\cup\gammap} - b\Ra^{-1}\mathcal{B} + a\mathcal{A},
    \end{align}
    for any $0<M,a<\infty$ and $0\leq b<(1+\max h-\min h)^{-1}$.
    First notice that by the definitions of $\mathcal{Q}$ and $\mathcal{B}$ 
    it follows that
    \begin{align}
        &
        M\Ra^2 + 2\langle \|\nabla \eta\|_2^2 \rangle-\mathcal{Q}
        \\
        &\qquad
        = \langle \|\nabla \eta\|_2^2 \rangle- \langle \|\nabla\varsigma\|_2^2\rangle - 2\langle\varsigma u\cdot\nabla \eta\rangle  -a\mathcal{A} - b(1+\max h-\min h) \Nu
        \\
        &\qquad
        = \langle \|\nabla \eta\|_2^2 \rangle- \langle \|\nabla\varsigma\|_2^2\rangle - 2\langle\varsigma u\cdot\nabla \eta\rangle - b(1+\max h-\min h) \Nu,
        \label{bauifsbafiu}
    \end{align}
    where in the last identity we used that $\mathcal{A}=0$ due to Lemma \ref{lemma_nonlin_vorticity_bound}. Using \eqref{bsahfbsahfhj} we can substitute the first three terms of the right-hand side of \eqref{bauifsbafiu} and find
    \begin{align}
        &(1 - b(1+\max h-\min h))\Nu
        \\
        &\qquad = \langle\|\nabla\eta\|_2^2\rangle - 2\langle \varsigma u\cdot\nabla\eta\rangle-\langle\|\nabla \varsigma\|_2^2\rangle - b(1+\max h-\min h)\Nu
        \\
        &\qquad = M\Ra^2 + 2\langle \|\nabla \eta\|_2^2 \rangle-\mathcal{Q}
        \\
        &\qquad \leq M\Ra^2 + C \delta^{-1} -\mathcal{Q},
        \label{ahusdfbafbhvcxubihxcvbiu}
    \end{align}
    where in the last inequality we used \eqref{ztvcubtcuzzvub}.
    The strategy now is as follows. As $b< (1+\max h-\min h)^{-1}$ if $\mathcal{Q}\geq 0$ we will get bounds for $\Nu$. Optimizing $M, a, \delta>0$ with respect to the best $\Ra$ exponent either if $\alpha$ and $\kappa$ are small or in general will yield the results.

    By Lemma \ref{lemma_nonlin_energy_bound} $\mathcal{B}\leq 0$ and therefore plugging in $\mathcal{A}$ one gets
    \begin{align}
        \mathcal{Q} &=M \Ra^2 + \langle \|\nabla\eta\|_2^2\rangle + \langle \|\nabla\varsigma\|_2^2\rangle + 2\langle\varsigma u\cdot\nabla \eta\rangle
        \\
        &\qquad + b\Ra^{-1}\langle\|\nabla u\|_2^2\rangle +b\Ra^{-1} \langle(2\alpha+\kappa)u_\tau^2\rangle_{\gammam\cup\gammap} - b\Ra^{-1}\mathcal{B} + a\mathcal{A}
        \\
        &\geq M \Ra^2 + \langle \|\nabla\eta\|_2^2\rangle + \langle \|\nabla\varsigma\|_2^2\rangle + 2\langle\varsigma u\cdot\nabla \eta\rangle
        \\
        &\qquad + b\Ra^{-1}\langle\|\nabla u\|_2^2\rangle +b\Ra^{-1} \langle(2\alpha+\kappa)u_\tau^2\rangle_{\gammam\cup\gammap} + a\langle \|\nabla \omega\|_2^2\rangle
        \\
        &\qquad + 2a\Pr^{-1} \langle (\alpha+\kappa)u\cdot (u\cdot\nabla)u\rangle_{\gammam\cup\gammap} + 2a\langle(\alpha+\kappa)u\cdot\nabla p\rangle_{\gammam\cup\gammap}
        \\
        &\qquad -2a\Ra \langle(\alpha+\kappa) u_\tau n_1 \rangle_{\gammam}- a\Ra \langle \omega \partial_1\theta\rangle
        \label{ubvcxhbiucnuvbihcnvbg}
    \end{align}
    Next we estimate the last four terms on the right-hand side of \eqref{ubvcxhbiucnuvbihcnvbg} individually, but before doing so we remark that by Lemma \ref{lemma_elliptic_regularity_periodic_domain}, Lemma \ref{lemma_nonlin_vorticity_bound} and Lemma \ref{lemma_nonlin_energy_bound}
    \begin{align}
        \|u\|_{W^{1,4}} &\lesssim \|\omega\|_4 + \left(1+\|\kappa\|_\infty^\frac{3}{2}\right) \|u\|_2
        \lesssim \|u_0\|_{W^{1,4}} + \|\alpha^{-1}\|_\infty\Ra,
        \quad \label{ubifbasduifh}
    \end{align}
    where we used that $\|\kappa\|_\infty,\|\alpha+\kappa\|_\infty\leq 1$ and $\|\alpha^{-1}\|_\infty\geq 1$.
    \begin{itemize}
        \item
            By Hölder's inequality and trace theorem
            \begin{align}
                \Pr^{-1} \left|\int_{\gammam\cup\gammap}(\alpha+\kappa) u\cdot(u\cdot\nabla) u \right| \lesssim \Pr^{-1} \|\alpha+\kappa\|_\infty \|u^2 \nabla u\|_{W^{1,1}},\qquad\quad
                \label{ibuafbdsf}
            \end{align}
            which can be bounded using Hölder's inequality and Sobolev embedding by
            \begin{align}
                \|u^2\nabla u\|_{W^{1,1}} &\lesssim \|u^2\nabla u\|_1 + \|u\nabla u\nabla u\|_1 + \|u^2 \nabla^2 u\|_1
                \\
                &\leq \|u\|_4^2 \|u\|_{H^1} + \|u\|_4 \|\nabla u\|_4\|\nabla u\|_2 + \|u\|_4^2 \|u\|_{H^2}
                \\
                &\lesssim \|u\|_{H^1}^2\|u\|_{H^2}
                \label{bzuhvcbcuvzh}
            \end{align}
            Due to \eqref{ibuafbdsf}, \eqref{bzuhvcbcuvzh} and the assumption $\|\alpha+\kappa\|_\infty\leq 1$, Young's inequality yields
            \begin{align}
                2\Pr^{-1} &\left|\int_{\gammam\cup\gammap}(\alpha+\kappa) u\cdot(u\cdot\nabla) u \right|
                \\
                &\qquad \lesssim \Pr^{-1} \|\alpha+\kappa\|_\infty  \|u\|_{H^2}\|u\|_{H^1}^2
                \\
                &\qquad \lesssim \epsilon \|u\|_{H^2}^2 + \epsilon^{-1} \Pr^{-2} \|u\|_{H^1}^4
                \\
                &\qquad \lesssim \epsilon \|u\|_{H^2}^2 + \epsilon^{-1} \Pr^{-2} (\|u_0\|_{W^{1,4}} + \|\alpha^{-1}\|_\infty\Ra)^2\|u\|_{H^1}^2
                \qquad\quad
                \label{ubvchbuicvbh}
            \end{align}
            for any $\epsilon>0$, where we used \eqref{ubifbasduifh} in the last estimate.
            
        \item
            The periodicity on the boundary, Hölder's inequality and trace theorem imply
            \begin{align}
                \left|\int_{\gammam\cup\gammap}(\alpha+\kappa)u\cdot\nabla p\right| &= \left|\int_{\gammam\cup\gammap}p \tau \cdot \nabla ((\alpha+\kappa)u_\tau) \right|
                \\
                &\lesssim \|\alpha+\kappa\|_{L^{\infty}} \|p\nabla u\|_{W^{1,1}} + \|\alpha+\kappa\|_{W^{1,\infty}} \|pu\|_{W^{1,1}}
                \\
                &\lesssim \left(\|\alpha+\kappa\|_\infty\|u\|_{H^2} + \|\alpha+\kappa\|_{W^{1,\infty}} \|u\|_{H^1}\right)\|p\|_{H^1}
                \label{bsdjgdfsgdfsghj}
            \end{align}
            The pressure bound, i.e. Lemma \ref{lemma_nonlin_pressure_bound} and Young's inequality yield
            \begin{align}
                &\left|\int_{\gammam\cup\gammap}(\alpha+\kappa)u\cdot\nabla p\right|
                \\
                &\qquad \lesssim \left(\|\alpha+\kappa\|_\infty\|u\|_{H^2} + \|\alpha+\kappa\|_{W^{1,\infty}} \|u\|_{H^1}\right)\|p\|_{H^1}
                \\
                &\qquad \lesssim \left(\|\alpha+\kappa\|_\infty\|u\|_{H^2} + \|\alpha+\kappa\|_{W^{1,\infty}} \|u\|_{H^1}\right)
                \\
                &\qquad\qquad\cdot \bigg[\left(\frac{1+\|\kappa\|_\infty}{\Pr} \|u\|_{W^{1,4}}+ \|\dot\alpha+\dot\kappa\|_\infty \right)\|u\|_{H^1}
                \\
                &\qquad\qquad\qquad + \|\alpha+\kappa\|_\infty\|u\|_{H^2} + \Ra \bigg]
                \\
                &\qquad \lesssim (\epsilon + \|\alpha+\kappa\|_\infty^2)\|u\|_{H^2}^2 + (1+\epsilon^{-1})\|\alpha+\kappa\|_{W^{1,\infty}}^2\Ra^2
                \\
                &\qquad\qquad + \left(\Pr^{-2}  (\|u_0\|_{W^{1,4}} + \|\alpha^{-1}\|_\infty\Ra)^2 + \|\alpha+\kappa\|_{W^{1,\infty}}^2 + 1\right) \|u\|_{H^1}^2,
                \label{uiasdfuiabsdf}
            \end{align}
            where we also used $\|\kappa\|_\infty\leq 1$ due to \eqref{dabsfkjadfsb_kappa_leq_1} and \eqref{ubifbasduifh}.

        \item
            Hölder's inequality and trace theorem as well as Young's inequality yield
            \begin{align}
                2\Ra \left|\int_{\gammam} (\alpha+\kappa) u_\tau n_1\right| \lesssim \|\alpha+\kappa\|_\infty \Ra \|u\|_{H^1} \lesssim \|\alpha+\kappa\|_\infty^2 \Ra^2 + \|u\|_{H^1}^2.
                \\[-10pt]
                \label{iubfasdbfiuasdf}
            \end{align}

        \item
            Note that by \eqref{nonlin_def_varsigma}, Hölder's and Young's inequality
            \begin{align}
                |a\Ra\langle\omega\partial_1 \theta\rangle|
                &= |a\Ra\langle\omega\partial_1 (\varsigma+\eta)\rangle|
                \\
                &\leq |\langle a\Ra\|\omega\|_2 (\|\nabla \varsigma\|_2 + \|\nabla\eta\|_2)\rangle|
                \\
                &\leq |\langle a^2\Ra^2\|\omega\|_2^2\rangle| + \frac{1}{4} |\langle(\|\nabla \varsigma\|_2 + \|\nabla\eta\|_2)^2\rangle|
                \\
                &\leq a^2\Ra^2\langle\|\omega\|_2^2\rangle + \frac{1}{2}\langle\|\nabla \varsigma\|_2^2\rangle + \frac{1}{2}\langle\|\nabla\eta\|_2^2\rangle.
                \label{cvhxbuzvhuczx}
            \end{align}
    \end{itemize}
    
    Taking the long-time average of \eqref{ubvchbuicvbh}, \eqref{uiasdfuiabsdf}, and \eqref{iubfasdbfiuasdf}, and plugging the results, as well as \eqref{cvhxbuzvhuczx} into \eqref{ubvcxhbiucnuvbihcnvbg}, one finds
    \begin{align}
        \mathcal{Q}
        &\geq M \Ra^2 + \langle \|\nabla\eta\|_2^2\rangle + \langle \|\nabla\varsigma\|_2^2\rangle + 2\langle\varsigma u\cdot\nabla \eta\rangle
        \\
        &\qquad + b\Ra^{-1}\langle\|\nabla u\|_2^2\rangle +b\Ra^{-1} \langle(2\alpha+\kappa)u_\tau^2\rangle_{\gammam\cup\gammap} + a\langle \|\nabla \omega\|_2^2\rangle
        \\
        &\qquad + 2a\Pr^{-1} \langle (\alpha+\kappa)u\cdot (u\cdot\nabla)u\rangle_{\gammam\cup\gammap} + 2a\langle(\alpha+\kappa)u\cdot\nabla p\rangle_{\gammam\cup\gammap}
        \\
        &\qquad -2a\Ra \langle(\alpha+\kappa) u_\tau n_1 \rangle_{\gammam}- a\Ra \langle \omega \partial_1\theta\rangle
        \\
        &\geq M \Ra^2 + \frac{1}{2} \langle \|\nabla\eta\|_2^2\rangle + \frac{1}{2}\langle \|\nabla\varsigma\|_2^2\rangle + 2\langle\varsigma u\cdot\nabla \eta\rangle
        \\
        &\qquad + b\Ra^{-1}\langle\|\nabla u\|_2^2\rangle +b\Ra^{-1} \langle(2\alpha+\kappa)u_\tau^2\rangle_{\gammam\cup\gammap} + a\langle \|\nabla \omega\|_2^2\rangle
        \\
        &\qquad - aC(\epsilon + \|\alpha+\kappa\|_\infty^2)\langle\|u\|_{H^2}^2\rangle - a \tilde C (1+\epsilon^{-1})\|\alpha+\kappa\|_{W^{1,\infty}}^2\Ra^2
        \\
        &\qquad - a C\left((1+\epsilon^{-1})\frac{\|u_0\|_{W^{1,4}}^2 + \|\alpha^{-1}\|_\infty^2\Ra^2}{\Pr^2} + \|\alpha+\kappa\|_{W^{1,\infty}}^2 + 1\right)
        \\
        &\qquad\qquad \cdot \langle\|u\|_{H^1}^2\rangle
        \\
        &\qquad - a^2\Ra^2 \langle\|\omega\|_2^2\rangle
    \end{align}
    Choosing $M=a\tilde C\|\alpha+\kappa\|_{W^{1,\infty}}^2$, using Lemma \ref{lemma_nonlin_varsigmaucdotnalbaeta} and $\frac{1}{2}\langle\|\nabla \eta\|_2^2\rangle\geq 0$
    \begin{align}
        \mathcal{Q}
        &\geq M \Ra^2 + \frac{1}{2} \langle \|\nabla\eta\|_2^2\rangle + \frac{1}{2}\langle \|\nabla\varsigma\|_2^2\rangle + 2\langle\varsigma u\cdot\nabla \eta\rangle
        \\
        &\qquad + b\Ra^{-1}\langle\|\nabla u\|_2^2\rangle +b\Ra^{-1} \langle(2\alpha+\kappa)u_\tau^2\rangle_{\gammam\cup\gammap} + a\langle \|\nabla \omega\|_2^2\rangle
        \\
        &\qquad - aC(\epsilon + \|\alpha+\kappa\|_\infty^2)\langle\|u\|_{H^2}^2\rangle - a\tilde C (1+\epsilon^{-1})\|\alpha+\kappa\|_{W^{1,\infty}}^2\Ra^2
        \\
        &\qquad - a C\left((1+\epsilon^{-1})\frac{\|u_0\|_{W^{1,4}}^2 + \|\alpha^{-1}\|_\infty^2\Ra^2}{\Pr^2} + \|\alpha+\kappa\|_{W^{1,\infty}}^2 + 1\right) 
        \\
        &\qquad\qquad \cdot \langle\|u\|_{H^1}^2\rangle
        \\
        &\qquad - a^2\Ra^2 \langle\|\omega\|_2^2\rangle
        \\
        &\geq b\Ra^{-1}\langle\|\nabla u\|_2^2\rangle +b\Ra^{-1} \langle(2\alpha+\kappa)u_\tau^2\rangle_{\gammam\cup\gammap} - a^2\Ra^2 \langle\|\omega\|_2^2\rangle
        \\
        &\qquad - aC(\epsilon + \|\alpha+\kappa\|_\infty^2)\langle\|u\|_{H^2}^2\rangle + a\langle \|\nabla \omega\|_2^2\rangle - \delta^6 (a\epsilon)^{-1} C\langle\|\partial_2 u\|_2^2\rangle
        \\
        &\qquad - a C\left((1+\epsilon^{-1})\frac{\|u_0\|_{W^{1,4}}^2 + \|\alpha^{-1}\|_\infty^2\Ra^2}{\Pr^2} + \|\alpha+\kappa\|_{W^{1,\infty}}^2 + 1\right)
        \\
        &\qquad\qquad \cdot \langle\|u\|_{H^1}^2\rangle.
        \label{ngnvbinbvuiovbn}
    \end{align}
    Next by Lemma \ref{lemma_elliptic_regularity_periodic_domain} and \eqref{dabsfkjadfsb_kappa_leq_1}
    \begin{align}
        \langle\|u\|_{H^2}^2\rangle
        &\lesssim \langle\|\omega\|_{H^1}^2\rangle +\|\kappa\|_\infty^2 \langle\|\omega\|_2^2\rangle + (1+\|\kappa\|_{W^{1,\infty}}+\|\kappa\|_\infty^2)^2\langle\|u\|_2^2\rangle
        \\
        &\lesssim \langle\|\nabla\omega\|_2^2\rangle + (1+\|\kappa\|_{W^{1,\infty}}^2)\langle\|u\|_{H^1}^2\rangle
    \end{align}
    and by Lemma \ref{lemma_nonlin_coercivity}
    \begin{align}
        \frac{b \|\alpha^{-1}\|_\infty^{-1}}{8\Ra} \langle \|u\|_{H^1}^2\rangle&\leq \frac{3b}{8\Ra}\langle\|\nabla u\|_2^2\rangle + \frac{b}{2\Ra}\langle(2\alpha+\kappa)u_\tau^2\rangle_{\gammam\cup\gammap},
    \end{align}
    which applied to \eqref{ngnvbinbvuiovbn} yields
    \begin{align}
        \mathcal{Q}
        &\geq \frac{5b}{8\Ra}\langle\|\nabla u\|_2^2\rangle +\frac{b}{2\Ra} \langle(2\alpha+\kappa)u_\tau^2\rangle_{\gammam\cup\gammap} - a^2\Ra^2 \langle\|\omega\|_2^2\rangle
        \\
        &\qquad + a\left(1-C_2(\epsilon + \|\alpha+\kappa\|_\infty^2)\right)\langle\|\nabla \omega\|_2^2\rangle - \delta^6 (a\epsilon)^{-1} C\langle\|\partial_2 u\|_2^2\rangle
        \\
        &\qquad+ \bigg[\frac{b\|\alpha^{-1}\|_\infty^{-1}}{8\Ra} - a C\bigg((1+\epsilon^{-1})\frac{\|u_0\|_{W^{1,4}}^2 + \|\alpha^{-1}\|_\infty^2\Ra^2}{\Pr^2}
        \\
        &\qquad\qquad+ (1+\epsilon)(\|\alpha\|_{W^{1,\infty}}^2+\|\kappa\|_{W^{1,\infty}}^4) + 1\bigg)\bigg] \langle\|u\|_{H^1}^2\rangle.        
    \end{align}
    Choosing $\epsilon = \frac{1}{2C_2}$ one has for $\|\alpha+\kappa\|_\infty^2\leq \frac{1}{2C_2}$
    \begin{align}
        \mathcal{Q}
        &\geq \left(\frac{5b}{8\Ra}-C a^{-1}\delta^6\right)\langle\|\nabla u\|_2^2\rangle +\frac{b}{2\Ra} \langle(2\alpha+\kappa)u_\tau^2\rangle_{\gammam\cup\gammap} - a^2\Ra^2 \langle\|\omega\|_2^2\rangle
        \\
        &\qquad+ \bigg[\frac{b\|\alpha^{-1}\|_\infty^{-1}}{8\Ra} - a C\bigg(\frac{\|u_0\|_{W^{1,4}}^2 + \|\alpha^{-1}\|_\infty^2\Ra^2}{\Pr^2}
        \\
        &\qquad\qquad+ \|\alpha\|_{W^{1,\infty}}^2+\|\kappa\|_{W^{1,\infty}}^4 + 1\bigg)\bigg] \langle\|u\|_{H^1}^2\rangle.     
        \label{cvxiuobcvzbviu}
    \end{align}
    Next, we distinguish between the conditions on $\kappa$. Note though that all the estimates follow a similar approach.
    \begin{itemize}
        \item Let $|\kappa|\leq \alpha$.
        
            Then by Lemma \eqref{lemma_grad_ids}
            \begin{align}
                \|\omega\|_2^2 &= \|\nabla u\|_2^2 - \int_{\gammam\cup\gammap} \kappa u_\tau^2 \leq \|\nabla u\|_2^2 + \int_{\gammam\cup\gammap} \alpha u_\tau^2
                \\
                &\leq \|\nabla u\|_2^2 + \int_{\gammam\cup\gammap} (2\alpha+\kappa) u_\tau^2,
            \end{align}
            which after taking the long-time average implies for \eqref{cvxiuobcvzbviu}
            \begin{align}
                \mathcal{Q}
                &\geq \left(\frac{b}{8\Ra}-C a^{-1}\delta^6\right)\langle\|\nabla u\|_2^2\rangle + \left( \frac{b}{2\Ra}- a^2\Ra^2 \right) \langle\|\omega\|_2^2\rangle
                \\
                &\qquad+ \bigg[\frac{b\|\alpha^{-1}\|_\infty^{-1}}{8\Ra} - a C\bigg(\frac{\|u_0\|_{W^{1,4}}^2 + \|\alpha^{-1}\|_\infty^2\Ra^2}{\Pr^2}
                \\
                &\qquad\qquad+ \|\alpha\|_{W^{1,\infty}}^2+\|\kappa\|_{W^{1,\infty}}^4 + 1\bigg)\bigg] \langle\|u\|_{H^1}^2\rangle.
                \label{bniadbsahfsd}
            \end{align}
            As $b$ has to fulfill $b<(1+\max h-\min h)^{-1}$ in order for the second round bracket on the right-hand side of \eqref{bniadbsahfsd} to be non-negative we set $a=a_0\Ra^{-\frac{3}{2}}$, which yields
            \begin{align}
                \mathcal{Q}
                &\geq \left(\frac{b}{8\Ra}-C a_0^{-1}\Ra^{\frac{3}{2}}\delta^6\right)\langle\|\nabla u\|_2^2\rangle + \left( \frac{b}{2}- a_0^2 \right) \Ra^{-1} \langle\|\omega\|_2^2\rangle
                \\
                &\qquad+ \bigg[\frac{b\|\alpha^{-1}\|_\infty^{-1}}{8} - a_0 \Ra^{-\frac{1}{2}} C\bigg(\frac{\|u_0\|_{W^{1,4}}^2 + \|\alpha^{-1}\|_\infty^2\Ra^2}{\Pr^2}
                \\
                &\qquad\qquad+ \|\alpha\|_{W^{1,\infty}}^2+\|\kappa\|_{W^{1,\infty}}^4 + 1\bigg)\bigg]\Ra^{-1} \langle\|u\|_{H^1}^2\rangle
            \end{align}
            and the assumption $\Pr\geq \|\alpha^{-1}\|_\infty^\frac{3}{2}\Ra^\frac{3}{4}$ yields
            \begin{align}
                \mathcal{Q}
                &\geq \left(\frac{b}{8\Ra}-C a_0^{-1}\Ra^{\frac{3}{2}}\delta^6\right)\langle\|\nabla u\|_2^2\rangle + \left( \frac{b}{2}- a_0^2 \right) \Ra^{-1} \langle\|\omega\|_2^2\rangle
                \\
                &\qquad+ \bigg[\frac{b}{8} - a_0 C\bigg(\|u_0\|_{W^{1,4}}^2\|\alpha^{-1}\|_\infty^{-2}\Ra^{-2} + 1
                \\
                &\qquad\qquad + \frac{\|\alpha\|_{W^{1,\infty}}^2+\|\kappa\|_{W^{1,\infty}}^4 + 1}{\Ra^{\frac{1}{2}}}\|\alpha^{-1}\|_\infty\bigg)\bigg]\|\alpha^{-1}\|_\infty^{-1}\Ra^{-1} \langle\|u\|_{H^1}^2\rangle
            \end{align}
            and since $\Ra^{-\frac{1}{2}}< \|\alpha^{-1}\|_\infty^{-1}<1$
            \begin{align}
                \mathcal{Q}
                &\geq \left(\frac{b}{8\Ra}-C a_0^{-1}\Ra^{\frac{3}{2}}\delta^6\right)\langle\|\nabla u\|_2^2\rangle + \left( \frac{b}{2}- a_0^2 \right) \Ra^{-1} \langle\|\omega\|_2^2\rangle
                \\
                &\qquad+ \bigg[\frac{b}{8} - a_0 C\big(\|u_0\|_{W^{1,4}}^2 + \|\alpha\|_{W^{1,\infty}}^2+\|\kappa\|_{W^{1,\infty}}^4 + 1\big)\bigg]
                \\
                &\qquad\qquad \cdot \|\alpha^{-1}\|_\infty^{-1}\Ra^{-1} \langle\|u\|_{H^1}^2\rangle
            \end{align}
            Without loss of generality, we can assume that $C\geq 1$ such that setting
            \begin{align}
                a_0 = \frac{b}{8C\left(\|u_0\|_{W^{1,4}}^2 + \|\alpha\|_{W^{1,\infty}}^2+\|\kappa\|_{W^{1,\infty}}^4 + 1\right)}
            \end{align}
            yields
            \begin{align}
                \mathcal{Q}
                &\geq \left(\frac{b}{8\Ra}-C a_0^{-1}\Ra^{\frac{3}{2}}\delta^6\right)\langle\|\nabla u\|_2^2\rangle.
            \end{align}
            Letting $\delta$ solve $\frac{b}{8\Ra}=C a_0^{-1}\Ra^{\frac{3}{2}}\delta^6$, i.e.
            \begin{align}
                \delta 
                = \left( \frac{a_0 b}{8C} \right)^\frac{1}{6}\Ra^{-\frac{5}{12}},
            \end{align}
            it holds $\mathcal{Q}\geq 0$ and therefore, setting $b=\frac{1}{2(1+\max h-\min h)}$,\eqref{ahusdfbafbhvcxubihxcvbiu} results in
            \begin{align}
                \Nu &\leq 2(1 - b(1+\max h-\min h))\Nu
                \\
                &\leq 2M\Ra^2 + C \delta^{-1} -2\mathcal{Q}
                \\
                &\lesssim a \|\alpha+\kappa\|_{W^{1,\infty}}^2\Ra^2 + a_0^{-\frac{1}{6}}b^{-\frac{1}{6}} \Ra^\frac{5}{12}
                \\
                &\lesssim a_0 \|\alpha+\kappa\|_{W^{1,\infty}}^2\Ra^\frac{1}{2} + a_0^{-\frac{1}{6}}b^{-\frac{1}{6}} \Ra^\frac{5}{12}
                \\
                &\lesssim \|\alpha+\kappa\|_{W^{1,\infty}}^2\Ra^\frac{1}{2} + \left(\|u_0\|_{W^{1,4}}^\frac{1}{3} + \|\alpha\|_{W^{1,\infty}}^\frac{1}{3}+\|\kappa\|_{W^{1,\infty}}^\frac{2}{3} + 1\right) \Ra^\frac{5}{12}
            \end{align}

        \item Let $|\kappa|\leq 2 \alpha+ \frac{1}{4\sqrt{1+(h')^2}}\sqrt{\alpha}$.

        Taking the long-time average of Lemma \ref{lemma_elliptic_regularity_periodic_domain}, using $\|\kappa\|_\infty\leq 1$,
        \begin{align}
            \langle \|\omega\|_2^2\rangle \leq \langle \|u\|_{H^1}^2\rangle,
        \end{align}
        which for \eqref{cvxiuobcvzbviu} yields
        \begin{align}
            \mathcal{Q}
            &\geq \left(\frac{5b}{8\Ra}-C a^{-1}\delta^6\right)\langle\|\nabla u\|_2^2\rangle
            \\
            &\qquad+ \bigg[\frac{b\|\alpha^{-1}\|_\infty^{-1}}{8\Ra} -a^2 \Ra^2 - a C\bigg(\frac{\|u_0\|_{W^{1,4}}^2 + \|\alpha^{-1}\|_\infty^2\Ra^2}{\Pr^2}
            \\
            &\qquad\qquad+ \|\alpha\|_{W^{1,\infty}}^2+\|\kappa\|_{W^{1,\infty}}^4 + 1\bigg)\bigg] \langle\|u\|_{H^1}^2\rangle
        \end{align}
        and by Lemma \ref{lemma_nonlin_coercivity}
        \begin{align}
            \mathcal{Q}
            &\geq \left(\frac{5b}{8\Ra}-C a^{-1}\delta^6\right)\langle\|\nabla u\|_2^2\rangle 
            \\
            &\qquad+ \bigg[\frac{b\|\alpha^{-1}\|_\infty^{-1}}{8\Ra} -a^2 \Ra^2 - a C\bigg(\frac{\|u_0\|_{W^{1,4}}^2 + \|\alpha^{-1}\|_\infty^2\Ra^2}{\Pr^2}
            \\
            &\qquad\qquad+ \|\alpha\|_{W^{1,\infty}}^2+\|\kappa\|_{W^{1,\infty}}^4 + 1\bigg)\bigg] \langle\|u\|_{H^1}^2\rangle
            \label{tgbvuznznvutnz}
        \end{align}
        \begin{itemize}
            \item 
                As before setting $a=a_0\Ra^{-\frac{3}{2}}$, the assumptions $\Ra^{-1}\leq \|\alpha^{-1}\|_\infty^{-1}$, $\Pr\geq \|\alpha^{-1}\|_\infty^\frac{5}{4}\Ra^\frac{3}{4}\geq 1$ yield
                \begin{align}
                    \mathcal{Q}
                    &\geq \left(\frac{5b}{8\Ra}-C a_0^{-1}\Ra^\frac{3}{2}\delta^6\right)\langle\|\nabla u\|_2^2\rangle 
                    \\
                    &\qquad+ \frac{1}{\Ra} \bigg[\frac{b}{8}\|\alpha^{-1}\|_\infty^{-1} - a_0^2
                    \\
                    &\qquad\qquad - a_0 \|\alpha^{-1}\|_\infty^{-\frac{1}{2}}C\left(\|u_0\|_{W^{1,4}}^2 + \|\alpha\|_{W^{1,\infty}}^2+\|\kappa\|_{W^{1,\infty}}^4 + 1\right)\bigg] 
                    \\
                    &\qquad\qquad \cdot \langle\|u\|_{H^1}^2\rangle
                \end{align}
                and, assuming without loss of generality $C\geq 1$, for
                \begin{align}
                    b &= \frac{1}{2(1+\max h-\min h)}
                    \\
                    a_0 &= \frac{b\|\alpha^{-1}\|_\infty^{-\frac{1}{2}}}{16C(\|u_0\|_{W^{1,4}}^2 + \|\alpha\|_{W^{1,\infty}}^2+\|\kappa\|_{W^{1,\infty}}^4 + 1)}
                    \\
                    \delta &= \left(\frac{5a_0b}{8C}\right)^\frac{1}{6}\Ra^{-\frac{5}{12}}
                \end{align}
                it holds $\mathcal{Q}\geq 0$ and therefore \eqref{ahusdfbafbhvcxubihxcvbiu} yields
                \begin{align}
                    \Nu
                    &\leq 2M\Ra^2 + C \delta^{-1} -2\mathcal{Q}
                    \\
                    &\lesssim a \|\alpha+\kappa\|_{W^{1,\infty}}^2\Ra^2 + a_0^{-\frac{1}{6}}b^{-\frac{1}{6}} \Ra^\frac{5}{12}
                    \\
                    &\lesssim a_0 \|\alpha+\kappa\|_{W^{1,\infty}}^2\Ra^\frac{1}{2} + a_0^{-\frac{1}{6}}b^{-\frac{1}{6}} \Ra^\frac{5}{12}
                    \\
                    &\lesssim \|\alpha^{-1}\|_\infty^{-\frac{1}{2}} \|\alpha+\kappa\|_{W^{1,\infty}}^2\Ra^\frac{1}{2}
                    \\
                    &\qquad+ \left(\|u_0\|_{W^{1,4}}^\frac{1}{3} + \|\alpha\|_{W^{1,\infty}}^\frac{1}{3}+\|\kappa\|_{W^{1,\infty}}^\frac{2}{3} + 1\right) \|\alpha^{-1}\|_\infty^{\frac{1}{12}}\Ra^\frac{5}{12}.
                \end{align}
            \item
                Setting instead $a=a_0\Ra^{-\frac{11}{7}}$ and using the assumption $\Pr\geq \Ra^\frac{5}{7}\geq 1$ in \eqref{tgbvuznznvutnz} yields
                \begin{align}
                    \mathcal{Q}
                    &\geq \left(\frac{5b}{8\Ra}-C a_0^{-1}\Ra^{\frac{11}{7}}\delta^6\right)\langle\|\nabla u\|_2^2\rangle 
                    \\
                    &\qquad+ \Ra^{-1} \bigg[\frac{b\|\alpha^{-1}\|_\infty^{-1}}{8} -a_0^2 \Ra^{-\frac{1}{7}} - a_0 \|\alpha^{-1}\|_\infty^2
                    \\
                    &\qquad\qquad - a_0 \Ra^{-\frac{4}{7}} C\left(\|u_0\|_{W^{1,4}}^2+ \|\alpha\|_{W^{1,\infty}}^2+\|\kappa\|_{W^{1,\infty}}^4 + 1\right)\bigg]
                    \\
                    &\qquad\qquad\qquad\cdot\langle\|u\|_{H^1}^2\rangle
                \end{align}
                and choosing
                \begin{align}
                    a_0 =\frac{\|\alpha^{-1}\|_\infty^{-1}b}{24C\left(\|u_0\|_{W^{1,4}}^2+\|\alpha^{-1}\|_\infty^2 + \|\alpha\|_{W^{1,\infty}}^2+\|\kappa\|_{W^{1,\infty}}^4 + 1\right)}
                \end{align}
                one gets
                \begin{align}
                    \mathcal{Q}
                    &\geq \left(\frac{5b}{8\Ra}-C a_0^{-1}\Ra^{\frac{11}{7}}\delta^6\right)\langle\|\nabla u\|_2^2\rangle,
                \end{align}
                which after setting $\delta = \left(\frac{5a_0b}{8C}\right)^\frac{1}{6}\Ra^{-\frac{3}{7}}$ implies $\mathcal{Q}\geq 0$. Therefore using \eqref{ahusdfbafbhvcxubihxcvbiu} results in
                \begin{align}
                    \Nu
                    &\leq 2M\Ra^2 + C \delta^{-1} -2\mathcal{Q}
                    \\
                    &\lesssim a\|\alpha+\kappa\|_{W^{1,\infty}}^2 \Ra^2 + C \delta^{-1}
                    \\
                    &\lesssim \left(a_0\|\alpha+\kappa\|_{W^{1,\infty}}^2 + a_0^{-\frac{1}{6}}\right)\Ra^{\frac{3}{7}}
                    \\
                    &\lesssim \Big( \|\alpha^{-1}\|_\infty^{-3}\|\alpha+\kappa\|_{W^{1,\infty}}^2 + \|\alpha^{-1}\|_\infty^\frac{1}{2}
                    \\
                    &\qquad + \|\alpha^{-1}\|_\infty^\frac{1}{6}\left(\|u_0\|_{W^{1,4}}^\frac{1}{3} + \|\alpha\|_{W^{1,\infty}}^\frac{1}{3}+\|\kappa\|_{W^{1,\infty}}^\frac{2}{3} + 1\right) \Big)\Ra^{\frac{3}{7}}
                \end{align}
        \end{itemize}
            
    \end{itemize}
\end{proof}

\newpage
\chapter{Thermally Non-Diffusive System}
\label{chapter:thermally_nd_system}

\section{The Model}

\begin{figure}
    \hrule
    \vspace{-\baselineskip}
    \begin{center}
        \includegraphics[width=\textwidth]{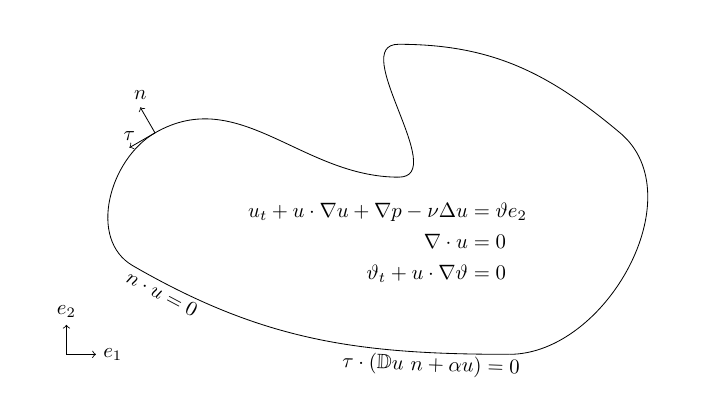}
        \vspace*{-10mm}
        \caption{Illustration of the thermally non-diffusive system.}
        \label{fig:nd_overview}
    \end{center}
    \hrule
\end{figure}

In the following, we study the system
\begin{alignat}{2}
    u_t + u\cdot \nabla u +\nabla p - \Delta u &=\theta e_2& \qquad\textnormal{ in }&\Omega
    \label{nd_navier_stokes}
    \\
    \theta_t + u\cdot \nabla \theta &=0 & \qquad\textnormal{ in }&\Omega
    \label{nd_transport_equation}
    \\
    \nabla \cdot u &= 0 & \qquad\textnormal{ in }&\Omega
    \label{nd_incompressible}
    \\
    n\cdot u &= 0& \qquad\textnormal{ on }&\partial\Omega
    \label{nd_no_penetration_bc}
    \\
    \tau \cdot (\D u \ n + \alpha u ) &= 0& \qquad\textnormal{ on }&\partial\Omega
    \label{nd_navier_slip}
    \\
    (u,\theta)(\cdot,0) &= (u_0,\theta_0) & \qquad\textnormal{ in }&\Omega,
    \label{nd_initial_data}
\end{alignat}
as motivated in Section \ref{Section:intro_nd_system}. The system is illustrated in Figure \ref{fig:nd_overview}.\noeqref{nd_initial_data}

The subsequent results and analysis are contained in \cite{bleitnerCarlsonNobili2023Large}. Although some proofs follow a similar strategy as in Chapter \ref{chapter:rayleigh-benard-convection}, we prove them here for the convenience of the reader and the sake of comprehensiveness.

\section{Main Results}
\label{section:nd_main_results}
The main focus in the study of this system is devoted to regularity results as well as the long-time asymptotics of solutions. The findings for the former are summarized in the following system.
\begin{theorem}[Regularity]
\label{theorem_nd_regularity}
    Let $\Omega$ be a bounded $C^{2,1}$ domain, $u_0\in H^2(\Omega)$, $\theta_0\in L^4(\Omega)$ and $0<\alpha\in W^{1,\infty}(\partial\Omega)$. Then
    \begin{align}
        u&\in L^\infty\left((0,\infty);H^2(\Omega)\right)\cap L^q\left((0,\infty);W^{1,q}(\Omega)\right)
        \\
        u_t&\in L^\infty\left((0,\infty);L^2(\Omega)\right)\cap L^2\left((0,\infty);H^1(\Omega)\right)
        \\
        p&\in L^\infty\left((0,\infty);H^1(\Omega)\right)
        \\
        \theta&\in L^\infty\left((0,\infty);L^4(\Omega)\right)
    \end{align}
    for any $2\leq q <\infty$.
    
    If additionally $\Omega$ is a simply connected $C^{3,1}$ domain, $\theta_0 \in W^{1,\tilde q}(\Omega)$ for some $\tilde q \in 2\mathbb{N}$ and $0<\alpha \in W^{2,\infty}(\partial\Omega)$. Then
    \begin{align}
        u&\in L^2\left((0,T);H^3(\Omega)\right) \cap L^\frac{2r}{r-2}\left((0,T);W^{2,r}(\Omega)\right)
        \\
        \theta&\in L^\infty\left((0,T);W^{1,q}(\Omega)\right)
    \end{align}
    for any $T>0$, $2\leq r <\infty$ and $1\leq q \leq \tilde q$.
\end{theorem}
These results are proven consecutively in Section \ref{section:nd_regularity_estimates}.

Theorem \ref{theorem_nd_regularity} improves the regularity results of \cite{hu2018Partially} given by
\begin{align}
    u \in L^\infty\left((0,T);H^1(\Omega)\right)\cap L^2\left((0,T);H^2(\Omega)\right)
\end{align}
for constant $\alpha>0$ in multiple ways. The approach used here allows spatially varying slip coefficients $\alpha$. The results given in Theorem \ref{theorem_nd_regularity} show higher regularity and the bounds hold uniform in time up to moderate order. In particular, these uniform in time estimates allow us to infer the long-time behavior of solutions as given in the following Theorem.

\begin{theorem}[Convergence]
    \label{theorem_nd_convergence}
    
    Let $\Omega$ be $C^{2,1}$, $u_0\in H^2(\Omega)$, $\theta_0 \in L^{4}(\Omega)$ and $0<\alpha \in W^{1,\infty}(\partial\Omega)$. Then for any $1\leq q<\infty$
    \begin{alignat}{2}
        \|u(t)&\|_{W^{1,q}}&&\to 0
        \label{nd_u_to_0_in_long_time}
        \\
        \|u_t(t)&\|_2&&\to 0
        \label{nd_ut_to_0_in_long_time}
        \\
        \|(\nabla p-\theta e_2)(t)&\|_{H^{-1}}&&\to 0
        \label{nd_nablap-thetae2_to_0_in_long_time}
    \end{alignat}
    for $t\to\infty$.
\end{theorem}

Theorem \ref{theorem_nd_convergence} is proven in Section \ref{section:nd_convergence}.

The long-time behavior described in the Theorem shows convergence to the hydrostatic equilibrium, where the buoyancy is balanced by the pressure gradient and the velocity field vanishes. This behavior is not surprising, as over time the diffusion operator is expected to dissipate fluctuations in the velocity field, while the Navier-slip boundary conditions create drag along the boundary, slowing the velocity field down. Due to the conservation of the temperature field, the pressure gradient has to balance the buoyancy forces.

\begin{figure}
    \hrule
    \begin{center}
        \includegraphics[width=0.5\textwidth]{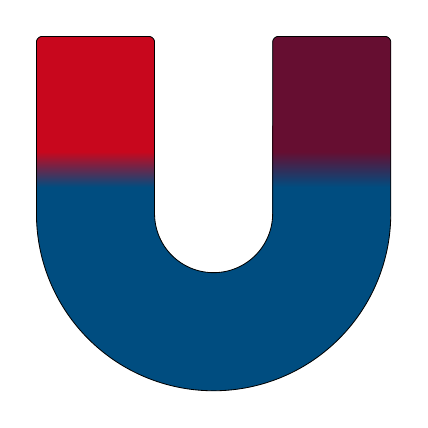}
        \vspace*{-5mm}
        \caption{Illustration of a hydrostatic equilibrium that is not uniformly stratified.}
        \label{fig:nd_x1_dependent_he}
    \end{center}
    \hrule
\end{figure}

Note that no specific profile for the hydrostatic equilibrium is given, and it could potentially vary in the horizontal direction. Such an $x_1$-dependent hydrostatic equilibrium is sketched in Figure \ref{fig:nd_x1_dependent_he}.

Finally, we want to remark that linear profiles for temperature field play an important role in this system. In fact, subtracting the hydrostatic equilibrium given by
\begin{align}
    u_\textnormal{he}&=0,& \theta_\textnormal{he} &= \beta x_2 + \gamma, & p_\textnormal{he}=\frac{1}{2}\beta x_2^2 + \gamma x_2 + \delta,
\end{align}
where $\beta>0$ and $\gamma,\delta\in \mathbb{R}$, from the solution, i.e. defining
\begin{align}
    \hat u &= u, &\hat \theta &= \theta-\beta x_2 -\gamma,& \hat p&=p-\frac{1}{2}\beta x_2^2 - \gamma x_2 - \delta,
\end{align}
the new variables satisfy
\begin{alignat}{2}
    \hat u_t + \hat u\cdot \nabla \hat u +\nabla \hat p - \Delta \hat u &=\hat \theta e_2& \qquad\textnormal{ in }&\Omega
    \label{bgfbfgi}
    \\
    \hat\theta_t + \hat u\cdot \nabla \hat\theta &=-\beta \hat u_2 & \qquad\textnormal{ in }&\Omega
    \label{njhgthjrberbher}
    \\
    \nabla \cdot \hat u &= 0 & \qquad\textnormal{ in }&\Omega
    \label{jggfdghfdgh}
    \\
    n\cdot \hat u &= 0& \qquad\textnormal{ on }&\partial\Omega
    \label{hnjrtkhnjrth}
    \\
    \tau \cdot (\D \hat u \ n + \alpha \hat u ) &= 0& \qquad\textnormal{ on }&\partial\Omega
    \label{bjnrkbjrgbnjf}
    \\
    (\hat u,\hat \theta)(\cdot,0) &= (u_0,\theta_0-\beta x_2 -\gamma) & \quad\textnormal{ in }&\Omega,
    \label{gnjrtgrtnjrnb}
\end{alignat}\noeqref{jggfdghfdgh}\noeqref{hnjrtkhnjrth}\noeqref{bjnrkbjrgbnjf}\noeqref{gnjrtgrtnjrnb}
When testing \eqref{bgfbfgi} with $\beta \hat u$ and \eqref{njhgthjrberbher} with $\hat \theta$ and adding both results together the right-hand sides cancel. This observation will be used in parts of the proof of Theorem \ref{theorem_nd_regularity}. Additionally, \cite{doering2018longTime} show that for $\beta<0$ such a solution in a spatially periodic setting is unstable. For $\beta>0$ and the, there imposed, free slip boundary conditions such a linear profile is linearly stable in $H^2\times L^2$ under additional assumptions on the domain or boundary conditions on the temperature field.

For reference, Theorem \ref{theorem_nd_convergence} shows the nonlinear stability of such a hydrostatic equilibrium in the class $W^{1,q}\times H^{-1}$.

\section{Regularity Estimates}
\label{section:nd_regularity_estimates}

Due to the absence of thermal diffusion, no gain in regularity for the scalar field can be expected. Note though that all $L^q$-norms of $\theta$ are conserved if $u$ is sufficiently smooth, i.e. for any $1 \leq q \leq \infty$
\begin{align}
    \|\theta(t)\|_q = \|\theta_0\|_q
    \label{nd_maximum_principle}
\end{align}
for all $t\geq 0$ if $\theta_0\in L^q(\Omega)$. In fact for $q \in 2\mathbb{N}$, the cases that will be used in the following analysis, testing \eqref{nd_transport_equation} with $\theta^{q-1}$ shows
\begin{align}
    \frac{d}{dt}\|\theta\|_q^q = -q\int_{\Omega} \theta^{q-1} u\cdot \nabla \theta = - \int_{\Omega} u\cdot \nabla \theta^q = - \int_{\partial\Omega} u\cdot n \theta^q + \int_{\Omega} \theta^q \nabla \cdot u = 0,
\end{align}
which implies \eqref{nd_maximum_principle} for these choices of $q$.

Having that the $L^2$-norm of $\theta$ is conserved, one can directly prove bounds for the energy of the fluid as given in the following Lemma.

\begin{lemma}[Energy Bound]
    \label{lemma_nd_energy_bound}
    Let $\Omega$ be $C^{1,1}$, $u_0,\theta_0\in L^2(\Omega)$ and $0<\alpha\in L^\infty(\partial\Omega)$. Then
    \begin{align}
        u\in L^\infty\left((0,\infty); L^2(\Omega)\right)
    \end{align}
    and it holds
    \begin{align}
        \|u\|_2 \lessapprox \|u_0\|_2 + \|\theta_0\|_2,
    \end{align}
    where the implicit constant depends on $\alpha$ and $\Omega$.
\end{lemma}

\begin{proof}
    Testing \eqref{nd_navier_stokes} with $u$ one finds
    \begin{align}
        \frac{1}{2}\frac{d}{dt} \|u\|_2^2 &= - \int_{\Omega} u \cdot (u\cdot \nabla) u - \int_{\Omega} u \cdot \nabla p + \int_{\Omega} u\cdot\Delta u + \int_{\Omega} \theta u_2
        \label{bzhubbuvgzcvb}
    \end{align}
    The first term on the right-hand side of \eqref{bzhubbuvgzcvb} vanishes due to \eqref{nd_incompressible} and \eqref{nd_no_penetration_bc} as
    \begin{align}
        2\int_{\Omega} u\cdot (u\cdot\nabla) u = \int_{\Omega} u\cdot \nabla |u|^2 = - \int_{\Omega} |u|^2 \nabla \cdot u + \int_{\partial\Omega} |u|^2 n\cdot u = 0\qquad
        \label{bcibvqwerqweczb}
    \end{align}
    and similarly for the second term on the right-hand side of \eqref{bzhubbuvgzcvb}
    \begin{align}
        \int_{\Omega} u\cdot \nabla p = \int_{\partial\Omega} u\cdot n p - \int_{\Omega} p \nabla \cdot u = 0.
        \label{tzuvbnzuvbtzn}
    \end{align}
    By Lemma \ref{lemma_int_Delta_u_v} and Lemma \ref{lemma_coercivity} we find
    \begin{align}
        -\int_{\Omega} u \cdot \Delta u = 2\|\D u\|_2^2 + 2 \int_{\partial\Omega} \alpha u_\tau^2 \gtrsim (1+\|\alpha^{-1}\|_\infty)^{-1} \|u\|_2^2.
        \label{wemwerwerqerwqrw}
    \end{align}
    Combining \eqref{bzhubbuvgzcvb}, \eqref{bcibvqwerqweczb}, \eqref{tzuvbnzuvbtzn}, \eqref{wemwerwerqerwqrw} and using Hölder's and Young's inequality there exists a constant $C>0$ such that
    \begin{align}
        \frac{1}{2}\frac{d}{dt} \|u\|_2^2 &+ C(1+\|\alpha^{-1}\|_\infty)^{-1}\|u\|_2^2
        \\
        &\leq \|\theta\|_2 \|u\|_2
        \\
        &\leq \frac{C}{2} (1+\|\alpha^{-1}\|_\infty)^{-1} \|u\|_2^2 + (2C)^{-1} (1+\|\alpha^{-1}\|_\infty) \|\theta\|_2^2
    \end{align}
    and therefore
    \begin{align}
        \frac{d}{dt} \|u\|_2^2 + C(1+\|\alpha^{-1}\|_\infty)^{-1}\|u\|_2^2
        &\lesssim (1+\|\alpha^{-1}\|_\infty) \|\theta_0\|_2^2,
        \label{rtenjtkrenj}
    \end{align}
    where we used \eqref{nd_maximum_principle},
    which after applying Grönwall's inequality yields
    \begin{align}
        \|u\|_2^2 \leq e^{-C(1+\|\alpha^{-1}\|_\infty)^{-1} t} \|u_0\|_2^2 + (1+\|\alpha^{-1}\|_\infty)^2\|\theta_0\|_2^2
    \end{align}
\end{proof}

\begin{lemma}
    \label{lemma_nd_L2H1_bound}
    Let $\Omega$ be $C^{1,1}$, $u_0,\theta_0\in L^2(\Omega)$ and $0<\alpha\in L^\infty(\partial\Omega)$. Then
    \begin{align}
        u\in L^2\left((0,\infty);H^1(\Omega)\right)
    \end{align}
    and
    \begin{align}
        \int_0^\infty \|u(s)\|_{H^1}^2\ ds \lessapprox\|u_0\|_2^2 + \|\theta_0\|_2^2 + 1,
    \end{align}
    where the implicit constant depends on $\alpha$ and $\Omega$.
\end{lemma}

\begin{proof}
    Defining
    \begin{align}
        \hat u &= u
        \\
        \hat \theta &= \theta - \beta x_2 - \gamma
        \\
        \hat p &= p - \frac{\beta}{2} x_2^2 - \gamma x_2 - \delta
    \end{align}
    for any $\beta,\gamma,\delta\in \mathbb{R}$ with $\beta > 0$, $(\hat u,\hat p,\hat\theta)$ solves \eqref{bgfbfgi}--\eqref{gnjrtgrtnjrnb}.
    Testing \eqref{bgfbfgi} with $\hat u$ and \eqref{njhgthjrberbher} with $\hat \theta$ and adding them together
    \begin{align}
        \frac{1}{2}\frac{d}{dt} \left(\beta\|u\|_2^2 + \|\hat\theta\|_2^2 \right) &= -\beta \int_{\Omega} u\cdot(u\cdot\nabla)u - \beta \int_{\Omega} u\cdot \nabla \hat p + \beta \int_{\Omega} u\cdot\Delta u
        \\
        &\qquad + \beta  \int_{\Omega} \hat \theta u_2 - \int_{\Omega} \hat\theta u\cdot\nabla\hat\theta - \int_{\Omega} \hat \theta \beta u_2
        \label{tzccvxzubtcb}
    \end{align}
    Note that the $\beta\hat\theta \hat u_2$ terms, originating from the right-hand sides of \eqref{bgfbfgi} and \eqref{njhgthjrberbher}, cancel. Similar to \eqref{bcibvqwerqweczb} and \eqref{tzuvbnzuvbtzn} the first and second term on the right-hand side of \eqref{tzccvxzubtcb} vanish, as does the fifth term since
    \begin{align}
        2\int_{\Omega} \hat\theta \hat u\cdot\nabla\hat\theta = \int_{\Omega} \hat u\cdot\nabla (\hat\theta^2) = \int_{\partial\Omega} \hat u\cdot n\hat\theta^2 - \int_{\Omega} \hat\theta^2 \nabla \cdot \hat u = 0,
    \end{align}
    implying
    \begin{align}
        \frac{1}{2}\frac{d}{dt} \left(\beta\|\hat u\|_2^2 + \|\hat\theta\|_2^2 \right) &= \beta \int_{\Omega} \hat u\cdot\Delta \hat u = - 2 \beta \|\D \hat u\|_2^2 - 2 \beta \int_{\partial\Omega} \alpha \hat u_\tau^2,
        \label{bczuvbzcu}
    \end{align}
    where in the last identity we used Lemma \ref{lemma_int_Delta_u_v}. It follows that
    \begin{align}
        \frac{1}{2}\frac{d}{dt} \big(\beta\|u\|_2^2 + \|&\theta-\beta x_2 -\gamma\|_2^2 \big) + 2 \beta \|\D u\|_2^2 + 2 \beta \int_{\partial\Omega} \alpha u_\tau^2
        \\
        & = \frac{1}{2}\frac{d}{dt} \left(\beta\|\hat u\|_2^2 + \|\hat \theta\|_2^2 \right) + 2 \beta \|\D \hat u\|_2^2 + 2 \beta \int_{\partial\Omega} \alpha \hat u_\tau^2
        \\
        & = 0
        \label{vbcvbcggbvhcj}
    \end{align}
    and integrating \eqref{vbcvbcggbvhcj} in time
    \begin{align}
        \beta \|u(t)\|_2^2 + \|\theta(t)-\beta x_2 -\gamma\|_2^2 &+ 4 \beta \int_0^t \left(\|\D u(s)\|_2^2 + \int_{\partial\Omega} \alpha u_\tau^2(s)\right) ds
        \\
        &\qquad\qquad\qquad = \beta \|u_0\|_2^2 + \|\theta_0-\beta x_2-\gamma\|_2^2.
    \end{align}
    As by Lemma \ref{lemma_coercivity} there exists a constant $C>0$ depending only on $\alpha$ and $\Omega$ such that
    \begin{align}
        C\|u\|_{H^1}^2 \leq \|\D u\|_2^2 + \int_{\partial\Omega} \alpha u_\tau^2
    \end{align}
    one finds
    \begin{align}
        C\beta\int_0^t \|u(s)\|_{H^1}^2 \ ds &\leq \beta\|u_0\|_2^2 + \|\theta_0-\beta x_2 -\gamma\|_2^2
        \\
        &\lesssim \beta\|u_0\|_2^2 + \|\theta_0\|_2^2 + \beta + |\gamma|
        \label{klvbxchjvc}
    \end{align}
    and since the right-hand side of \eqref{klvbxchjvc} is independent of time, this bound holds uniform in time, concluding the regularity result and choosing $\beta=1$ and $\gamma=0$ yields the bound.
\end{proof}

Next we prove $L^\infty\left((0,\infty); W^{1,q}(\Omega)\right)$ bounds for $u$. The proof follows the same strategy as the one of Lemma \ref{lemma_vorticity_bound}.

Analogous to before we first derive a vorticity formulation. The vorticity, defined by $\omega = \nabla^\perp \cdot u = -\partial_2 u_1 + \partial_1 u_2$ fulfills
\begin{alignat}{2}
    \omega_t + u\cdot\nabla \omega - \Delta \omega &= \partial_1 \theta & \qquad \textnormal{ in }&\Omega
    \label{nd_vorticity_pde}
    \\
    \omega &= -2(\alpha+\kappa)u_\tau & \qquad \textnormal{ on }&\partial\Omega.
    \label{nd_vorticity_bc}
\end{alignat}

\eqref{nd_vorticity_pde} follows immediately when applying $\nabla^\perp\cdot$ to \eqref{nd_navier_stokes} and using
\begin{align}
    \nabla^\perp\cdot \left((u\cdot\nabla)u\right) = u\cdot \nabla \omega + \omega \nabla \cdot u
\end{align}
and \eqref{nd_incompressible}. In order to prove the boundary condition \eqref{nd_vorticity_bc}, note that by \eqref{nd_navier_slip} and \eqref{n_tau_grad_u}
\begin{align}
    -2(\alpha+\kappa)u_\tau &= 2 \tau \cdot \D u \ n - 2 n\cdot (\tau\cdot \nabla) u
    \\
    &= (\tau_j n_i - n_j \tau_i)\partial_i u_j
    \\
    &= (\tau_1 n_2 - n_1 \tau_2)\partial_2 u_1 + (\tau_2 n_1 - n_2 \tau_1)\partial_1 u_2
    \\
    &= -\partial_2 u_1 + \partial_1 u_2
    \\
    &= \omega
\end{align}
where in the second to last estimate we used $\tau = n^\perp$ and $|n|=1$.

\begin{lemma}
    Let $\Omega$ be $C^{1,1}$, $2< \tilde q \in 2\mathbb{N}$, $u_0\in W^{1,\tilde q}(\Omega)$, $\theta_0 \in L^{\tilde q}(\Omega)$ and $0<\alpha \in L^\infty(\partial\Omega)$. Then
    \begin{align}
        u \in L^\infty\left((0,\infty);W^{1,q}(\Omega)\right)
    \end{align}
    for any $1\leq q \leq \tilde q$ and
    \begin{align}
        \|u\|_{W^{1,q}}\lessapprox \|u\|_{W^{1,\tilde q}}+ \|\theta_0\|_{\tilde q},
    \end{align}
    where the implicit constant depends only on $q$, $\tilde q$, $\alpha$ and $\Omega$.
    \label{lemma_nd_W1q_bound}
\end{lemma}

\begin{proof}
    Fix some arbitrary $\tmax>0$, define
    \begin{align}
        \Lambda = 2 \|(\alpha+\kappa)u_\tau \|_{L^\infty\left([0,\tmax]\times \partial\Omega\right)}
    \end{align}
    and let $\omegatildepm$ solve
    \begin{alignat}{2}
        \omegatildepm_t + u\cdot \nabla\omegatildepm - \Delta \omegatildepm &= \partial_1 \theta \qquad &\textnormal{ in }&\Omega
        \\
        \omegatildepm &= \pm \Lambda \qquad &\textnormal{ on }&\partial\Omega
        \\
        \omegatildepm(\cdot,0)&= \pm |\omega_0|\qquad &\textnormal{ in }&\Omega.
    \end{alignat}
    Then $\omegabarpm$, defined by
    \begin{align}
        \omegabarpm = \omega - \omegatildepm,
    \end{align}
    solves
    \begin{alignat}{2}
        \omegabarpm_t + u\cdot \nabla \omegabarpm - \Delta \omegabarpm &= 0 \qquad &\textnormal{ in }&\Omega
        \\
        \omegatildepm &= -2(\alpha+\kappa) u_\tau \mp \Lambda \qquad &\textnormal{ on }&\partial\Omega
        \\
        \omegabarpm(\cdot,0)&= \omega_0 \mp |\omega_0|\qquad &\textnormal{ in }&\Omega.
    \end{alignat}
    The initial and boundary values of $\omegabarp$ and $\omegabarm$ non-positive, respectively non-negative the maximum principle implies $\omegabarp\leq 0 $ and $\omegabarm\geq 0$, and therefore, $\omega$ can be bounded by
    \begin{align}
        |\omega|\leq |\omegatildep|+|\omegatildem|.
        \label{cvhubcvuhbui}
    \end{align}
    We define
    \begin{align}
        \omegahatpm = \omegatildepm \mp \Lambda
        \label{uvbchbihvcub}
    \end{align}
    in order to remove the boundary condition. The following analysis will hold the same way for $\omegahatp$ as it does for $\omegahatm$. Therefore, we focus on $\omegahatm$ and omit the $^-$ to simplify the notation. Note that $\hat \omega$ solves
    \begin{alignat}{2}
        \hat \omega + u\cdot\nabla\hat \omega - \Delta \hat\omega &= \partial_1\theta \qquad &\textnormal{ in }&\Omega
        \label{cbhuxixucvbh}
        \\
        \hat \omega &= 0 \qquad &\textnormal{ on }&\partial\Omega
        \label{nd_hatomega_bc}
        \\
        \hat\omega(\cdot,0) &= -|\omega_0|+\Lambda \qquad &\textnormal{ in }&\Omega.\label{sadfjsdafsdfi}
    \end{alignat}\noeqref{sadfjsdafsdfi}%
    Testing \eqref{cbhuxixucvbh} with $\hat \omega^{{\tilde q}-1}$ yields
    \begin{align}
        \frac{1}{{\tilde q}}\frac{d}{dt}\|\hat\omega\|_{\tilde q}^{\tilde q} = - \int_{\Omega} \hat\omega^{{\tilde q}-1} u\cdot \nabla\hat\omega + \int_{\Omega} \hat\omega^{{\tilde q}-1} \Delta \hat \omega + \int_{\Omega} \hat\omega^{{\tilde q}-1} \partial_1\theta
        \label{cvxbzhuicvbixcvb}
    \end{align}
    and, after integration by parts, the first term on the right-hand side of \eqref{cvxbzhuicvbixcvb} vanishes as
    \begin{align}
        {\tilde q}\int_{\Omega} \hat\omega^{{\tilde q}-1} u\cdot\nabla\hat\omega = \int_{\Omega} u\cdot\nabla (\hat\omega^{{\tilde q}}) = -\int_{\Omega} \hat\omega^{\tilde q} \nabla \cdot u + \int_{\partial\Omega} \hat \omega^{\tilde q} n\cdot u = 0.\quad
        \label{vcuixcvbz}
    \end{align}
    Note that for the second term on the right-hand side of \eqref{cvxbzhuicvbixcvb}, integration by parts using \eqref{nd_hatomega_bc} yields
    \begin{align}
        \int_{\Omega} \hat\omega^{{\tilde q}-1} \Delta\hat\omega
        &= - ({\tilde q}-1) \int_{\Omega} \hat\omega^{{\tilde q}-2} |\nabla \hat\omega|^2
        = - ({\tilde q}-1) \|\hat\omega^\frac{{\tilde q}-2}{2} \nabla \hat\omega\|_2^2.
        \label{vgxzhvbhjcvxbhcj}
    \end{align}
    Again using integration by parts and \eqref{nd_hatomega_bc} the last term on the right-hand side can be written as
    \begin{align}
        \int_{\Omega} \hat\omega^{{\tilde q}-1} \partial_1 \theta
        &= - ({\tilde q}-1) \int_{\Omega} \theta \hat\omega^{{\tilde q}-2} \partial_1 \hat\omega,
    \end{align}
    which, after applying Hölder's and Young's inequality, yields
    \begin{align}
        \int_{\Omega} \hat\omega^{{\tilde q}-1} \partial_1 \theta
        &\leq ({\tilde q}-1) \|\hat\omega^{\frac{{\tilde q}-2}{2}} \nabla \hat\omega\|_2 \|\hat\omega^\frac{{\tilde q}-2}{2}\theta\|_2
        \\
        &\leq \frac{{\tilde q}-1}{2} \|\hat\omega^{\frac{{\tilde q}-2}{2}} \nabla \hat\omega\|_2^2 + \frac{{\tilde q}-1}{2} \|\hat\omega^\frac{{\tilde q}-2}{2}\theta\|_2^2
        \label{uzvcubcv}
    \end{align}
    Combining \eqref{cvxbzhuicvbixcvb}, \eqref{vcuixcvbz}, \eqref{vgxzhvbhjcvxbhcj} and \eqref{uzvcubcv} we obtain
    \begin{align}
        \frac{1}{{\tilde q}}\frac{d}{dt}\|\hat\omega\|_{\tilde q}^{\tilde q} + \frac{{\tilde q}-1}{2} \|\hat\omega^{\frac{{\tilde q}-2}{2}} \nabla \hat\omega\|_2^2
        &\leq \frac{{\tilde q}-1}{2} \|\hat\omega^\frac{{\tilde q}-2}{2}\theta\|_2^2
        \\
        &= \frac{{\tilde q}-1}{2} \|\hat\omega^{{\tilde q}-2}\|_{\frac{{\tilde q}}{{\tilde q}-2}} \|\theta^2\|_{\frac{{\tilde q}}{2}}
        \\
        &= \frac{{\tilde q}-1}{2} \|\hat\omega\|_{\tilde q}^{{\tilde q}-2}\|\theta\|_{\tilde q}^2
    \end{align}
    and using
    \begin{align}
        \|\hat\omega^{\frac{{\tilde q}-2}{{\tilde q}}}\nabla\hat\omega\|_2^2 = \int_\Omega \left(\hat\omega^\frac{{\tilde q}-2}{2}\nabla\hat\omega\right)^2 = \int_\Omega \left(\tfrac{2}{{\tilde q}}\nabla \left(\hat\omega^\frac{{\tilde q}}{2}\right)\right)^2 = \frac{4}{{\tilde q}^2} \|\nabla(\hat\omega^\frac{{\tilde q}}{2})\|_2^2
    \end{align}
    one gets
    \begin{align}
        \frac{1}{{\tilde q}}\frac{d}{dt}\|\hat\omega\|_{\tilde q}^{\tilde q} + \frac{2({\tilde q}-1)}{{\tilde q}^2} \|\nabla (\hat\omega^{\frac{{\tilde q}}{2}})\|_2^2
        &\leq \frac{{\tilde q}-1}{2} \|\hat\omega\|_{\tilde q}^{{\tilde q}-2}\|\theta\|_{\tilde q}^2.
    \end{align}
    Since $\hat\omega=0$ on $\partial\Omega$, Poincaré's inequality yields
    \begin{align}
        \|\hat \omega\|_{\tilde q}^{\tilde q} = \|\hat \omega^\frac{{\tilde q}}{2}\|_2^2 \lesssim \|\nabla (\hat\omega^\frac{{\tilde q}}{2})\|_2^2
    \end{align}
    and therefore, there exists a constant $C>0$ such that
    \begin{align}
        \frac{1}{2}\|\hat\omega\|_{\tilde q}^{{\tilde q}-2}\frac{d}{dt}\|\hat\omega\|_{\tilde q}^2 + C\frac{{\tilde q}-1}{{\tilde q}^2} \|\hat\omega\|_{\tilde q}^{\tilde q}
        &=
        \frac{1}{{\tilde q}}\frac{d}{dt}\|\hat\omega\|_{\tilde q}^{\tilde q} + C\frac{{\tilde q}-1}{{\tilde q}^2} \|\hat\omega\|_{\tilde q}^{\tilde q}
        \\
        &\leq \frac{{\tilde q}-1}{2} \|\hat\omega\|_{\tilde q}^{{\tilde q}-2} \|\theta\|_{\tilde q}^2.
        \label{zubxicvbhuxcb}
    \end{align}
    Dividing \eqref{zubxicvbhuxcb} by $\|\hat\omega\|_{\tilde q}^{{\tilde q}-2}$ implies
    \begin{align}
        \frac{d}{dt} \|\hat\omega\|_{\tilde q}^2 + C\frac{{\tilde q}-1}{{\tilde q}^2}\|\hat\omega\|_{\tilde q}^2 \leq ({\tilde q}-1) \|\theta\|_{\tilde q}^2 =({\tilde q}-1) \|\theta_0\|_{\tilde q}^2,
        \label{hsafsuidffhs}
    \end{align}
    where in the last identity we used \eqref{nd_maximum_principle}. Applying Grönwall's inequality to \eqref{hsafsuidffhs} yields
    \begin{align}
        \|\hat\omega(\cdot,t)\|_{\tilde q}^2 \leq \|\hat\omega_0\|_{\tilde q}^2 + C {\tilde q}^2 \|\theta_0\|_{\tilde q}^2
        \label{cuixbcvxjbklcj}
    \end{align}
    for some constant $C>0$ only depending on $\Omega$ and $\tilde q$.
    Next, we estimate $\Lambda$. One has
    \begin{align}
        \Lambda &= 2\|(\alpha+\kappa)u_\tau\|_{L^\infty([0,\tmax]\times \partial\Omega)}
        \\
        &\lesssim \|\alpha+\kappa\|_\infty \|u\|_{L^\infty\left((0,\tmax);L^\infty(\partial\Omega)\right)}
        \\
        &\lessapprox \|u\|_{L^\infty\left((0,\tmax);L^\infty(\Omega)\right)},
        \label{hcvubixhbxucbh}
    \end{align}
    where the $\lessapprox$ hides the additional dependency of the implicit constant on $\|\alpha+\kappa\|_\infty$.
    For details the reader is refered to Satz 7.1.26 in \cite{emmrichGewoehnlicheUndOperatorDifferentialgleichungen}, 
    respectively Exercise 4.1 in \cite{troeltzschOptimaleSteuerungpartiellerDifferentialgleichungen}. 
    Using Gagliardo-Nirenberg interpolation 
    and Young's inequality one can further estimate \eqref{hcvubixhbxucbh} by
     \begin{align}
        \Lambda
        &\lessapprox \|u\|_{L^\infty_t;L^\infty_x}
        \\
        &\lessapprox \|\nabla u\|_{L^\infty_t;L^{\tilde q}_x}^{\frac{{\tilde q}}{2({\tilde q}-1)}}\|u\|_{L^\infty_t;;L^2_x}^\frac{{\tilde q}-2}{2({\tilde q}-1)} + \|u\|_{L^\infty_t;L^2_x}
        \\
        &\lesssim \epsilon \|\nabla u\|_{L^\infty_t;L^{\tilde q}_x} + \left(1+\epsilon^{\frac{{\tilde q}}{2-{\tilde q}}}\right) \|u\|_{L^\infty_t;L^2_x}
        \\
        &\lessapprox \epsilon \|\omega\|_{L^\infty_t;L^{\tilde q}_x} + \left(1+\epsilon+\epsilon^{\frac{{\tilde q}}{2-{\tilde q}}}\right) \|u\|_{L^\infty_t;L^2_x}
        \label{uihcvubichvub}
     \end{align}
     for any $\epsilon>0$, where $L^\infty_t;L^r_x=L^\infty\left((0,T);L^r(\Omega)\right)$ and in the last estimate we used Lemma \ref{lemma_elliptic_regularity_bounded_domain}. Combining \eqref{cvhubcvuhbui}, \eqref{uvbchbihvcub}, \eqref{cuixbcvxjbklcj} and \eqref{uihcvubichvub} yields
     \begin{align}
         \|\omega\|_{L^\infty_t;L^{\tilde q}_x} &\leq \|\omegatildem\|_{L^\infty_t;L^{\tilde q}_x}+\|\omegatildep\|_{L^\infty_t;L^{\tilde q}_x}
         \\
         &\lesssim \|\omegahatp\|_{L^\infty_t;L^{\tilde q}_x}+\|\omegahatm\|_{L^\infty_t;L^{\tilde q}_x}+\Lambda
         \\
         &\lesssim \|\omegahatp_0\|_{\tilde q} + \|\omegahatm_0\|_{\tilde q} + \|\theta_0\|_{\tilde q} + \Lambda
         \\
         &\lesssim \|\omega_0\|_{\tilde q} + \|\theta_0\|_{\tilde q} + \Lambda
         \\
         &\lesssim \|\omega_0\|_{\tilde q} + \|\theta_0\|_{\tilde q} + \Lambda
         \\
         &\lesssim \|\omega_0\|_{\tilde q} + \|\theta_0\|_{\tilde q} + \Lambda
         \\
         &\lessapprox \|\omega_0\|_{\tilde q} + \|\theta_0\|_{\tilde q} + \epsilon \|\omega\|_{L^\infty_t;L^{\tilde q}_x} + \left(1+\epsilon+\epsilon^{\frac{{\tilde q}}{2-{\tilde q}}}\right) \|u\|_{L^\infty_t;L^2_x}
     \end{align}
     and choosing $\epsilon$ sufficiently small, one can compensate the second term on the right-hand side implying
     \begin{align}
         \|\omega\|_{L^\infty_t;L^{\tilde q}_x}
         &\lessapprox \|\omega_0\|_{\tilde q} + \|\theta_0\|_{\tilde q} + \|u\|_{L^\infty_t;L^2_x}.
     \end{align}
     Note that by Lemma \ref{lemma_elliptic_regularity_bounded_domain}, the full Sobolev norm of $u$ is bounded by
     \begin{align}
         \|u\|_{L^\infty_t;W^{1,{\tilde q}}_x} 
         &\lessapprox \|\omega\|_{L^\infty_t;L^{{\tilde q}}_x} + \|u\|_{L^\infty_t;L^2_x}
         \\
         &\lessapprox \|\omega_0\|_{\tilde q} + \|\theta_0\|_{\tilde q} + \|u\|_{L^\infty_t;L^2_x}.
     \end{align}
     Then, using Hölder's inequality and the energy bound, i.e. Lemma \ref{lemma_nd_energy_bound}, one finds
     \begin{align}
         \|u\|_{L^\infty_t;W^{1,q}_x} &\lesssim \|u\|_{L^\infty_t;W^{1,{\tilde q}}_x}
         \lessapprox \|\omega_0\|_{\tilde q} + \|\theta_0\|_{\tilde q} + \|u_0\|_2 + \|\theta_0\|_2
         \\
         &\lesssim \|u_0\|_{W^{1,{\tilde q}}} + \|\theta_0\|_{\tilde q},
     \end{align}
     where in the last estimate we additionally used the definition of the vorticity. Finally, since the right-hand side is independent of $\tmax$, it holds uniform in time, proving the claim.
\end{proof}

Next, we want to derive bounds for the pressure. Instead of first deriving a Poisson equation and boundary conditions for the pressure as in Section 2.2 of \cite{bleitnerCarlsonNobili2023Large} and Lemma \ref{lemma_nonlin_pressure_bound}, here we use the similar but more direct approach of Lemma \ref{lemma_pressure_bound}.

\begin{lemma}[Pressure Bound]
    \label{lemma_nd_pressure_bound}
    Let $\Omega$ be $C^{2,1}$, $u_0\in W^{1,4}(\Omega)$, $\theta_0 \in L^{4}(\Omega)$ and $0<\alpha \in W^{1,\infty}(\partial\Omega)$. Then
    \begin{align}
        p\in L^\infty\left((0,\infty);H^1(\Omega)\right).
    \end{align}
\end{lemma}

\begin{proof}
    Testing \eqref{nd_navier_stokes} with $\nabla p$ yields
    \begin{align}
        \|\nabla p\|_2^2 = - \int_{\Omega} \nabla p \cdot u_t - \int_{\Omega} \nabla p \cdot (u\cdot\nabla)u + \int_{\Omega} \nabla p \cdot \Delta u + \int_{\Omega} \partial_2 p \theta.
        \label{cvnibcvnub}
    \end{align}
    The first term on the right-hand side of \eqref{cvnibcvnub} vanishes under partial integration since
    \begin{align}
        \int_{\Omega} \nabla p \cdot u_t = \int_{\partial\Omega} p n\cdot u_t - \int_{\Omega} p \nabla \cdot u_t = 0,
        \label{bbnvcbjkcvnb}
    \end{align}
    where in the last identity we used $n \cdot u_t = 0$ and $\nabla \cdot u_t=0$, which both follow from taking the time derivative of \eqref{nd_no_penetration_bc} and \eqref{nd_incompressible}, respectively.
    
    For the third term on the right-hand side of \eqref{cvnibcvnub}, \eqref{id_Delta_u_is_nablaPerp_omega} and integration by parts yield
    \begin{align}
        \int_{\Omega} \nabla p\cdot\Delta u
        &= \int_{\Omega} \nabla p\cdot\nabla^\perp\omega
        \\
        &= \int_{\partial\Omega} \nabla p \cdot n^\perp \omega - \int_{\Omega} \nabla^\perp\cdot\nabla p \omega
    \end{align}
    The term in the bulk vanishes since $\nabla^\perp\cdot\nabla p = 0$ and, using $\tau = n^\perp$ and \eqref{nd_vorticity_bc}, the boundary term can be written as
    \begin{align}
        \nabla p \cdot n^\perp \omega = -2 (\alpha+\kappa)\nabla p \cdot \tau u_\tau = -2(\alpha+\kappa) u\cdot\nabla p.
    \end{align}
    Combining these identities implies
    \begin{align}
        \int_{\Omega} \nabla p \cdot \Delta u = - 2\int_{\partial\Omega} (\alpha+\kappa) u\cdot \nabla p,
    \end{align}
    which, together with \eqref{cvnibcvnub} and \eqref{bbnvcbjkcvnb}, yields
    \begin{align}
        \|\nabla p\|_2^2 = - \int_{\Omega} \nabla p \cdot (u\cdot\nabla)u - 2 \int_{\partial\Omega} (\alpha+\kappa) u\cdot \nabla p + \int_{\Omega} \partial_2 p \theta.
        \label{vnbicvnib}
    \end{align}
    Note that as $p$ is only defined up to a constant, we can choose it such that $p$ has a vanishing average and Poincaré's inequality holds. Therefore \eqref{vnbicvnib} implies
    \begin{align}
        \|p\|_{H^1}^2 &= \|\nabla p\|_2^2 + \| p \|_2^2
        \\
        &\lesssim \|\nabla p\|_2^2 
        \\
        &= - \int_{\Omega} \nabla p \cdot (u\cdot\nabla)u - 2 \int_{\partial\Omega} (\alpha+\kappa) u\cdot \nabla p + \int_{\Omega} \partial_2 p \theta.
        \label{cvnbicvnub}
    \end{align}
    The first and third terms on the right-hand side of \eqref{cvnbicvnub} can be bounded using Hölder's inequality, and the second one using Lemma \ref{lemma_u_cdot_nabla_p_on_boundary}, implying
    \begin{align}
        \|p\|_{H^1}^2
        &\lessapprox \|\nabla p\|_2 \|u\|_4 \|u\|_{W^{1,4}} + \|p\|_{H^1} \|u\|_{H^1} + \|\theta\|_2 \|\nabla p\|_2
        \\
        &\leq (\|u\|_{W^{1,4}}^2 + \|u\|_{H^1} +\|\theta\|_2)\|p\|_{H^1}
    \end{align}
    and dividing by $\|p\|_{H^1}$ one has
    \begin{align}
        \|p\|_{H^1}
        &\lessapprox (1+\|u\|_{W^{1,4}})\|u\|_{W^{1,4}} +\|\theta\|_2.
    \end{align}
    Finally using Lemma \ref{lemma_nd_W1q_bound}, \eqref{nd_maximum_principle} and Hölder's inequality
    \begin{align}
        \|p\|_{H^1} \lessapprox (1+\|u_0\|_{W^{1,4}}+\|\theta_0\|_4) (\|u_0\|_{W^{1,4}}+\|\theta_0\|_4) \lesssim  1+\|u_0\|_{W^{1,4}}^2+\|\theta_0\|_4^2,
    \end{align}
    where we used Young's inequality in the last estimate.
\end{proof}

The following Lemma corresponds to Lemma 2.5 in \cite{bleitnerCarlsonNobili2023Large}.
\begin{lemma}[Time Derivative Bound]
    \label{lemma_nd_time_derivative_bound}
    Let $\Omega$ be $C^{1,1}$, $u_0\in H^2(\Omega)$, $\theta_0 \in L^{4}(\Omega)$ and $0<\alpha \in L^\infty(\partial\Omega)$. Then
    \begin{align}
        u_t\in L^\infty\left((0,\infty);L^2(\Omega)\right)\cap L^2\left((0,\infty);H^1(\Omega)\right).
    \end{align}
\end{lemma}

\begin{proof}
    Differentiating \eqref{nd_navier_stokes}, \eqref{nd_incompressible}, \eqref{nd_no_penetration_bc} and \eqref{nd_navier_slip} and plugging in \eqref{nd_transport_equation} yields the system
    \begin{alignat}{2}
        u_{tt} + u_t\cdot\nabla u + u\cdot\nabla u_t + \nabla p_t - \Delta u_t &= - u\cdot\nabla \theta e_2 &\qquad \textnormal{ in }&\Omega
        \label{nd_ut_pde}
        \\
        \nabla \cdot u_t &= 0 &\qquad \textnormal{ in }&\Omega
        \label{nd_ut_incompressible}
        \\
        \tau \cdot (\D u_t + \alpha u_t) &= 0 &\qquad \textnormal{ on }&\partial\Omega\label{nd_ut_nav_slip}
        \\
        n \cdot u_t &= 0 &\qquad \textnormal{ on }&\partial\Omega.
        \label{nd_ut_no_pen}
    \end{alignat}\noeqref{nd_ut_nav_slip}%
    Testing \eqref{nd_ut_pde} with $u_t$ one finds
    \begin{align}
        \frac{1}{2}\frac{d}{dt}\|u_t\|_2^2 &= -\int_{\Omega} u_t\cdot (u_t\cdot\nabla)u - \int_{\Omega} u_t\cdot (u\cdot \nabla) u_t - \int_{\Omega} u_t \cdot \nabla p_t  + \int_{\Omega} u_t\cdot \Delta u_t
        \\
        &\qquad - \int_{\Omega} (u_t\cdot e_2)u\cdot\nabla \theta.
        \label{invijxcnvsjn}
    \end{align}
    Note that the second and third terms on the right-hand side of \eqref{invijxcnvsjn} vanish due to
    \begin{align}
        2\int_{\Omega} u_t\cdot (u\cdot\nabla) u_t = \int_{\Omega} u\cdot\nabla (u_t\cdot u_t) = \int_{\partial\Omega} u\cdot n u_t^2 - \int_{\Omega} u_t^2 \nabla \cdot u = 0\qquad
        \label{ernjternt}
    \end{align}
    and
    \begin{align}
        \int_{\Omega} u_t\cdot \nabla p_t = \int_{\partial\Omega} u_t\cdot n p_t - \int_{\partial\Omega} p_t \nabla \cdot u_t = 0.
        \label{erntretkj}
    \end{align}
    Note that Lemma \ref{lemma_int_Delta_u_v} also applies to $u_t$, implying
    \begin{align}
        - \int_{\Omega} u_t\cdot\Delta u_t = 2\|\D u_t\|_2^2 + 2 \int_{\partial\Omega} \alpha u_t^2
        \label{ioerjtioerj}
    \end{align}
    covering the fourth term on the right-hand side of \eqref{invijxcnvsjn}.
    Combining \eqref{invijxcnvsjn}, \eqref{ernjternt}, \eqref{erntretkj} and \eqref{ioerjtioerj}, and using integration by parts, one gets
    \begin{align}
        \frac{1}{2}\frac{d}{dt}\|u_t\|_2^2 &+ 2 \|\D u_t\|_2^2 + 2 \int_{\partial\Omega}\alpha u_\tau^2
        \\
        &= -\int_{\Omega} u_t\cdot (u_t\cdot\nabla)u - \int_{\Omega} (u_t\cdot e_2)u\cdot\nabla \theta
        \\
        &= -\int_{\Omega} u_t\cdot (u_t\cdot\nabla)u - \int_{\partial\Omega} (u_t\cdot e_2)u\cdot n \theta + \int_{\Omega} (u_t\cdot e_2)\theta \nabla \cdot u
        \\
        &\qquad + \int_{\Omega} \theta u\cdot\nabla(u_t\cdot e_2)
        \\
        &= -\int_{\Omega} u_t\cdot (u_t\cdot\nabla)u + \int_{\Omega} \theta u\cdot\nabla(u_t\cdot e_2),
    \end{align}
    where in the last identity we used \eqref{nd_no_penetration_bc} and \eqref{nd_incompressible}.
    Lemma \ref{lemma_coercivity} also applies to $u_t$, implying the existence of a constant $C>0$ depending on $\Omega$ such that
    \begin{align}
        \frac{1}{2}\frac{d}{dt} \|u_t\|_2^2 + C\|u_t\|_{H^1}^2 &\leq \int_{\Omega} |u_t\cdot (u_t\cdot\nabla)u | + \int_{\Omega}| \theta u\cdot\nabla(u_t\cdot e_2)|.
    \end{align}
    Using Hölder's and Ladyzhenskaya's inequality, we obtain
    \begin{align}
        \frac{1}{2}\frac{d}{dt} \|u_t\|_2^2 + C\|u_t\|_{H^1}^2 
        &\lesssim \|u_t\|_4^2 \|u\|_{H^1} + \|\theta\|_4 \|u\|_4 \|\nabla u_t\|_2
        \\
        &\lesssim \|\nabla u_t\|_2\|u_t\|_2 \|u\|_{H^1} + \|\theta\|_4\|u_t\|_2 \|u\|_{H^1}
        \\
        &\lesssim \epsilon \|\nabla u_t\|_2^2 + \epsilon^{-1} \|u_t\|_2^2 \|u\|_{H^1}^2 + \epsilon \|u_t\|_2^2 + \epsilon^{-1}\|\theta\|_4^2\|u\|_{H^1}^2
        \label{nertnekrt}
    \end{align}
    for any $\epsilon>0$, where in the last estimate we used Young's inequality. Choosing $\epsilon$ sufficiently small one can compensate the first and third term on the right-hand side of \eqref{nertnekrt}, implying
    \begin{align}
        \frac{d}{dt}\|u_t\|_2^2 + C \|u_t\|_{H^1}^2 &\lesssim \|u_t\|_2^2 \|u\|_{H^1}^2 + \|\theta\|_4^2 \|u\|_{H^1}^2
        \\
        &= \|u_t\|_2^2 \|u\|_{H^1}^2 + \|\theta_0\|_4^2 \|u\|_{H^1}^2,
        \label{jerjrewnbjkr}
    \end{align}
    where we used \eqref{nd_maximum_principle}. Grönwall's inequality and Lemma \ref{lemma_nd_L2H1_bound} now yield
    \begin{align}
        \|u_t(t)\|_2^2 &\lessapprox \|u_t(0)\|_2^2 e^{C\int_0^t\|u(s)\|_{H^1}^2 ds} + \|\theta_0\|_4^2 \int_0^t e^{C\int_s^t \|u(\lambda)\|_{H^1}^2d\lambda}\|u(s)\|_{H^1}\ ds
        \\
        &\lessapprox \left(\|u_t(0)\|_2^2+ (\|u_0\|_2^2+\|\theta_0\|_2^2+1)\|\theta_0\|_4^2\right) e^{C(\|u_0\|_2^2+\|\theta_0\|_2^2+1)}.\label{njkbnjcbnj}
    \end{align}
    For smooth solutions, testing \eqref{nd_navier_stokes} with $u_t$ one has
    \begin{align}
        \|u_t\|_2^2 \leq \|u_t\|_2 \|u\cdot\nabla u\|_2 + \|u_t\|_2 \|\Delta u\|_2 + \|u_t\|_2 \|\theta\|_2,
        \label{inbxucvhxb}
    \end{align}
    where again the pressure term vanishes due to \eqref{nd_ut_incompressible} and \eqref{nd_ut_no_pen}. Dividing \eqref{inbxucvhxb} by $\|u_t\|_2$ and using Hölder's and Ladyzhenskaya's inequality, it follows, that
    \begin{align}
        \|u_t\|_2 \leq \|u\|_4\|\|\nabla u\|_4 + \|u\|_{H^2} + \|\theta_0\|_2 \lesssim (1+\|u\|_{H^1})\|u\|_{H^2} + \|\theta\|_2,
        \label{njkretnjek}
    \end{align}
    implying $\|u_t(0)\|_2 \leq (1+\|u_0\|_{H^1})\|u_0\|_{H^2}+\|\theta_0\|_2$. Therefore the right-hand side of \eqref{njkbnjcbnj} is uniformly bounded in time, proving
    \begin{align}
        u_t\in L^\infty\left((0,\infty);L^2(\Omega)\right).
    \end{align}
    Plugging \eqref{njkretnjek} into \eqref{jerjrewnbjkr} yields
    \begin{align}
        \frac{d}{dt}\|u_t\|_2^2 + C \|u_t\|_{H^1}^2 &\leq C_0 \|u\|_{H^1}^2
        \label{njkerntjj}
    \end{align}
    with a constant $C_0>0$ depending on $\|u_0\|_{H^2}, \|\theta_0\|_4$, $\alpha$ and $\Omega$. Integrating \eqref{njkerntjj} in time results in
    \begin{align}
        \|u_t(t)\|_2^2 + C \int_0^t \|u_t(s)\|_{H^1}^2\ ds \leq C_0  \int_0^t \|u\|_{H^1}^2\ ds \leq C_0,
    \end{align}
    where in the last estimate we used Lemma \ref{lemma_nd_L2H1_bound}. As the right-hand side is independent of $t$, this bound holds uniformly in time, implying
    \begin{align}
        u_t \in L^2\left((0,\infty);H^1(\Omega)\right).
    \end{align}
\end{proof}

We are now able to prove uniform in time $H^2$-bounds for the velocity.

\begin{lemma}[{\texorpdfstring{$H^2$}{H2}-Bounds}]
\label{lemma_nd_LinftyH2_bound}
    Let $\Omega$ be $C^{2,1}$, $u_0\in H^2(\Omega)$, $\theta_0 \in L^{4}(\Omega)$ and $0<\alpha \in W^{1,\infty}(\partial\Omega)$. Then
    \begin{align}
        u\in L^\infty\left((0,\infty);H^2(\Omega)\right)\cap L^q\left((0,\infty);W^{1,q}(\Omega)\right)
    \end{align}
    for any $2\leq q <\infty$.
\end{lemma}

\begin{proof}
    Using \eqref{id_Delta_u_is_nablaPerp_omega}, taking the $L^2$-norm of \eqref{nd_navier_stokes} and applying Hölder's inequality we find
    \begin{align}
        \|\nabla \omega\|_2 &= \|\Delta u\|_2 = \|u_t + u\cdot \nabla u +\nabla p -\theta e_2\|_2
        \\
        &\leq \|u_t\|_2 + \|u\|_{W^{1,4}}^2 + \|\nabla p \|_2 + \|\theta\|_2
        \label{enjrtkkenrjt}
    \end{align}
    Combining Lemma \ref{lemma_elliptic_regularity_bounded_domain} with \eqref{enjrtkkenrjt}, we get
    \begin{align}
        \|u\|_{H^2} &\lessapprox \|\nabla \omega\|_2 + \|u\|_{H^1}
        \\
        &\lesssim \|u_t\|_2 + \|u\|_{W^{1,4}}^2 + \|u\|_{H^1} + \|\nabla p \|_2 + \|\theta\|_2.
        \label{njrekternjt}
    \end{align}
    By Lemma \ref{lemma_nd_time_derivative_bound}, Lemma \ref{lemma_nd_pressure_bound}, Lemma \ref{lemma_nd_W1q_bound} and \eqref{nd_maximum_principle}, the right-hand side of \eqref{njrekternjt} is uniformly bounded in time, implying
    \begin{align}
        u\in L^\infty\left((0,\infty);H^2(\Omega)\right).\label{njekwrn}
    \end{align}
    Additionally, by Gagliardo-Nirenberg interpolation, we obtain
    \begin{align}
        \|u\|_{W^{1,q}}^q \lesssim \|u\|_{H^2}^{q-2}\|u\|_{H^1}^2
    \end{align}
    for $2\leq q<\infty$, which, in view of \eqref{njekwrn} and Lemma \ref{lemma_nd_L2H1_bound}, yields
    \begin{align}
        u\in L^q\left((0,\infty);W^{1,q}(\Omega)\right).
    \end{align}
\end{proof}

\begin{lemma}[\texorpdfstring{$H^3$}{H3}-Bounds]
    Let $\Omega$ be $C^{3,1}$, $u_0\in H^2(\Omega)$, $\theta_0 \in W^{1,\tilde q}(\Omega)$ for some $\tilde q \in 2\mathbb{N}$ and $0<\alpha \in W^{2,\infty}(\partial\Omega)$. Then
    \begin{align}
        u&\in L^2\left((0,T);H^3(\Omega)\right) \cap L^\frac{2r}{r-2}\left((0,T);W^{2,r}(\Omega)\right)
        \\
        \theta&\in L^\infty\left((0,T);W^{1,q}(\Omega)\right)
    \end{align}
    for any $T>0$, $2\leq r <\infty$ and $1\leq q \leq \tilde q$.
\end{lemma}

\begin{proof}
    Testing \eqref{nd_vorticity_pde} with $\Delta\omega$ one has
    \begin{align}
        \int_\Omega \omega_t \Delta \omega + \int_{\Omega} u\cdot\nabla \omega \Delta \omega  -\|\Delta\omega\|_2^2 = \int_{\Omega}\partial_1\theta\Delta\omega.
        \label{nkrentkertjn}
    \end{align}
    For the first term on the left-hand side integration by parts yields
    \begin{align}
        \int_{\Omega} \omega_t \Delta \omega &= \int_{\partial\Omega} \omega_t n\cdot\nabla \omega - \int_\Omega \nabla \omega_t \cdot\nabla \omega
        \\
        &= - 2\int_{\partial\Omega} (\alpha +\kappa) u_t\cdot \tau n\cdot\nabla \omega - \frac{1}{2}\frac{d}{dt} \|\nabla\omega\|_2^2,
        \label{nrtjkernjter}
    \end{align}
    where in the last identity we have used that taking the derivative of \eqref{vorticity_bc} with respect to time yields $\omega_t = -2(\alpha+\kappa)u_t\cdot \tau$ on $\partial\Omega$. Combining \eqref{nkrentkertjn} and \eqref{nrtjkernjter}, using the trace theorem, as well as Hölder's and Ladyzhenskaya's inequality, one gets
    \begin{align}
        \frac{1}{2}\frac{d}{dt}\|\nabla\omega\|_2^2 +\|\Delta \omega\|_2^2
        &= -2\int_{\partial\Omega} (\alpha +\kappa) u_t\cdot \tau n\cdot\nabla \omega + \int_{\Omega} u\cdot\nabla\omega \Delta \omega - \int_{\Omega} \partial_1 \theta\Delta \omega
        \\[-12pt]
        &\lessapprox \|u_t\|_{H^1}\|\nabla\omega\|_{H^1} + \|\Delta \omega\|_2 \|\nabla \omega\|_4 \|u\|_4 + \|\nabla\theta\|_2 \|\Delta\omega\|_2
        \\
        &\lesssim \|u_t\|_{H^1}\|u\|_{H^3} + \|u\|_{H^3}^\frac{3}{2} \|u\|_{H^2}^\frac{1}{2} \|u\|_{H^1} + \|\nabla\theta\|_2 \|u\|_{H^3}
        \\
        &\lesssim \epsilon \|u\|_{H^3}^2 + \epsilon^{-1} \|u_t\|_{H^1}^2 + \epsilon^{-3} \|u\|_{H^2}^2 \|u\|_{H^1}^{4} + \epsilon^{-1}\|\nabla\theta\|_2^2\\[-2pt]
        \label{nberjtrb}
    \end{align}
    for any $\epsilon>0$ where in the last estimate we used Young's inequality. By Lemma \ref{lemma_elliptic_regularity_bounded_domain} we obtain
    \begin{align}
        \|u\|_{H^3}^2 \lessapprox \|\Delta \omega\|_2^2 + \|u\|_{H^2}^2,
    \end{align}
    which combined with \eqref{nberjtrb} yields
    \begin{align}
        \frac{1}{2}\frac{d}{dt}\|\nabla\omega\|_2^2 +\|\Delta \omega\|_2^2
        &\lessapprox \epsilon \|\Delta \omega\|_2^2 + \epsilon \|u\|_{H^2}^2 + \epsilon^{-1} \|u_t\|_{H^1}^2
        \\
        &\qquad + \epsilon^{-3} \|u\|_{H^2}^2 \|u\|_{H^1}^{4} + \epsilon^{-1}\|\nabla\theta\|_2^2
    \end{align}
    and therefore, choosing $\epsilon$ sufficiently small in order to compensate for the $\Delta\omega$-term one has
    \begin{align}
        \frac{d}{dt}\|\nabla\omega\|_2^2 +\|\Delta \omega\|_2^2
        &\lessapprox \|u_t\|_{H^1}^2 + (1+\|u\|_{H^1}^4) \|u\|_{H^2}^2 + \|\nabla\theta\|_2^2.
        \label{erntjkertnj}
    \end{align}
    Next, we focus on $\nabla \theta$.

    Taking the gradient of \eqref{nd_transport_equation} and testing it with $|\nabla \theta|^{{\tilde q}-2}\nabla\theta$
    \begin{align}
        \frac{1}{{\tilde q}}\frac{d}{dt} \|\nabla \theta\|_{\tilde q}^{\tilde q}
        = - \int_{\Omega} |\nabla\theta|^{{\tilde q}-2} \nabla \theta \cdot(\nabla\theta\cdot\nabla)u - \int_{\Omega} |\nabla\theta|^{{\tilde q}-2}\nabla \theta \cdot (u\cdot\nabla)\nabla\theta
        \\[-10pt]\label{njerktnejrt}
    \end{align}
    Note that second term on the right-hand side of \eqref{njerktnejrt} vanishes due to \eqref{nd_incompressible} and \eqref{nd_no_penetration_bc} as integration by parts yields
    \begin{align}
        \int_{\Omega} |\nabla\theta|^{{\tilde q}-2}\nabla \theta \cdot (u\cdot\nabla)\nabla\theta
        &= \frac{1}{2} \int_{\Omega} |\nabla\theta|^{{\tilde q}-2}u\cdot \nabla |\nabla\theta|^2 = \frac{1}{{\tilde q}} \int_{\Omega} u\cdot \nabla |\nabla\theta|^{{\tilde q}}
        \\
        &= \frac{1}{{\tilde q}} \int_{\partial\Omega} u\cdot n |\nabla\theta|^{{\tilde q}} - \frac{1}{{\tilde q}} \int_{\Omega} |\nabla\theta|^{{\tilde q}} \nabla\cdot u = 0\quad
        \label{nertjker}
    \end{align}
    Plugging \eqref{nertjker} into \eqref{njerktnejrt} and using Hölder's inequality one finds
    \begin{align}
        \frac{1}{{\tilde q}}\frac{d}{dt} \|\nabla \theta\|_{\tilde q}^{\tilde q}
        &= - \int_{\Omega} |\nabla\theta|^{{\tilde q}-2} \nabla \theta \cdot(\nabla\theta\cdot\nabla)u
        \leq \|\nabla\theta\|_{\tilde q}^{\tilde q} \|\nabla u\|_\infty.
        \label{renjtb}
    \end{align}
    To estimate the velocity term note that by the logarithmic Sobolev inequality (see Lemma 2.2 in \cite{doering2018longTime} for a proof) we obtain
    \begin{align}
        \|\nabla u\|_\infty \lessapprox (1+\|\nabla u\|_{H^1})\log(1+\|\nabla \Delta u\|_2)
    \end{align}
    and therefore, \eqref{renjtb} and Lemma \ref{lemma_elliptic_regularity_bounded_domain} yield
    \begin{align}
        \frac{1}{{\tilde q}}\frac{d}{dt} \|\nabla \theta\|_{\tilde q}^{\tilde q} &\leq \|\nabla \theta\|_{\tilde q}^{\tilde q} \|\nabla u\|_\infty
        \\
        &\lessapprox \|\nabla\theta\|_{\tilde q}^{\tilde q} (1+\|\nabla u\|_{H^1})\log(1+\|u\|_{H^3})
        \\
        &\lessapprox \|\nabla\theta\|_{\tilde q}^{\tilde q}(1+\|u\|_{H^2})\log(1+\|\Delta \omega\|_2 + \|u\|_{H^2}).
        \label{njkbnrjetb}
    \end{align}
    Note that
    \begin{align}
        \log(1+\|\Delta \omega\|_2 + \|u\|_{H^2})
        &\leq \log\left((1+\|\Delta \omega\|_2)(1+ \|u\|_{H^2})\right)
        \\
        &= \log(1+\|\Delta \omega\|_2)+\log(1+ \|u\|_{H^2})
        \\
        &\leq \log(1+\|\Delta \omega\|_2)+\|u\|_{H^2},
    \end{align}
    which applied to \eqref{njkbnrjetb} yields
    \begin{align}
        \frac{1}{{\tilde q}}\frac{d}{dt} \|\nabla \theta\|_{\tilde q}^{\tilde q}
        &\lessapprox \|\nabla\theta\|_{\tilde q}^{\tilde q}(1+\|u\|_{H^2})\log(1+\|\Delta \omega\|_2 + \|u\|_{H^2})
        \\
        &\leq \|\nabla\theta\|_{\tilde q}^{\tilde q}(1+\|u\|_{H^2}^2) +  \|\nabla\theta\|_{\tilde q}^{\tilde q}(1+\|u\|_{H^2})\log(1+\|\Delta \omega\|_2),\\[-5pt]
        \label{erjteb}
    \end{align}
    where we additionally applied Young's inequality. Adding \eqref{erjteb} to \eqref{erntjkertnj} implies
    \begin{align}
        \frac{d}{dt}&\left(\|\nabla\omega\|_2^2 + {\tilde q}^{-1}\|\nabla \theta\|_{\tilde q}^{\tilde q}+1\right) +\|\Delta \omega\|_2^2
        \\
        &\qquad =\frac{d}{dt}\left(\|\nabla\omega\|_2^2 + {\tilde q}^{-1} \|\nabla \theta\|_{\tilde q}^{\tilde q}\right) +\|\Delta \omega\|_2^2
        \\
        &\qquad\lessapprox \|u_t\|_{H^1}^2 + \|u\|_{H^1}^4 \|u\|_{H^2}^2 + \|u\|_{H^2}^2 + \|\nabla\theta\|_2^2 + \|\nabla\theta\|_{\tilde q}^{\tilde q}(1+\|u\|_{H^2}^2)
        \\
        &\qquad\qquad + \|\nabla\theta\|_{\tilde q}^{\tilde q}(1+\|u\|_{H^2})\log(1+\|\Delta \omega\|_2)
        \\
        &\qquad\lesssim \|u_t\|_{H^1}^2 + \|u\|_{H^1}^4 \|u\|_{H^2}^2 + \|u\|_{H^2}^2 + 1 + \|\nabla\theta\|_{\tilde q}^{\tilde q}(1+\|u\|_{H^2}^2)
        \\
        &\qquad\qquad + \|\nabla\theta\|_{\tilde q}^{\tilde q}(1+\|u\|_{H^2})\log(1+\|\Delta \omega\|_2),
        \label{ndjkfgndjfkg}
    \end{align}
    where in the last estimate we used Hölder's inequality. Let us further define
    \begin{align}
        Y(t) &= \|\nabla\omega\|_2^2 + \frac{1}{{\tilde q}} \|\nabla \theta\|_{\tilde q}^{\tilde q}+1
        \\
        Z(t) &= \|\Delta \omega\|_2^2
        \\
        A(t) &= C(1+{\tilde q})(\|u_t\|_{H^1}^2 + \|u\|_{H^1}^4 \|u\|_{H^2}^2 + \|u\|_{H^2}^2 + 1)
        \\
        B(t) &= C(1+{\tilde q})(1+\|u\|_{H^2}),
    \end{align}
    where $C$ is the implicit constant of \eqref{ndjkfgndjfkg}, the bound
    \begin{align}
        \frac{d}{dt} Y(t) + Z(t) \leq A(t)Y(t) + B(t)Y(t)\log(1+Z(t))
    \end{align}
    holds by \eqref{ndjkfgndjfkg}. Note that by Lemma \ref{lemma_nd_time_derivative_bound} and Lemma \ref{lemma_nd_LinftyH2_bound} one has $A\in L^1(0,T)$ for any $T>0$. Further Lemma \ref{lemma_nd_LinftyH2_bound} implies $B\in L^2(0,T)$. Therefore, the logarithmic Grönwall's inequality (see Lemma 2.3 in \cite{doering2018longTime}) implies
    \begin{align}
        Y(t)&<\infty
        \\
        \int_0^t Z(s)\ ds &< \infty
    \end{align}
    for any $0\leq t\leq T$, which in turn yields
    \begin{align}
        \nabla\theta &\in L^\infty\left((0,T);L^{\tilde q}(\Omega)\right),
        \label{nthjzbtrzbjrt}
        \\
        \Delta \omega &\in L^2\left((0,T);L^2(\Omega)\right).
        \label{vubhcvubhczv}
    \end{align}
    By \eqref{nd_maximum_principle}, it holds
    \begin{align}
        \|\theta\|_{W^{1,{\tilde q}}}^{\tilde q} = \|\theta\|_{\tilde q}^{\tilde q} + \|\nabla\theta\|_{\tilde q}^{\tilde q} = \|\theta_0\|_{\tilde q}^{\tilde q} + \|\nabla\theta\|_{\tilde q}^{\tilde q},
    \end{align}
    which combined with \eqref{nthjzbtrzbjrt} results in
    \begin{align}
        \theta\in L^\infty\left((0,T);W^{1,{\tilde q}}(\Omega)\right)
    \end{align}
    for any $T>0$ and due to Hölder's inequality also for any $1\leq q \leq {\tilde q}$. Similarly, by Lemma \ref{lemma_elliptic_regularity_bounded_domain}
    \begin{align}
        \|u\|_{H^3}^2 \lessapprox \|\Delta \omega\|_2^2+\|u\|_{H^2}^2,
    \end{align}
    which, by \eqref{vubhcvubhczv} and Lemma \ref{lemma_nd_LinftyH2_bound}, yields
    \begin{align}
        u\in L^2\left((0,T);H^3(\Omega)\right)
        \label{bhjtbrthzjtbr}
    \end{align}
    for any $T>0$. Finally, Gagliardo-Nirenberg interpolation implies
    \begin{align}
        \| u\|_{W^{2,r}}^\frac{2r}{r-2} \lesssim \|u\|_{H^3}^2 \|u\|_{H^2}^\frac{4}{r-2}
    \end{align}
    for any $2<r<\infty$, which, by \eqref{bhjtbrthzjtbr} and Lemma \ref{lemma_nd_LinftyH2_bound}, results in
    \begin{align}
        u\in L^\frac{2r}{r-2}\left((0,T);W^{2,r}(\Omega)\right)
    \end{align}
    for any $2\leq r <\infty$, where the limit case is covered by Lemma \ref{lemma_nd_LinftyH2_bound}.
\end{proof}

\section{Convergence}
\label{section:nd_convergence}

In order to show the convergence to the hydrostatic equilibrium we will extensively use the following Lemma, which is a slight variation of Lemma 3.1 in \cite{doering2018longTime}.

\begin{lemma}
    Let $f\geq 0$, $f\in L^1(0,\infty)$, and $f'(t)\leq C$ for some constant $C> 0$ and all $t\geq 0$. Then
    \begin{align}
        f(t)\to 0
    \end{align}
    as $t\to \infty$.
    \label{lemma_nd_fGeqZero_fLone_fPrimeLeqC_convergence}
\end{lemma}

For the proof of this Lemma, we need the following technical result, showing that for sufficiently large $t$, the function has to approach $0$ in every time interval.

\begin{lemma}
    Let $f\geq 0$, $f\in L^1(0,\infty)$ and $\epsilon,\delta>0$. Then there exists some $T_0>0$ such that for all $T\geq T_0$ there exists some $t\in [T,T+\delta]$ such that
    \begin{align}
        f(t)<\epsilon.
    \end{align}
    \label{lemma_nd_technical_helper_lemma}
\end{lemma}

\begin{proof}
    We prove this Lemma by contradiction. Assume there exists a sequence $T_k\to\infty$ such that for all  $t\in [T_k, T_k+\delta]$
    \begin{align}
        f(t)\geq \epsilon.
    \end{align}
    Choosing a subsequence $T_{k_i}$ such that the intervals do not overlap one finds
    \begin{align}
        \infty > \int_0^\infty f(t) \ dt \geq \sum_{i=0}^\infty \int_{T_{k_i}}^{T_{k_i}+\delta} f(t) \ dt \geq \sum_{i=0}^\infty \int_{T_{k_i}}^{T_{k_i}+\delta} \epsilon \ dt = \sum_{i=0}^\infty \epsilon\delta = \infty,
    \end{align}
    a contradiction.
\end{proof}

With this, we can prove Lemma \ref{lemma_nd_fGeqZero_fLone_fPrimeLeqC_convergence}.
\begin{proof}[Proof of Lemma \ref{lemma_nd_fGeqZero_fLone_fPrimeLeqC_convergence}]\\
    We will show that for any $\epsilon>0$ we can find $T>0$ such that for any $t>T$ one has $f(t)<\epsilon$.
    By Lemma \ref{lemma_nd_technical_helper_lemma}, there exists some $T_0>0$ such that for every $\tilde T\geq T_0$ there exists some $\tilde t\in [\tilde T,\tilde T+\frac{\epsilon}{4C}]$ such that
    \begin{align}
        f(\tilde t)< \frac{\epsilon}{2}.
    \end{align}
    Then for $t\in [\tilde T+\frac{\epsilon}{4C},\tilde T+\frac{\epsilon}{2C}]$
    \begin{align}
        f(t) &= f(\tilde t) + \int_{\tilde t}^{t} f'(t)\ dt 
        \leq f(\tilde t) + \int_{\tilde T}^{\tilde T+\frac{\epsilon}{2C}} C\ dt
        < \frac{\epsilon}{2}+C\frac{\epsilon}{2C}\leq \epsilon.
    \end{align}
    So for $t\geq T= T_0 + \frac{\epsilon}{4C}$ one has $f(t)< \epsilon$ concluding the proof.
\end{proof}

Now, we show that the appropriate norms fulfill the assumptions of Lemma \ref{lemma_nd_fGeqZero_fLone_fPrimeLeqC_convergence} and therefore decay in time. We start with the energy.

\begin{lemma}
    Let $\Omega$ be $C^{1,1}$, $u_0,\theta_0\in L^2(\Omega)$ and $0<\alpha\in L^\infty(\partial\Omega)$. Then
    \begin{align}
        \|u(t)\|_2\to 0
    \end{align}
    as $t\to\infty$. Additionally, for any $\beta>0$ and $\gamma\in \mathbb{R}$ there exists a constant $C\geq 0$ satisfying $C^2\leq \beta\|u_0\|_2^2 + \|\theta_0-\beta x_2-\gamma\|_2^2$ such that
    \begin{align}
        \|\theta(t)-\beta x_2-\gamma\|_2 \to C
    \end{align}
    for $t\to\infty$.
    \label{lemma_nd_u_to_0}
\end{lemma}

\begin{proof}
    By \eqref{rtenjtkrenj}
    \begin{align}
        \frac{d}{dt} \|u\|_2^2
        &\lesssim (1+\|\alpha^{-1}\|_\infty) \|\theta_0\|_2^2,
    \end{align}
    and by Lemma \ref{lemma_nd_L2H1_bound} $\|u(t)\|_2^2\in L^1(0,\infty)$. Therefore, Lemma \ref{lemma_nd_fGeqZero_fLone_fPrimeLeqC_convergence} directly implies
    \begin{align}
        \|u(t)\|_2^2\to 0\label{nrejtkernjt}
    \end{align}
    as $t\to\infty$. Next, by \eqref{vbcvbcggbvhcj}
    \begin{align}
        \frac{1}{2}\frac{d}{dt} &\left(\beta\|u\|_2^2 + \|\theta-\beta x_2 -\gamma\|_2^2 \right)
        \\
        &\qquad \leq \frac{1}{2}\frac{d}{dt} \left(\beta\|u\|_2^2 + \|\theta-\beta x_2 -\gamma\|_2^2 \right) + 2 \beta \|\D u\|_2^2 + 2 \beta \int_{\partial\Omega} \alpha u_\tau^2 = 0
    \end{align}
    for any $\beta >0$ and $\gamma\in \mathbb{R}$, implying that
    \begin{align}
        \beta\|u(t)\|_2^2 + \|\theta(t)-\beta x_2 -\gamma\|_2^2 
        \label{njtkertjnke}
    \end{align}
    is a non-increasing, non-negative function of time. By \eqref{nrejtkernjt}, it holds $\|u(t)\|_2^2\to 0$ for $t\to \infty$. Therefore, one obtains
    \begin{align}
        \|\theta(t)-\beta x_2-\gamma\|_2^2\to C^2
    \end{align}
    for $t\to\infty$ and some constant $C\geq 0$. As \eqref{njtkertjnke} is non-increasing in time, the constant can be bounded by $C^2\leq \beta\|u_0\|_2^2 + \|\theta_0-\beta x_2-\gamma\|_2^2$.
\end{proof}

Next, we show the decay of the vorticity.

\begin{lemma}
    Let $\Omega$ be $C^{2,1}$, $u_0\in H^2(\Omega)$, $\theta_0 \in L^{4}(\Omega)$ and $0<\alpha \in W^{1,\infty}(\partial\Omega)$. Then for any $1\leq q<\infty$
    \begin{align}
        \|u(t)\|_{W^{1,q}} \to 0
    \end{align}
    as $t\to \infty$.
    \label{lemma_nd_uW1q_decay}
\end{lemma}

\begin{proof}
    Testing \eqref{nd_vorticity_pde} with $\omega$ one has
    \begin{align}
        \frac{1}{2}\frac{d}{dt} \|\omega\|_2^2 = - \int_{\Omega} \omega u\cdot\nabla \omega + \int_{\Omega} \omega\Delta\omega + \int_{\Omega} \omega \partial_1 \theta.
        \label{rktjertnbjbgkn}
    \end{align}
    Note that the first term on the right-hand side vanishes due to \eqref{nd_no_penetration_bc} and \eqref{nd_incompressible} as
    \begin{align}
        2\int_{\Omega} \omega u\cdot\nabla \omega = \int_{\Omega} u\cdot \nabla (\omega^2) = \int_{\partial\Omega} \omega^2 u\cdot n - \int_{\Omega} \omega^2 \nabla\cdot u = 0.
        \label{rtjkbbrdsfasdfdsa}
    \end{align}
    For the second term on the right-hand side of \eqref{rktjertnbjbgkn}, integration by parts and \eqref{nd_vorticity_bc} yield
    \begin{align}
        \int_{\Omega} \omega \Delta \omega &= \int_{\partial\Omega} \omega n\cdot\nabla \omega - \|\nabla\omega\|_2^2
        = -2 \int_{\partial\Omega} (\alpha+\kappa)u_\tau n\cdot\nabla \omega - \|\nabla\omega\|_2^2\qquad
        \label{jbetjhber}
    \end{align}
    In order to handle the boundary term, note that \eqref{id_Delta_u_is_nablaPerp_omega} and tracing \eqref{nd_navier_stokes} along the boundary implies
    \begin{align}
        n\cdot\nabla \omega = n^\perp\cdot \nabla^\perp \omega = \tau\cdot\Delta u = \tau\cdot u_t + \tau \cdot(u\cdot\nabla) u+\tau \cdot \nabla p - \theta \tau\cdot e_2.\qquad
        \label{rtznrthjzbtrhjz}
    \end{align}
    Combining \eqref{jbetjhber} and \eqref{rtznrthjzbtrhjz}, one has
    \begin{align}
        \int_{\Omega} \omega \Delta \omega &= - \|\nabla\omega\|_2^2 -2 \int_{\partial\Omega} (\alpha+\kappa)u_\tau n\cdot\nabla \omega
        \\
        &= - \|\nabla\omega\|_2^2 -2 \int_{\partial\Omega} (\alpha+\kappa)u \cdot u_t -2 \int_{\partial\Omega} (\alpha+\kappa)u \cdot (u\cdot\nabla) u
        \\
        &\qquad -2 \int_{\partial\Omega} (\alpha+\kappa)u \cdot \nabla p + 2 \int_{\partial\Omega} (\alpha+\kappa)u_\tau \tau_2 \theta
        \\
        &= - \|\nabla\omega\|_2^2 - \frac{d}{dt} \int_{\partial\Omega} (\alpha+\kappa) u_\tau^2 -2 \int_{\partial\Omega} (\alpha+\kappa)u \cdot (u\cdot\nabla) u
        \\
        &\qquad -2 \int_{\partial\Omega} (\alpha+\kappa)u \cdot \nabla p -  \int_{\partial\Omega} \omega \tau_2 \theta
        \label{betrhjreh}
    \end{align}
    For the third term on the right-hand side of \eqref{rktjertnbjbgkn} integration by parts yields
    \begin{align}
        \int_{\Omega} \omega\partial_1 \theta = \int_{\partial\Omega} n_1 \omega \theta - \int_{\Omega} \theta \partial_1 \omega.
        \label{brhtbtrhrbh}
    \end{align}
    Combining \eqref{rktjertnbjbgkn}, \eqref{rtjkbbrdsfasdfdsa}, \eqref{betrhjreh} and \eqref{brhtbtrhrbh}, it holds
    \begin{align}
        \frac{1}{2}\frac{d}{dt} &\left(\|\omega\|_2^2 + 2\int_{\partial\Omega} (\alpha+\kappa) u_\tau^2\right) + \|\nabla\omega\|_2^2
        \\
        &\qquad =  -2 \int_{\partial\Omega} (\alpha+\kappa)u \cdot (u\cdot\nabla) u - 2 \int_{\partial\Omega} (\alpha+\kappa)u \cdot \nabla p -  \int_{\partial\Omega} \omega \tau_2 \theta 
        \\
        &\qquad\qquad + \int_{\partial\Omega} n_1 \omega \theta - \int_{\Omega}  \theta \partial_1 \omega
        \\
        &\qquad =  -  \int_{\partial\Omega} (\alpha+\kappa)u \cdot \nabla (u\cdot u) - 2 \int_{\partial\Omega} (\alpha+\kappa)u \cdot \nabla p - \int_{\Omega} \theta \partial_1 \omega,\qquad
        \label{erjtetbjh}
    \end{align}
    where in the last identity we used that the boundary terms cancel as $\tau = n^\perp$ and therefore $\tau_2 = n_1$. By Lemma \ref{lemma_u_cdot_nabla_p_on_boundary} and an application of Hölder's and Young's inequality
    \begin{align}
        \left|\int_{\partial\Omega} (\alpha+\kappa)u \cdot \nabla (u\cdot u)\right| &\lessapprox \|u\|_{W^{1,3}} \|u^2\|_{W^{1,\frac{3}{2}}} \lesssim \|u\|_{W^{1,3}}^3
        \label{rekjterb}
    \end{align}
    and
    \begin{align}
        \left|\int_{\partial\Omega} (\alpha+\kappa)u \cdot \nabla p\right| &\lessapprox \|u\|_{H^1}\|p\|_{H^1}.
        \label{brthjbre}
    \end{align}
    The last term on the right-hand side can be estimated by Hölder's and Young's inequality as
    \begin{align}
        \left|\int_{\Omega} \theta \partial_1\omega \right|\leq \|\theta\|_2\|\partial_1\omega\|_2\leq \epsilon \|\nabla\omega\|_2^2 + \epsilon^{-1}\|\theta\|_2^2 = \epsilon \|\nabla\omega\|_2^2 + \epsilon^{-1}\|\theta_0\|_2^2,\qquad
        \label{bejktbetjh}
    \end{align}
    where we also used \eqref{nd_maximum_principle}. Additionally, note that by Lemma \ref{lemma_grad_ids}
    \begin{align}
        \|\omega\|_2^2 + 2\int_{\partial\Omega}\kappa u_\tau^2 = 2\|\D u\|_2^2.
        \label{rejntkbrhj}
    \end{align}
    Combining \eqref{erjtetbjh}, \eqref{rekjterb}, \eqref{brthjbre}, \eqref{bejktbetjh} and \eqref{rejntkbrhj} yields
    \begin{align}
        \frac{d}{dt} &\left(\|\D u\|_2^2 + \int_{\partial\Omega} \alpha u_\tau^2\right) + \|\nabla\omega\|_2^2
        \\
        &\qquad \qquad \lessapprox \|u\|_{W^{1,3}}^3 +\|u\|_{H^1}\|p\|_{H^1} + \epsilon \|\nabla \omega\|_2^2 + \epsilon^{-1} \|\theta_0\|_2^2
    \end{align}
    and choosing $\epsilon$ sufficiently small one has
    \begin{align}
        \frac{d}{dt} \left(\|\D u\|_2^2 + \int_{\partial\Omega} \alpha u_\tau^2\right)
        &\leq \frac{d}{dt} \left(\|\D u\|_2^2 + \int_{\partial\Omega} \alpha u_\tau^2\right) + \frac{1}{2}\|\nabla\omega\|_2^2
        \\
        &\lessapprox \|u\|_{W^{1,3}}^3 +\|u\|_{H^1}\|p\|_{H^1} + \|\theta_0\|_2^2
    \end{align}
    Note that by Lemma \ref{lemma_nd_W1q_bound} and Lemma \ref{lemma_nd_pressure_bound}, the right-hand side is uniformly bounded in time and by trace theorem and Lemma \ref{lemma_nd_L2H1_bound}, it holds
    \begin{align}
        \|\D u\|_2^2 +\int_{\partial\Omega} \alpha u_\tau^2 \lesssim (1+\|\alpha\|_\infty) \|u\|_{H^1}^2 \in L^1(0,\infty).
    \end{align}
    Therefore, $\|\D u\|_2^2 + \int_{\partial\Omega}\alpha u_\tau^2$ fulfills the assumptions of Lemma \ref{lemma_nd_fGeqZero_fLone_fPrimeLeqC_convergence} implying
    \begin{align}
        \|\D u(t)\|_2^2 + \int_{\partial\Omega} \alpha u_\tau^2(t) \to 0
    \end{align}
    as $t\to \infty$ and due to Lemma \ref{lemma_coercivity}
    \begin{align}
        \|u(t)\|_{H^1}^2 \lessapprox \|\D u(t)\|_2^2 + \int_{\partial\Omega} \alpha u_\tau^2(t) \to 0
    \end{align}
    for $t\to\infty$. Finally as $\|u\|_{H^2}$ is uniformly bounded in time due to Lemma \ref{lemma_nd_LinftyH2_bound}, Gagliardo-Nirenberg interpolation yields
    \begin{align}
        \|u(t)\|_{W^{1,q}}^q \lesssim \|u(t)\|_{H^2}^{q-2} \|u(t)\|_{H^1}^2\to 0
    \end{align}
    as $t\to\infty$ for any $2 \leq q< \infty$. The, Hölder's inequality implies the decay for any $1\leq q \leq 2$.
\end{proof}

Next, we show that the time derivative of the $L^2$-norm of $u_t$ is uniformly bounded in time, which together with the previous regularity estimates proves its decay.

\begin{lemma}
    Let $\Omega$ be $C^{1,1}$, $u_0\in H^2(\Omega)$, $\theta_0 \in L^{4}(\Omega)$ and $0<\alpha \in L^\infty(\partial\Omega)$. Then
    \begin{align}
        \|u_t(t)\|_2 \to 0
    \end{align}
    as $t\to\infty$.
    \label{lemma_nd_ut_decay}
\end{lemma}

\begin{proof}
    By \eqref{jerjrewnbjkr}
    \begin{align}
        \frac{d}{dt}\|u_t\|_2^2 &\leq \frac{d}{dt}\|u_t\|_2^2 + C \|u_t\|_{H^1}^2 
        \\
        &\lesssim\|u_t\|_2^2 \|u\|_{H^1}^2 + \|\theta_0\|_4^2 \|u\|_{H^1}^2
        \label{jnbwejrbwejhrwb}
    \end{align}
    and by Lemma \ref{lemma_nd_time_derivative_bound} and Lemma \ref{lemma_nd_W1q_bound} the right-hand side of \eqref{jnbwejrbwejhrwb} is uniformly bounded in time. Additionally, by Lemma \ref{lemma_nd_time_derivative_bound} $\|u_t(t)\|_2^2\in L^1(0,\infty)$, so it fulfills the assumptions of \eqref{lemma_nd_fGeqZero_fLone_fPrimeLeqC_convergence}, implying
    \begin{align}
        \|u_t(t)\|_2^2 \to 0
    \end{align}
    as $t\to\infty$.
\end{proof}

As $u$ decays, one can now separate the other terms in \eqref{nd_navier_stokes} and show their decay. Note that we do not have decay estimates for $\|u\|_{H^2}$. Therefore, the diffusion term is covered by Lemma \ref{lemma_int_Delta_u_v}. This has the disadvantage of only achieving weak convergence.

\begin{lemma}
    Let $\Omega$ be $C^{2,1}$, $u_0\in H^2(\Omega)$, $\theta_0 \in L^{4}(\Omega)$ and $0<\alpha \in W^{1,\infty}(\partial\Omega)$. Then
    \begin{align}
        \|(\nabla p-\theta e_2)(t)\|_{H^{-1}}\to 0
    \end{align}
    for $t\to\infty$.
\end{lemma}

\begin{proof}
    Let $v\in H^1$, then by \eqref{nd_navier_stokes}
    \begin{align}
        \langle \nabla p-\theta e_2,v\rangle &= \langle -u_t-u\cdot\nabla u +\Delta u,v\rangle = -\langle u_t,v\rangle - \langle u\cdot\nabla u,v\rangle + \langle \Delta u,v\rangle
        \\[-3pt]
        \label{renjtkerjnt}
    \end{align}
    Note that the first and second term on the right-hand side can be estimated, using Hölder's inequality, by
    \begin{align}
        |\langle u_t,v\rangle|\leq \|u_t\|_2 \|v\|_2
        \label{nerjtkren}
    \end{align}
    and
    \begin{align}
        |\langle u\cdot\nabla u,v\rangle &\leq \|u\|_4\|u\|_{H^1}\|v\|_4 \leq \|u\|_{H^1}^2\|v\|_{H^1},
        \label{ntjktrn}
    \end{align}
    where, in the last estimate, we also used Ladyzhenskaya's inequality.
    Due to Lemma \ref{lemma_int_Delta_u_v}, the third term on the right-hand side of \eqref{renjtkerjnt} satisfies
    \begin{align}
        |\langle \Delta u,v\rangle| \leq 2\|\D u\|_2\|\D v\|_2 + 2\int_{\partial\Omega} |\alpha u\cdot v| \lesssim (1+\|\alpha\|_\infty) \|u\|_{H^1}\|v\|_{H^1},\qquad
        \label{nerjtkrej}
    \end{align}
    where in the last inequality we used the trace theorem. Combining \eqref{renjtkerjnt}, \eqref{nerjtkren}, \eqref{ntjktrn} and \eqref{nerjtkrej} yields
    \begin{align}
        \|\nabla p-\theta e_2\|_{H^{-1}}
        &= \sup_{\|v\|_{H^1}=1} \langle \nabla p -\theta e_2, v \rangle
        \\
        &= \sup_{\|v\|_{H^1}=1} \left( -\langle u_t,v\rangle - \langle u\cdot\nabla u,v\rangle +\langle\Delta u,v\rangle \right)
        \\
        &\lessapprox \sup_{\|v\|_{H^1}=1} \left( (\|u_t\|_2 + \|u\|_{H^1}^2 + \|u\|_{H^1}) \|v\|_{H^1} \right)
        \\
        &= \|u_t\|_2 + \|u\|_{H^1}^2 + \|u\|_{H^1}.
        \label{berjtbhetrjbher}
    \end{align}
    Due to Lemma \ref{lemma_nd_ut_decay} and Lemma \ref{lemma_nd_uW1q_decay}, the right-hand side of \eqref{berjtbhetrjbher} vanishes in the long-time limit, implying
    \begin{align}
        \|(\nabla p-\theta e_2)(t)\|_{H^{-1}}\to 0
    \end{align}
    for $t\to\infty$.
\end{proof}

\chapter*{References}
\addcontentsline{toc}{chapter}{References}
\printbibliography[heading=none]
\end{document}